\documentclass[12pt]{amsart}
\usepackage{latexsym,amssymb,amsmath,amsthm}
\usepackage{comment,graphicx,color}
\usepackage{morefloats}
\usepackage[section]{placeins}
\definecolor{red1}{rgb}{1,0.9,0.9}
\definecolor{blue1}{rgb}{0.9,0.9,1}
\definecolor{green1}{rgb}{0.9,1,0.9}
\definecolor{yellow1}{rgb}{1,1,0.9}
\definecolor{yellow2}{rgb}{1,1,0.8}

\def\question#1{ \vspace{2mm} \begin{center} \fcolorbox{green1}{green1}{ \parbox{11.2cm}{{\bf Question:} #1}} \vspace{2mm} \end{center} }
\def\conjecture#1{ \vspace{2mm} \begin{center} \fcolorbox{green1}{green1}{ \parbox{11.2cm}{{\bf Conjecture:} #1}} \vspace{2mm} \end{center} }
\def\resultlemma#1{ \vspace{2mm} \begin{center} \fcolorbox{yellow1}{yellow1}{ \parbox{11.2cm}{{\bf Lemma:} #1}} \vspace{2mm} \end{center} }
\def\resulttheorem#1{ \vspace{2mm} \begin{center} \fcolorbox{yellow1}{yellow1}{ \parbox{11.2cm}{{\bf Theorem:} #1}} \vspace{2mm} \end{center} }
\def\resultcorollary#1{ \vspace{2mm} \begin{center} \fcolorbox{yellow1}{yellow1}{ \parbox{11.2cm}{{\bf Corollary:} #1}} \vspace{2mm} \end{center} }
\def\definition#1{ \vspace{2mm} \begin{center} \fcolorbox{red1}{red1}{ \parbox{11.2cm}{{\bf Definition:} #1}} \vspace{2mm} \end{center} }

\theoremstyle{definition}

\def\Dcal{\mathcal{D}} 
\def\Wcal{\mathcal{W}} 
\def\Ecal{\mathcal{E}} 
\def\Ccal{\mathcal{C}} 
\def\Gcal{\mathcal{G}} 
\def\Scal{\mathcal{S}} 
\def\Bcal{\mathcal{B}}
\def\Vcal{\mathcal{V}}
\def\Pcal{\mathcal{P}}
\def\Hcal{\mathcal{H}}

\newcommand{\R}{\mathbb{R}}
\newcommand{\Z}{\mathbb{Z}}

\title{Coloring graphs using topology}
\author{Oliver Knill}
\date{December 21, 2014}
\address{
        Department of Mathematics \\
        Harvard University \\
        Cambridge, MA, 02138
        }

\begin{document}
\maketitle

\begin{abstract}
Higher dimensional graphs can be used to color $2$-dimensional geometric graphs $G$. 
If $G$ the boundary of a three dimensional graph $H$ for example, we can refine the interior 
until it is colorable with $4$ colors. The later goal is achieved if all interior edge degrees
are even. Using a refinement process which cuts the interior along surfaces $S$
we can adapt the degrees along the boundary $S$. More efficient is a self-cobordism of $G$ with itself
with a host graph discretizing the product of $G$ with an interval.  
It follows from the fact that Euler curvature is zero everywhere for three dimensional geometric graphs,
that the odd degree edge set $O$ is a cycle and so a boundary if $H$ is simply connected. A reduction 
to minimal coloring would imply the four color theorem. The method is expected to give a reason ``why 4 colors
suffice" and suggests that every two dimensional geometric graph of arbitrary degree and orientation
can be colored by $5$ colors: since the projective plane can not be a boundary of a $3$-dimensional graph
and because for higher genus surfaces, the interior $H$ is not simply connected,
we need in general to embed a surface into a 4-dimensional simply connected graph in order to color it.
This explains the appearance of the chromatic number $5$ for higher degree or non-orientable situations,
a number we believe to be the upper limit.
For every surface type, we construct examples with chromatic number $3,4$ or $5$, where the construction of 
surfaces with chromatic number $5$ is based on a method of Fisk. 
We have implemented and illustrated all the topological aspects described in this paper on a computer. 
So far we still need human guidance or simulated annealing to do the refinements in 
the higher dimensional host graph. 
\end{abstract}

\section{Spheres}

A finite simple graph $G=(E,V)$ is {\bf 2-dimensional geometric},
if every unit sphere $S(x)$ in $G$ is a circular graph $C_n$ with $n \geq 4$ vertices.
Let $\Gcal_2$ denote the class of these graphs. 
From the five platonic solids, only the icosahedron or octahedron are in $\Gcal_2$. 
By definition a graph $G \in \Gcal_2$ does not contain tetrahedral subgraphs. 
If $v=|V|$ is the number of vertices, $e=|E|$ is the 
number of edges and $f$ is the number of triangles in $G$, then
the {\bf Euler characteristic} of $G \in \Gcal_2$ is
defined by {\bf Euler's formula} $\chi(G) = v-e+f$. Let $\Omega_0$ be the set of functions 
on $V$, let $\Omega_1$ denote the set of functions on $E$ and let $\Omega_2$ be the set of 
functions on triangles. $\Omega_k$ can be realized as anti-symmetric functions in $k+1$ variables. 
Given an ordering of vertices on $E$ and $F$, the {\bf gradient} $d_0: \Omega_0 \to \Omega_1$ is
defined as $df(a,b)=f(b)-f(a)$ and the {\bf curl} is the linear map 
$d_1: \Omega_1 \to \Omega_2$ as $dF((a,b,c)) = F((b,c))-F(a,c))+F((b,c))$. To see the structure better, 
define $d_k=0$ for $k \geq 2$. The {\bf scalar Laplacian} $L_0=d_0^* d_0$ and {\bf 1-form Laplacian} 
$L_1=d_1^* d_1 + d_2 d_2^*$ and {\bf 2-form Laplacian} $L_2=d_2^* d_2 + d_1 d_1^*$ are symmetric 
matrices which do not depend on the basis defined by the chosen ordering of edges and triangles. 
The {\bf cohomology groups} of $G$ are defined as $H^k(G)={\rm ker}(d_k)/{\rm im}(d_{k-1})$ with
$d_{-1}=0$. The zeroth {\bf Betti number} $b_0={\rm dim}(H^0(G))={\rm ker}(d_0)$ counts the number 
of connectivity components of $G$, the first Betti number 
$b_1 = {\rm dim}(H^1(G)) = {\rm dim}({\rm ker}(d_1)/{\rm im}(d_0))$ counts twice the number of holes in 
the discrete surface $G$ and $b_2= {\rm dim}(H^2(G)) = {\rm dim}(\Omega_2)/{\rm im}(d_1)$ counts the
dimension of volume forms on $G$. 
If {\bf Poincar\'e duality} $b_2={\rm dim}(H^2(G))={\rm dim}({\rm ker}(d_2)/{\rm im}(d_1)) 
= b_0$ holds, then $G$ is called {\bf orientable}. The {\bf Euler-Poincar\'e formula}, which follows from 
the rank-nullity theorem in linear algebra, tells that $\chi(G) = b_0-b_1+b_2$. The Euler characteristic is 
$2-b_1$ if $G \in \Gcal_3$ is orientable and connected. In that case, the {\bf genus} $g=b_0/2 =1-\chi(G)/2$ of $G$ 
determines the topological type of $G$. For genus $g=0$, the graph is a sphere with $\chi(G)=2$. 
A graph is called {\bf 4-connected}, if removing $1$,$2$ or $3$ vertices leaves $G$ connected. 
A graph $G=(V,E)$ is {\bf maximal planar}, if adding an additional edge while keeping $V$ destroys 
planarity. One can characterize $\Scal_2$ as follows:

\resultlemma{ $\Scal_2 = \{ G \; | \; G$ is 4-connected and maximal planar $\}$.  }

The appendix contains a proof. The class of 4-connected, maximal planar graphs has been studied by Whitney
and Tutte who showed that such graphs are Hamiltonian \cite{SaatyKainen}. 
Related is the {\bf Steinitz theorem} in polyhedral combinatorics which tells that every 
3-connected planar graph is the graph of a convex polyhedron. 
If we take such a graph and make it maximal, we get a convex polyhedron. 
It can still have triangular ``necks" however. We need 4-connectivity to be in $\Gcal_2$. \\

Given a finite simple graph $G$. If $v_i$ is the number of subgraphs $K_{i+1}$ 
then the {\bf Euler polynomial} $e(x)=\sum_i v_i x^i$ and {\bf Poincar\'e polynomial} $p(x)=\sum_i b_i x^i$ 
satisfy the {\bf Euler-Poincar\'e formula}
$e(-1)=p(-1)$. The matrices $d_i: \Omega_i \to \Omega_{i+1}$ define {\bf exterior derivative} $d$
which is a lower triangular $e(1) \times e(1)$ matrix satisfying $d^2=0$. These are large matrices
in general. For the octahedron already, the smallest graph in $\Gcal_2$, we deal 
with $21 \times 21$ matrices as $e(1)=6+12+8=21$. 
While the {\bf Dirac operator} $D=d+d^*$  \cite{KnillILAS,knillmckeansinger,CauchyBinetKnill}
depends on the chosen ordering of simplices, the {\bf form Laplacian} $L=D^2$ 
containing the block diagonal $L_k$ does not. The Dirac matrix $D$ is useful: 
{\bf the wave equation} $u_{tt} = -L u$ on forms 
for example has the explicit solution $u(t) = \cos(D t) u(0) + \sin(D t) D^{-1} u'(0)$
if the initial velocity is in the orthogonal complement of the kernel of the symmetric matrix $D$. 
The kernel of $L_k$ is the space of {\bf harmonic $k$-forms}. 
Its dimension is $b_k$ by the Hodge theorem. The {\bf heat evolution} $e^{-Lt} f$ converges for $t \to \infty$
to harmonic forms. 
It gives a concrete {\bf Hurewizc homomorphism} from the fundamental group $\pi_1(G)$
to the first cohomology group $H^1(G)$: 
a path $\gamma$ defines naturally a \{0-1\}-valued function $F$ in $\Omega_1$, which 
is discrete version of a {\bf de Rham current}. Applying the heat flow on $F$ gives in the limit $t \to \infty$
an element in $H^1(G)$. Since the initial embedding map is a homomorphism and the heat flow is linear,
this is a homomorphism. Other results like the {\bf McKean-Singer formula}  \cite{McKeanSinger}
${\rm str}(e^{-t L}) = \chi(G)$ relate for $t=0$ the definition of $\chi(G)$ 
with $\chi(G)=b_0-b_1+\dots$ in the limit $t \to \infty$. 
The discrete McKean-Singer result \cite{knillmckeansinger} 
can be used to construct form-isospectral non-isomorphic graphs \cite{knillmckeansinger}.
Calculus and in particular the {\bf Green-Stokes theorem} on $\Scal_2$ is known 
since Kirkhoff \cite{SunadaFundamental}.  The matrices $d_i$ and cohomology were introduced by Poincar\'e.
The work of Veblen played an important role \cite{BarnettEarlyWritings}, especially \cite{Veblen1912} which
dealt with the four color problem.
To formulate Stokes, take a simple closed path $C$ bounding a surface $S$ in $G$.
For $F \in \Omega_1$, we have $\int_C F =\int_S d_1 F$
which is {\bf Stokes theorem}. The left hand side sums up the values of $F$ on the edges, the later 
sums the curls in the interior.  Writing each curl as a line integral along the triangle 
leaves only the boundary as the interior terms cancel.  
Given a path $C$ starting at $a \in V$ and $b \in V$, the {\bf fundamental theorem of line integrals} 
gives $\int_C df = f(b)-f(a)$. 

\section{Sphere coloring} 

Let $\Pcal$ denote the set of planar graphs. Denote by $\Ccal_c$ the set of graphs which can be vertex
colored with $c$ or less colors. The argument will be to color each $4$-connectivity component in 
$\Ccal_4 \cap \Pcal$ separately and glue them together after possibly permuting the colors in each component and
because maximality assures that critical sets, where the $4$-connectivity components break must be triangles.

\resulttheorem{ $\Scal_2 \subset \Ccal_4$ is equivalent to $\Pcal \subset \Ccal_4$.  }

\begin{proof}
(i) Assume first that $\Scal_2 \subset \Ccal_4$. Given $G \in \Pcal$. Denote by $\overline{G}$ 
its maximal planar completion. We will show $\overline{G} \in \Ccal_4$ implying  $G \in \Ccal_4$. 
Let $G_i$ denote the 4-connectivity components of $\overline{G}$. Each of them is maximal planar.
By the lemma, each is in $\Scal_2$. By assumption, it can be colored by $4$ colors. Now,
glue together the individual components one by one. 
It follows from maximal planarity that the vertices
along we have cut, form a triangle, so that the coloring gives each of its vertices a
different color.
After coloring each component and doing a color permutation, the pieces can be 
glued back together. Now $\overline{G}$ is colored with $4$ colors. \\
(ii) To show the reverse, assume that $\Pcal \subset \Ccal_4$. Given $G \in \Scal_2$.
By the lemma, it is a 4-connected maximal planar graph. Especially, it is in $\Pcal$, 
As by assumption $\Pcal \subset \Ccal_4$, we know that $G \in \Ccal_4$. We have 
shown $\Scal_2 \subset \Ccal_4$. 
\end{proof} 

The fact that the four color problem can be studied by looking at discrete spheres
continues what Cayley or Tait stated, when introducing {\bf cubic maps},
maps which have triangulated graphs as their dual graph. A graph is {\bf cubic}, if all vertex degrees
are 3. Let $\Hcal$ denote the class of Hamiltonian graphs. Whitney showed $\Hcal \cap \Pcal \subset \Ccal_4$
is equivalent to $\Pcal \subset \Ccal_4$. Whitney also showed that a 4-connected, maximal planar graph is 
Hamiltonian. As by the lemma, 4-connected maximal planar graphs coincide with $\Scal_2$.
Whitney might have known the lemma proven here, but we have derived it without the detour 
over Hamiltonian graphs. The reduction to 4-connectivity is super-seeded by
``splicing" of graphs to a larger graph which according to \cite{Stromquist} 
{\it "has been simultaneously discovered by almost every other researcher in the field"}. 
The same remark might hold for the lemma. All graphs $G$ in $\Gcal_2$ have cubic dual graphs $\hat{G}$.
The assumption to have circular unit spheres forces graphs in $\Gcal_2$ to be $4$-connected. 
In differential geometric language, a circular unit sphere at $x$ means 
{\bf positive radius of injectivity} at $x$. 
Gluing two octahedra along a triangle is an example of a maximal planar graph which is 
{\bf 4-disconnected}. This graph has unit spheres which contain a star graph $S_3$,
the geodesic circle $\exp_1(x) = S(x)$ of radius $1$ already contains a {\bf cut locus}. 
One can check quickly whether a given graph $G$ is 
in $\Scal_2$ or not: find all unit spheres $S(x)$, determine whether they are cyclic graphs 
$C_{n(x)}$ with $n(x) \geq 4$, then compute the Euler characteristic with the 
{\bf planar Gauss-Bonnet formula}
$\chi(G) = \sum_{x \in V} (1-n(x)/6)$ and see whether it is $2$. The planar Gauss-Bonnet formula
must have been known since the 19th century while the general Gauss-Bonnet-Chern case is only more
recent \cite{cherngaussbonnet}. The reduction to 4-connected graphs holds only for simply connected graphs. 
Stacking tetrahedra along triangles in a circular matter leads to difficulties which are beyond the 
complexity to color the individual tetrahedra  as we get an action of the fundamental group on the 
permutation group of the four colors. 

\section{Discrete surfaces}

Since the four color theorem is equivalent to $\Scal_2 \subset \Ccal_4$, one can ask whether
every orientable $G \in \Gcal_2$ is in $\Ccal_4$. Fisk has shown that the answer is no. Under a minimal condition
for the minimal length (removed in \cite{Thomassen}), toroidal graphs in $\Gcal_2$ are 
in $\Ccal_5$ \cite{AlbertsonStromquist}
and there are examples of any genus not in $\Ccal_4$: one possibility is to have
two odd degree vertices exist which are adjacent \cite{Fisk1978}. A reformulation of a question
raised by Albertson-Stromquist (\cite{AlbertsonStromquist} question 1) 
asks whether every toroidal $G \in \Gcal_2$ is in $\Ccal_5$. In the light of {\bf Fisk theory}
and \cite{Thomassen}, one can be more bold and strengthen the question to a 
stronger conjecture

\conjecture{ $\Gcal_2 \subset \Ccal_5$. }

\cite{AlbertsonStromquist} gave an affirmative answer for the genus 1 case if the length of the 
shortest homotopically nontrivial circle is $8$ or more. They ask the question which
\cite{Thomassen} answers adding  more confidence that the conjecture should hold.
In a proposal to a HCRP project of Summer 2014 with Jenny Nitishinskaya,
we explored graphs $\Gcal_2 \subset \Gcal_3$ and the hypothesis $\Gcal_2 \subset \Gcal_4$ 
with geometric methods. 
The reason for working on this was that at that time, foliation methods
looked like a good way to color the graph. 
The projective plane found by Jenny forced us to assume orientability. 
We learned about Fisk theory in December 2014. The relevant work is
\cite{AlbertsonStromquist,Fisk1978,Fisk1977a,Fisk1977b,Fisk1980}.
As an orientable $G \in \Ccal_3$ can be embedded in a 4-dimensional sphere $H \in \Scal_4$ with
a simply connected complement $H \setminus G$, we believe that the complement of that surface can 
be chopped up to be in $\Ccal_5$. This embedding and refinement procedure using embeddings is 
the main point of the present story. It is because of the embedding that the above
claim $\Gcal_{2} \subset \Ccal_5$ is a reasonable conjecture. \\

A future proof of the four color theorem could even strengthen it by 
answering affirmatively the following question:

\question{ Is every $G \in \Scal_2$ a subgraph of some $H \in \Scal_3 \cap \Ccal_4$? }

We believe even that similarly, as any $G \in \Scal_1$ is a subgraph of a refinement 
$H \in \Scal_2 \cap \Scal_3$ of the octahedron, any $G \in \Scal_2$ could
be realized as a subgraph of a 
refinement $H \in \Scal_3 \cap \Ccal_4$ of the 16-cell, one of the 4D platonic solids
and more generally that any $G \in \Scal_d$ can be realized as a subgraph of a refinement 
$H \in \Scal_{d+1} \cap \Ccal$ of the $(d+1)$-dimensional cross polytop.  \\
Refinements can be done in various ways without changing the interior 
degrees of edges: we can cut edges or glue in 16-cells or other spheres in 
$\Scal_3$. The embedding question could be compared to the {\bf Nash embedding problem} 
in Riemannian geometry, where the question is when a given manifold
can be embedded isometrically in a higher dimensional Euclidean space. 
But since the refinements allow to tune
the curvature of the host graph pretty arbitrarily, the discrete question could be easier and
the 16 cell could be {\bf universal} in the sense that any $G \in \Scal_2$ can be embedded 
into a suitable refinement and so be in $\Ccal_4$. \\

While for the sphere, the four color theorem shows that there 
is a dichotomy of chromatic numbers $c=3$ or $c=4$ 
and where $\Scal_2 \cap \Ccal_3$ consists of exactly the Eulerian graphs, 
there are at least 3 cases for any surface type. None is empty: 

\resulttheorem{
For $c=3,4,5$ and any genus $g>0$ and orientation $o$,
there is $G \in \Gcal_2$ with chromatic number $c$ and type $(g,o)$.
}

\begin{proof} 
Because the sphere is excluded, the only simply connected case is the 
projective plane. Concrete implementations for $c=3,4,5$ are given in
the figure. The case $c=5$ has also been found first by Fisk. 
In the torus case, the situation $c=3$ is obtained for the flat torus $T_{3k,3k}$
with integer $k \geq 2$. All other flat tori $T_{n,n}$ are in $\Ccal_4$. 
Fisk has given examples with $c=5$: the simplest way to get them is to start
with $T_{4k,4k}$ with $k \geq 1$, then draw a geodesic $\gamma$ of type $(2,1)$ of
length $8k$. Now subdivide along $k$ edges so that the geodesic could be drawn
as a straight line of slope $2$ in a universal cover. Now Dehn twist the 
torus along a fundamental cycle by $1$. This has the effect that the 
curve $\gamma$ now opens up and that the initial point $A$ and end point $B$
are neighbors. The vertices $A,B$ are now the only ones with odd degree.
An odd degree vertex already forces $c=4$ at least but due to the fact that 
along $\gamma$, all vertex degrees are even, the holonomy transport along  the ``Fisk chain"
$\gamma$ does not allow the colorings of $A$ and $B$ to match. 
Therefore, $c=5$. The construction works for arbitrary large graphs. \\
A general surface of genus $g$ can be obtained in various ways.
One is the {\bf branched cover} construction. An other is the 
{\bf connected sum} construction:  take a large example of a Fisk graph which 
si flat except near the curve $\gamma$ producing a fisssure. 
Now cut out a $g-1$ squares with enough distance between each other and also
bounded away from the Fisk geodesic. Now glue the tori together along
the boundaries of these squares. In the non-orientable case, cut an other hole
and glue in a projective plane with chromatic number $c=3$. 
This can be done in such a way that all vertex degrees are even and 
so that the {\bf Fisk cut} connects two vertices
have odd degree. Without Fisk cut and with the right geometry so that
all fundamental cycles have length which are a multiple of $3$, then 
$c=3$. If the graph is Eulerian and one of the lengths does not match, we
get $c=4$. If we take such a graph and additionally cut along a Fisk 
cycle, we get $c=5$. 
\end{proof} 

As the conjecture strenghtheining the Albertson-Stromquist question 
indicates, we expect this to be a trichotomy for all surfaces. 
Unlike on the sphere, the 
Eulerian property does not distinguish between $c=3$ or $c=4$ in the higher genus case.
It would be nice to have
simple necessary and sufficient conditions which decide in which chromatic
class a graph is. See \cite{KnillNitishinskaya}. \\

The Fisk construction should work also in higher dimension higher dimensions. 
Lets sketch how to  get examples of $3$-dimensional
tori in $\Gcal_3$ with chromatic number $c=6$: 
tessellate the standard flat $3$ torus $R^3/Z^3$ with $10^3=1000$ solid octahedra which got
triangulation $G$ by placing a central vertex in each center. The graph is in $\Gcal_3$
and has the Betti vector $(1,3,3,1)$. 
Because the shortest geodesics has a length which is not a multiple of $4$, this 
graph can not be colored with $4$ colors but needs $5$. 
Now build a closed geodesic $\gamma$ of type 
$(1,2,0)$ and use it to cut edges as we do when constructing refinements.
The cutting produces two parallel closed lines which are parallel 
to $\gamma$ along which the degree is odd. If this construction is done with $5$ parallel 
geodesics, all odd degree vertices can be removed. We have now obtained a new refined 
three dimensional torus. Now produce a three dimensional Dehn twist 
by cutting  a two dimensional torus and translate perpendicular to the line of curves. 
After gluing back, there are 5 curves $\gamma_i$ with initial points $A_i$ and
end points $B_i$ which are pairwise neighboring. The same Fisk argument shows that we can 
not color the unit ball around $A_i$ or $B_i$ with $5$ colors as one coloring can be 
``parallel transported" along $\gamma_i$ from $A_i$ to $B_i$ producing a mismatch.  \\

So, what singles out the sphere? It is the only surface type which is
the boundary of a simply connected three dimensional manifold. It is this
topological connection which makes us so confident that the above conjecture
is true. 

\section{More remarks} 

Unlike the coloring questions for higher dimensional maps, where one needs 
arbitrary many colors in 3 or more dimensions \cite{Tietze,RobinWilson}, the coloring of geometric
graphs in $\Scal_d$ is interesting and to a large degree unexplored. 
We believe the chromatic numbers to depend very much on the topology of the graph 
and the possible embeddings. 
We have defined  the class $\Scal_d$ graph theoretically which is not a 
triangularization approach: our spheres $G$ in $\Scal_d$ are nice in the sense 
that every unit sphere in $G$ are in $\Scal_{d-1}$ which are again nice in the sense that
every unit sphere is again nice etc. 
This inductive definition gives confidence that not some pathological 
counter example exists. The intuition whether the chromatic numbers depend 
on the topology of the graph and the possible embeddings needs to be explored more. \\

Back to two-dimensions we have again to stress that the {\bf 7 color theorem} dealing 
with graph embedding of $G=K_7$ into the classical torus is a different question. 
A variant of geometric graphs have been introduced 
under the name {\bf locally Euclidean} and especially {\bf locally planar} in two dimensions.
In particular relevant is the work of Fisk 
who wrote his thesis in 1972 under the guidance of Gian-Carlo Rota. 
It seems not have made a big impact yet, presumably because Fisk theory of 1977
was published only a year after the Appel-Haken proof. 

Unlike Fisk theory, we remain in the discrete and in graph theory and avoid the notion of 
``triangularization". Unless one is a finitist, this is to a large
degree just a change of language but it avoids some difficulties as the {\bf Hauptvermutung} has
illustrated. It is important to us to have entirely graph theoretical notions everywhere. 
It is easier to implement graph theoretical methods in the computer because
these structures are hardwired into computer algebra systems and because the language is easier.
One has seen that already when switching from {\bf map colorings} to 
{\bf graph colorings}, which must have happened primarily at the time of Birkhoff. 
That transition from topology to discrete mathematics also illustrates the connections between
the fields which present since the K\"onigsberg Bridge problem which was the seed for graph theory. \\

Discrete {\bf projective planes} are especially interesting from a chromatic perspective: 
there are examples in $\Ccal_3$ or $\Ccal_4$ as well as examples (already found by Fisk), not in 
$\Ccal_4$. It could be related to the fact that the projective plane is a nontrivial element in the 
cobordism group defined by $\Gcal_2$.  Orientable surfaces $G$ on the other hand are all 
{\bf 0-cobordant} but the filling $H$ satisfying $G=\delta H$ is not simply connected. 
So, if orientability and simply connectedness both are obstructions to $\Ccal_4$ what about
the Klein bottle? It is neither orientable, nor simply connected. But it is 0-cobordant as the
boundary of the {\bf solid Klein bottle}. It is reasonable therefore to suggest that any
Klein bottle in $\Gcal_2$ is also in $\Ccal_5$. The torus appears also naturally as $f^{-1}(x)$ in a 
{\bf Hopf fibration} $f:S^3 \to S^2$, a fact which could be used to color it naturally.  \\

Let $\Ecal$ denote the class of {\bf Eulerian graphs}, graphs which admit Eulerian circuits, 
closed paths hitting every edge exactly once. 
By the Euler-Hierholzer theorem they agree with graphs for which all vertices have even degree. 
Kempe discovered and Heawood proved that $\Pcal \cap \Ecal \subset \Ccal_3$. See also
\cite{Fisk1977a}. The corresponding question for not necessarily simply connected 
graphs in $\Gcal_2$ is more subtle, as global properties start to matter. 
There are discrete toral graphs like the flat $T_{6,6}$ which are in $\Ccal_3$ but others 
like $T_{5,5}$ which are not. 
It is not only the minimal lengths of geodesics on the torus which matter. 
The flat toral graphs like $T_{4,4}$ have chromatic number $4$ while other 
implementations with 16 vertices and nowhere vanishing curvature 
(obtained by diagonal flips from $T_{4,4}$) have chromatic number $3$. Non-simply 
connected graphs feature holonomy constraints for coloring which do not occur in 
the simply connected case. \\

Given a closed simple circuit $C \subset G$, where $G \in \Scal_2$. If $C$ has no self intersections
and length $\geq 4$, then by the Jordan curve theorem, it divides $G$ into 
two connectivity components $C_+,C_-$, from which we can assume that $C_+$ is nonempty. 
The graph so constructed from the union of $C$ and $C_+$ is called a {\bf disc}. If $G$ is a disc
and $x$ is a vertex in $G$, then its unit sphere $S(x)$ is either a circular graph $C_n$
with $n \geq 4$ or an interval graph $I_n$ with $n \geq 2$. The former are called 
{\bf interior points} of $G$, the later are called {\bf boundary points} of $G$. 
Let $\Bcal_2$ denote the set of all discs.  If $G \in \Bcal_2$, then the boundary is a 
single circular graph $C_n$ with $n \geq 4$.  
On the other hand, removing a single vertex from $G \in \Scal_2$ produces a 
graph in $\Bcal_2$. We will below give a new definition of $\Scal_2$ and $\Bcal_2$ which generalizes
to higher dimensions. 
The theory of reducibility and minimal configurations shows that $\Bcal_2 \subset \Ccal_4$ is equivalent
to the four color problem. We still would like to 
see directly that $\Bcal_2 \subset \Ccal_4$ is equivalent to $\Scal_2 \subset \Ccal_4$. 
One direction is clear: if $\Scal_2 \subset \Ccal_4$ and $G \in \Bcal_2$ is given, we can complete it
to a sphere, color it and get a coloring of $G$.
Assume on the other hand that $\Bcal_2 \subset \Ccal_4$. Can we conclude $\Scal_2 \subset \Ccal_4$?
Given $G \in \Scal$. Removing an edge $e=(a,b)$ produces an element in $\Bcal_2$ which can be colored.
If there is no Kempe chain from $a$ to $b$, we can recolor one
of the colors and reconnect. Otherwise, by planarity, there is no Kempe chain with the other two vertices.
This classical Kempe argument shows that if there is an example in $\Scal_2$ which
is not 4 colorable, a diagonal flip renders it 4 colorable. A similar Kempe
argument shows that a counter example in $\Scal_2$ has no degree $4$ vertices
as removing such a vertex allows after a Kempe recoloring to fill in the vertex again.
As flipping a diagonal in a wheel graph $W_5$ which by Gauss-Bonnet is present in $\Scal_2$ 
produces a wheel graph $W_4$, there is an other reason why flipping diagonals can enable colorings.
Cayley knew already that a ``minimal criminal" (= irreducible) counter example to the four color theorem 
is a triangularization.
The 4 color theorem implies that the boundary of $G \in \Bcal$ of an irreducible graph is in $\Ccal_3$.
Disks with 4 colors therefore satisfy a {\bf maximum principle}: there are colorings so that the
{\bf maximal chromatic number of a closed circuit} happens in the interior. 
This maximum principle has been realized by Tutte who called an example of a disc $G$ for
which the boundary has chromatic number $4$ a {\bf chromatic obstacle}. 

\section{Cobordisms}

A graph is called {\bf 3-dimensional geometric}, if every unit sphere is in $\Scal_2$. 
Such graphs play the role of $3$-dimensional manifolds. If $t$ is the number of tetrahedra,
then the Euler characteristic of such a graph is $\chi(G) = v-e+f-t$. 
A graph for which every unit sphere
is either in $\Scal_2$ or $\Bcal_2$ is called a {\bf $3$-dimensional geometric graph with boundary}.
The Euler characteristic of $G \in \Bcal_3$ is defined as $\chi(G) = v-e+f-t$, where 
$v$ is the number of vertices $K_1$ , $e$ is the number of edges $K_2$, $f$ is the number of triangles $K_3$ 
and $t$ the number of tetrahedra $K_4$ in $G$.
Euler-Poincar\'e assures that $\chi(G) = b_0-b_1+b_2-b_3$. For connected, orientable case $G \in \Bcal_3$, 
one has $b_0=b_3=1$. We will see below that $b_1=b_2$ because $\chi(G)=0$. \\

The {\bf curvature} of a vertex $x \in G$ with $G \in \Gcal_2$ is defined as $K(x) = 1-V_0/2+V_1/3$,
where $V_i$ are the number of $K_{i+1}$ subgraphs in $S(x)$. As $S(x)$ is a discrete circle 
with Euler characteristic zero, we have $V_0=V_1={\rm deg}(x)$ and $K(x) = 1-{\rm deg}(x)/6$. 
From the {\bf Euler handshaking lemma} $\sum_x V_0(x)=2v_1$ follows 
immediately the {\bf Gauss-Bonnet theorem}
$\sum_{x \in V} K(x) = \chi(G)$ for surfaces $G \in \Gcal_2$. In particular, The Euler formula
$\sum_{x \in V} K(x) = 2$ holds for $G \in \Scal_2$.  \\

For a graph $G \in \Gcal_3$  the {\bf curvature} of a vertex is defined as 
$$  K(x) = 1 - \frac{V_0}{2} + \frac{V_1}{3} - \frac{V_2}{4}  \; . $$

\resultlemma{ Every interior vertex of $G \in \Gcal_3$ has zero curvature.  }
\begin{proof}
Since $S(x)$ is in $\Scal$, we have $V_2(x) - V_1(x) + V_0(x) = 2$.
The handshaking result $3 V_2(x) = 2 V_1(x)$ holds because two triangles
meet and "handshake" at an edge. Three times the number of triangles double
counts the number of edges in $S(x)$. Expressing now $K(x)$ in terms of one 
variable establishes it.
\end{proof}

With the assumption $V_{-1}=1$, the curvature of a vertex in a finite
simple graph is defined in \cite{cherngaussbonnet} as
$$  K(x) = \sum_{k=0}^{\infty} (-1)^k \frac{V_{k-1}(x)}{k+1}  \; $$

\resulttheorem{
$\sum_{x \in V} K(x) = \chi(G)$
}
\begin{proof}
The handshaking lemma generalizes to $\sum_{x \in V} V_k(x)=v_{k+1} (k+2)$
as any of the $v_k$ subgraphs $K_{k+2}$ contribute a count to each of its 
$(k+2)$ vertices. The Euler handshaking lemma is in the special case $k=0$
sometimes called the "first theorem" of graph theory.
\end{proof}

Sometimes the handshaking result is called the "fundamental theorem of graph
theory". We have just seen that it is essentially equivalent to Gauss-Bonnet.  \\

Curvature has an other relation with colorings.
the expectation of the index $i_f(x)$ defined in \cite{poincarehopf} has first
shown in \cite{indexexpectation} to be curvature.  In 
\cite{colorcurvature} we remarked that one can take the expectation also over
the space of all colorings. 

\resultcorollary{
If $G$ in $\Gcal_3$ has no boundary, then $\chi(G)=0$.
}

We have shown in full generality that $K(x)$ 
is constant zero for every vertex $x$ in an odd-dimensional geometric graph $G$ \cite{indexformula}.
We established this result using 
{\bf integral geometry} by writing curvature as the expectation of indices $i_f(x)$ satisfying
the Poincar\'e-Hopf theorem 
$$  \sum_{x \in V} i_f(x) = \chi(G) $$ 
which after taking expectation gives Gauss-Bonnet. \\

Two graphs $G_1,G_2 \in \Gcal_2$ are called {\bf cobordant}, if there is graph $H \in \Gcal_3$
such that the disjoint union $G_1 \cup G_2$ agrees with the boundary of $H$. 
We worked with this in the summer 2014 and learned in December 2014 about the work of Fisk 
\cite{Fisk1980} who works with triangularizations and ``homotopies" in which case 
$H$ is a triangularization of $G \times [0,1]$. 

\resultlemma{
Cobordism is an equivalence relation in $\Scal_2$.
}

\begin{proof}
Transitivity and symmetry are clear. For reflexivity, we have to construct from $G=(V,E)$ 
a graph $H \in \Gcal_3$ which has the disjoint union $G \cup G$ as a boundary.  
To do so, place the completion $\overline{G}$ of the {\bf dual graph} 
$\hat{G}=(\hat{V},\hat{E})$ between. The vertices of
$H$ is the union of $V_i$ and $\hat{V}$. The edges of $H$ is the union of $E_i$ and $\hat{E}$
as well as connections $(v,w)$ where $v$ is a vertex in $G \cup G$ and $w$ is 
a vertex in $\hat{G}$. 
\end{proof}

Theorem 77 in \cite{Fisk1980} shows this for $\Gcal_2$ under some holonomy assumption.
As in the continuum, the cobordism classes form a group where addition is the disjoint union.

\resulttheorem{
The cobordism group of $\Gcal_2$ is nontrivial. 
}

Actually, since graphs in $\Gcal_d$ define in a 
natural way topological manifolds with the same topological properties, the cobordism groups are
the same as in the continuum. It is $Z_2$ for $\Gcal_2$. 

\resultlemma{
The projective plane is not zero cobordant.
}

This follows from the fact that $\chi(G)=1$ for the projective plane and because:

\resultlemma{
If $H \in \Gcal_3$ has the boundary $G$, then $\chi(G)$ is even.  }

\begin{proof}
Glue $G \cup G$ along
$\delta H$ to get a graph $H' \in \Gcal_2$ without boundary.
The later has Euler characteristic $0$.
By the intersection property of $\chi$ we have
$2 \chi(H) - \chi(G) = 0$ showing that $\chi(G)$ is even. 
\end{proof}

This classically known result was one of the triggers to look at cobordisms for coloring
after having seen projective planes in $\Gcal_2$ with chromatic number $5$.  \\

Both homotopy and cobordisms are rough equivalence relations for graphs: 
cobordisms does preserve dimension but not Euler characteristic $\chi$ nor cohomology;
homotopy does not preserve dimension but leaves $\chi$ and cohomology invariant.
We have introduced in \cite{KnillTopology} a topological equivalence relation for graphs
based on dimension and homotopy which captures topological features. 
It is based on {\bf \v{C}ech constructions},  preserves $\chi$, homotopy classes, 
cohomology and cobordism classes and is therefore a natural discrete analogue of what 
homeomorphisms are in the continuum. It is based on {\bf inductive dimension} for graphs 
which was introduced in \cite{elemente11} inspired from the Menger-Uryson dimension 
and studied for example in \cite{randomgraph}. See also 
\cite{KnillWolframDemo1,KnillWolframDemo2} for example code. 
It uses {\bf homotopy} which is inspired from a definition given by 
Whitehead \cite{Whitehead} and introduced by Ivashchenko \cite{I94} to graph theory. 

\section{Minimal colorings}

The {\bf edge degree} ${\rm deg}(e)$ of an edge $e=(a,b)$ of a graph $G \in \Gcal_3$ 
is defined as the number of 
elements in the circular graph $S(a) \cap S(b)$. It is the number 
of tetrahedra which hinge at $e$. We sometimes just call it {\bf degree}. The edge degree
is related to {\bf Ricci curvature} in the continuum because as the later, it is a quantity 
defined by a direction. The relation with curvature is also clear as it is related to the
scalar curvature on the $2$-dimensional unit sphere. The sum of the Ricci curvatures 
of all edges at a vertex $x$ is a type of {\bf scalar curvature} for $3$-dimensional graphs. 
While there is no Euler curvature for $3$-dimensional manifolds (the Euler curvature involves
a Pfaffian only defined in even dimensions (see e.g. \cite{Cycon}), there is an Euler curvature in odd dimensions
but it is zero everywhere \cite{indexformula}. \\

For a graph $G \in \Gcal_3$ with boundary,
an edge is called {\bf interior} if $S(a) \cap S(b)$ is a circular graph. For boundary edges
the graph $S(a) \cap S(b)$ is an interval graph. It can be possible that
two boundary points $a,b$ define an edge $(a,b)$ which is interior. \\

The {\bf fundamental group} $\pi_1(G)$ of a finite simple graph $G=(V,E)$ is defined
as  the set of equivalence classes of closed paths 
$$ \gamma = \{x_0,x_1,\dots ,x_n=x_0 \; \} $$ 
in $G$, where $a_i \in V$ and $(a_i,a_{i+1}) \in E$.
In the literature, graphs often are considered $1$-dimensional simplicial complexes so that
in such a setup only trees are simply connected. With the just given definition, a nice triangularization
of a manifold $M$ has the same fundamental group as $M$.
A deformation step of a path within $G$ consist of one of the transformations
$S,S^{-1},T,T^{-1}$, where $S$ shortens the leg $(abc)$ in a triangle with the
third edge $(ac)$ and $T$ removes a backtracking step $(aba)$ to $(a)$ or $(aa)$ to $(a)$.
$G$ is called simply connected if $\pi_1(G)$ is trivial. A more natural definition of curve
homotopy can be given by looking a graph representing the embedding $C_n \to G$. Two curves are
homotopic if the corresponding embedding graphs are homotopic. \\

Given $G \in \Gcal_3$, the {\bf dual graph} $\hat{G}$ is the graph whose
vertices are the tetrahedra in $G$ and where two tetrahedra are connected, if they
share a triangle. The dual graph can be defined in the union of all graphs $\Gcal_n$
where vertices are the maximal simplices in $G$. 
It is always triangle free and contains less information than $G$:
the circular graph $C_4$ and the wheel graph $W_4$ have the same dual graph.
Assume for simplicity that $G$ is connected. Then also $\hat{G}$ is connected and
$\chi(\hat{G}) = \hat{v}-\hat{e} = 1-b_1$ where $b_1={\rm dim}(H^1(\hat{G}))$ is the
{\bf genus} of $\hat{G}$. As $\hat{G}$ is $1$-dimensional, the genus 
is the number of shortest closed loops which generate the fundamental group of $\hat{G}$. 
For an octahedron $G$ for example, where $\hat{G}$ is the cube, the genus of $\hat{G}$ is 
$b_1=g=5$ as $\chi(\hat{G})=8-12=-4$ for the cube because of the lack of triangles.
Indeed, if $\gamma_1,\dots ,\gamma_5$ are counter clockwise oriented loops around the faces,
then $\gamma_1+ \cdots \gamma_5 = \gamma_6$ is a closed loop around the 6th face. \\

The cube is not geometric as all unit spheres are three point
graphs $P_3$ of dimension $0$. We can complete the cube however by capping its 6
faces and embed it in a geometric graph $\overline{G} \in \Scal_2$ which has $v=8+6=14$ vertices, $e=12+24=36$ edges
and $f=4 \cdot 6=24$ triangles. Now $\chi(G) =v=e+f=2$ is the Euler formula and $\overline{G}$ is 
simply connected, satisfies $\chi(G)=b_0-b_1+b_2=2-b_1 = 2$ as $b_1=0$ by Hurewicz.
Most graph theory literature treats graphs as $1$-dimensional simplicial complexes.
A more geometrically view leads to a closer relation to the continuum (see \cite{KnillBaltimore} for an overview). 
An illustration is the {\bf Brouwer fixed point theorem} \cite{brouwergraph}:

\resulttheorem{
For any graph endomorphism $T$ of $G \in \Bcal_d$ there exists a complete subgraph of $G$ which is fixed by $T$.}

An other consequence is sphere automorphism result like

\resulttheorem{
A orientation preserving graph endomorphism $T$ of $G \in \Scal_2$ or
an orientation reversing graph endomorphism $T$ of $G \in \Scal_3$ has two complete subgraphs which are fixed by $T$.}

The reason for the second statement is that the {\bf Lefshetz number} of an orientation reversing automorphism is $2$.
Also duality holds in general for geometric graphs but then the dual graph has first to be completed to become geometric.
As cross polytopes $G$, the smallest graphs in $\Scal_d$ illustrates, the dual graph $\hat{G}$, a hypercube is not geometric
as it has never triangles. But there is a way to fill in $d$ dimensional fillings in each chamber to get a geometric graph.
We have used such completions to explore the hypothesis that for odd dimensional graphs, the Euler curvature identically
vanishes. For orientable geometric graphs, there is a {\bf Poincar\'e duality} relating the Betti numbers of the dual graph
with the Betti numbers of the graph itself. This is essentially what Whitney did when inventing the concept of manifold. 
In two and three dimensions, $\hat{G}$ at least allows to reconstruct the data $v_i(G)$: 

\resultlemma{
The graph $\hat{G}$ determines all numbers $v_i(G)$ for $G \in \Gcal_3$. 
}
\begin{proof}
Use $\hat{G}$ as a blueprint to stack together tetrahedra to construct
examples of graphs which have $\hat{G}$ as dual graphs.
The genus of $\hat{G}$ is the number of edges of $G$. 
The number $\hat{v}$ of vertices of $\hat{G}$ is the number of tetrahedra in $G$. 
The number $\hat{e}$ of edges of $\hat{G}$ is the number of triangles in $G$. 
Because $\chi(G)=0$ the number of vertices of $G$ are determined. 
\end{proof}

By Poincar\'e duality, which implies $b_0-b_1+b_2-b_3=0$, we know that $b_1(G)=b_2(G)$ if $G$ is 
connected and orientable. The example of homology spheres shows that $b_1(G)=0$ does not force $G$ to be a 
3-sphere. But we can ask what is the class of graphs $G \in \Gcal_3$ which have the same dual 
graph $\hat{G}$. If $\pi_1(G)$ is trivial for $G \in \Gcal_3$ then the Perelmann-Poincar\'e theorem 
assures that $G$ is in $\Scal_3$. P-P implies the corresponding
discrete version $|\pi_1(G)=1|  \Rightarrow G \in \Scal_3$ if $G \in \Gcal_3$ as $G$ 
naturally defines a 3-manifold. A discrete P-P theorem could be easy to prove by 
{\bf edge contraction}: take a union of two balls $B(x) \cup B(y)$ with $e=(x,y) \in E$,
pull together the two vertices to one vertex $z$. This reduces the number of vertices. 
Now continue this process while we have geometric graphs. 
Such a reduction stops if every $B(x) \cup B(y)$ contains two boundary points
$u,v$ for which $(u,v) \in E$. But then we have a cross polytop, the smallest
element in $\Scal_3$. As in the continuum with the Ricci flow, one presumably has to deal
with obstacles where reductions are no more possible because a unit sphere $S(x)$ 
stops to be in $\Scal_2$, but regularizing the discrete flow in the discrete is expected to
be much easier than in the continuum. Back to the question whether $\hat{G}$ determine the graph $G$: 

\question{
If $G_1,G_2 \in \Scal_2$ have isomorphic $\hat{G}_1,\hat{G}_2$. Is $G_1$ isomorphic to $G_2$?
}

The following result generalizes the Kempe-Heawood minimal coloring result 
to three dimensions:

\resulttheorem{
If $G \in \Gcal_3$ is simply connected and every interior edge has
even degree then $G \in \Ccal_4$. 
}

\begin{proof}
Start with a tetrahedron in $G$ and color it with one of the $24=4!$ possible colorings. 
For any neighboring tetrahedron, the colors are now determined, as the interface
with the first tetrahedron fixes three colors. Now continue coloring. 
An obstruction would manifest from a small circle around an edge as the
dual graph of $G$ is a triangle-free graph whose fundamental group is generated
by circular graphs belonging to chains of tetrahedra hinging at one edge. 
\end{proof} 

Kempe mentioned the $2$-dimensional version without proof and Heawood was the first to give the argument.
Kempe therefore should be credited for its discovery and Heawood for the proof \cite{Barnette}.
It is a rather general argument that in the simply connected case, one can
minimally color $G$ with a minimal number of colors 
if the dual graph $\hat{G}$ is bipartite, in arbitrary dimensions. \\

The Perelman-Poincar\'e theorem assures that simply connected $G \in \Gcal_3$ is in $\Scal_3$. One could
therefore classify $\Scal_3$ as the simply connected graphs in $\Gcal_3$. One the same note,
one can define $\Scal_3$ graph theoretically as the class of graphs which after removing a vertex
are in $\Bcal_3$ and $\Bcal_3$ as the class of $3$-dimensional graphs with boundary which are
contractible. One could also define the class $\Scal_3$ Morse-Reeb theoretically as graphs 
admitting functions with two critical points and none with one, using homotopy
as a subset of graphs with Lusternik-Schnirelman  \cite{josellisknill} category 2 or using an analogue of 
Zippin's result as graphs in $\Scal_2$ satisfy the Jordan Curve theorem. 
In the non-simply connected case, the situation is more interesting and difficult. 
We understand the $2$-dimensional case. In the $3$-dimensional case, torsion effects matter in the 
sense that not only the parity of the lengths of minimal cycles in the fundamental group matter. 
This later possibility was pointed out by Jenny Nitishinskaya.
Torsion effects could also be responsible for the existence of non-isometric graphs in $\Scal_3$ with 
the same dual graph, provided that they exist. \\

Let $G \in \Gcal_3$ be a connected graph which can be minimally colored with $4$ colors. 
Fixing the color on a tetrahedron now determines the color of all other tetrahedra:
just follow paths on the dual graph $\hat{G}$ to color along chains of simplices. 
There are two type of paths: some 
of the closed paths in $\hat{G}$ belong to the set of tetrahedra attached to an edge
$e=(a,b)$. If $f$ is the colors on a,b are fixed, then on the circle
$S(a) \cap S(b)$ the other two colors have to be used which reflects that 
the circle has even length ${\rm deg}(e)$. If $G$ has a nontrivial fundamental group, then
chains $C: t_0,t_1, \dots ,t_n=t_0$ of tetrahedra are possible whose union generate subgraphs which are
not contractible. Let $S_4$ denote the group of permutations with 4 elements. The coloring
along the chain $C$ now induces a holonomy permutation in that permutation group 
when going along the chain. Only if all the holonomy conditions are satisfied, the case 
$G \in \Ccal_4$ is possible. Already in $\Gcal_2$, not only the length of the closed path 
matter in the non-simply connected case \cite{KnillNitishinskaya}.

\section{Coloring the boundary}

If a minimal colorings with $d+1$ colors in $\Gcal_d$ is possible, it is straightforward to construct the
coloring. Especially, for graphs $G \in \Ccal_4 \cap \Gcal_3$, the coloring procedure is trivial if it exists
because the color of one tetrahedra determines the coloring of all its neighbors. \\

The idea is now to color a graph $G$ in $\Scal_2$ by realizing it as the boundary of 
a graph $H \in \Bcal_3 \in \Ccal_4$ or alternatively to build a cobordism $G \to G'$ 
like a self-cobordism with a filling graph $H \in \Bcal_3 \cap \Ccal_4$. In both cases, now modify 
the interior to make the graph and so its boundary colorable with 4 colors. \\

Let $O(G) = \{ e \in E \; | \; {\rm deg}(e) \; {\rm odd} \; \}$ 
denote the set of {\bf odd degree interior edges} of a graph $G=(V,E) \in \Bcal_3$. 
The functional 
$$  \phi(G) = \sum_{e \in O} {\rm dist}(e,\delta(G))  $$
counts the sum of all geodesic distances of odd degree edges to the boundary,
where ${\rm dist}(e,A)$ for an edge $e=(a,b)$ is the maximum of the geodesic distances
${\rm dist}(a,A)$ and ${\rm dist}(b,A)$. \\

The {\bf cobordism coloring algorithm} for $G$ chops up the graph $H$ for a cobordism $H(G,G')$ without 
changing the boundary $G$. One possibility is a $0$-cobordism $H(G,\emptyset)$. 
The chopping aims to get rid of all interior odd edges leading to a minimal coloring with $4$ colors. 
One refinement step is to cut interior edges $e=(a,b)$ and fill in a wheel graph with $S(a) \cap S(b)$ as boundary. 
Given $e=(a,b) \in O$, divide it with a new vertex $v$ and add edges $(a,v),(b,v)$ 
and $(v,x_i)$ where $x_i$ are the vertices of $S(a) \cap S(b)$. Call this an {\bf edge cut}
and call a sequence of edge cuts a {\bf surface cut}. We will justify the name below. 
Additionally, we can refine simplicies also by filling in 16-cells or 
completed 3D-cubes. All these refinements do not change the interior degrees. \\

Cobordisms have been introduced by Poincar\'e. While Fisk \cite{Fisk1980} introduces
cobordism in a chromatic setup, the definition relies on Euclidean structures as all
graphs are assumed to be triangularizations.  He also looked at the language of chains:
a {\bf surface} $S$ in a finite simple graph $G$ is a 
set of triangles $t_j$ in $G$ together with 
an orientation having the property that adjacent boundaries cancel. 
A {\bf curve} $C$ in a graph $G$ is a set of edges $e_j$ in $G$ together with 
an orientation such that adjacent boundaries cancel. 
Since we are interested only in the parity of the degree of an edge, we use the group 
$Z_2$ and build finite formal sum $S=\sum_i a_i t_i$ with $a_i \in \Z_2$. 
These {\bf 2-chains} are finite formal sum $C=\sum_i a_i e_i$ with $a_i \in \Z_2$ is a {\bf 1-chain}. 
The boundary of a 2-chain is $\delta S = \sum_i a_i \delta t_i$, 
where $\delta (x_0,x_1,x_2) = (x_1,x_2) - (x_0,x_2) + (x_0,x_1)$ modulo $2$.
The boundary of a $1$-chain is $\delta C = \sum_i a_i \delta e_j$, where
$\delta (x_0,x_1) = (x_1)-(x_0)$. Since the boundary of a boundary is zero,
the {\bf first homology group} $H_1(G,\Z_2)$ as the group of 1-chains, the kernel of $\delta$
modulo the group of boundaries of $2$-chains.
A {\bf surface} is a 2-chain $S$ with $|a_j|=1$ for which $\delta S$ is a curve. 
A {\bf curve} is a 1-chain $S$ with $|a_j|=1$ for which $\delta S$ is a discrete 
point set. A {\bf discrete point set} is a $0$-chain $P$ for which $|a_j|=1$. 
As the fundamental group is trivial, the first homologies are trivial so that
every closed curve is the boundary of a $2$-chain. 
Because $O$ is closed, it is a union of closed curves $C_j$. Each 
is the boundary of a surface $S_j$. 
Write $S=\sum_j S_j$ modulo $2$ as $\sum_j a_j t_j$. Its boundary modulo 2 is $O$. 
The equation $\delta S = O$ without projection onto $Z_2$  defines an orientation of the $t_j$.
It is this {\bf homological triviality} which adds
confidence that one can get rid of $O$ by refinements. The case of higher genus surfaces, where
cohomology starts to matter indicates a strong topological connection. \\

With Kirchhoff we can interpret $O$ as a {\bf current} $j$ flowing through the network, 
An other analogy is a  {\bf spin network} with the group $Z_2$. 
{\bf Conservation of current} implies that the sum of all currents at a node is zero. This can
be rewritten as ${\rm div}(\vec{j})= d_0^* j=0$. In physics, this is called a {\bf Coulomb gauge}.
{\bf Electromagnetism} on a graph reduces to linear algebra: take a $1$-form $A$, define the electromagnetic 
field $F=dA$. Now,$dF=0$ and $d^* F=j$ can in a Coulomb gauge for $j$ be written as the Poisson equation 
$LA=d^* dA=j$ which leads to the electromagnetic potential $A=L^{-1} j$. If $G \in \Gcal_4$ (a space-time
manifold), the values of $F$ on the 6 faces of $K_5$ subgraphs are the electromagnetic field 
$(E_1,E_2,E_3,B_1,B_2,B_3)$. But we don't need a semi-Riemannian structure at all as
the electromagnetic formalism makes sense for any finite simple graph
as do the {\bf Maxwell equations} $dF=0,d^*F=j$ which are just geometry. \\

The four color theorem would follow from the following: 

\conjecture{
Given $G \in \Bcal_3$ with $\phi(G)>0$, there is a surface-cut which lowers the value of $\phi$. 
}

We have implemented this in the form of {\bf simulated annealing} or by brute force going through
edge cuttings and see whether $\phi$ can be lowered and doing so if possible. 
What is the effect of an edge cut? 

\resultlemma{
An edge cut at $e=(a,b) \in O$ increases the degrees at interior edges in $S(a) \cap S(b)$ by one.
}
\begin{proof}
The degree of the two parts $(a,v),(b,v)$ is the same as the degree of $e$. 
The degrees of the new vertices $(v,x_i)$ are all $4$ and so even.
The degrees of edges in $S(a) \cap S(b)$ gets increased by $1$. 
\end{proof}

We believe that we can chop up the interior because the odd degree edges can be simplified with 
$1$-dimensional curves $S(a) \cap S(b)$ for edges $(a,b)$. The reason is:

\resultlemma{
Given an interior vertex $v$. The sum of the edge degrees over all edges
attached to $v$ is even.  }

\begin{proof} 
The edge degree of $e = (v,w)$ is the vertex degree of $w$ in the unit sphere $S(x) \in \Scal_2$.
By the Euler handshaking lemma applied to this sphere, this is twice the number of edges in $S(x)$. 
\end{proof} 

\resultlemma{
If $H$ is simply connected, then the edges in $O$ can be oriented
to be the boundary of a surface $S$.  }
\begin{proof}
The degree condition assures that $O$ is a cycle. Since the homology group $H_1(G,Z_2)$
is trivial for $H$ as $H_1(G,\Z)$ is trivial.  So, $O$ is a boundary $O = \delta S$.
\end{proof}

To reduce the functional $\phi$, chose a closed curve $\gamma$, fill it with a 
$2$-dimensional surface in such a way that the boundary of this surface consist of 
$\gamma$ and boundary parts and such that the surface can be partitioned into wheel
graphs after splitting. If this is possible, the parity of $\gamma$ has changed
and the curve $\gamma$ is no more part of $O$. Building the surface might need some
additional cutting first. Continue splitting surfaces like that
until no element in $O$ is left. As experiments show, sometimes, we have to make
some preparation cuts or go slowly, moving the odd edges closer to the boundary. 
In our experience, its good to cut from the interior and move outwards because if the 
outside is chopped up early we can not see inside and need more work to dig us out.
By induction, we can clean out the degrees up to a neighorhood of the newly added vertex. 
But this only works if the outer boundary of $G$ is small. In general it is known
that an addition of a vertex has global consequences of the colors elsewhere.
This also manifests here. Changing the graph a little can produce completely different
tetrahedral partition in the interior which will be used to compute the coloring. 

\section{Higher dimensional spheres}

Here is a more formal approach to the notion of sphere and geometric graph. 
As before, we often identify $W \subset V$ with the maximal subgraph of $G$ containing $W$. 
The unit sphere $S(x)$ of a vertex $x$ for example is the graph generated by all vertices connected to $x$. Graph theorists 
sometimes call the unit sphere $S(x)$ the {\bf link} of $x$ and the unit ball $B(x) = S(x) \cup \{x\}$ the {\bf star} or
{\bf first neighborhood} of $x$.
The one point graph $K_1$ is contractible. Inductively, $G=(V,E)$ is {\bf contractible}, if there is $x \in V$ such that both
the unit sphere $S(x)$ and $V \setminus \{x \;\}$ are contractible. Examples are trees, complete graphs $K_n$ or wheel graphs.
All unit balls $B(x) = \{ y \; | \; d(x,y) \leq 1 \; \}$ in $G$ are contractible.
Circular graphs $C_n, n \geq 4$, the octahedron or icosahedron are examples of non-contractible graphs.
The Euler characteristic $\chi(G)$ is the super count $\sum_{k=0}^{\infty} (-1)^k v_k$
of the number $v_k$ of $K_{k+1}$ subgraphs of $G$.
Given subgraphs $A,B$ of $G$, the counting formula $\chi(A \cup B) = \chi(A) + \chi(B) - \chi(A \cap B)$ holds.
By induction, any contractible graph $G$ has Euler characteristic $\chi(G)=1$ because
$\chi(G) = \chi(G \setminus \{x\}) + \chi(B(x)) - \chi(S(x))$
and $\chi(B(x))=\chi(S(x))$ by the assumption that $B(x) \setminus \{x\} = S(x)$ is contractible.
Also the definition of {\bf geometric graph}, {\bf sphere} and {\bf ball} is inductive:
let $\Gcal_0$ denote the set of graphs without edges,
set $\Scal_0 \subset \Gcal_0$ to be the set of graphs with two vertices and no edges and denote
by $\Bcal_0 \subset \Gcal_0$ the set of graphs with one vertex.
Define inductively the set of {\bf geometric graphs} $\Gcal_d = \{ G \; | \;  S(x) \in \Scal_{d-1}$
or $S(x) \in \Bcal_{d-1} \; \forall x \in V \; \}$. For $G \in \Gcal_d$, the {\bf boundary} is
$\delta G = \{ x \; | \; S(x) \in \Bcal_{d-1} \; \}$
and $G \setminus \delta G$ is the {\bf interior} which we ask to be nonempty. The set of {\bf balls}
$\Bcal_d \subset \Gcal_d$ is defined as the set of contractible graphs in $\Gcal_d$ for which the
boundary is in $\Scal_{d-1}$. The set of {\bf spheres} $\Scal_d \subset \Gcal_d$ is the set of
non-contractible graphs for which removing an arbitrary vertex produces a graph in $\Bcal_d$.
Removing a vertex for which $S(x)$ is contractible has an inverse operation: add a vertex and connect
it to a contractible subgraph. These two operations generate {\bf homotopy}. Homotopy is an equivalence relation but
contractibility is not: there are contractible $G$ which can be contracted to the dunce hat $H$ while
$H$ is not contractible. From the inclusion-exclusion formula follows that $\chi$ is preserved by homotopy.
By induction, spheres $\Scal_d$ have Euler characteristic $2+2 (-1)^d$. They are never contractible.
Every graph in $\Bcal_d$ is contractible and has Euler characteristic $1$.
Any $d$-dimensional geometric graph $G$ naturally defines a topological manifold $M$
for which $G$ is a triangularization: just realize the discrete unit balls $B(x)$ as marks to be filled with 
Euclidean space and use them as charts.
We have $\chi(G)=\chi(M)$ and if $G_i$ are homotopic, their manifolds $M_i$ are classically homotopic.
The {\bf chromatic number} of a graph is the maximal integer $c$ for which a locally injective $Z_c$-valued
function exists on $V$. If $G$ can be colored by $c$ or less colors, it is called {\bf $c$-colorable}.
Let $C(G)$ be the set of colorings and $f \in C(G)$, the {\bf Poincar\'e-Hopf}
formula $\sum_{v \in V} i_f(x) = \chi(G)$ \cite{poincarehopf} holds, where $i_f(x) = 1-\chi(S^-(x))$ and $S^-(x)$ is part
of the sphere where $f(y)<f(x)$ assuming the usual ordering on $Z_c$. This follows by induction.
If we put the uniform measure on $C(G)$, then the expectation $E[i_f(x)]$ is the curvature $K(x)$ of a vertex $x$
\cite{colorcurvature}.
This is an integral geometric derivation of the
{\bf Gauss-Bonnet formula} $\sum_{x \in V} K(x) = \chi(X)$. A graph homomorphism $h:H \to G$ is called an 
{\bf embedding} of $H$ to $G$. 
The {\bf embedding graph} $G_H$ has as vertices the disjoint union of vertices of $H$ and $G$.
The edges of $G_H$ are the union of the edges of $H$ and $G$ together with all pairs $(x,h(x))$. Two embeddings are
called {\bf homotopic}, if their embedding graphs are homotopic. An embedding of $C_n, n \geq 4$
is called a {\bf closed loop} in $G$. The set of homotopy classes of all closed loops in $G$ is called the {\bf fundamental
group} $\pi_1(G)$ of $G$ as it carries a group operation: identify $C_n,C_m$ along a single vertex $v$ and blow this
up to a $0$-dimensional sphere $S_0 = \{a,b\}$. After defining $h(a)=h(b)=h(v)$, this leads to an embedding
of $C_{n+m+1}$. The fundamental group of $G$ agrees with the fundamental group of its topological filling
$M$. One can define higher homotopy groups $\pi_k(G)$ in the same way. As in the continuum, they are Abelian
for $k \geq 2$ and agree with the corresponding groups of their topological manifold realizations $M$. \\

Given a graph $G$ and a coloring $f$, we can look at {\bf level surfaces} $\{ f=c \}$. 
If $G$ is geometric, then a completion of $\{ f=c \}$ is geometric. 
Since we are interested in $\Scal_2$, lets look at this:

\resultlemma{
Given $G \in \Scal_2$, and $f$ a coloring, then the graph $H=\{ f=c \}$ is in $\Gcal_1$, a finite set of 
circular graphs. }

\begin{proof}
If $G=(V,E)$, denote $H=(W,F)$ the level curve. Let $T$ be the set of triangles in $G$.
Given $t=(a,b,c) \in T$, lets look at an edge where $f$ changes from $f<c$ to $f>c$. Such an edge
is a vertex in $H$. As $f(a),f(b),f(c)$ are different, there is exactly one other edge which belongs
to the vertex set of $H$. The triangle $t$ now represents the connection between these two vertices.
As every edge in $G$ has exactly two neighboring triangles, every vertex in $H$ has exactly two 
neighbors. The graph $H$ contains no triangles as otherwise, we had a tetrahedron in $G$. 
We see that every connectivity component of $H$ is a circular graph $C_n$ with $n \geq 4$. 
\end{proof} 

More generally, we can look at a collection  coloring functions
$f_1,\cdots ,f_k$ on $G \in \Gcal_d$ and define a {\bf discrete algebraic variety} 
$H=\{ f_1=c_1, \dots ,f_k=c_k \; \}$ as 
follows: its vertices are the $k$-dimensional simplices in $G$ for which all $f_j$ change sign. 
Two such simplices are connected if they intersect in a $(k-1)$-dimensional simplex. For $k=1$, we
know that these {\bf hyper-surfaces} can be completed to be geometric. For $G \in \Gcal_3$ for 
example, it is a $2$-dimensional graph with triangular or square faces. Surfaces of codimension $2$
play a role in \cite{eveneuler} as the intersection of $f=c$ with a unit sphere $(x)$ is a graph $B_f(x)$ 
whose Euler characteristic determines the {\bf symmetric index} $j_f(x) = (i_f(x)+i_{-f}(x))/2$ of $f$. The Euler
characteristic of a $4$-dimensional graph can now be written by Poincar\'e-Hopf 
as an expectation over the Euler characteristic of all $2$-dimensional $2$-dimensional subgraphs in a 
well defined sense and justify that Euler characteristic
is in spirit a Hilbert action: the index average result shows that $\chi(G)$ is a 
sum of an exotic scalar curvature obtained by adding up in a natural way all sectional 
curvatures through a point. The argument was that because Hilbert action is important in physics, 
that Euler characteristic should be taken seriously as a functional in graphs 
(see also \cite{KnillFunctional,randomgraph}. 
Graphs have played an important role in gravitational theories, as a computation tool or as spin networks
\cite{Rovelli}. \\

Lets think of $G$ as the host graph like the projective space in algebraic geometry. Then 
the graph $H=\{ f_1=c_1, \dots ,f_k=c_k \; \}$ is the analogue of an {\bf algebraic variety} $X$
and the points where the unit sphere $S(x)$ fails to be a sphere is a {\bf singular point}. 
This is analogue to the classical situation, where at singular points $x$ a small sphere $S_r(x) \cap X$
is no more a sphere. At a double point for example, where a curve crosses, the unit sphere $S_r(x) \cap X$
consists of 4 points and not of two. Already the hyper-surface case shows that we have to complete the 
surface to make it geometric. The condition of a function $f$ to {\bf be a coloring} is analogue to 
{\bf having no singular point}, the local injectivity condition being analogue to the non-vanishing of 
the gradient. 

\question{
Under which conditions can a discrete algebraic variety $V$ in a graph $G \in \Scal_d$ be 
completed to be in $\Vcal_k$ or $\Gcal_k$?  }

As in the continuum, the local injectivity condition assures
that discrete algebraic varieties are variety like. In the discrete they have a 
geometric completion if the host space $G$ is geometric. For a hyper surface $f=c$
the local injectivity condition avoids singularities. \\
 
To define what ``variety like" we need a notion of dimension. 
The definition of {\bf inductive dimension} starts with the assumption that the empty graph 
$G=\{\}$ has dimension $-1$. The dimension of a vertex $x$ in a graph is defined as
${\rm dim}(S(x))+1$. The dimension of the graph is the average over all dimensions on $V$. 
For any graph $G$, the dimension ${\rm dim}(G)$ is a rational number. 
Now we define a class $\Vcal_d$ of {\bf $d$-dimensional varieties}. It is larger than the 
class $\Gcal_d$ of geometric graphs. Let $\Vcal_{-1} = \{ \emptyset \}$. Set
$\Vcal_0=\Gcal_0$. For $d \geq 1$, $G \in \Vcal_d$ if every unit sphere $S(x)$ is in $\Vcal_{d-1}$.
The set $\sigma(G)$ of {\bf singularities} of $G$ is defined as the graph generated by 
vertices $x$ for which $S(x)$ is not in $\Scal_{d-1}$
nor in $\Bcal_{d-1}$. We require that both $\sigma(G)$ and $\delta(G)$ are in 
$\Vcal_{d-1} \cup \Vcal_{d-2} \cup \cdots \cup \Vcal_{0}$.
A graph $G \in \Vcal_d$ is geometric if $\sigma(G)$ is empty and if 
$\delta(G)$ is empty or in $\Gcal_{d-1}$.
A graph $G \in \Gcal$ is {\bf closed} if its boundary is empty. A closed geometric
graph is a variety for which both $\sigma(G)$ and $\delta(G)$ are empty. \\

A $1$-dimensional graph for example unit spheres with one, two or more discrete points. 
A $1$-dimensional variety is a $1$ dimensional graph for which the set
of vertices of degree larger than two are isolated. The star graph $S_3$ is the smallest variety
which is not a geometric graph. The figure $8$ graph is an example of a variety which
is not simply connected. A tree obtained by gluing two star graphs $S_3$ long an edge
is a $1$-dimensional graph but not a variety since the singularity set is $1$-dimensional.
Various chemical molecules are $1$-dimensional varieties, discretizations of classical 
varieties are often varieties as graphs. The following is related to Gr\"otzsch's theorem
establishing that all planar $1$-dimensional graphs are in $\Ccal_3$ but its much easier
as we do not have to deal with general planar graphs but with $1$-dimensional graphs 
for which the singularity set is untangled: 

\resultlemma{
$\Vcal_1 \subset \Ccal_3$. A graph $G \in \Vcal_1$ is in $\Ccal_2$ if and only 
if there are no odd-length cycles. }

\begin{proof}
The singularity set by definition is isolated. We can color each
of its vertices with the color $A$. Now color the connections between
any two singularities. If an odd number of vertices is between two
nodes, we can color with a second color $B$ only. In general, we need
three colors.
For the second statement see \cite{TuckerGross} Theorem 1.5.2: fix
a vertex and color the vertices with even distance with one color and
the others with the other.
\end{proof}

We believe that for coloring purposes, singularities 
do not matter as they are of lower dimension: 

\question{ Is $\Vcal_2 \subset \Ccal_5$?  }

\section{Outlook}

We have not seen any obstacles for chopping up the interior of a graph $H \in \Bcal_3$ 
with boundary $\delta H \in S_2$ in such a way that all interior edges have
even degree. In this respect we can ``see" why the four color theorem
is true. The proofs in \cite{AppelHaken,RSS,Gonthier} establish proofs which some
consider less insightful. But of course, the proof that a refinement can always be done
still needs to be given. We are in the ``Heesch stage" where we use the computer as a guide
to do the cutting. The goal is to have a deterministic constructive cutting algorithm which 
always works. The proof of the four color theorem would not only become {\bf visible}, it also 
would be {\bf constructive}. \\
 
Due to the fact that the obstacle set $O$ is the boundary of a surface, we expect this to be 
doable even in an elegant way.  We also hope that the four color theorem one day will become {\bf ``teachable"} 
meaning that a detailed proof can be given in the classroom in one session. 
But we have to remember that many predictions regarding the 4 color
problem had been wrong, even statements by leading mathematicians of their time have turned out to be
false. The history of the four color problem has shown, that many surprises can await
\cite{Ore,Tietze,Heesch,FritschFritsch,SaatyKainen,Soifer,ChartrandZhang2,RobinWilson}.
Our own experiments do not indicate any impediments. They ignite hope for a constructive approach.
The geometric calculation procedure trivially works for $1$-dimensional graph: just
chop up the interior of a disc to make all interior vertex degrees even. Two subdivisions are enough 
in every connected component. \\

One can ask for example, whether the chopping up procedure introduced here can be rephrased with ``discharging
procedures" introduced by Heesch, the architect of the proof of the four color theorem. As the topology 
of three dimensional space is much richer and much, much more complicated than the topology of the plane, 
we suspect however that this is not the case and that the ``computations" done by the tetrahedra in the 
inside could turn out to be more powerful than the ``field theoretical" planar methods introduced by Heesch. \\

Anyway, we hope to have shown that the story of graph coloring remains interesting. In particular,
because coloring questions of {\bf higher genus surfaces} in $\Gcal_2$ and higher dimensional graphs
like $\Scal_d$ are not settled. Already the question whether genus $1$ surface in $\Gcal_2$ 
is in $\Ccal_4$ is unanswered. \\

For the classical 4 color theorem, not everybody is happy with a machine proof. 
In general, mathematicians lament the {\bf lack of insight} to see
why four colors are enough. Much has been written about this, also in cross-disciplinary areas:
see for example \cite{RotaBeauty,MacKenzie,Ulianov,RuelleBrain,Kraken,StewartConcepts,Heintz2000,
Halmos1990,KrantzPudding,Soifer}. Oscar Lanford III, 
a pioneer in computer assisted methods, who proved the Feigenbaum conjectures about the existence
of solutions to the functional equations $g(x)=-\lambda^{-1} g(g(-\lambda x))$
with the help rigorous estimates provided by a computer comments on this in his articles.
The Computer-Human debate will probably always remain. More agreement is the fact that
{\bf high complexity} which can be dissatisfying, both in mathematics as well as in computer programs.
Complex proofs or programs have higher the risk of a mistake as proofreading
entangled the arguments or coding is challenging. This is unrelated 
whether a proof or program has been written in a computer assisted way.
\cite{Lanford84} explains how  the computer task is not complex as it essentially establishes bounds
(which can be done by computing with integers) which ensure that some Euclidean domain is mapped into itself
or that a contraction condition is satisfied: fixed point theorems like the Schauder or Banach
assure that the fixed point exists. The proof of the four color theorem proof on the other hand is still 
complex. In view of the history of the theorem, one wondered whether  the reduction analysis water tight 
and whether it does cover all cases. It is work of several decades with most theoretical groundbreaking work done
by Birkhoff or Heesch. In the case of the Kepler conjecture, the complexity is even higher 
\cite{HalesFerguson}. Its tempting to work on the later problem too with the computer as the
machine can generate millions of packings with the aim to look for a counter example \cite{Kni95d} or
denser packings in higher dimensions, where the answers are still unsettled. \\

Human proofs sometimes need some smoothing until they are bullet proof.
There are important results in mathematics which were born with flawed
proofs or gaps. The {\bf fundamental theorem of algebra} is an example:
Gauss own proof in 1799 is not complete from todays standards. And much before, Leibniz even
believed in 1702 that the fundamental theorem of algebra is wrong, claiming that $x^4+t^4$
could not be written as the product of two real algebraic factors \cite{Carrera}. 
But we have now various elegant proofs
like that a degree $n$ polynomial $f$ defines a Riemannian metric $g=|f|^{-2/n} |dz|^2$
on the Riemann sphere having curvature $n^{-1} \Delta \log|f|$. Curvature would be zero everywhere 
if ${\rm Re} \log(f)$ would be harmonic meaning $f$ had no roots \cite{AlmiraRomero} contradicting so
Gauss-Bonnet for spheres. In the four color theorem, 
11 years were needed to refute the published and acclaimed proof of Kempe 
which leading experts of the time had accepted. An other famous story explaining how a theorem evolves is
described in \cite{lakatos}: the Euler formula had to be adapted as higher genus surfaces were involved
which were needed already for Kepler polytopes.
In that case, the culprit was the notion of ``polyhedron" which had been foggy at first \cite{Richeson}. 
We restrict therefore mostly to convex polytopes \cite{gruenbaum}. \\

Chromatic graph theory remains fascinating because of the 
many outstanding open problems, in particular $NP$ questions as well as the {\bf Hadwiger conjecture} asking 
whether any graph with chromatic number $c(G)=k$ has $K_k$ as a minor. 
The long historical struggle as well as its many connections to other fields are other reason why the subject is so attractive.
The story of the 4-color problem has been told many times \cite{Mayer1982}: \\

According to Richard Baltzer, it was M\"obius in 1840 who spread a precursor problem (\cite{BallCoxeter,FritschFritsch, RobinWilson}
given to him by Benjamin Gotthold Weiske. It is now called {\bf M\"obius-Weiske puzzle} \cite{Soifer}. The actual problem was first posed
when Francis Guthrie \cite{MaritzMouton}, after thinking about it with his brother Frederick, who communicated it to his teacher
Augustus de Morgan, a former teacher of Francis, in 1852. The later communicated it to William Hamilton, who however 
could not be warmed up to the problem.
Arthur Cayley in 1878 helped to spread the problem by putting in first in print and giving it the status 
of an  open problem. Even Cayley, who worked in graph theory did not use the language of graphs yet
to describe the problem \cite{Crilly}. Alfred Kempe published a proof in 1879 which turned out
to be incomplete. The gap was noticed by Heawood 11 years later in 1890 \cite{Heawood90}.
There was also an unsuccessful attempt by Peter Tait in 1880.
\cite{Soifer} calls this time a {\bf ``victorian comedy of errors"}.
After considerable work including Charles Pierce, George Birkhoff, Oswald Veblen, Philip Franklin, 
Hassler Whitney, Hugo Hadwiger, Leonard Brooks,  William Tutte, Yoshio Shimamoto,
Heinrich Heesch and Karl D\"urre (rumors of proofs by Heesch and Br\"odel appeared in 1947 \cite{Tietze}, a proof claim
of Shimamoto is from 1971) and pushing the limit
of reducibility without computers like in Walter Stromquist's 1975 thesis \cite{Stromquist}, and more exhaustive computer
analysis like Frank Allaire and Edwart Swart,
a computer assisted proof of the 4-color theorem was obtained by Ken Appel and Wolfgang Haken
in 1976. A new proof appeared in 1997 by Neil Robertson, Daniel Sanders, Paul Seymour, and Robin Thomas. 
and the computer programs by these four authors are still publicly available today. 
Geometric coloring theory was pioneered by Steve Fisk \cite{Fisk1977a,Fisk1977b,Fisk1980}, 
(Harvard PhD who after postoc at MIT was at Bowdoin college until his death in 2010)  just about at the 
time when the four color theorem was solved. The concept of {\bf locally planar graphs} is due Walter
Stromquist \cite{StromquistPlanar} who also proved the toroidal 5 color theorem together with Michael Albertson 
(who worked at Smith college until his death in 2009) \cite{AlbertsonStromquist}. 
Work on the 4 color theorem continued. Goerge Gonthier produced in 2004 a fully machine-checked proof of the 
four-color theorem \cite{RobinWilson}. \\

The hope is that one day, various different teachable proofs will exist for the four color theorem. 
We believe that the theory pioneered by Fisk could play an important role. 
The proof of the fundamental theorem of algebra mentioned above uses Riemannian geometry and potential theory.
It illustrates how {\bf different fields of mathematics} can help to understand a subject. We hope to have 
given hope that ideas from geometric topology  like Fisk theory and discrete differential geometry 
could lead to a proof of the four color theorem one day. \\

Geometric coloring approaches are less studied for higher dimensional geometric graphs. 
We expect that the chromatic number of a graph $G$ depends on topological embedding realizations
of $G$ into higher dimensional graphs which can be minimally colored. \\

As many interesting problems are left and in particular since
the general coloring problem is {\bf NP-complete}, we expect graph coloring problems 
to remain interesting for a long time to come. One the same note, the Hamiltonian path problem is NP complete.
Graph coloring and Hamiltonian path problems have always been very close. It is no accident that some pioneers
were both concerned with colorings and Hamiltonian paths, most notably Whitney, who started his mathematical career as
a graduate student of Birkhoff in the field of chromatic graph theory. Tait approached the four color problem with
Hamiltonian path methods. The {\bf Whitney-Tutte theorem} \cite{Whitney1931,BM,Ore} tells that any 
4-connected planar graph is Hamiltonian. It follows from the easier direction of the lemma that

\resultcorollary{
All graphs in $\Scal_2$ are Hamiltonian.
}

Tait had conjectured that any 3-connected planar cubic graph has a 
Hamiltonian cycle. While this would have implied the 4-color theorem,
counter examples were found by Tutte. We can now ask: 

\question{Are all graphs in $\Gcal_2$ Hamiltonian? 
Are all graphs in $\Scal_d$ Hamiltonian? }

The Fisk torus in $\Gcal_2 \setminus \Scal_2$ and the projective plane outside $\Ccal_4$ 
are Hamiltonian despite the fact that it has chromatic number $5$. We checked also with other graphs
like Kleinbottles. Both geometric platonic solids in $\Scal_3$, the 16-cell as well as the 600-cell are Hamiltonian. \\

Constructive approaches to the 4-color theorem appeared early on. Kempe's attempt was constructive:
an {\bf ab-Kempe chain} with colors a,b attached to vertex is a maximal
connected component only using 2 colors. There is a natural involution
switching colors of such a Kempe chain. The Kempe algorithm first looks in $G=G_0$ for a vertex
of degree 5 or less, then in the remaining graph $G_1$ for a vertex of degree 5 or less
etc. This produces an ordering $v_1,\dots v_n$ of the $n$ vertices and a filtration $G_i$ of
subgraphs. The vertices $v_i$ are now colored recursive in a ``greedy way": define $f(v_n)=0$, then
$f(v_{n-1})=0$ if not adjacent or $v_{n-1}=1$ if adjacent etc. If a coloring is not
possible any more, perform a Kempe chain switch in a predefined way, then move on to the next vertex.
If the algorithm gets through to $v_0$, the coloring is done.
While the Kempe algorithm does not always produce a coloring, it often does, in particular,
if the maximal degree is less than 9 \cite{GethnerSpringer}.
In \cite{involve} the failure rate of Kempe's algorithm was investigated.  \\

Let $\Gcal_d$ be the set of geometric $d$-dimensional graphs, graphs for
which the unit spheres are in $\Scal_{d-1}$, the class of $(d-1)$-dimensional
spheres. In view of the {\bf Whitney embedding theorem}, we expect: 

\question{
Is it true that $\Scal_d \subset \Ccal_{d+2}$ and $\Gcal_d \subset \Ccal_{2d+1}$?
}

More generally, one can ask whether every $G \in \Gcal_d$ zero-cobordant with a simply connected $H$ 
is in $\Ccal_{d+2}$. To study this, one has to study to adjust the degree of a $(d-2)$-dimensional 
simplex $x=(x_1,\dots,x_{d-1})$ as the number of vertices of the 
intersection $C$ of all unit spheres $S(x_i)$ which is a circular graph. 
If all these degrees are even in $G \in \Scal_{d}$, then $G \in \Scal_{d} \cap \Ccal_{d+2}$. 
To adjust the degrees, one has to split $(d-1)$-dimensional surfaces by adding a new 
vertex in the center of $x$ and connect to all $x_i$ as well as to all vertices of $C$.
As for $d=1$ or $d=2$, there should be no impediments in the simply connected case. 
We have implemented the chopping procedure on the computer, but still 
rely on human assistance in doing the cutting as checking through all possible cuttings is
too costly. We believe that a systematic cutting procedure can be developed but we 
have not done so yet. Of course one has to be open to refine the graphs differently.
Like a $2$-dimensional graph can be refined by filling in an octahedron inside a triangle,
we can fill in the 3-cross polytop into a tetrahedron. This does not increase the set $O$. 

\section{The projective plane}

The projective plane is the generator of the cobordism group of $\Gcal_2$ and so chromatically
interesting. We don't know whether every projective plane is in $\Gcal_5$ but believe this to be the case. 
The projective plane can be obtained by gluing a disc to a M\"obius strip. To illustrate some of the
geometric notions, lets look closer at the Moebius strip. \\

The {\bf Moebius strip} $G$ is a non-orientable graph.
The smallest version has no interior points and is obtained from $C_7 = Z_7=\{0, \dots ,6\}$ 
by adding additional edges $(x,x+3)$ modulo $Z_7$.  
Having $v=v_0=7$ vertices and $e=v_1=16$ edges and $f=v_2=7$ triangles, we have
$\chi(G)=7-14+7=0$. Its chromatic polynomial is $f(x)=(x-3)(x-2)(x-1)x(x^3-8x^2+25x-29)$
so that $c(G)=4$ and $f(4)=7 \cdot 4!$ reflects the $7$ fold symmetry besides the color 
permutation symmetry. If we add an other vertex $v$ to $C_7$ and connect $v$ to every 
point in $C_7$, we get a discrete disc, a wheel graph, a graph in $\Bcal_2$. 
which is a $2$-dimensional geometric graph with boundary. 
Gluing the wheel graph and the Moebius strip together at $C_7$
gives the discrete {\bf projective plane} $P^2$. It is in $\Gcal_2$ and has chromatic number $5$.
The chromatic polynomial is $f(x) = x(x-1)(x-2)(x-3)(x-4)$
$(708386-1540046x+1556883x^2-968364x^3+411914x^4-125510x^5+27769x^6$ 
$- 4403x^7 + 478x^8 - 32x^9 + x^{10})$ so that $G$ admits $f(5)=2^3 \cdot 21 \cdot 17 \cdot 5!$ 
colorings with $5$ colors. \\

To see the cohomology connection while not having to deal with too large matrices, we go back to
the M\"obius strip $G$. Given an orientation of the $14$ edges of $G$, 
the {\bf gradient} is a $14 \times 7$ matrix 
$$ d_0 = 
\left[
\begin{array}{ccccccc}
 -1 & 0 & 0 & 0 & 1 & 0 & 0 \\
 -1 & 0 & 0 & 0 & 0 & 0 & 1 \\
 -1 & 1 & 0 & 0 & 0 & 0 & 0 \\
 -1 & 0 & 0 & 1 & 0 & 0 & 0 \\
 0 & -1 & 0 & 0 & 0 & 1 & 0 \\
 0 & -1 & 1 & 0 & 0 & 0 & 0 \\
 0 & -1 & 0 & 0 & 1 & 0 & 0 \\
 0 & 0 & 1 & -1 & 0 & 0 & 0 \\
 0 & 0 & 0 & -1 & 1 & 0 & 0 \\
 0 & 0 & 0 & -1 & 0 & 0 & 1 \\
 0 & 0 & -1 & 0 & 0 & 0 & 1 \\
 0 & 0 & -1 & 0 & 0 & 1 & 0 \\
 0 & 0 & 0 & 0 & -1 & 1 & 0 \\
 0 & 0 & 0 & 0 & 0 & -1 & 1 \\
\end{array} \right] $$
Given also an orientation of the triangles, we get a map, which assigns 
to a function on the edges a function on the triangles. This linear map
$d_1$ is the curl: 
$$  d_1 = 
\left[
\begin{array}{cccccccccccccc}
 1 & 0 & -1 & 0 & 0 & 0 & -1 & 0 & 0 & 0 & 0 & 0 & 0 & 0 \\
 1 & 0 & 0 & -1 & 0 & 0 & 0 & 0 & -1 & 0 & 0 & 0 & 0 & 0 \\
 0 & 1 & 0 & -1 & 0 & 0 & 0 & 0 & 0 & -1 & 0 & 0 & 0 & 0 \\
 0 & 0 & 0 & 0 & 1 & -1 & 0 & 0 & 0 & 0 & 0 & -1 & 0 & 0 \\
 0 & 0 & 0 & 0 & 1 & 0 & -1 & 0 & 0 & 0 & 0 & 0 & -1 & 0 \\
 0 & 0 & 0 & 0 & 0 & 0 & 0 & -1 & 0 & 1 & -1 & 0 & 0 & 0 \\
 0 & 0 & 0 & 0 & 0 & 0 & 0 & 0 & 0 & 0 & 1 & -1 & 0 & -1 \\
\end{array}
\right]
 $$
We can check that ${\rm curl}({\rm grad}(f)=0$. The matrix
$d_0^*$ is called the {\bf divergence}. As in calculus, ${\rm div}({\rm grad}(f)$ 
is the scalar Laplacian of $G$. It is a $7 \times 7$ matrix
$$ L_0 = \left[
\begin{array}{ccccccc}
 4 & -1 & 0 & -1 & -1 & 0 & -1 \\
 -1 & 4 & -1 & 0 & -1 & -1 & 0 \\
 0 & -1 & 4 & -1 & 0 & -1 & -1 \\
 -1 & 0 & -1 & 4 & -1 & 0 & -1 \\
 -1 & -1 & 0 & -1 & 4 & -1 & 0 \\
 0 & -1 & -1 & 0 & -1 & 4 & -1 \\
 -1 & 0 & -1 & -1 & 0 & -1 & 4 \\
\end{array}
\right] $$
It agrees with $L=B-A$, where $B$ is the diagonal matrix
containing the degree and $A$ is the adjacency matrix.
The eigenvalues  of the Laplacian are 
{6.24698, 6.24698, 4.55496, 4.55496, 3.19806, 3.19806, 0}
There is one eigenvalue $0$, belonging to the constant function.
Since ${\rm ker}(L_0)$ is $1$-dimensional, there is only one connected
component. The matrix $L_1 = d_0 d_0^* +  d_1^* d_1$ is
$$L_1 =\left[
\begin{array}{cccccccccccccc}
 4 & 1 & 0 & 0 & 0 & 0 & 0 & 0 & 0 & 0 & 0 & 0 & -1 & 0 \\
 1 & 3 & 1 & 0 & 0 & 0 & 0 & 0 & 0 & 0 & 1 & 0 & 0 & 1 \\
 0 & 1 & 3 & 1 & -1 & -1 & 0 & 0 & 0 & 0 & 0 & 0 & 0 & 0 \\
 0 & 0 & 1 & 4 & 0 & 0 & 0 & -1 & 0 & 0 & 0 & 0 & 0 & 0 \\
 0 & 0 & -1 & 0 & 4 & 0 & 0 & 0 & 0 & 0 & 0 & 0 & 0 & -1 \\
 0 & 0 & -1 & 0 & 0 & 3 & 1 & 1 & 0 & 0 & -1 & 0 & 0 & 0 \\
 0 & 0 & 0 & 0 & 0 & 1 & 4 & 0 & 1 & 0 & 0 & 0 & 0 & 0 \\
 0 & 0 & 0 & -1 & 0 & 1 & 0 & 3 & 1 & 0 & 0 & -1 & 0 & 0 \\
 0 & 0 & 0 & 0 & 0 & 0 & 1 & 1 & 3 & 1 & 0 & 0 & -1 & 0 \\
 0 & 0 & 0 & 0 & 0 & 0 & 0 & 0 & 1 & 4 & 0 & 0 & 0 & 1 \\
 0 & 1 & 0 & 0 & 0 & -1 & 0 & 0 & 0 & 0 & 4 & 0 & 0 & 0 \\
 0 & 0 & 0 & 0 & 0 & 0 & 0 & -1 & 0 & 0 & 0 & 4 & 1 & 0 \\
 -1 & 0 & 0 & 0 & 0 & 0 & 0 & 0 & -1 & 0 & 0 & 1 & 3 & -1 \\
 0 & 1 & 0 & 0 & -1 & 0 & 0 & 0 & 0 & 1 & 0 & 0 & -1 & 3 \\
\end{array}
\right]  \; . $$
It has one zero eigenvector 
$$ v=[1, -2, 2, -1, 1, 2, -1, -2, 2, -1, 1, -1, 2, 2]^T  \; . $$ 
This reflects the fact that
the graph is not simply connected. Actually, we can construct this eigenvector using the 
{\bf Hurewicz map}: take a non-contractible closed curve $\gamma$ in $G$ corresponding to a
nontrivial element in the fundamental group $\pi_1(G)$, take a $1$-form $F \in \Omega_1$ 
which is constant $1$ on the edges of $\gamma$ and $0$ on the rest. Then apply the heat flow 
$e^{-t L_1} F$. Because all nonzero eigenvalues of $L_1$ are positive, the heat flow kills the
corresponding eigenvectors and projects the original 1-form onto the kernel of $L_1$ which is as a harmonic
element in $\Omega_1$ a representative of a nontrivial cohomology class. 
Finally, look at $L_2 = d_1 d_1^*$ (note that as $G$ is $2$-dimensional, $d_2$ is zero). 
$$ L_2 = \left[
\begin{array}{cccccccc}
 3 & -1 & 0 & 1 & 0 & 0 & 0 & 0 \\
 -1 & 3 & 0 & 0 & 0 & 0 & 1 & 0 \\
 0 & 0 & 3 & 1 & 0 & 1 & 0 & 0 \\
 1 & 0 & 1 & 3 & 0 & 0 & 0 & 0 \\
 0 & 0 & 0 & 0 & 3 & 1 & 0 & 1 \\
 0 & 0 & 1 & 0 & 1 & 3 & 0 & 0 \\
 0 & 1 & 0 & 0 & 0 & 0 & 3 & 1 \\
 0 & 0 & 0 & 0 & 1 & 0 & 1 & 3 \\
\end{array}
\right] $$
which is invertible. This reflects that $G$ is not orientable. 
The Betti numbers of the Moebius strip are the dimensions of the
kernels of $L_k$. We have $\vec{b}  = (1,1,0)$ so that the Poincar\'e polynomial
is $p(x) = 1-x$.  By Euler-Poincar\'e, the 
number $p(1) = b_0-b_1=0$ agrees with the Euler characteristic.  \\

The projective plane $P^2$ has Euler characteristic $1$ as $b_0=1,b_1=0$ but $P^2$ is
not simply connected as a path $\gamma$ from a vertex to its antipode is closed but not 
contractible. We see why the Hurewitz homomorphism is trivial: applying the heat flow
to $F=1_{\gamma}$ converges to $0$ as the form-Laplacian $L_1$ has a trivial kernel. \\

A double cover of $P^2$ the sphere $S^2$ is simply connected. This means $\gamma + \gamma=0$
when seen as a chain and $H_1(P^2,\Z)=Z_2$.  As there are no branch points for
this cover, the {\bf Riemann-Hurwitz formula} establishes $\chi(S^2) = 2 \chi(P^2)$. 
Every orientable graph $G$ of genus $g$ can also in the discrete be seen as a double cover of the 2 sphere
with $2g+2$ branch points. The Riemann-Hurwitz formula gives $\chi(G) = 2 \chi(S) - (2g+2)
=2-2g$. We could also glue together two projective planes to a Klein bottle $K$ with 2
branch points so that $\chi(K) = 2 \chi(P)-2=0$. In the context of Riemann-Hurwitz one has to
mention the obvious: \\
if $H$ is a cover of $G$ then $c(G) \geq c(H)$.
Indeed, it is true in general that if $G$ is the image of a graph homomorphism $\phi: H \to G$, then
a coloring of $G$ can be lifted to a coloring of the cover $H$.  \\

As we can construct higher genus surfaces as branched covers of the sphere, we get from the four color theorem
$\Scal \subset \Ccal_4$ that these higher genus surfaces are in $\Ccal_4$. Examples of tori can be obtained for example
by taking two copies of a spheres, cut out two discs in both and glue them to a doughnut. Glueing like this two doughnuts
gives a genus $2$ surface etc. Of course, we do not know about all projective planes nor Klein bottles but
only about those having a nontrivial involution in the automorphism group allowing to take the quotient.
While we know the torus to be in $\Ccal_5$ only in general we are ambivalent about
the Klein bottle, as doing refinements in discrete solid Klein bottles could be more subtle. 
Its most likely that every Klein bottle in $\Gcal_2$ is also in $\Ccal_5$.

\section{Illustrations}

The following figures try to illustrate the story. The pictures also should demonstrate that one can deal
pretty well visually with geometric graphs using a computer. Computer algebra systems like Mathematica have
graph theory built in. The computer can generate from a given graph various types of embeddings, both planar
as well as in space. Since our modification and cutting steps cahnge the embedding drastically, 
we had built a data structure, where both the graph as well as the location of the vertices is processed
simultaneously. The list of thumbnails of 36 mathematicians mentioned is not complete. We have been in particular
sorry not to find a picture of the late Steve Fisk, whose work from the 70ies is very close to what is done
here. 

\begin{figure}[h]
\scalebox{0.29}{\includegraphics{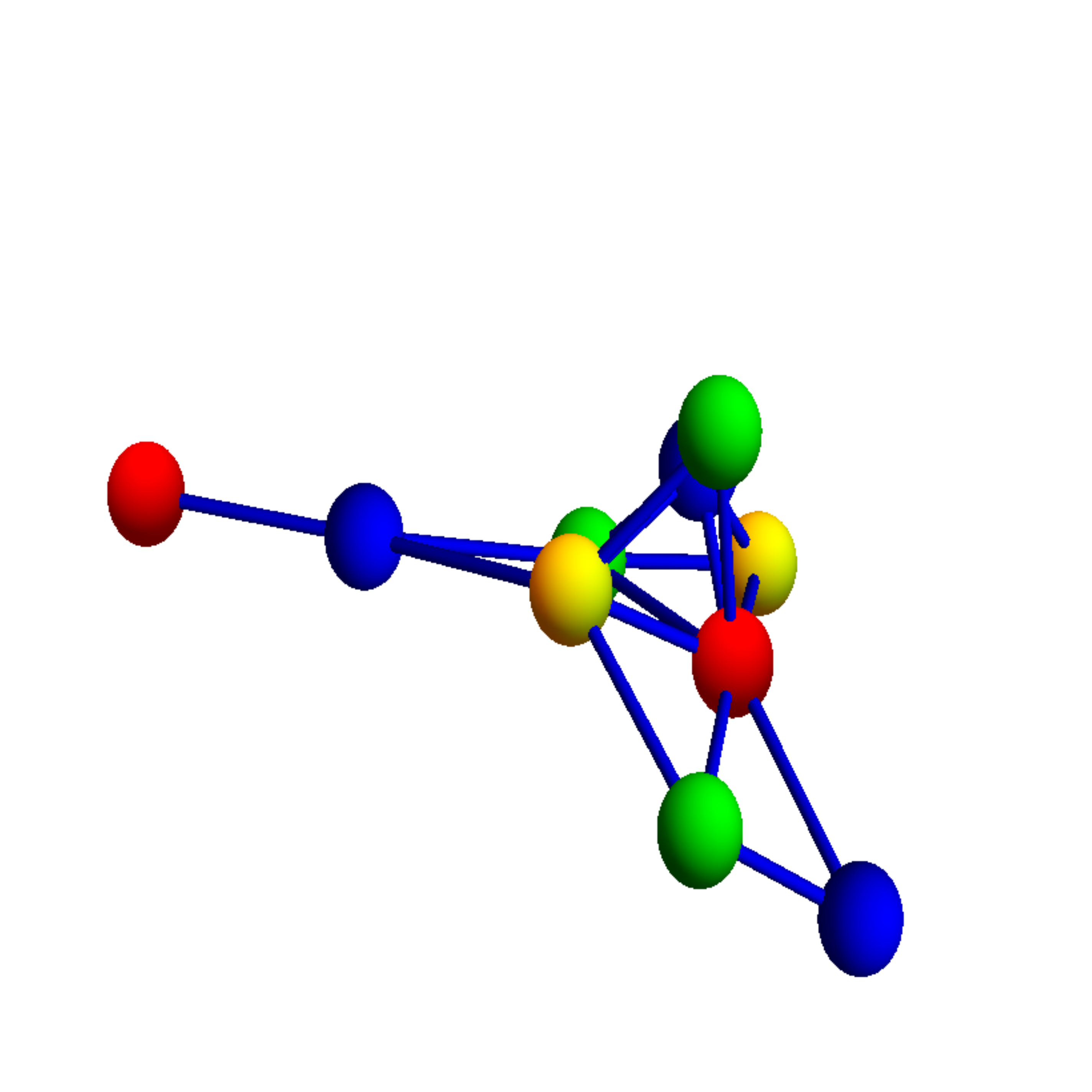}}
\scalebox{0.29}{\includegraphics{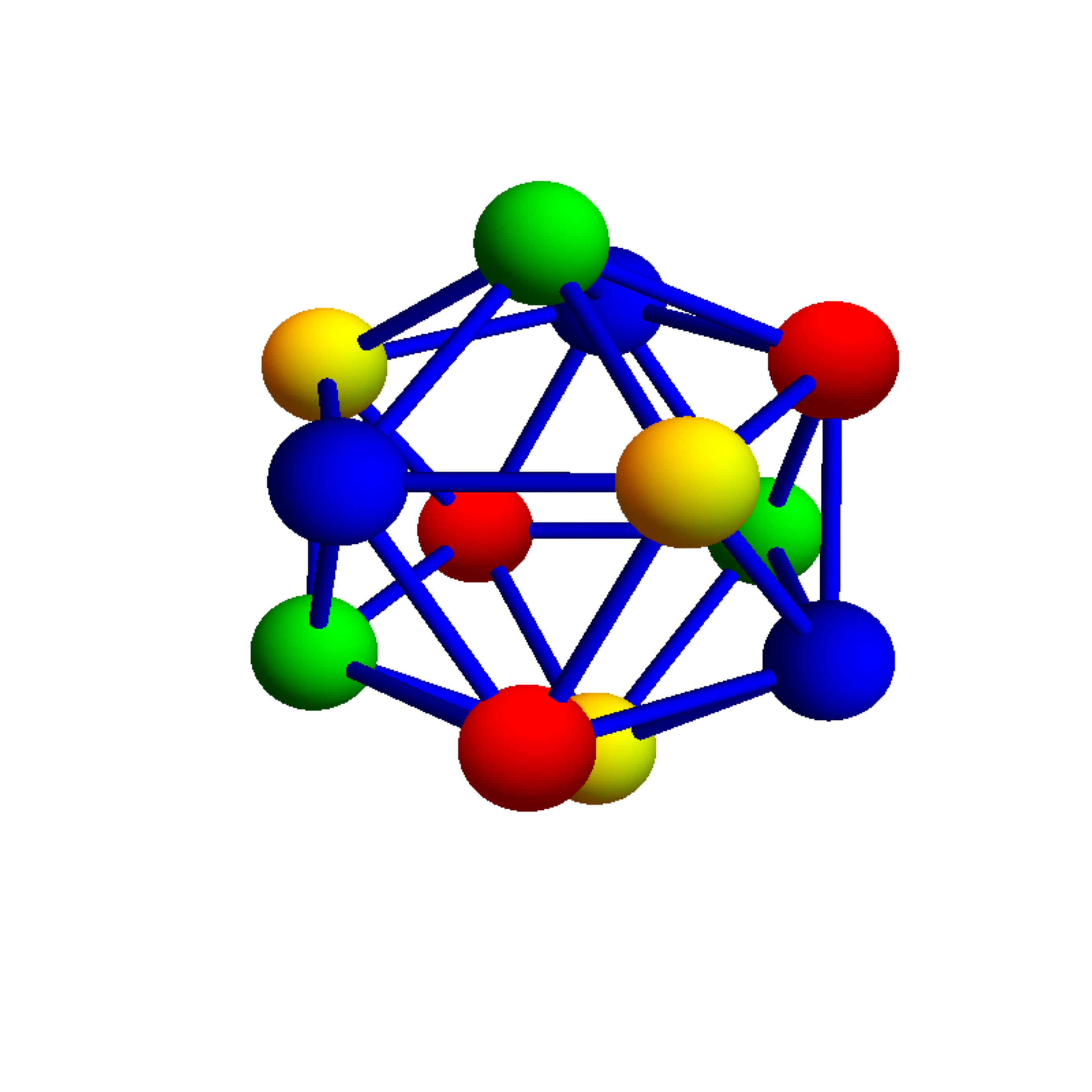}}
\caption{
The first figure shows a random graph in $\Ccal_4$. It is not geometric. 
The icosahedron $G$ is in $\Scal_2$. It is not in $\Ccal_3$ as 
it is not Eulerian.  
The chromatic polynomial of $G$ is
$c(G) =  (x-3) (x-2) (x-1)x$ $(x^8-24 x^7+260 x^6-1670 x^5+6999 x^4-19698 x^3$
$+36408 x^2-40240 x+20170)$ which evaluates to $240=10 \cdot 4!$ at $x=4$. 
The icosahedron is the largest two dimensional graph with positive curvature. 
}
\end{figure}

\begin{figure}[h]
\scalebox{0.22}{\includegraphics{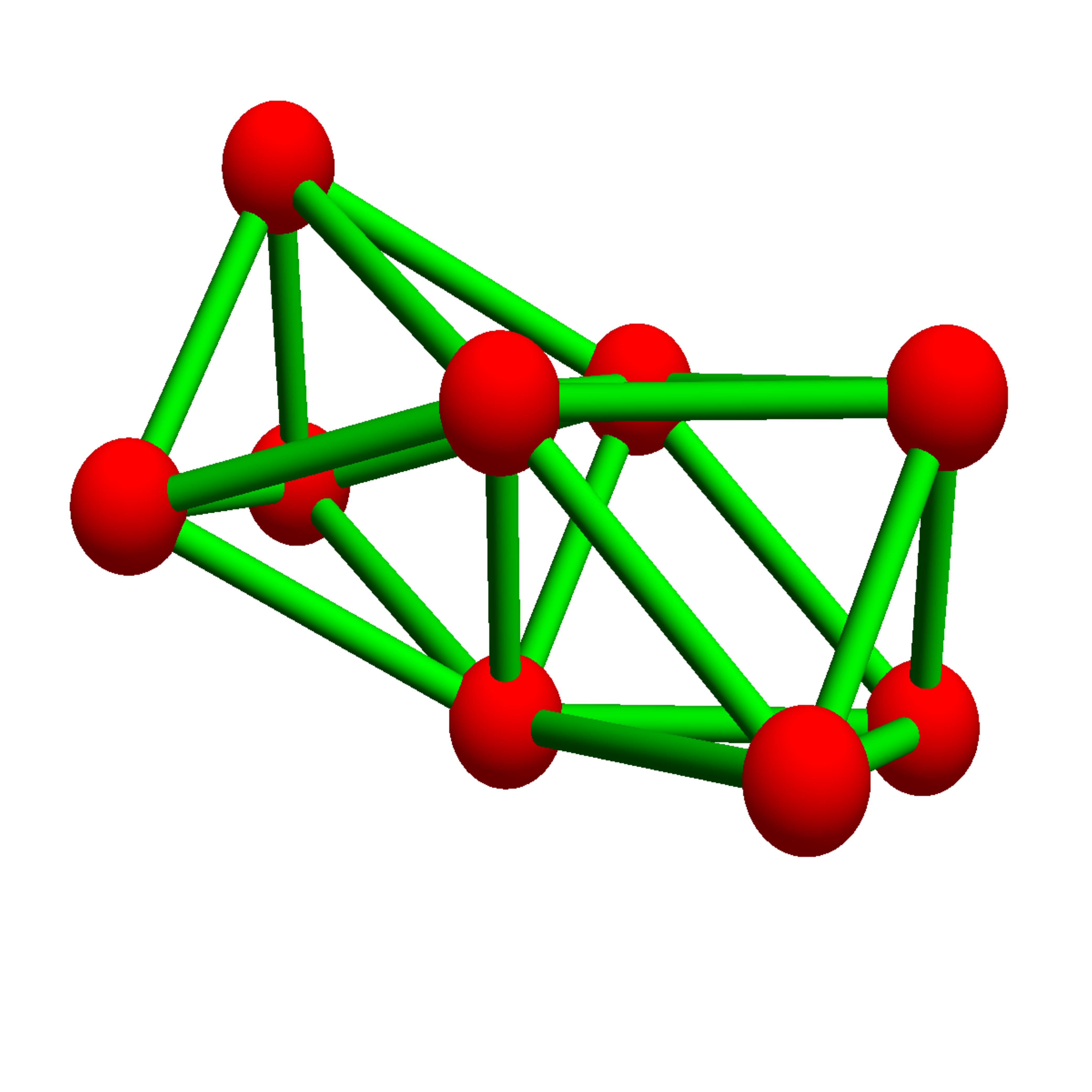}} 
\scalebox{0.22}{\includegraphics{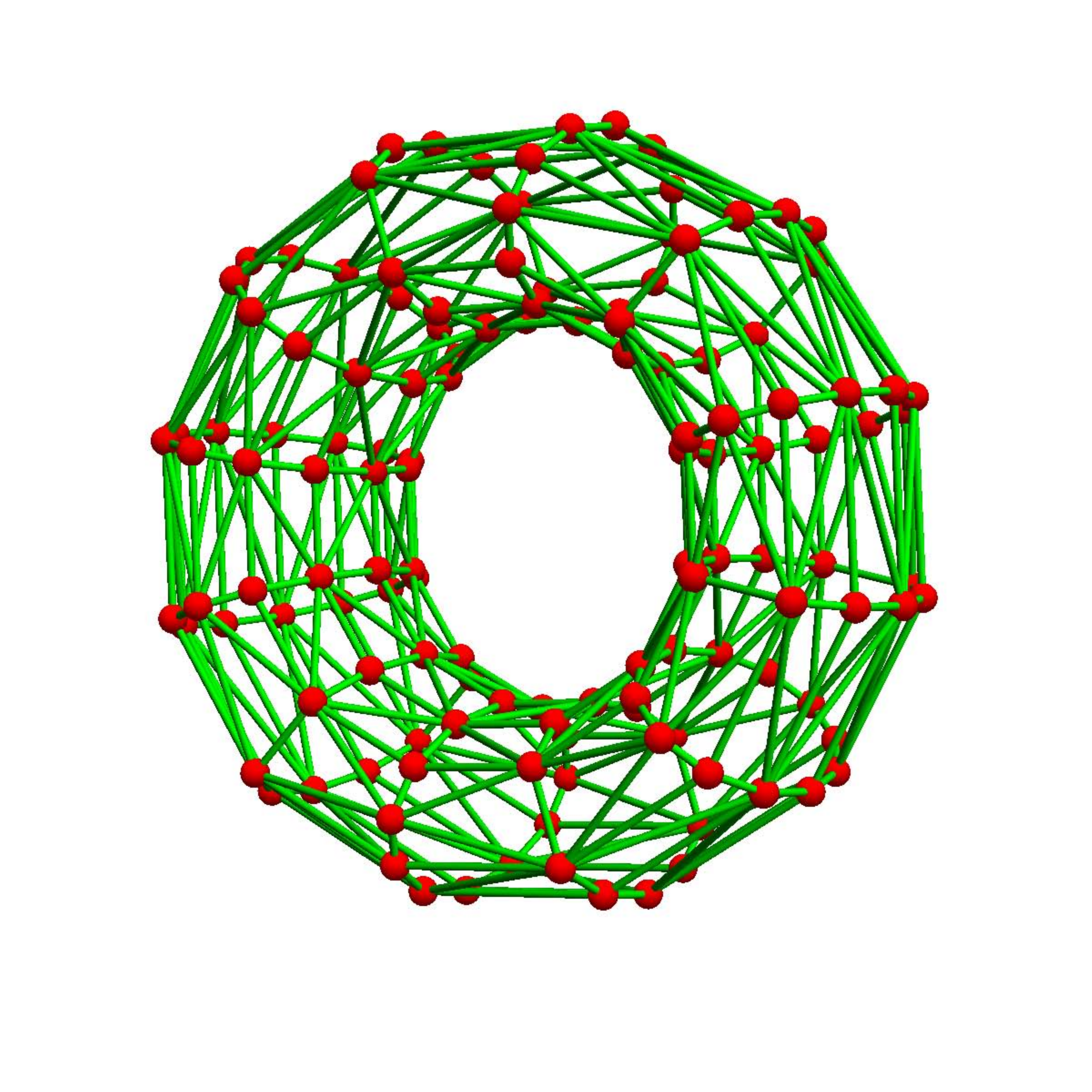}} 
\caption{
The twin octahedron is in $\Pcal$ but not in $\Scal$: not all unit spheres are circular. 
The graph is not 4-connected. It can be broken into two octahedra along a triangle. 
The torus below is in $\Gcal_2$ but not in $\Scal_2$ as the graph
is not simply connected. Alternatively, punching a hole into it keeps it not 
simply connected.}
\end{figure}

\begin{figure}[h]
\parbox{14.8cm}{
\parbox{6.6cm}{\scalebox{0.15}{\includegraphics{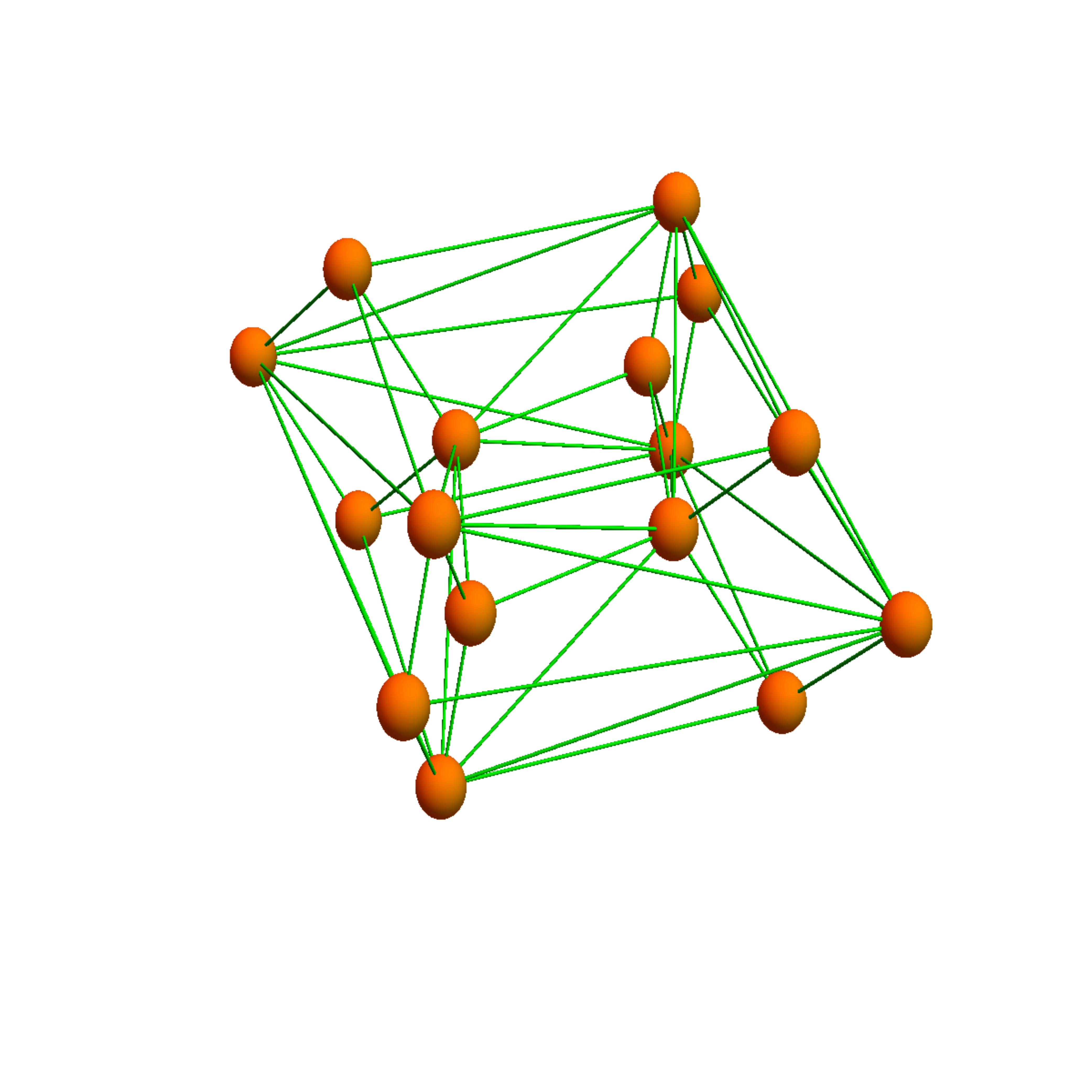}} }
\parbox{6.6cm}{\scalebox{0.15}{\includegraphics{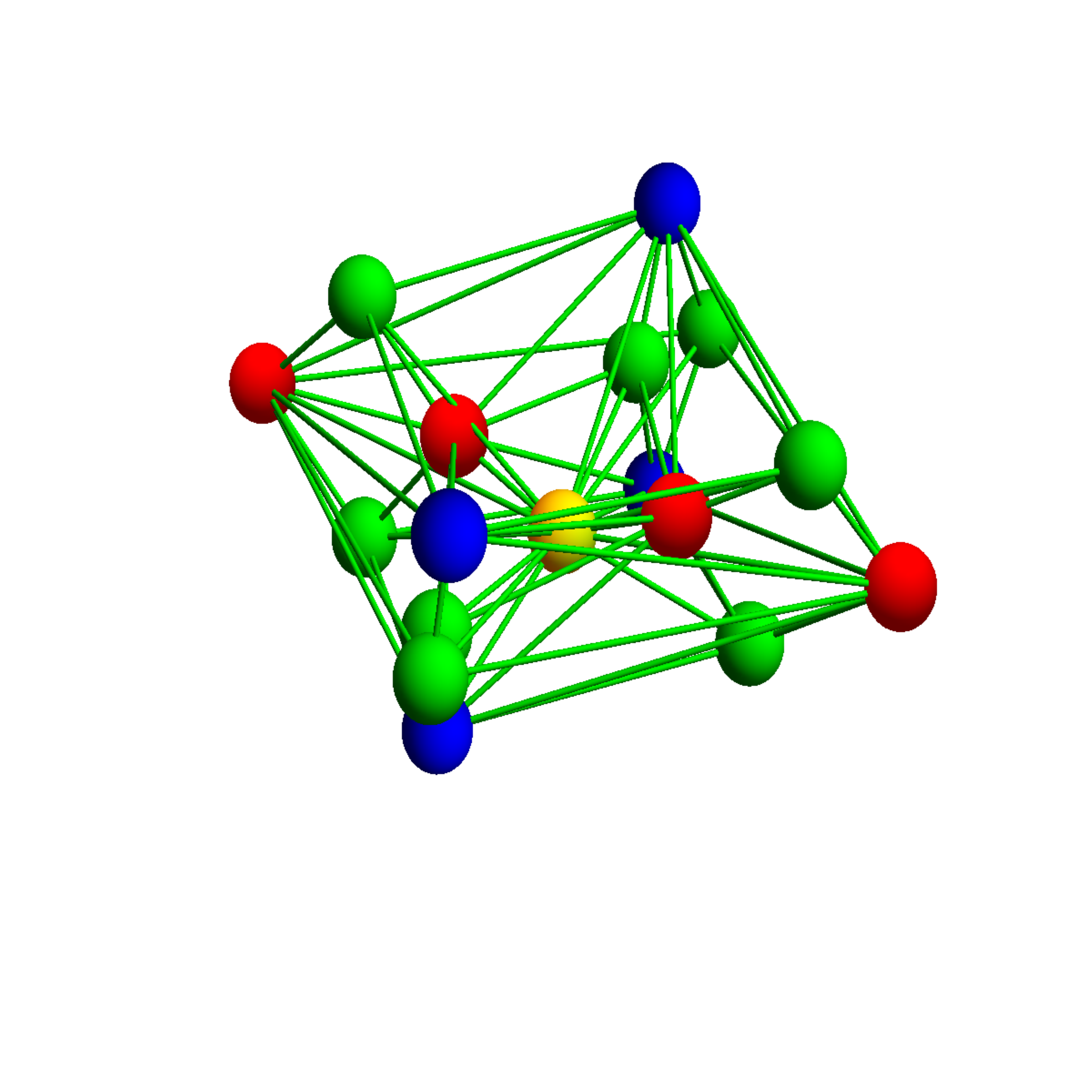}} }
}
\caption{
This torus in $\Gcal_2 \cap \Ccal_3$ resembles the embedding of the torus in $3$-dimensional 
space. Its curvatures take values in $\{ 1/3,-1/3 \}$. To the right we see the coloring. It was obtained
by seeing it as the unit sphere of a $3$-dimensional graph $H$ which has a $G$ as boundary.
While $H$ is contractible it is not geometric as there is a vertex for which $S(x)=G$ has
Euler characteristic $0$. But it easily can be refined to become geometric. 
The torus $G$ is then rendered 0-cobordant. But the three dimensional filling is not 
simply connected. }
\end{figure}

\begin{figure}[h]
\parbox{6.2cm}{\scalebox{0.12}{\includegraphics{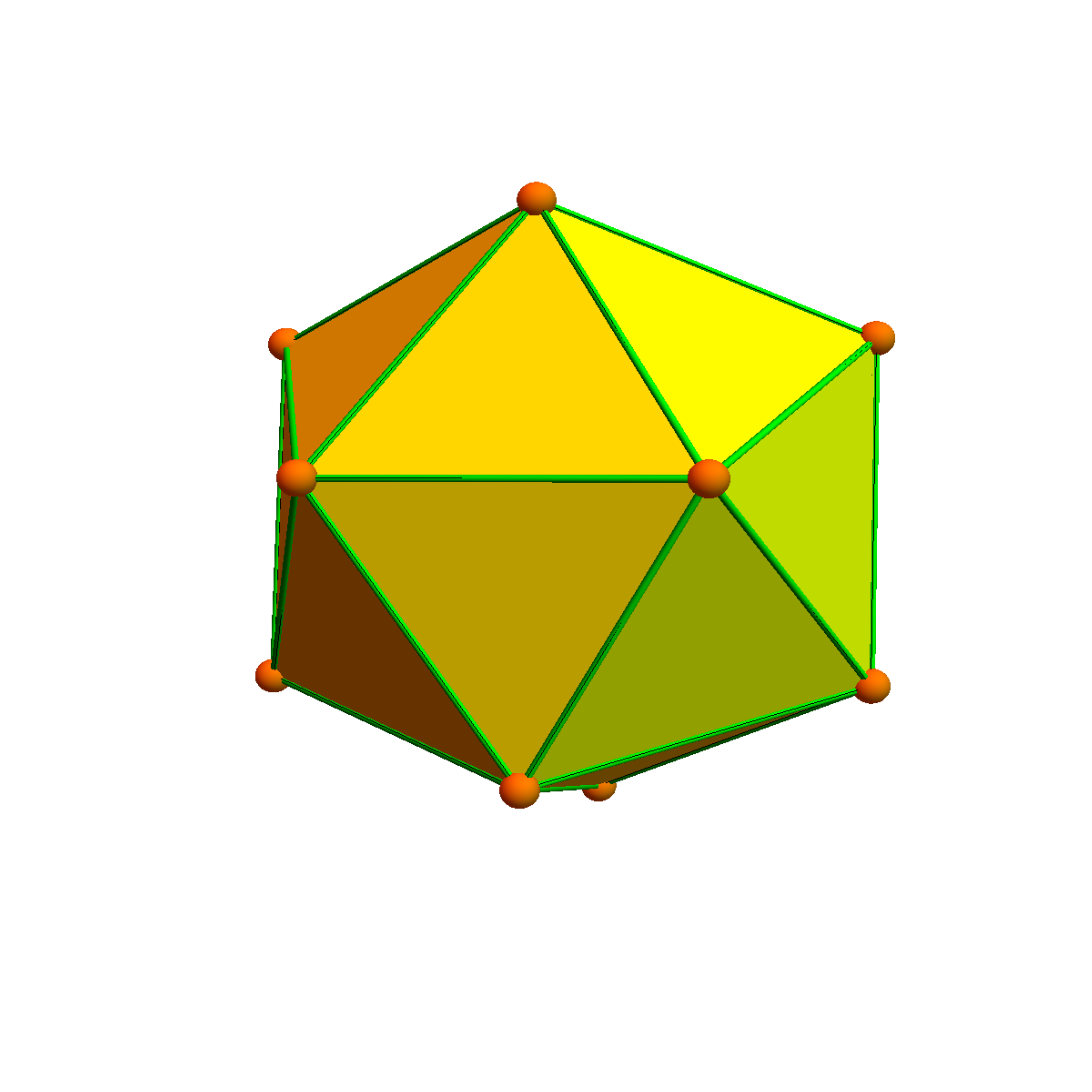}} }
\parbox{6.2cm}{\scalebox{0.12}{\includegraphics{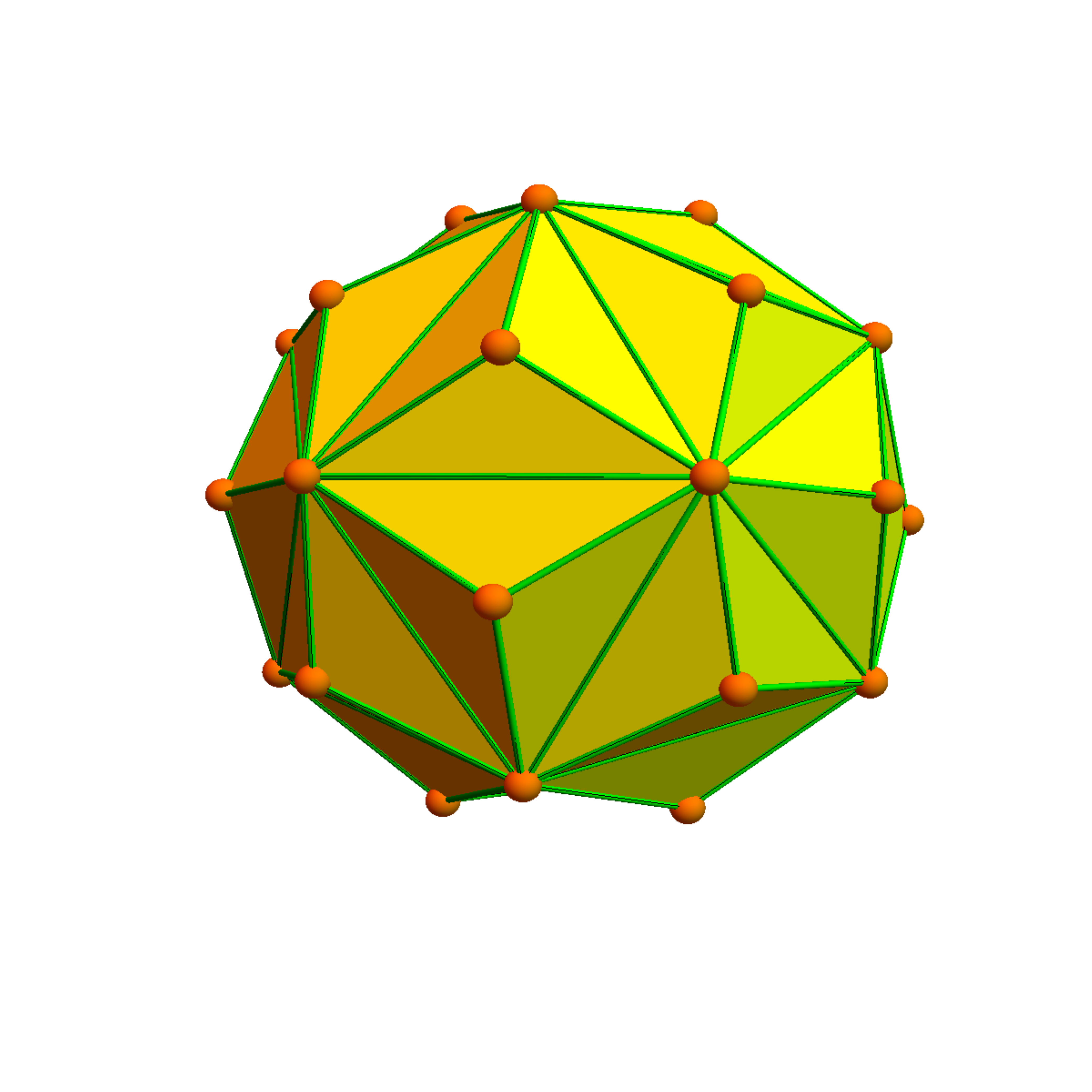}} }
\caption{
The tetrahedral completion of the icosahedron in $\Scal_2$ is obtained by 
capping each face. It is still a planar graph but it is no more $4$-connected:
cutting along a triangle isolates a single vertex. 
}
\end{figure}

\begin{figure}[h]
\scalebox{0.1}{\includegraphics{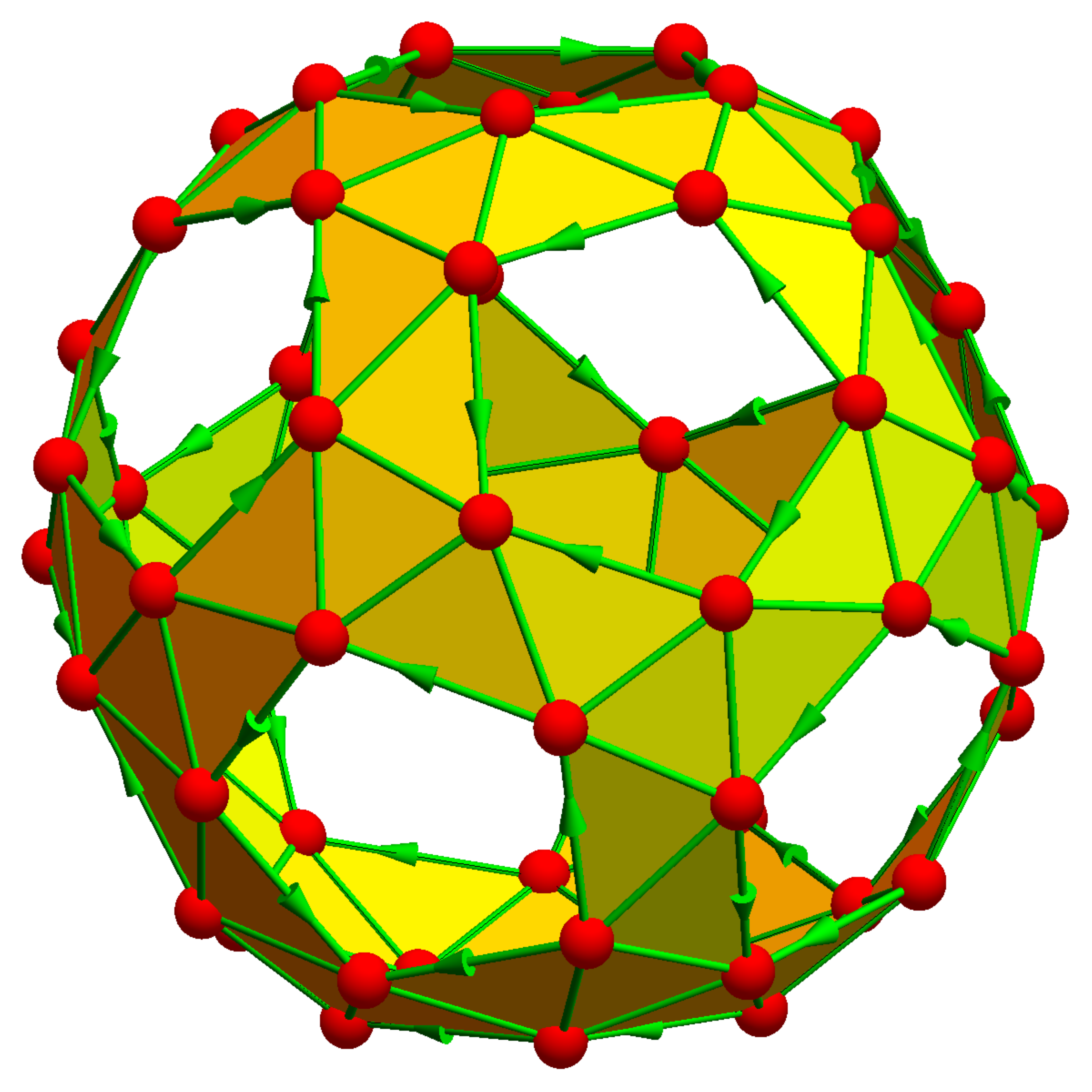}}
\scalebox{0.15}{\includegraphics{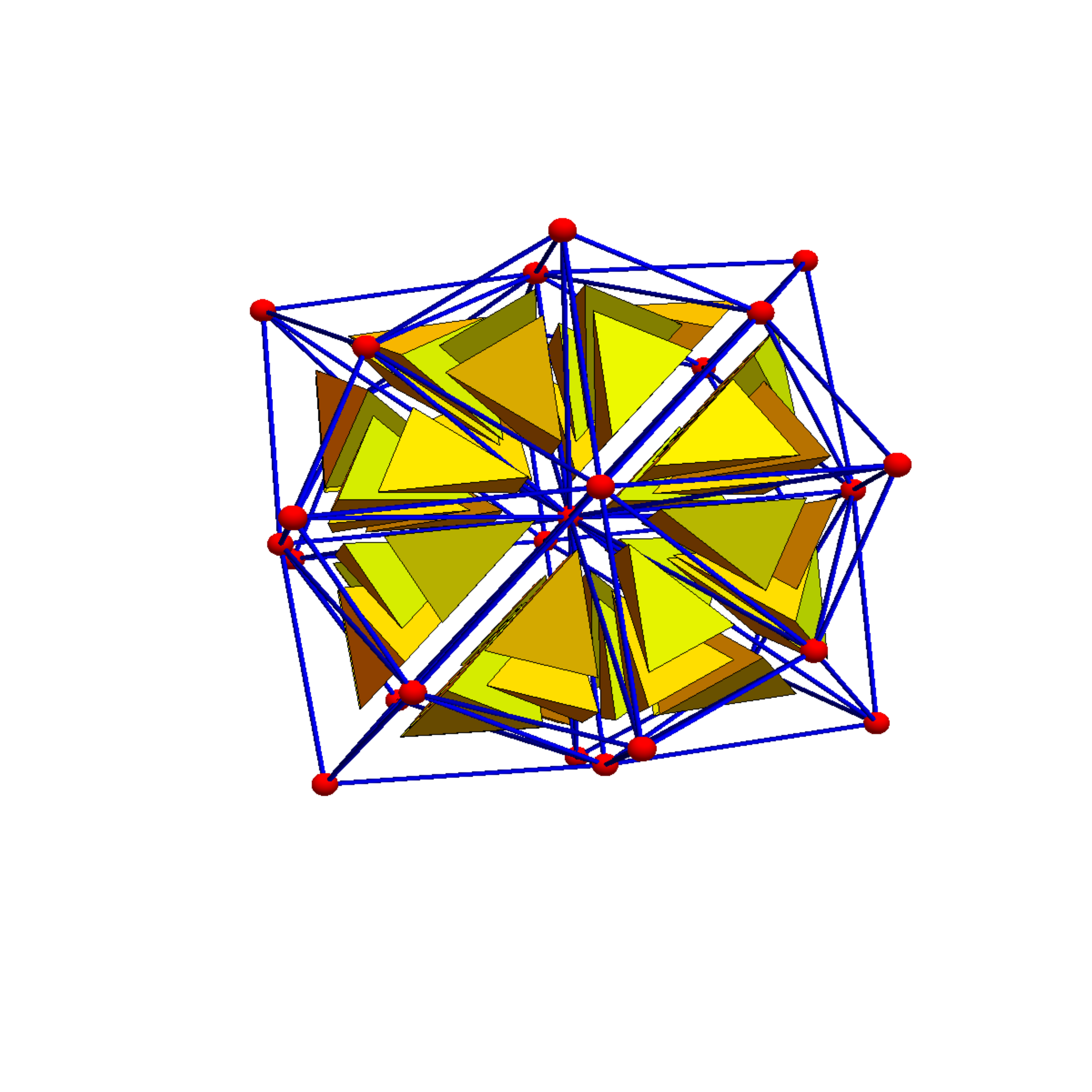}}
\caption{
A planar Archimedian polyhedron which is not in $\Scal_2$ as it is not maximal. It is a discrete surface
with boundary. A vector field $F$ is a color-valued function on $E$.
The line integral along the boundary of $S$ is equal to the total curl of $F$, summing over all
triangles. As $G$ is not simply connected, there are many fields on $G$ which are not gradients. 
To the right we see an illustration of the divergence theorem. Adding up the divergences over all 
tetrahedra gives the flux integral along the boundary. 
See also \cite{ArchimedesFunctions}. }
\end{figure}

\begin{figure}[h]
\parbox{6.2cm}{\scalebox{0.22}{\includegraphics{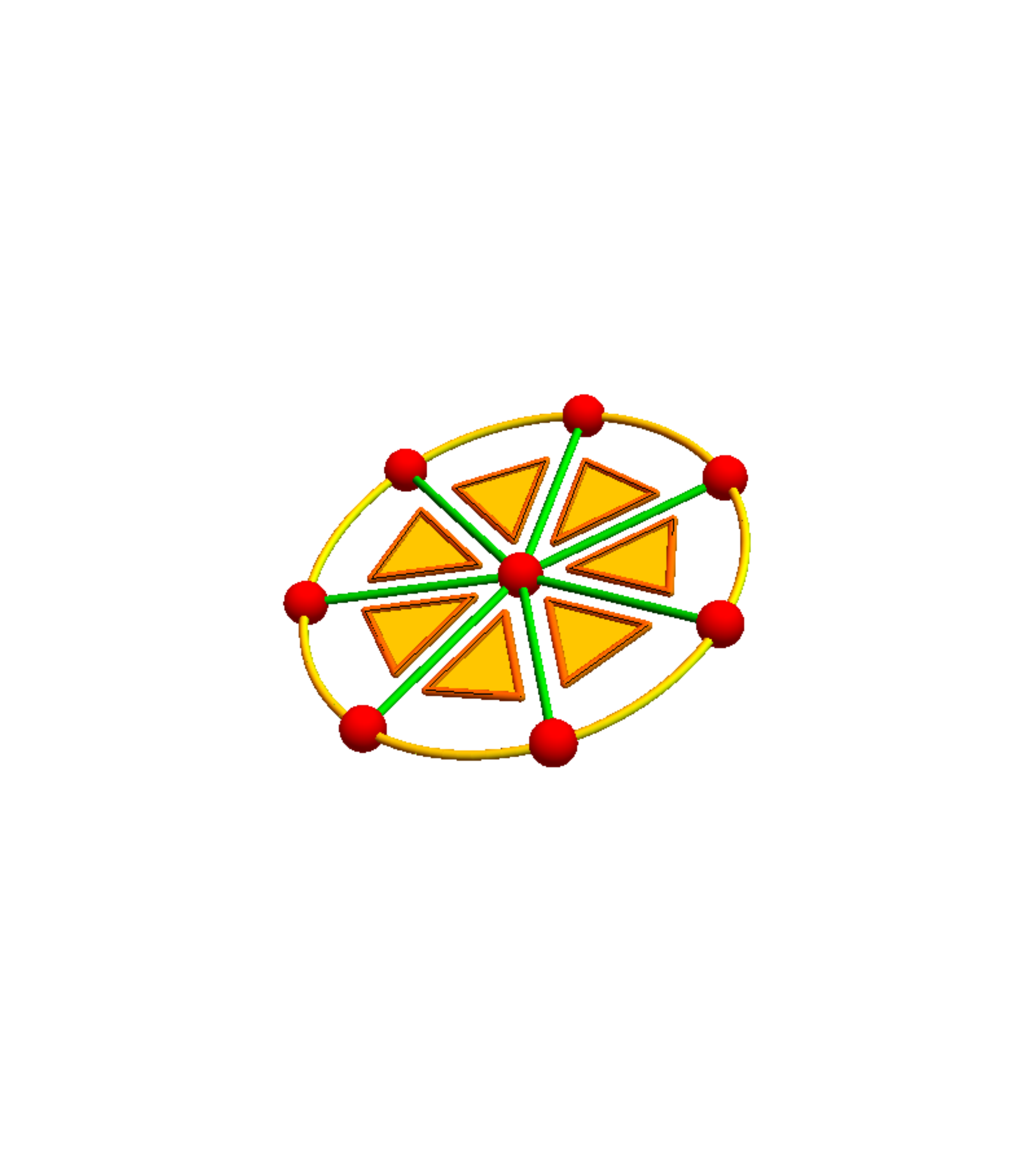}}}
\parbox{6.2cm}{\scalebox{0.22}{\includegraphics{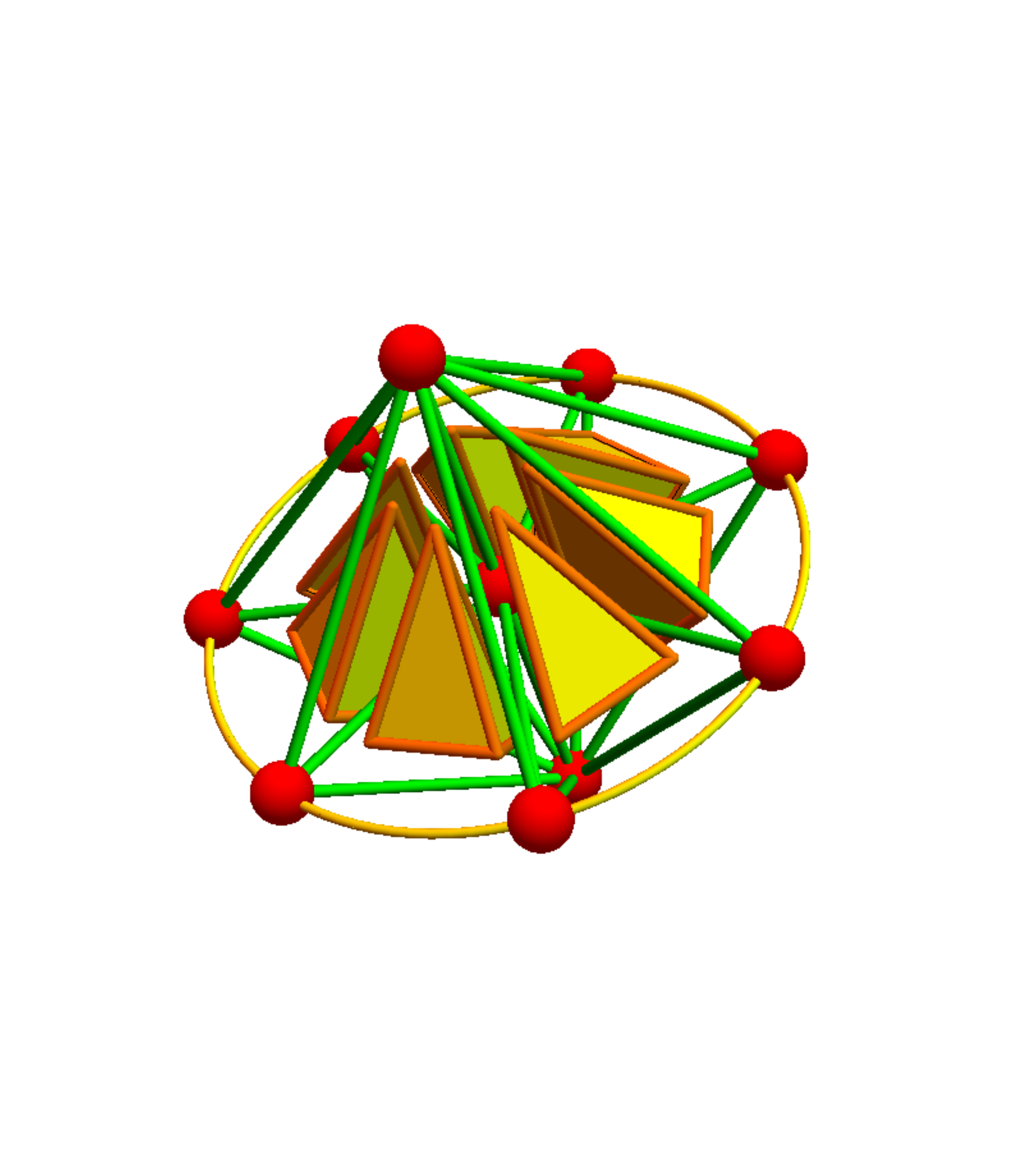}}}
\caption{
The vertex degree of a vertex of $G \in \Gcal_2$ and the 
edge degree of an edge of $G \in \Gcal_3$.
}
\end{figure}

\begin{figure}[h]
\parbox{6.2cm}{\scalebox{0.2}{\includegraphics{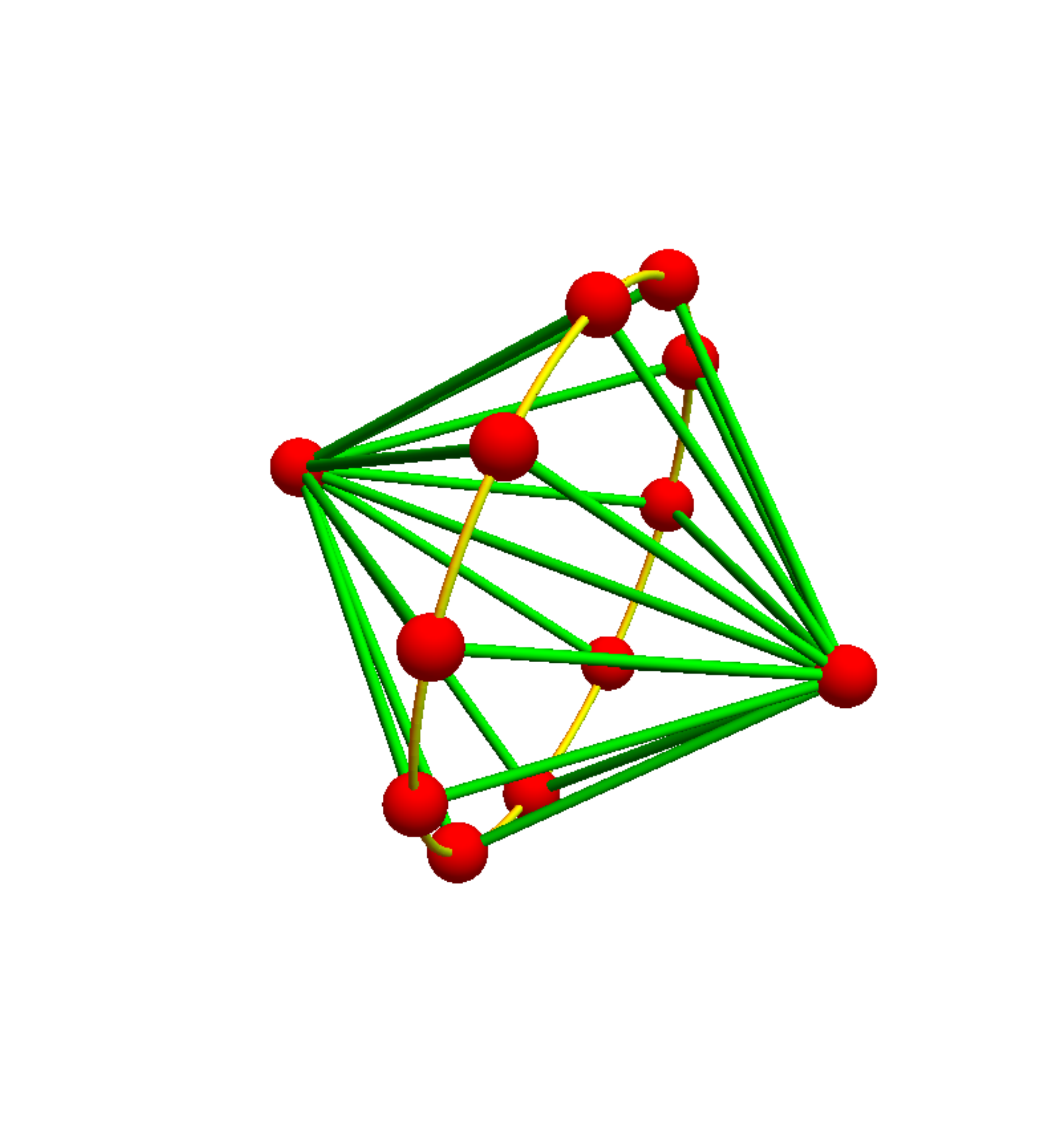}} }
\parbox{6.2cm}{\scalebox{0.2}{\includegraphics{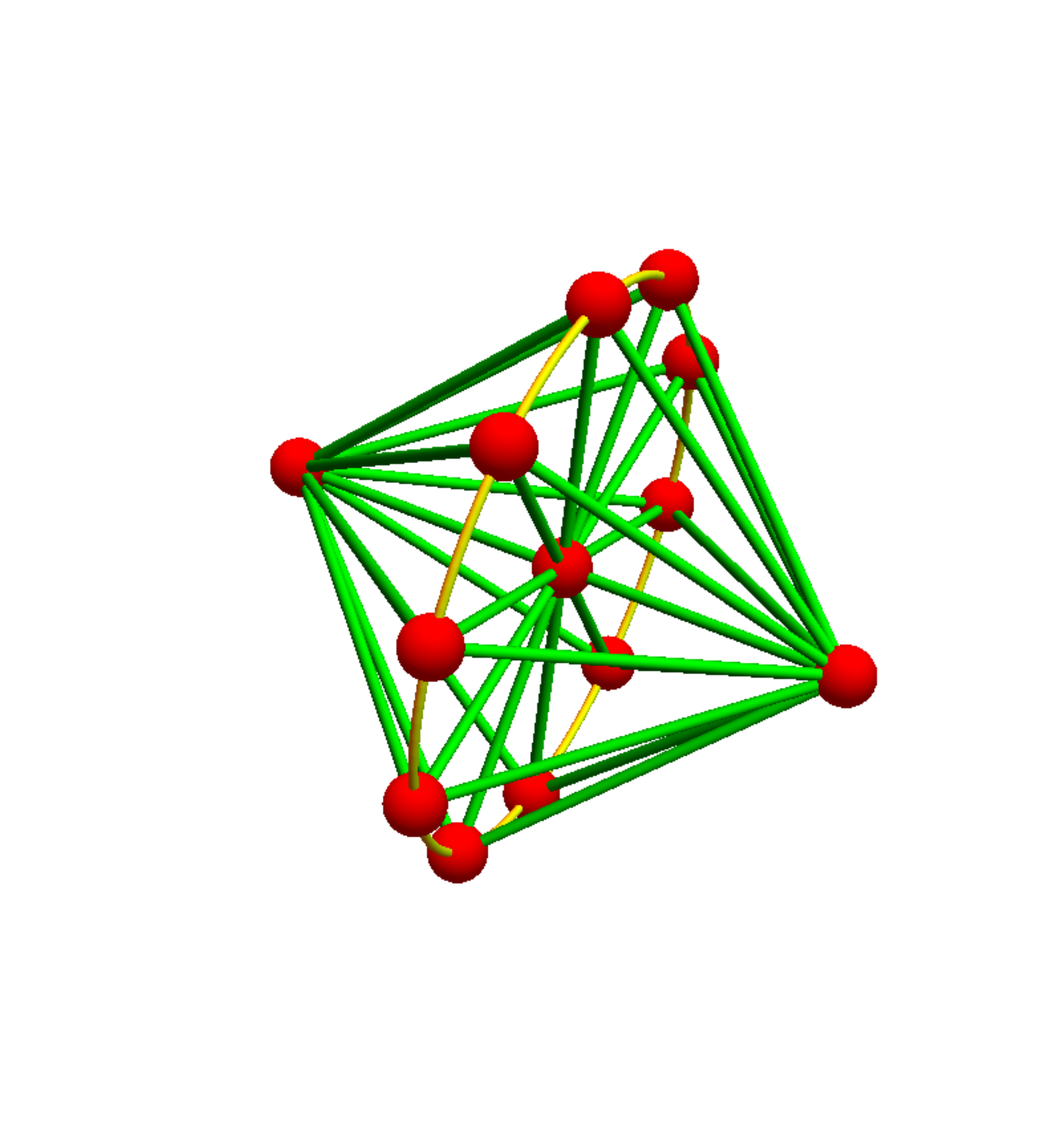}} }
\caption{
Cutting an edge $e=(a,b)$ in a $3$-dimensional graph doubles the number of tetrahedra
attached to $e$. It adds a surface in the form of a wheel graph and 
changes the degree of each edge in the circular graph $S(a) \cap S(b)$ by $1$. 
}
\end{figure}

\begin{figure}[h]
\parbox{8.2cm}{\scalebox{0.2}{\includegraphics{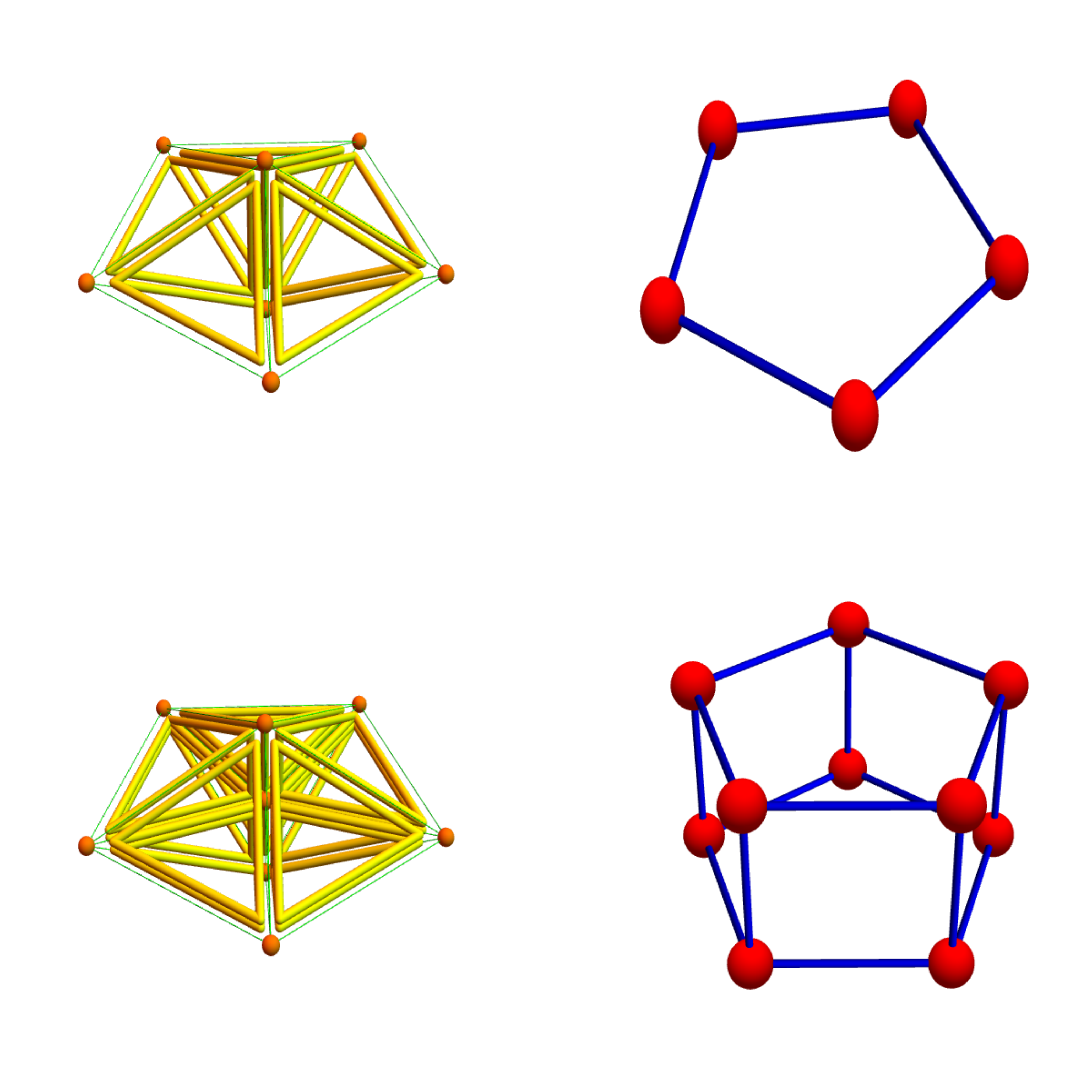}} }
\caption{
The procedure of cutting an  edge has the effect that in the dual graph a cycle becomes dihedral.
}
\end{figure}

\begin{figure}[h]
\scalebox{0.21}{\includegraphics{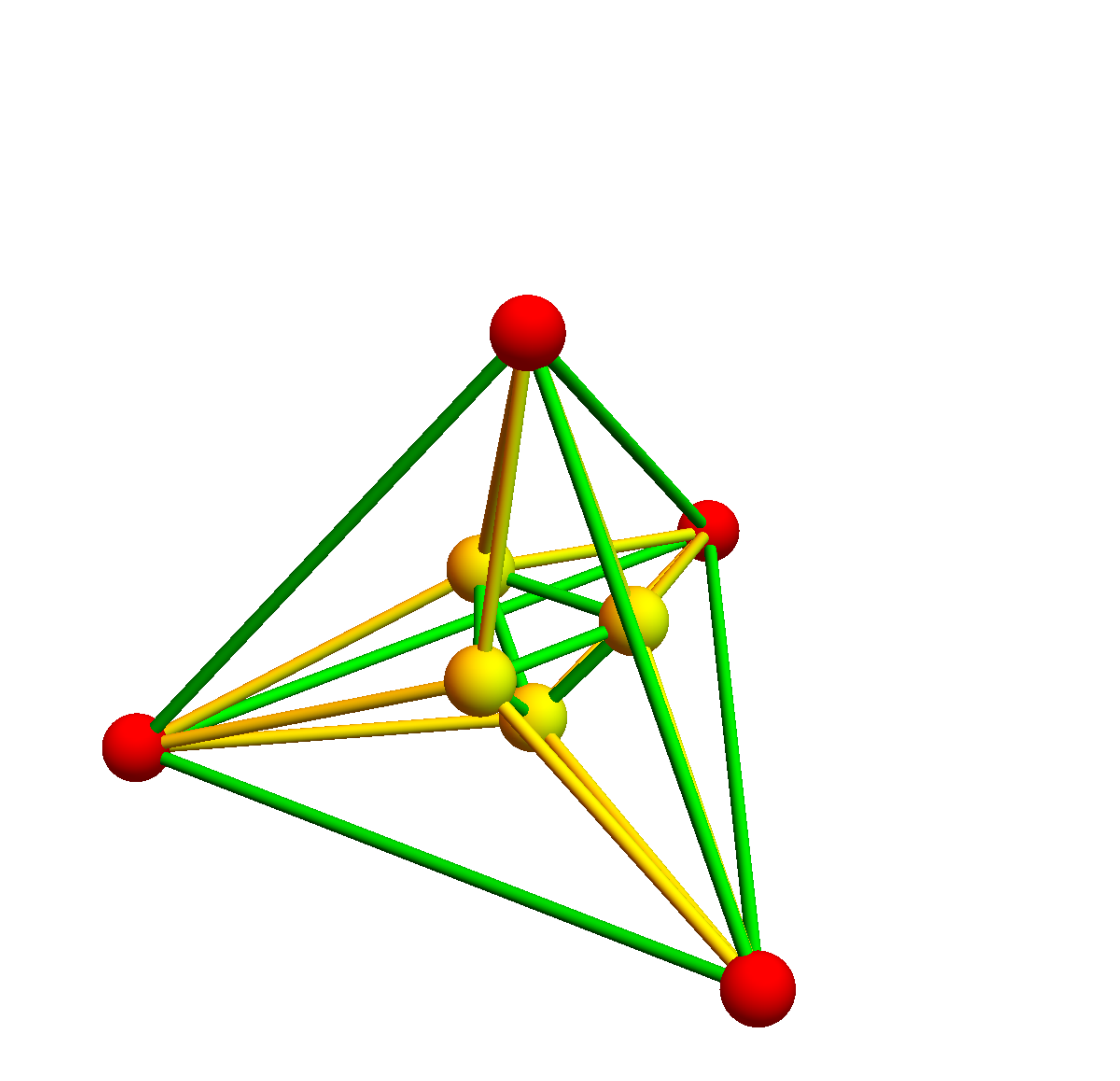}} 
\scalebox{0.21}{\includegraphics{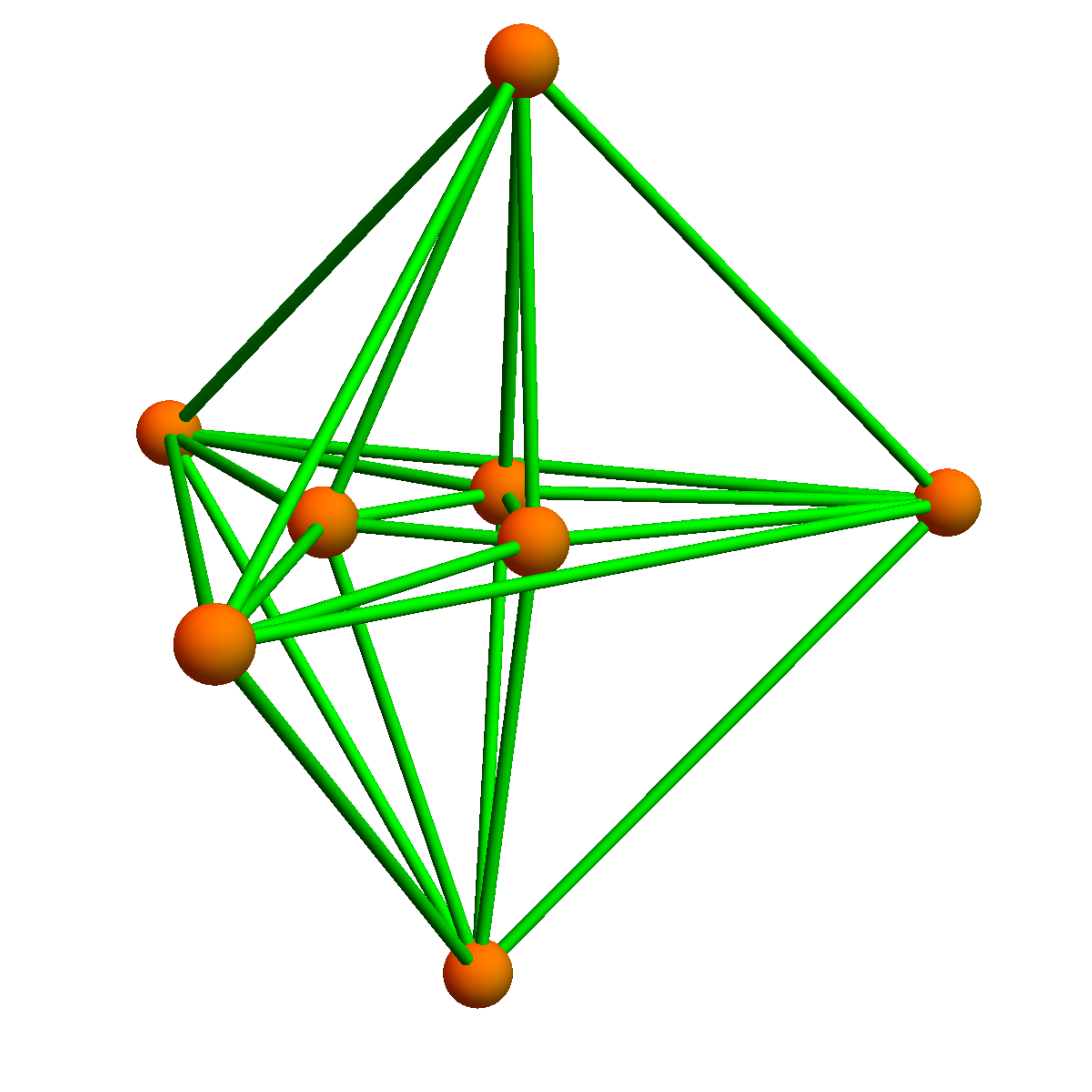}} 
\caption{
Appollonian network refinements which look geometric at first but are not valid. 
In the left case the unit sphere contains truncated tetrahedra which are not 
two-dimensional discrete spheres. In the right case, a two dimensional triangle
was refined and the rest cut up. The central sphere is not geometric. 
Also barycentric network refinements are not possible 
because we can not divide up the boundary. 
}
\end{figure}

\begin{figure}[h]
\scalebox{0.22}{\includegraphics{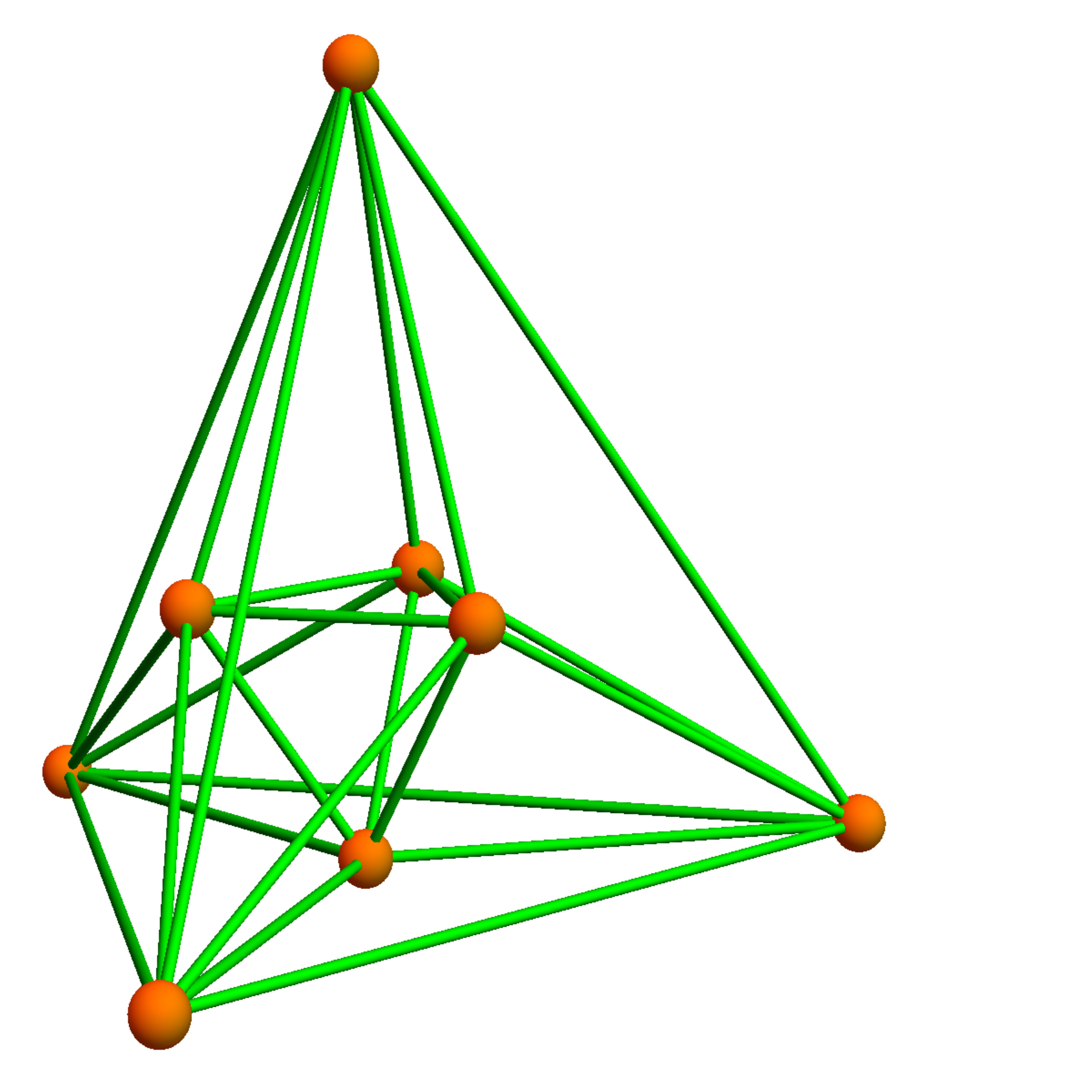}} 
\scalebox{0.22}{\includegraphics{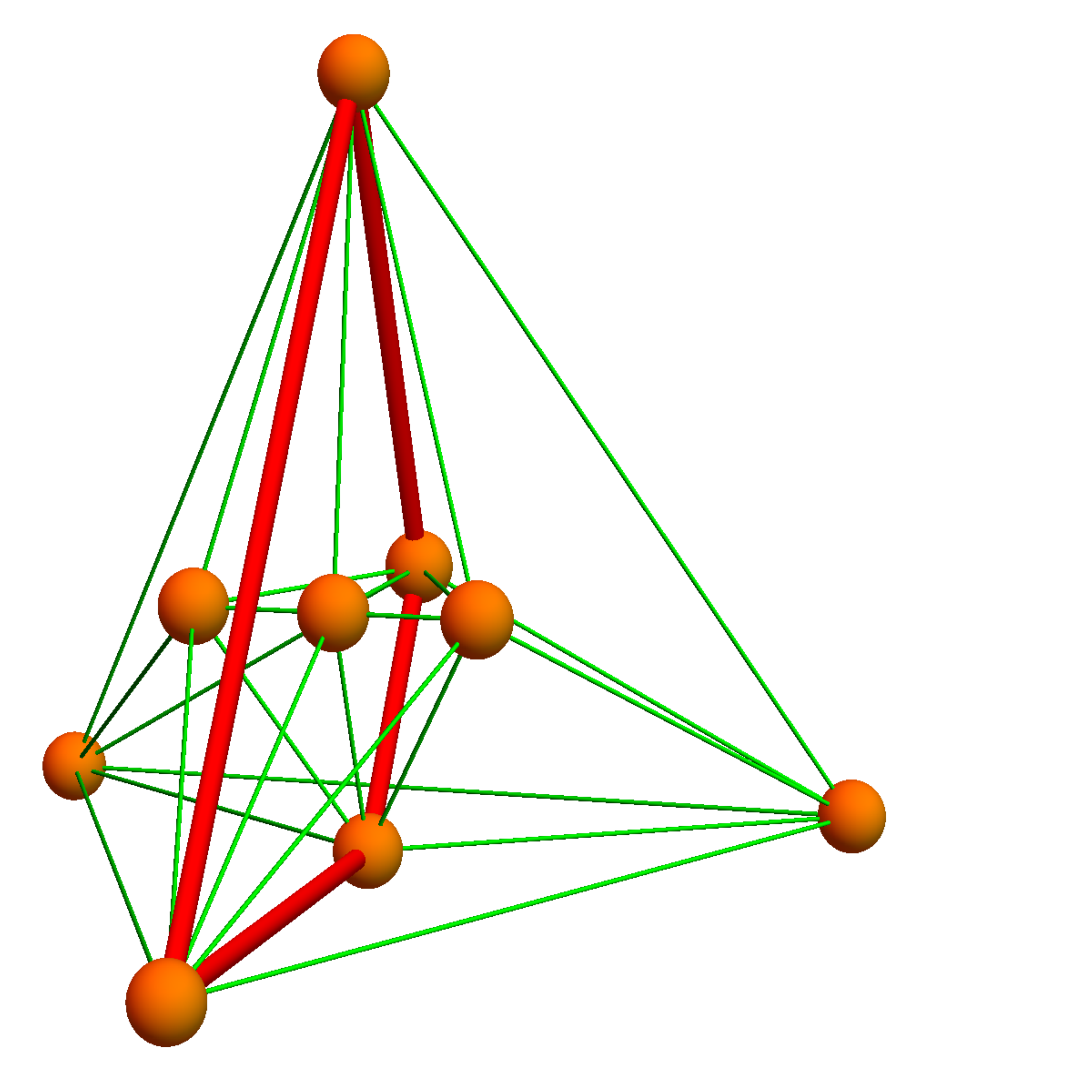}} 
\caption{
Replacing a tetrahedron by a 16-cell produces a refinement 
which does not affect the interior edge degrees. We see that graphs in 
$\Scal_3$ = three dimensional polytopes can be useful for very concrete things. 
In the figure below, an edge cut has been applied leading to a cycle of odd degree edges. 
}
\end{figure}

\begin{figure}[h]
\scalebox{0.24}{\includegraphics{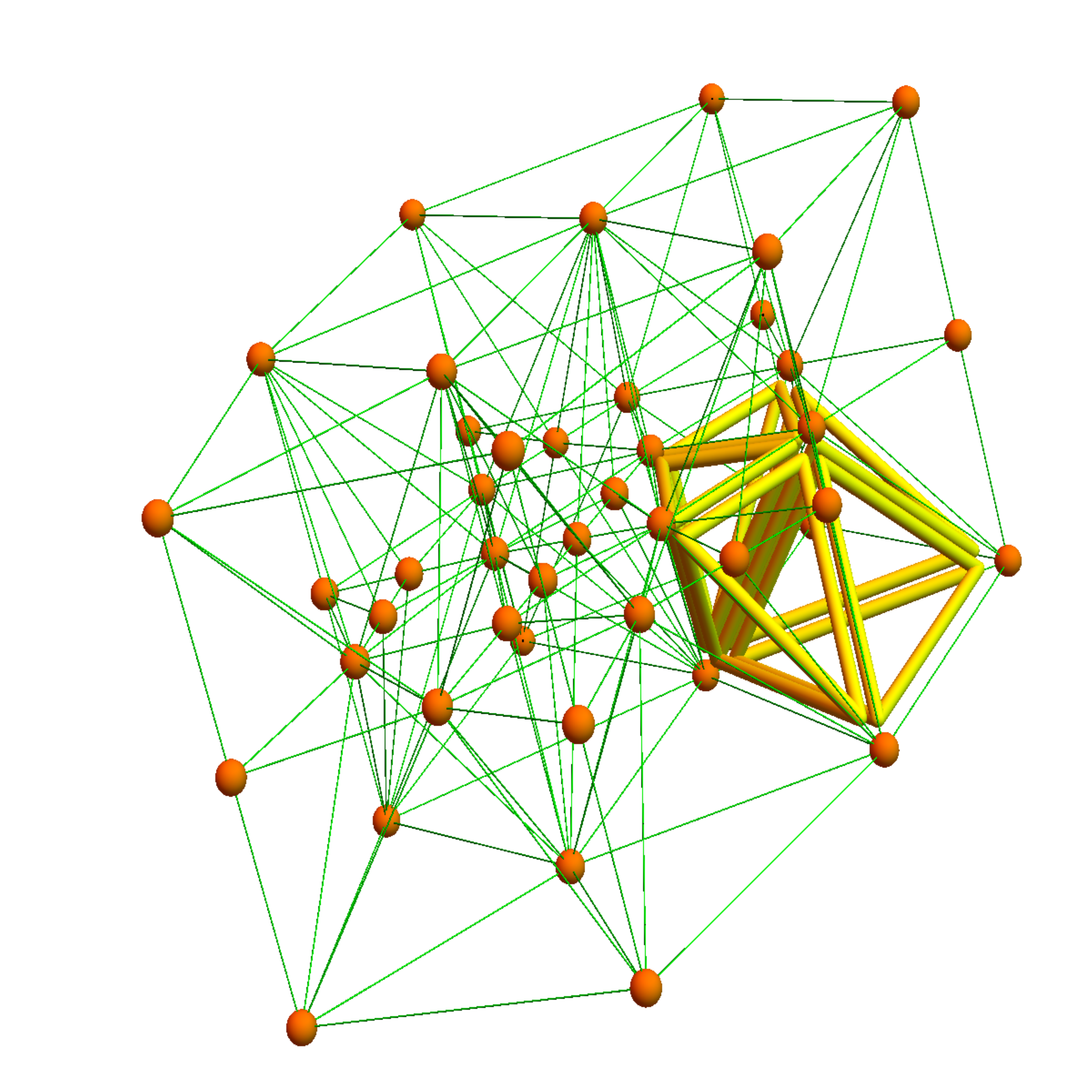}} 
\scalebox{0.24}{\includegraphics{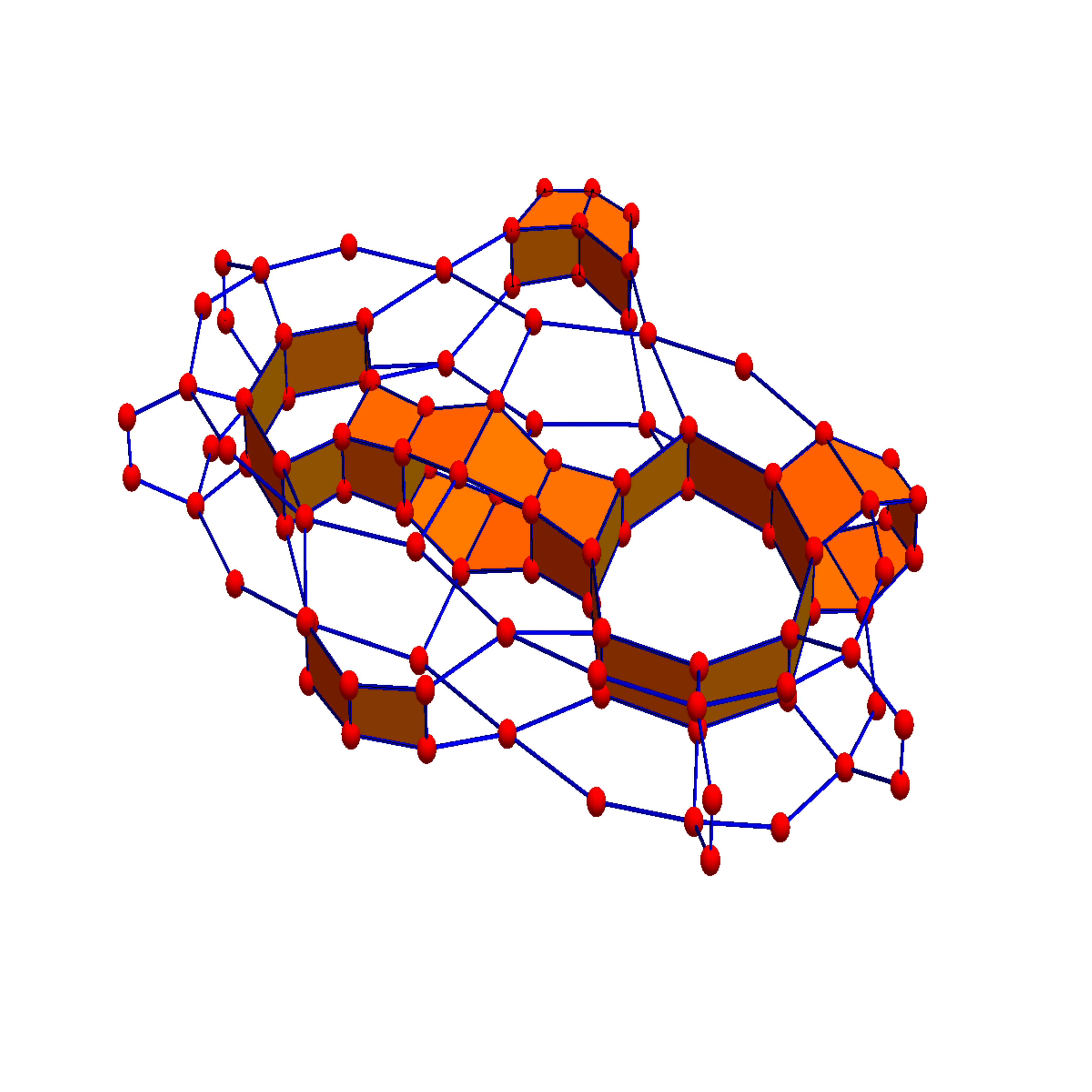}} 
\caption{
A graph in $\Gcal_3$ with an edge $e$ in the odd degree vertex set 
$O$. The edge $e$ connects two boundary vertices.
The dual graph $\hat{G}$ to the right shows there
a cycle of length $5$. }
\end{figure}

\begin{figure}[h]
\parbox{16.8cm}{
\parbox{7.2cm}{\scalebox{0.22}{\includegraphics{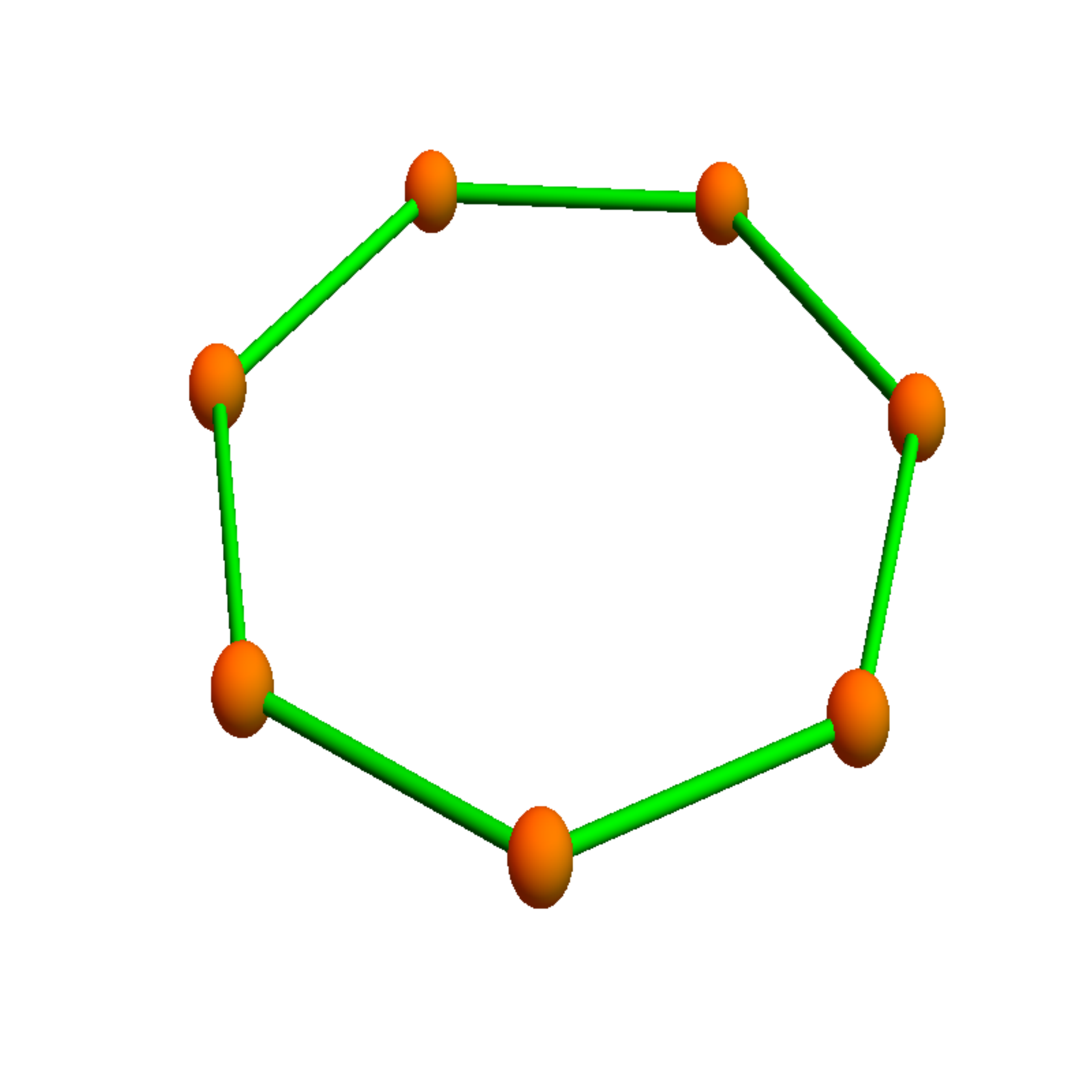}} }
\parbox{7.2cm}{\scalebox{0.22}{\includegraphics{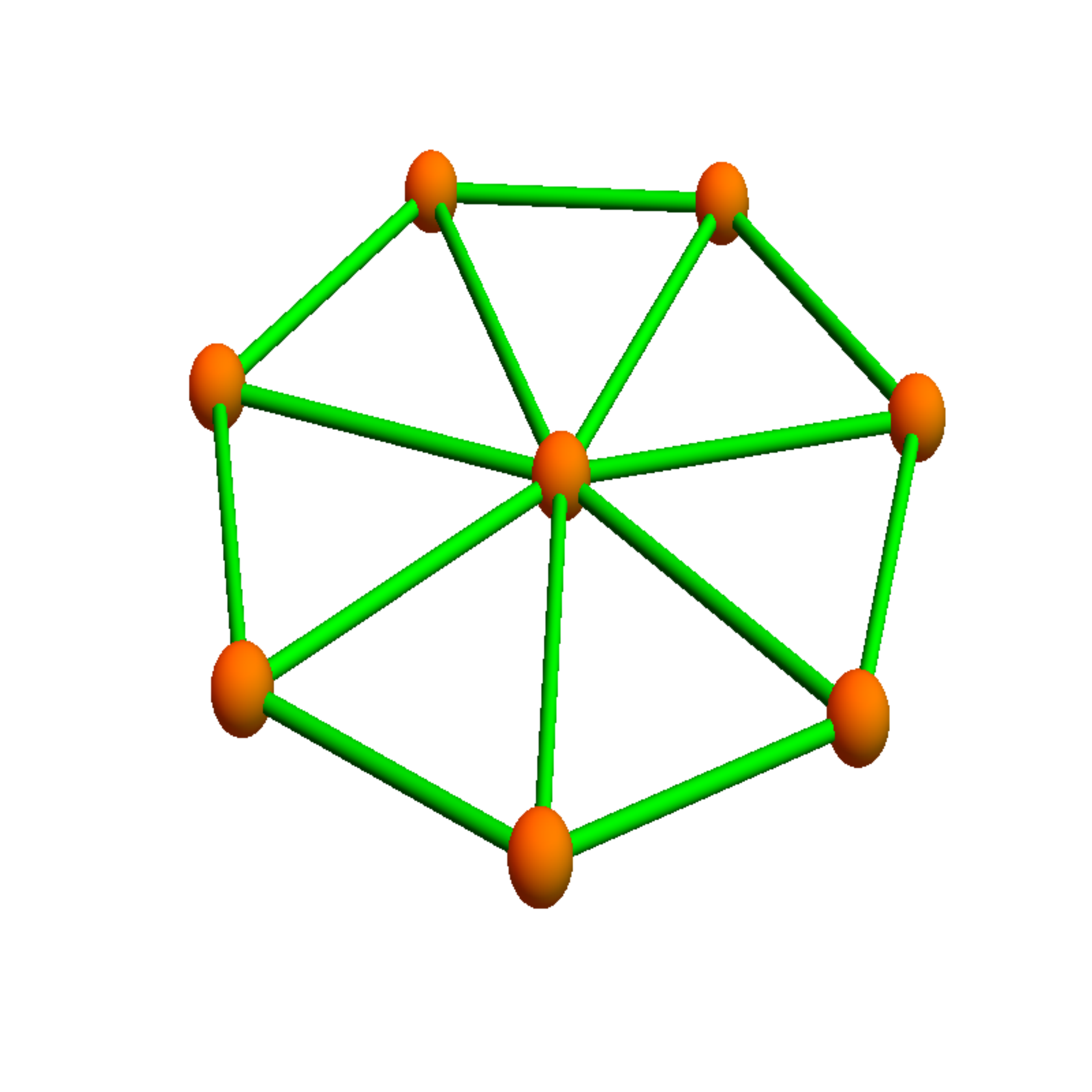}} }
}
\parbox{16.8cm}{
\parbox{7.2cm}{\scalebox{0.22}{\includegraphics{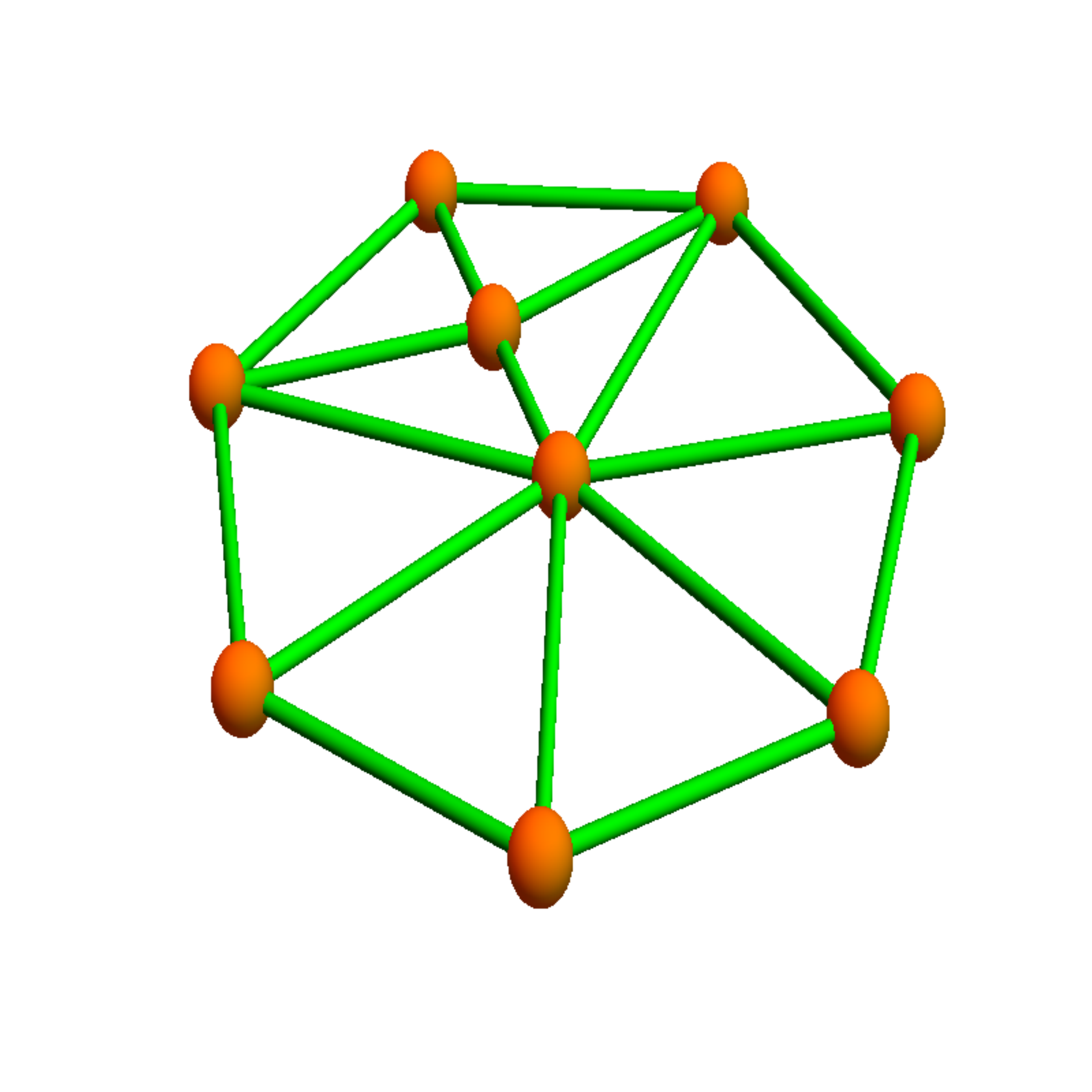}} }
\parbox{7.2cm}{\scalebox{0.22}{\includegraphics{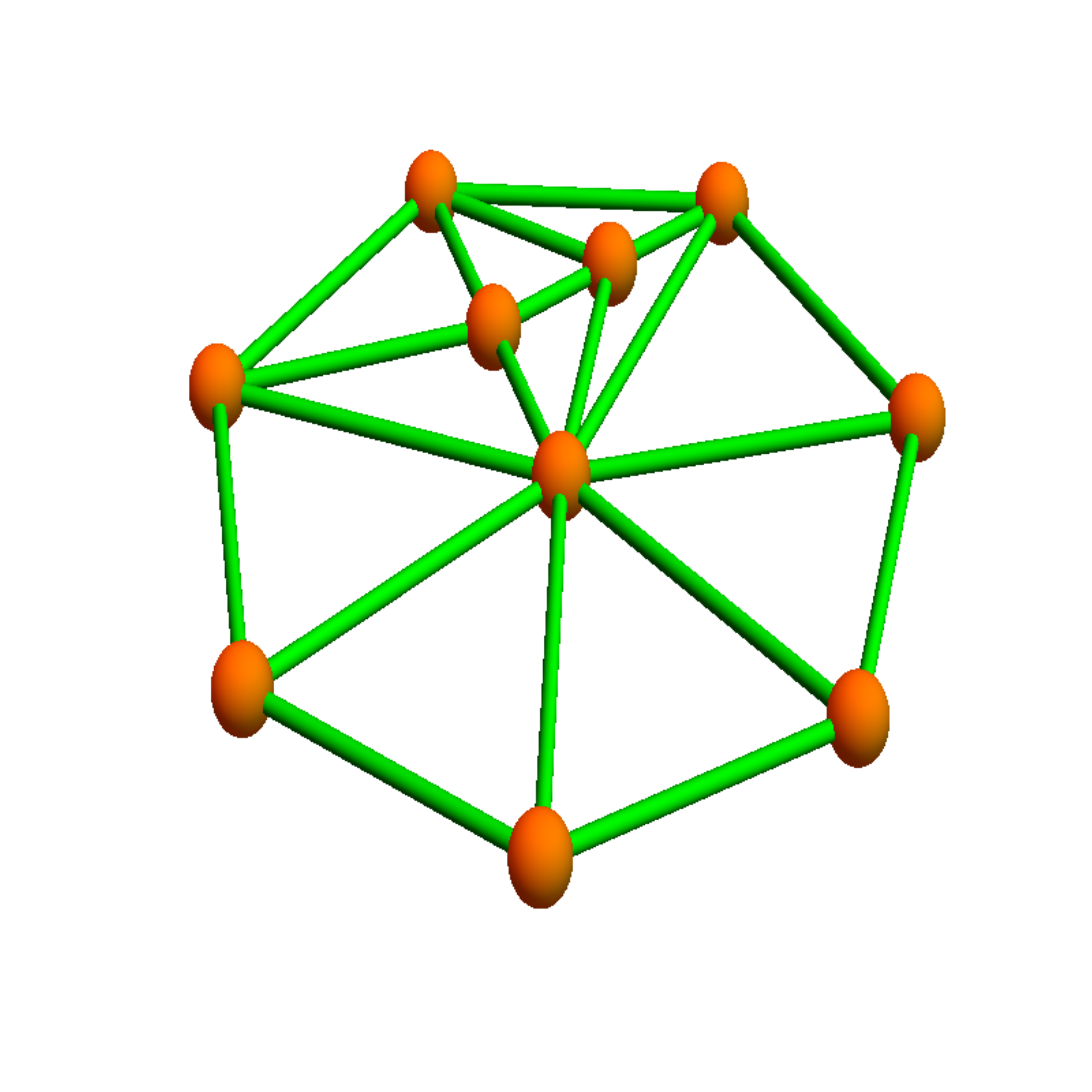}} }
}
\parbox{16.8cm}{
\parbox{7.2cm}{\scalebox{0.22}{\includegraphics{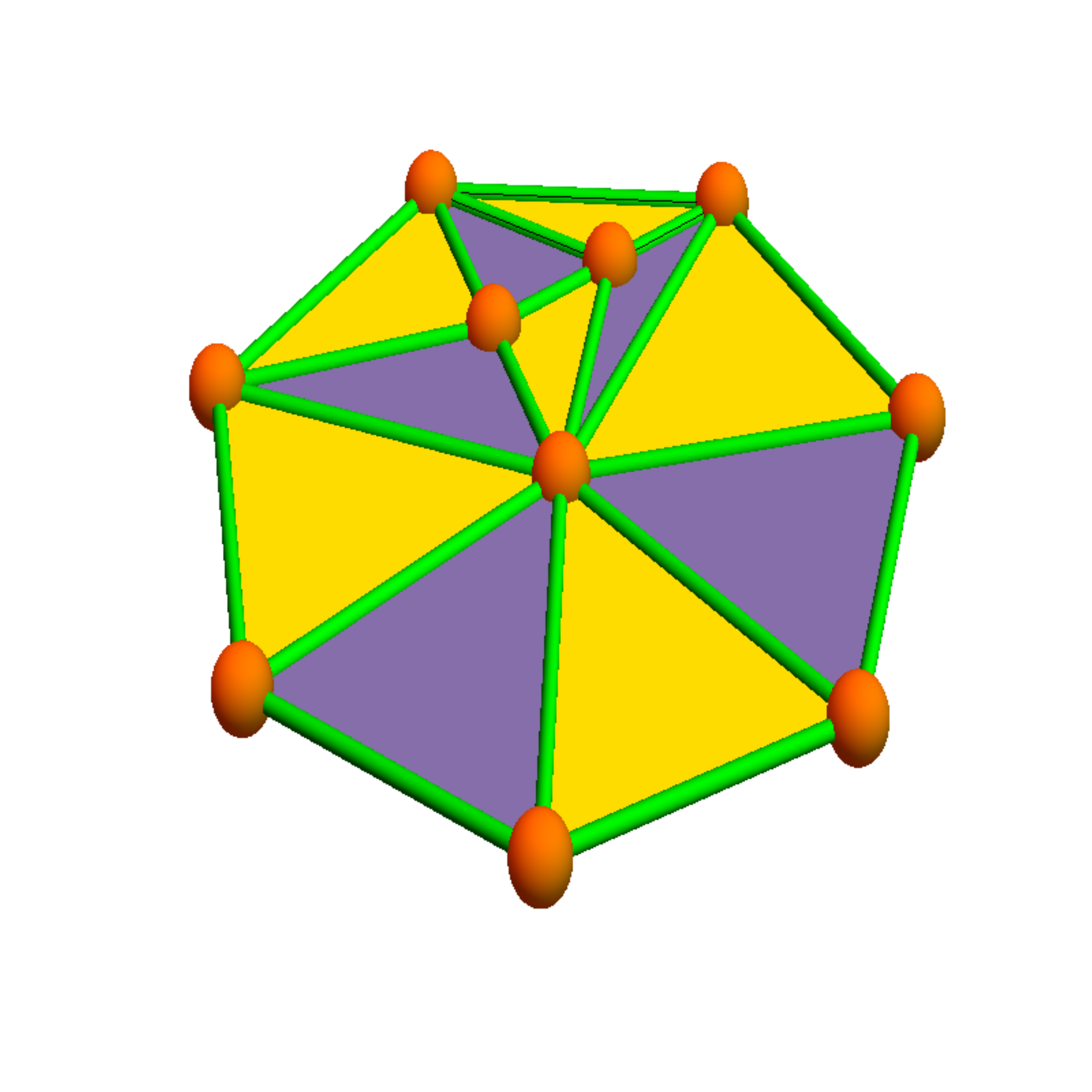}} }
\parbox{7.2cm}{\scalebox{0.22}{\includegraphics{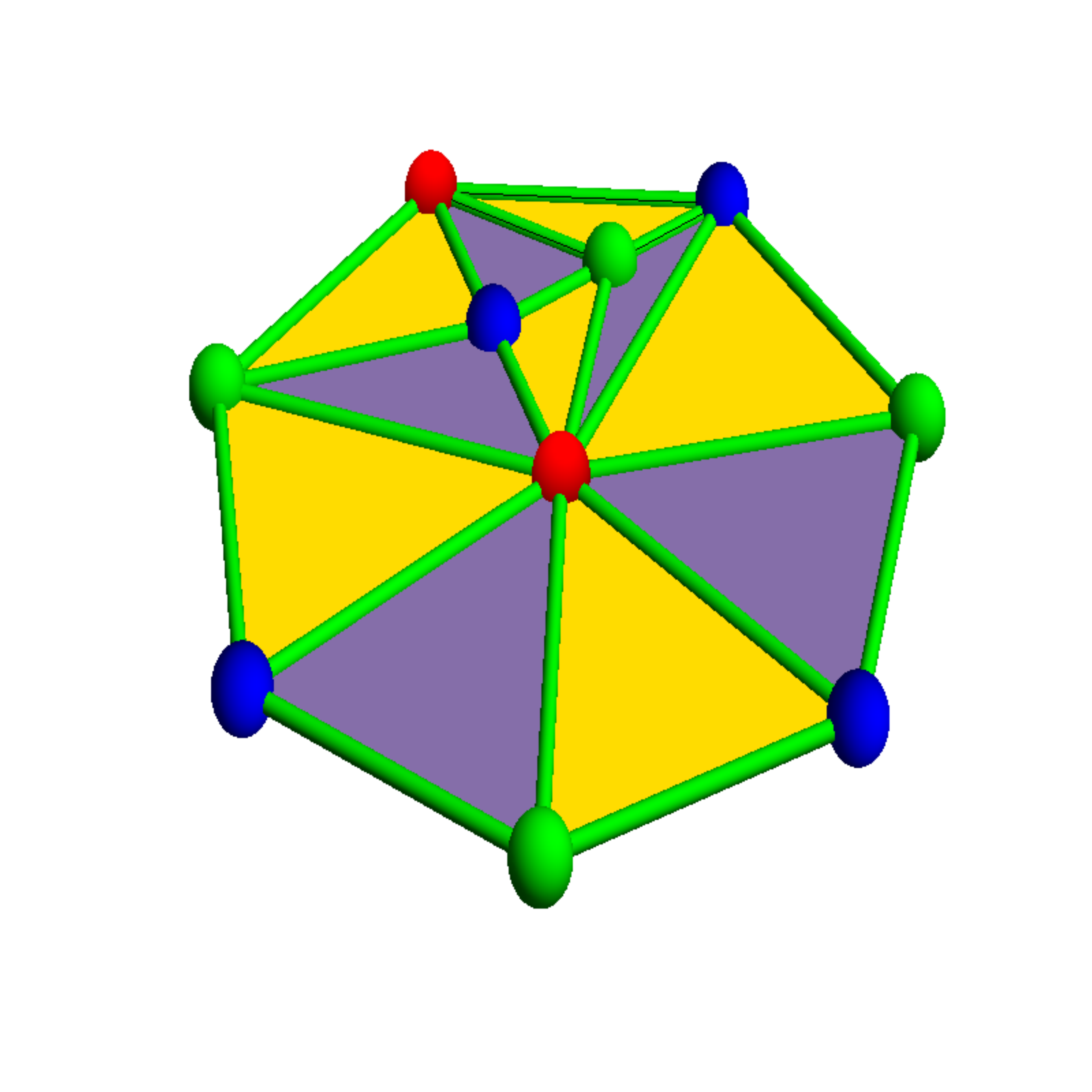}} }
}
\caption{
Coloring $C_n$ using two dimensions. For odd $n$,we refine the triangularization.
It appears silly to color the trivial case $C_n$ but it illustrates the main
idea introduced here. 
}
\end{figure}

\begin{figure}[h]
\parbox{16.8cm}{
\parbox{7.2cm}{\scalebox{0.17}{\includegraphics{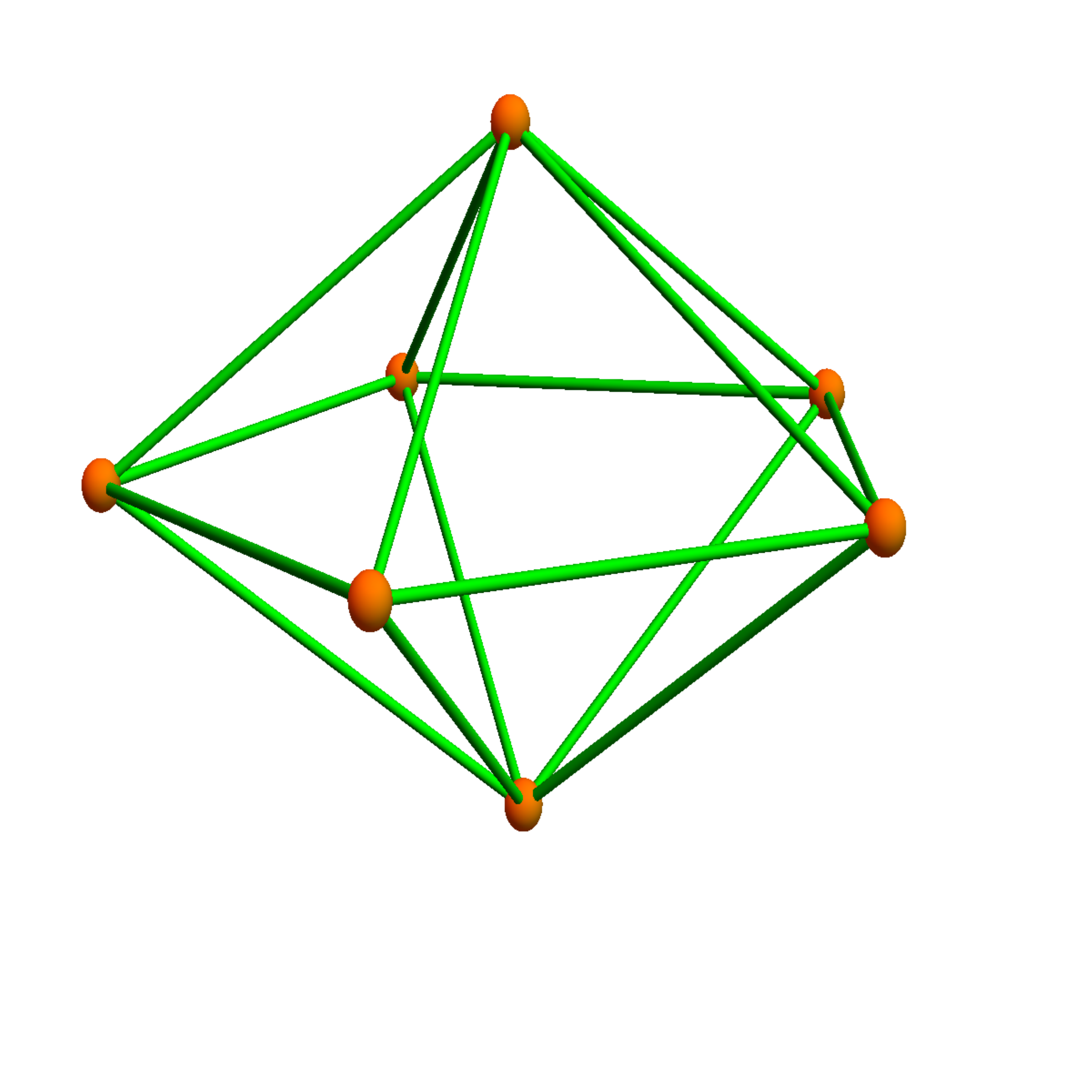}} }
\parbox{7.2cm}{\scalebox{0.17}{\includegraphics{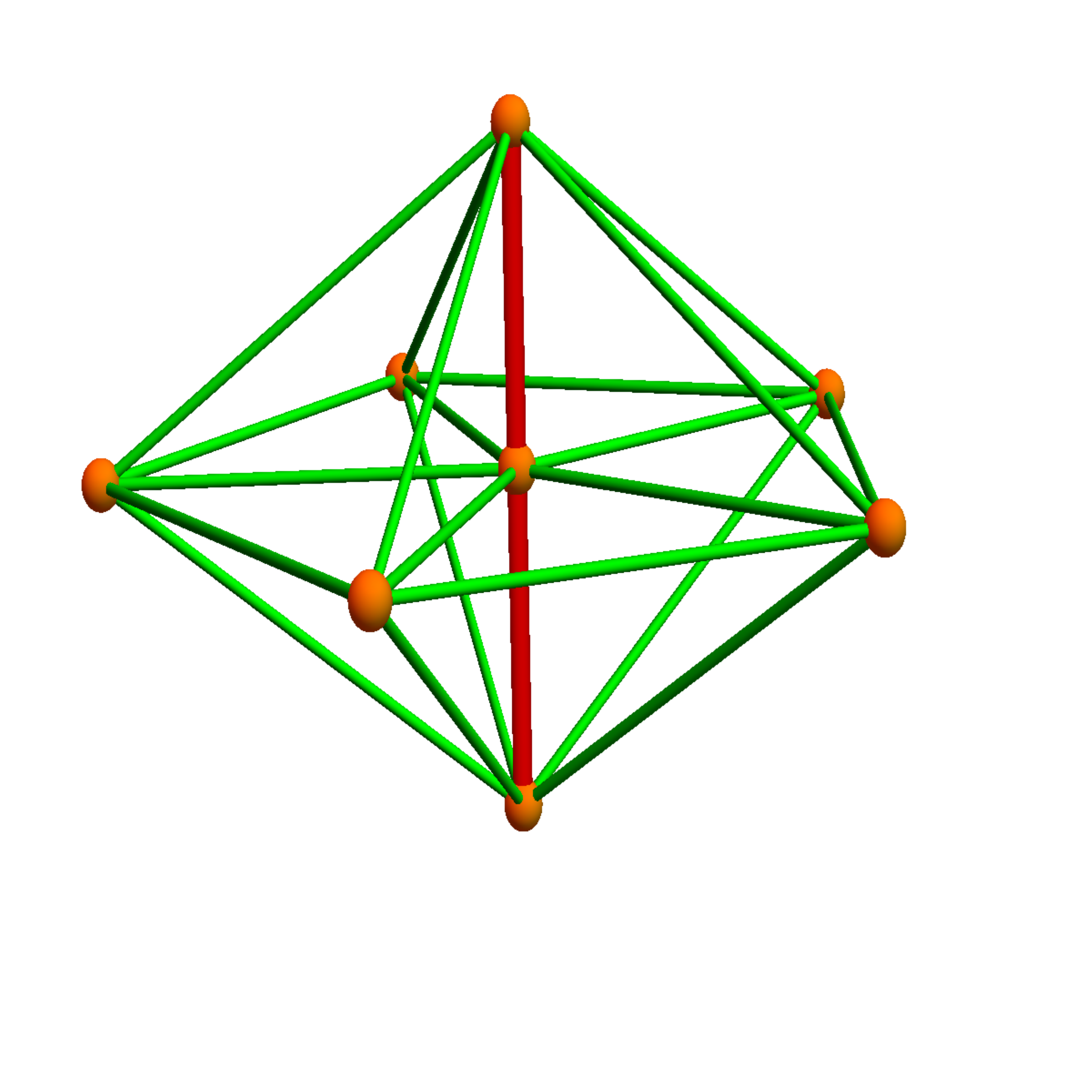}} }
}
\parbox{16.8cm}{
\parbox{7.2cm}{\scalebox{0.17}{\includegraphics{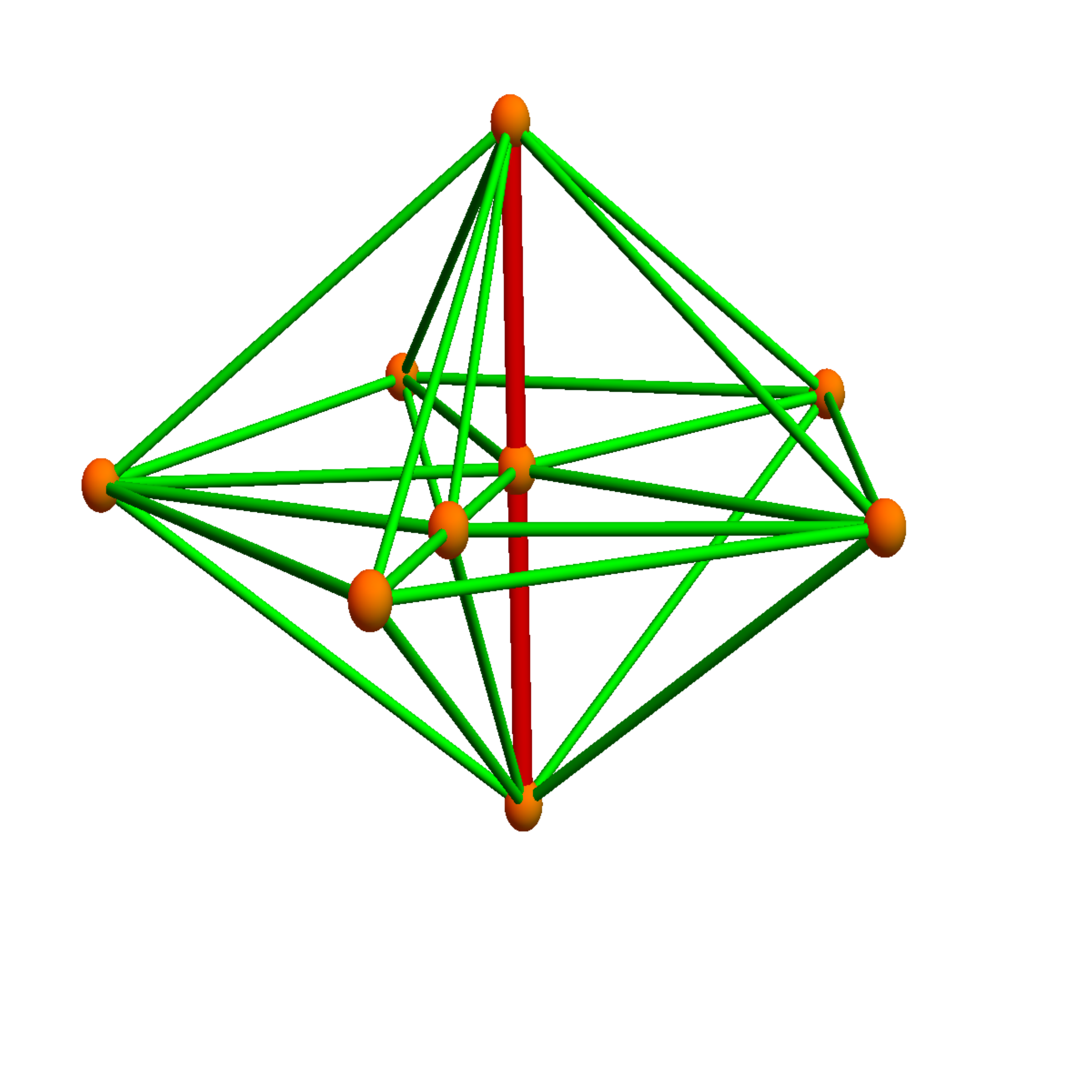}} }
\parbox{7.2cm}{\scalebox{0.17}{\includegraphics{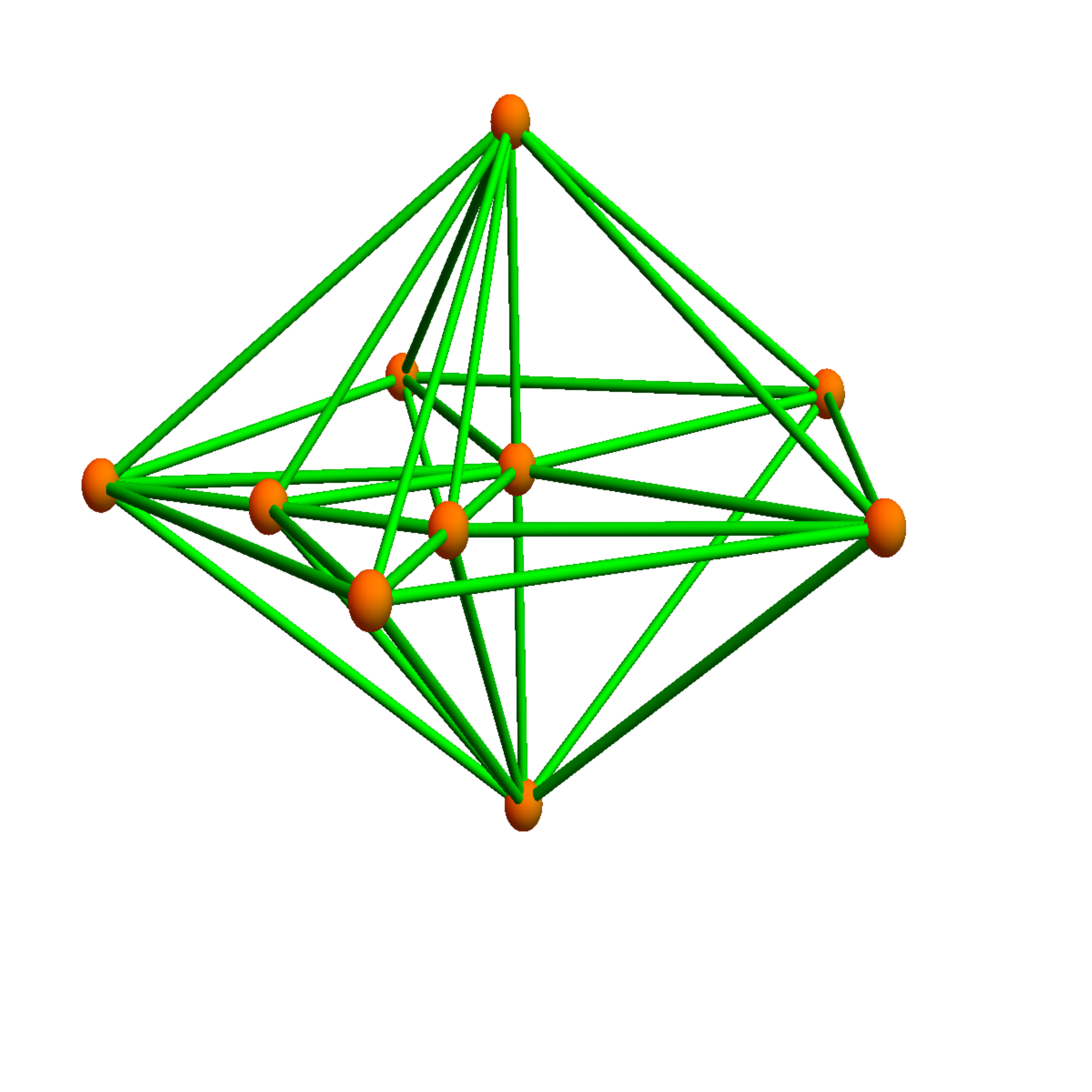}} }
}
\parbox{16.8cm}{
\parbox{7.2cm}{\scalebox{0.17}{\includegraphics{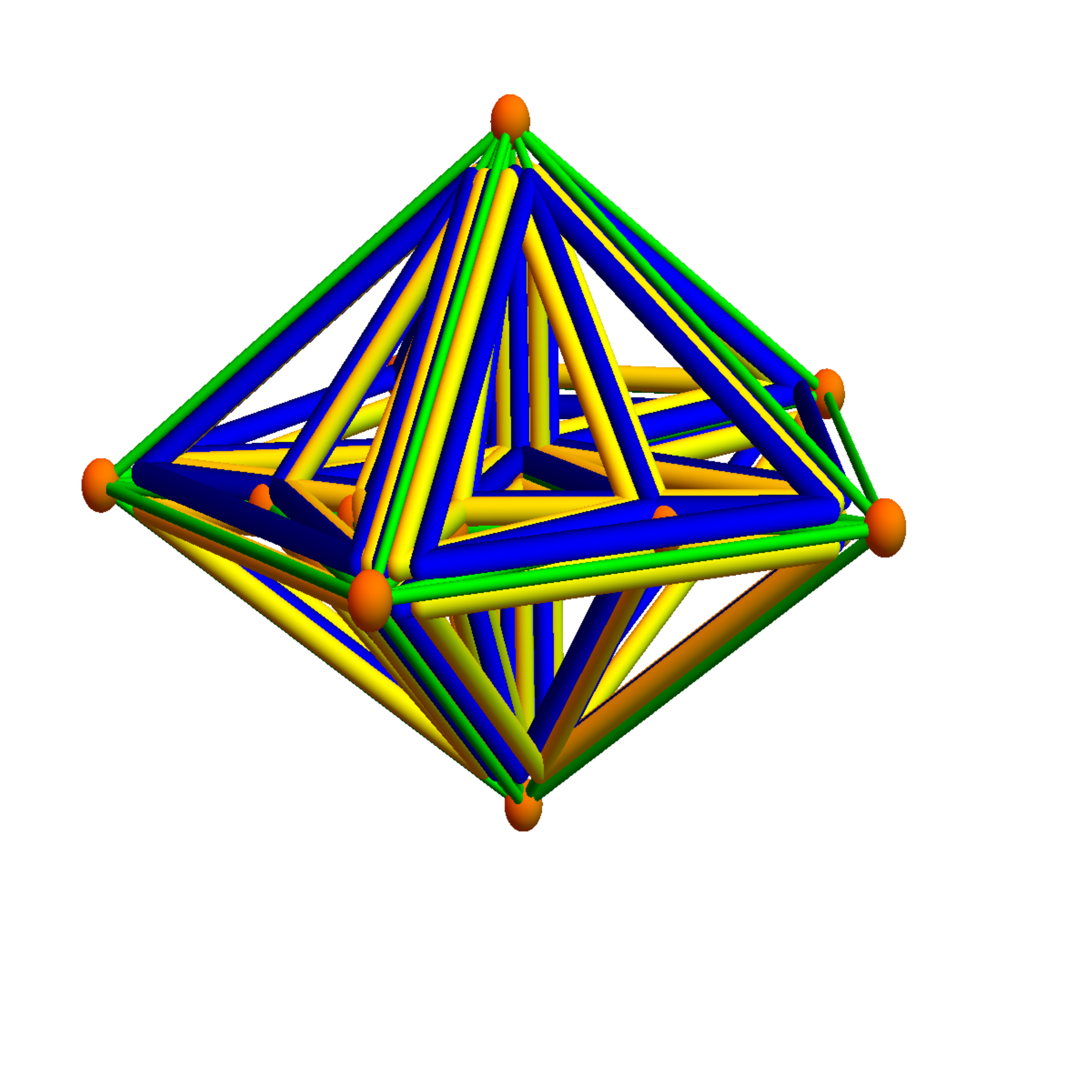}} }
\parbox{7.2cm}{\scalebox{0.17}{\includegraphics{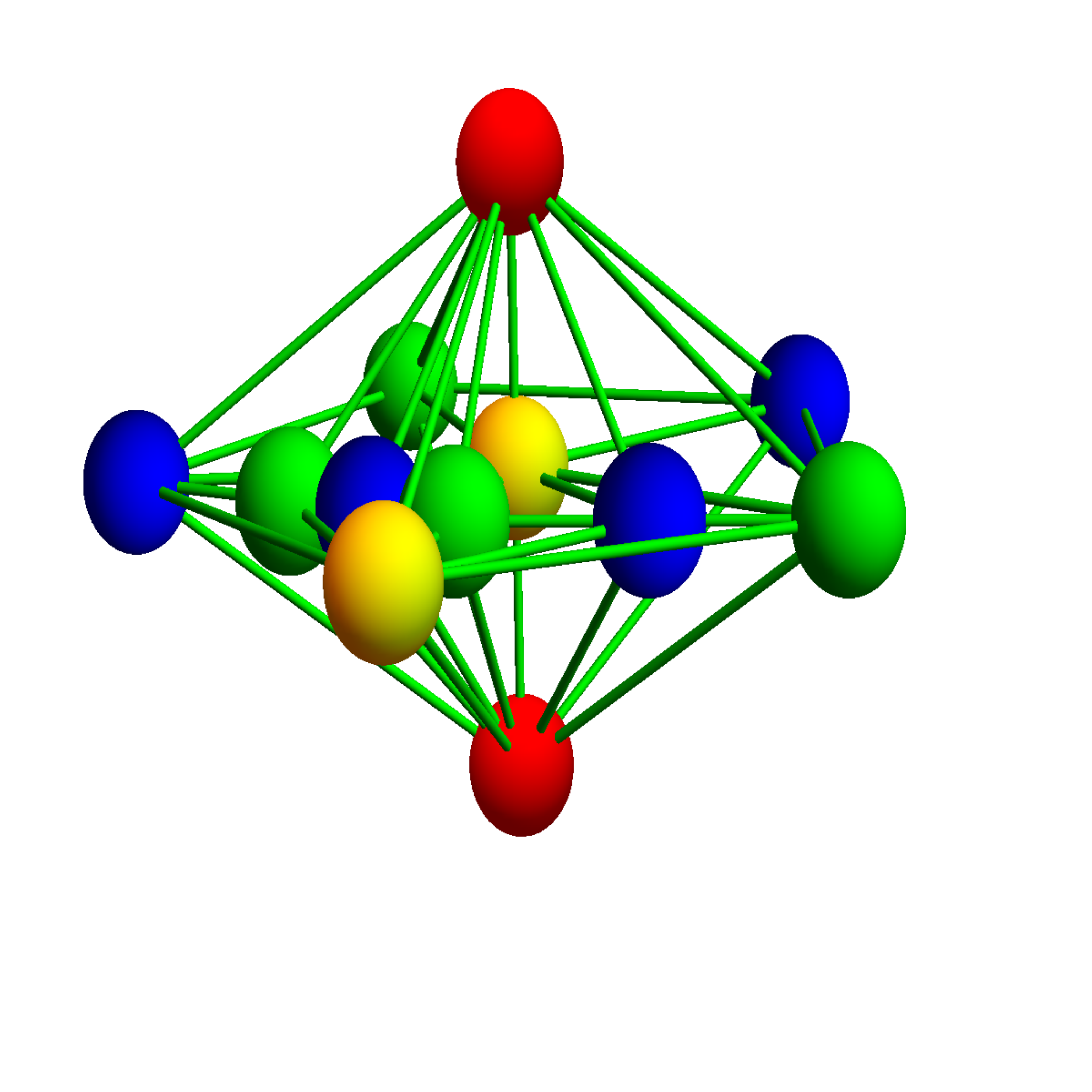}} }
}
\caption{
Coloring $G \in \Scal_2$ by seeing it as a boundary of a ball. 
}
\end{figure}

\begin{figure}[h]
\scalebox{0.22}{\includegraphics{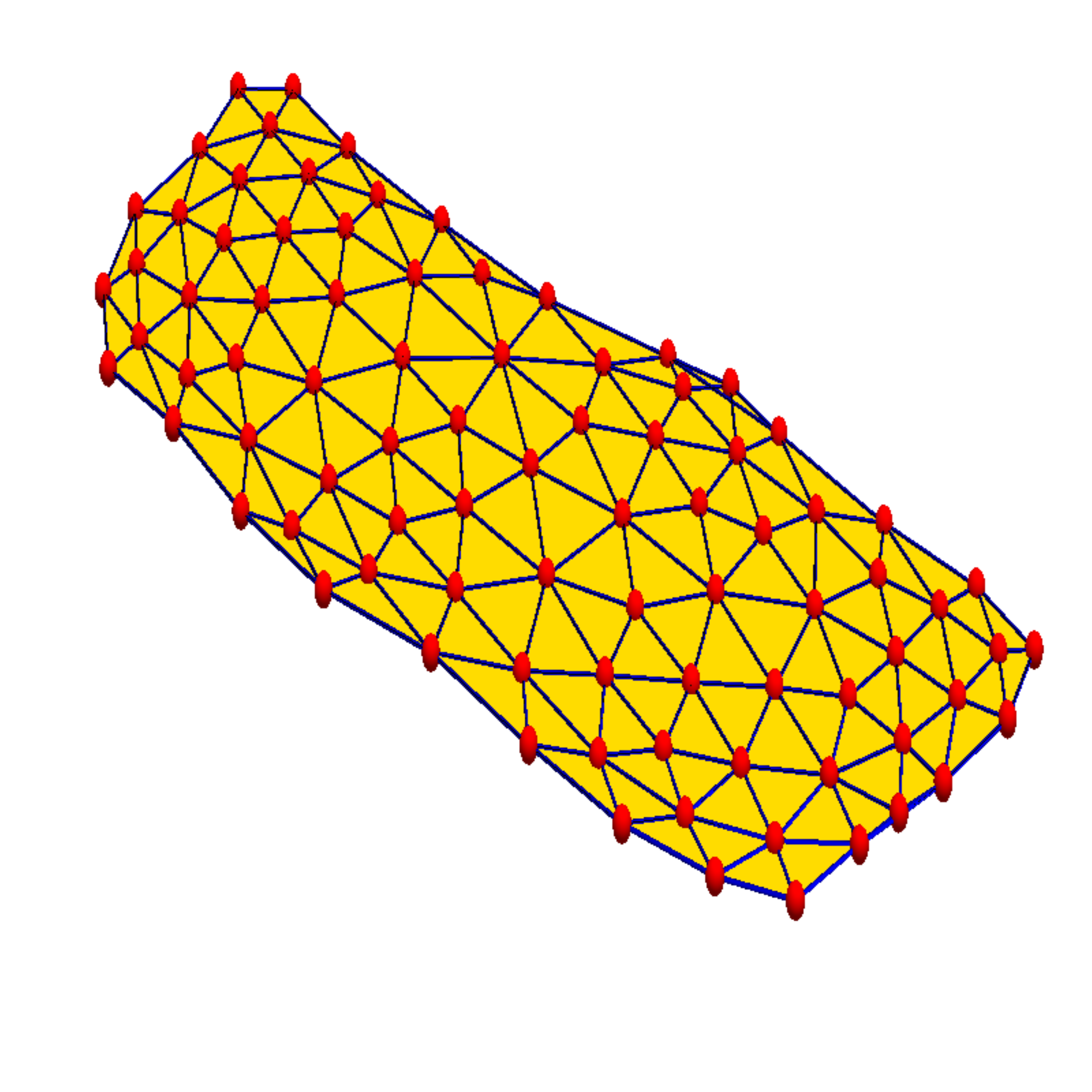}}
\scalebox{0.22}{\includegraphics{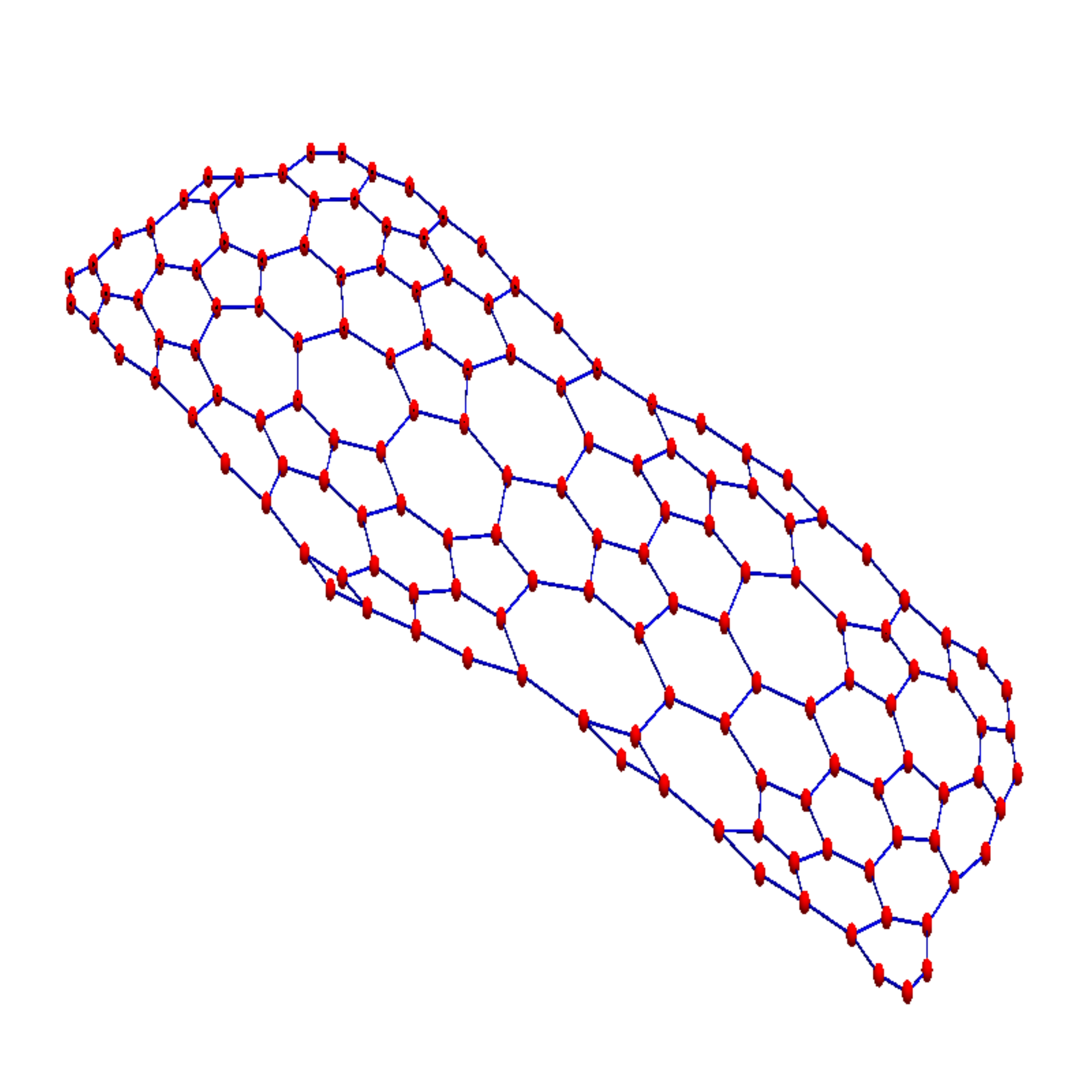}}
\caption{
A graph in $G \in \Bcal_2$ and its  dual graph $\hat{G}$
whose vertices are the triangles of $G$.
\label{orientable}
}
\end{figure}

\begin{figure}[h]
\parbox{15cm}{
\parbox{7.2cm}{\scalebox{0.15}{\includegraphics{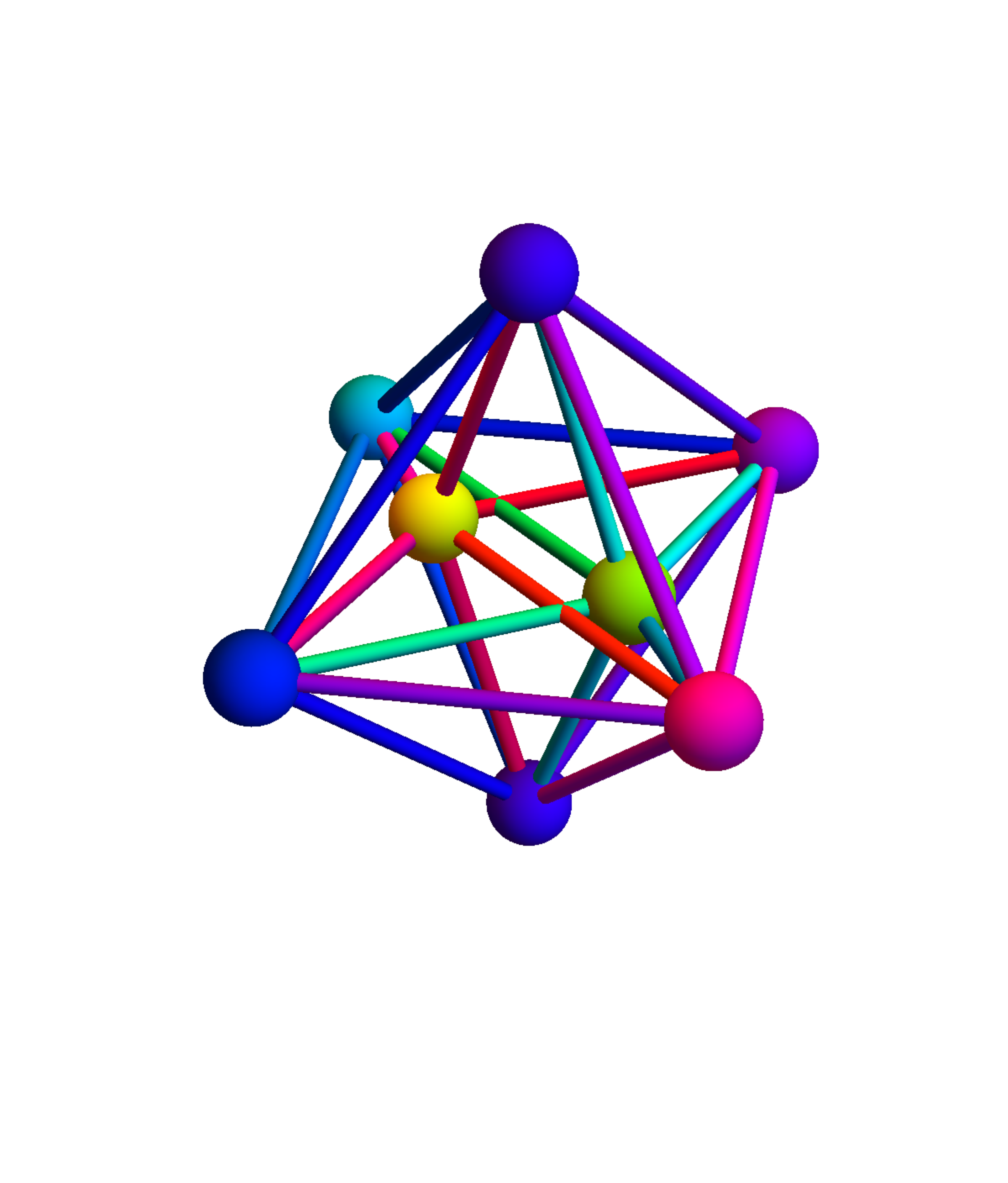}}}
\parbox{7.2cm}{\scalebox{0.15}{\includegraphics{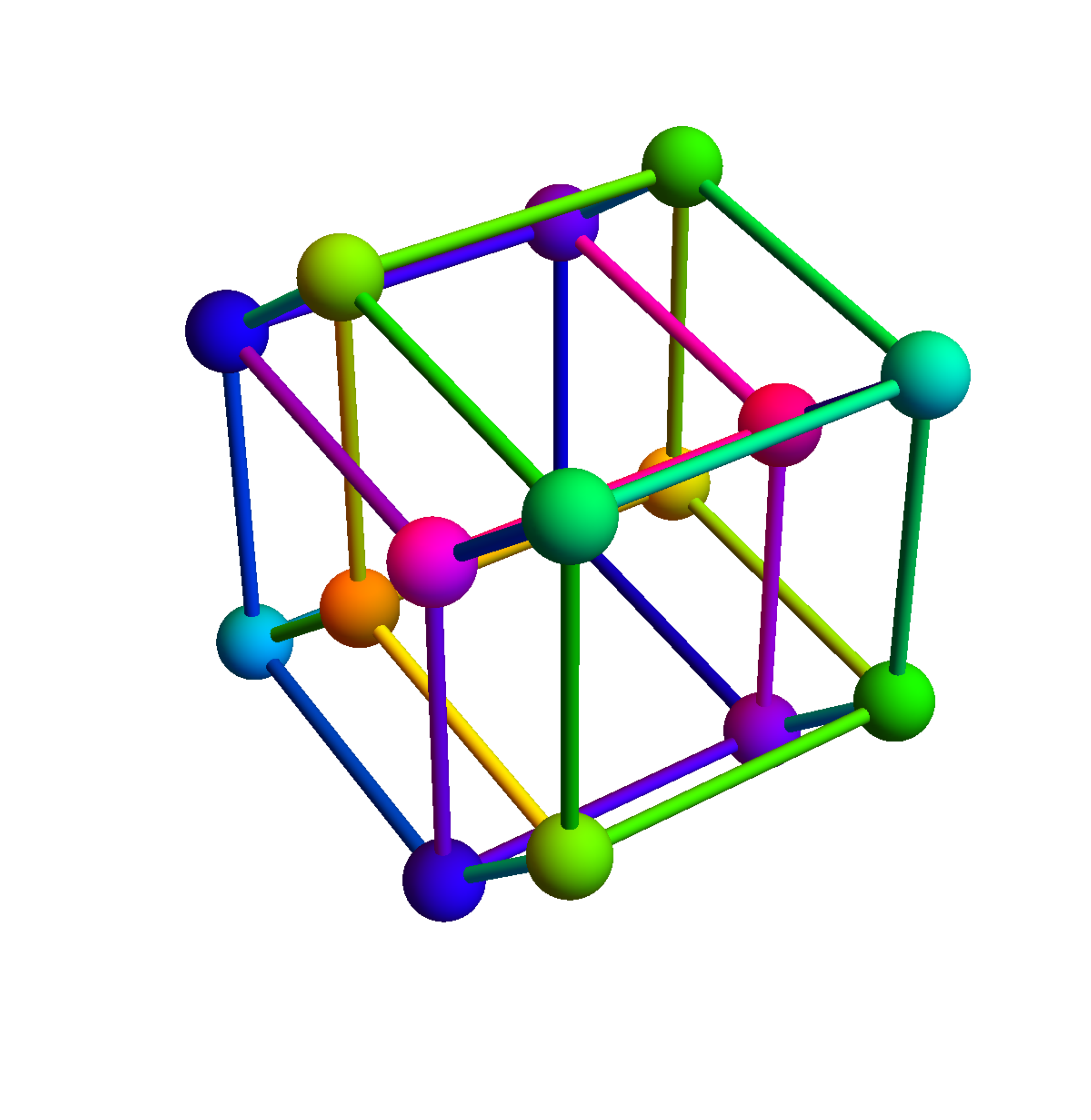}}}
}
\caption{
The 16 cell $G \in \Scal_3$ and its dual $\hat{G}$, the {\bf tesseract}.
Both are platonic solids in 4 dimensions. We look at the first as
an element in $\Scal_3$, a sphere. 
\label{orientable}
}
\end{figure}

\begin{figure}[h]
\scalebox{0.22}{\includegraphics{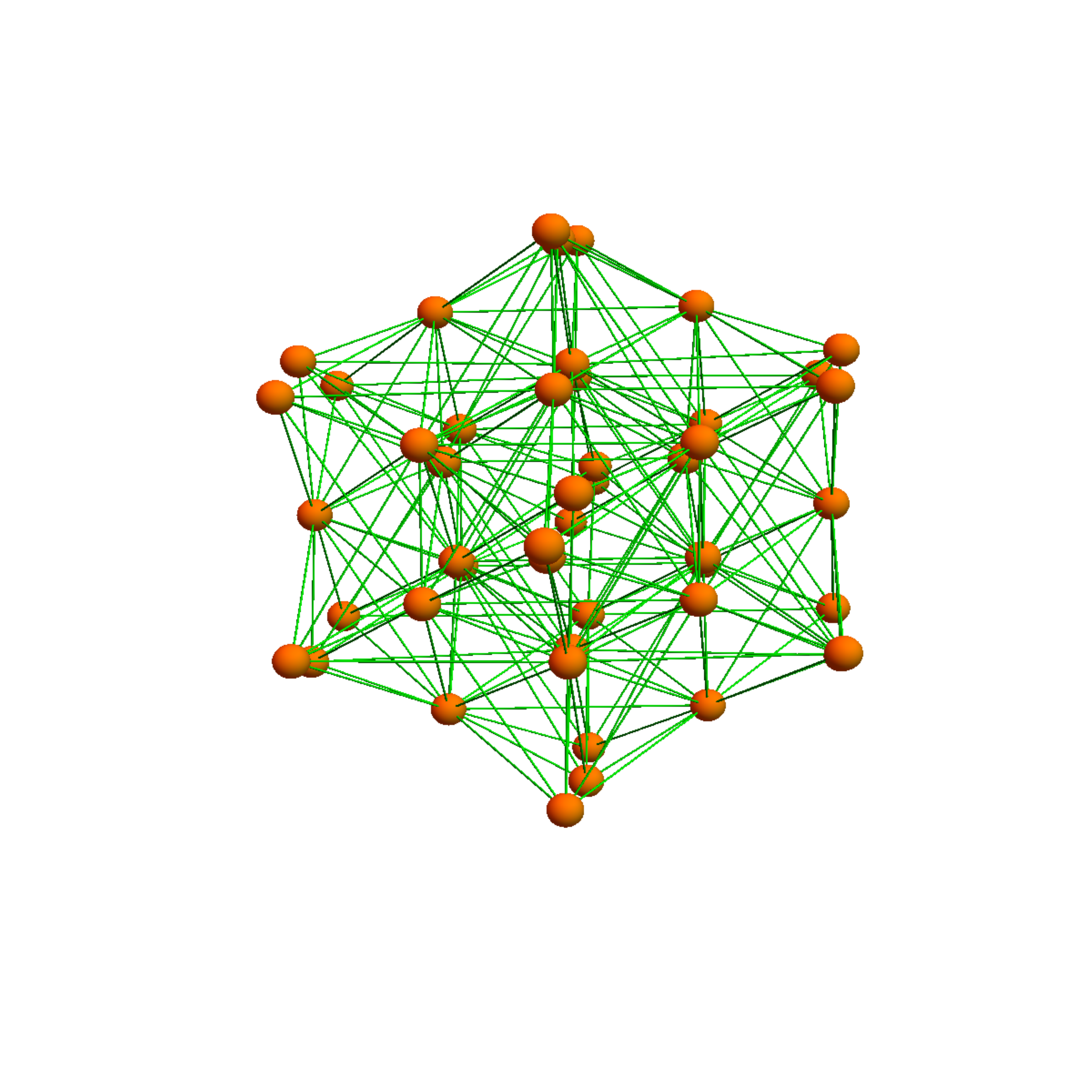}}
\caption{
The caped 2D cube is Eulerian and so in $\Ccal_3$. 
The tesseract is a polytop. It can be completed to become a 
three dimensional sphere $G \in \Scal_3$. It has the property that
all edge degrees are even so that it is in $\Ccal_4$. 
It could therefore be used for refinements: replace a tetrahedron with such a cube. 
}
\end{figure}      

\begin{figure}[h]
\parbox{15cm}{
\parbox{7.2cm}{\scalebox{0.12}{\includegraphics{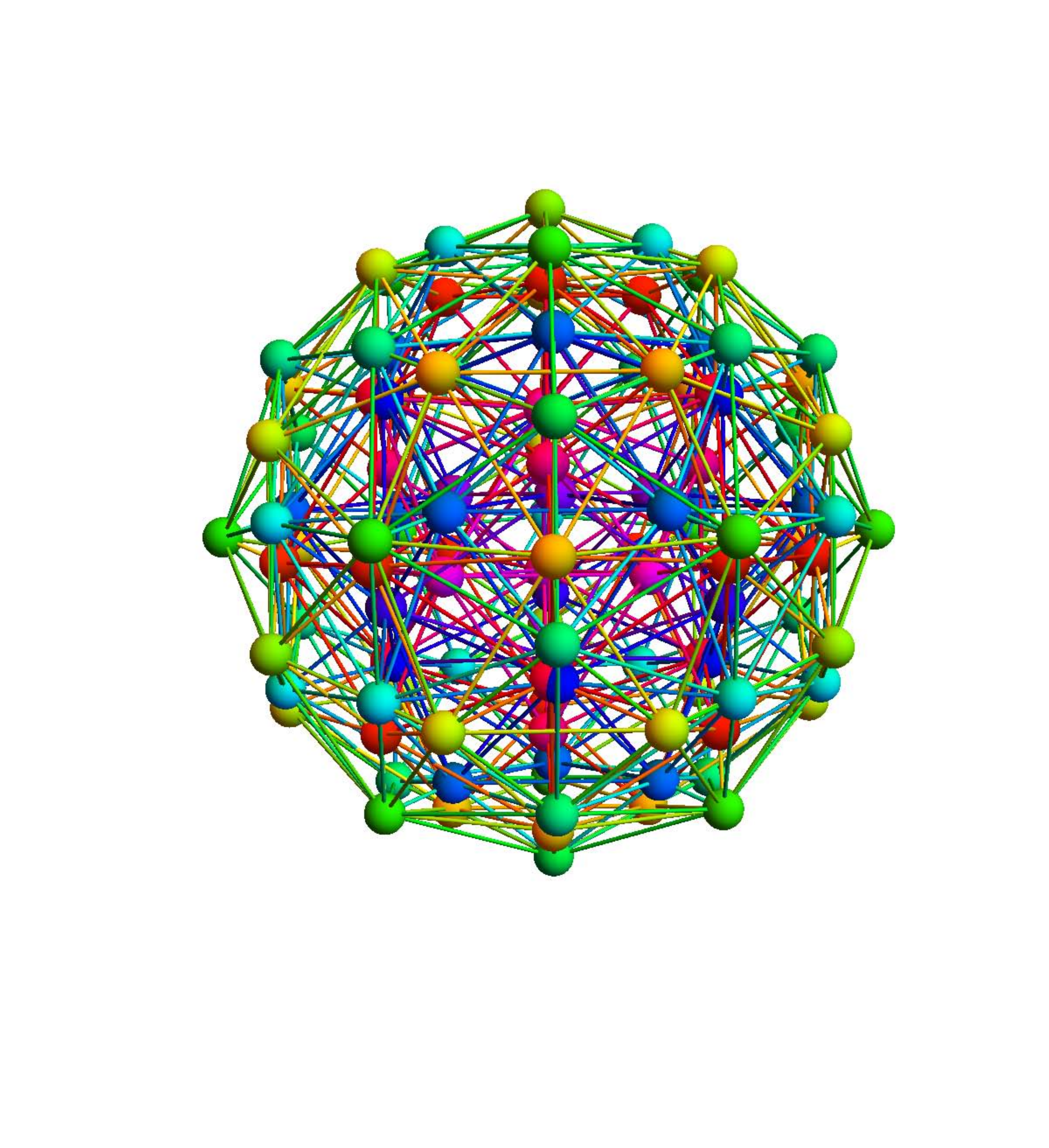}}}
\parbox{7.2cm}{\scalebox{0.12}{\includegraphics{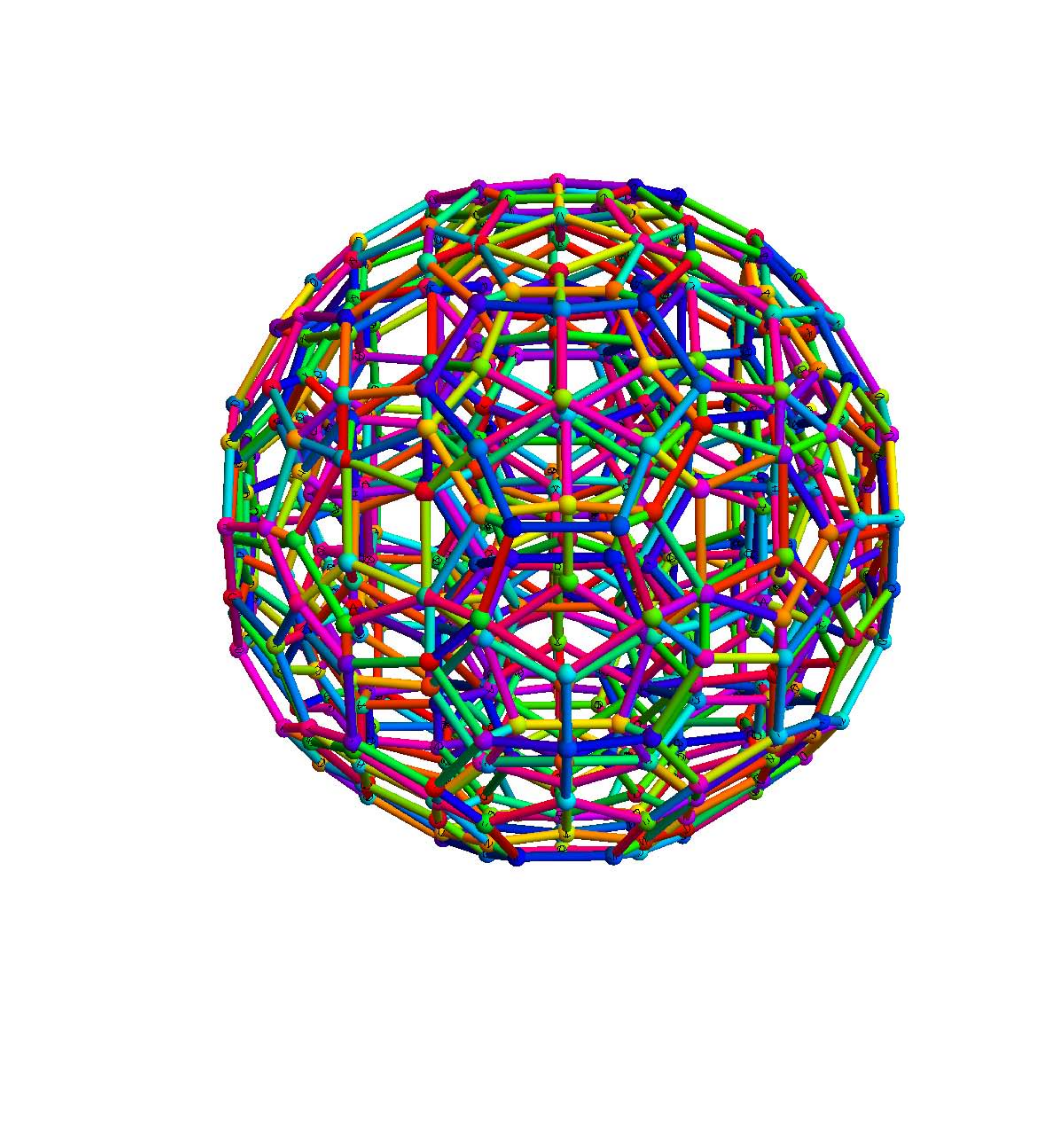}}}
}
\caption{
The 600 cell is a graph $G \in \Scal_3$. Its dual graph is 
the 120 cell $\hat{G}$. 
\label{orientable}
}
\end{figure}

\begin{figure}[h]
\scalebox{0.4}{\includegraphics{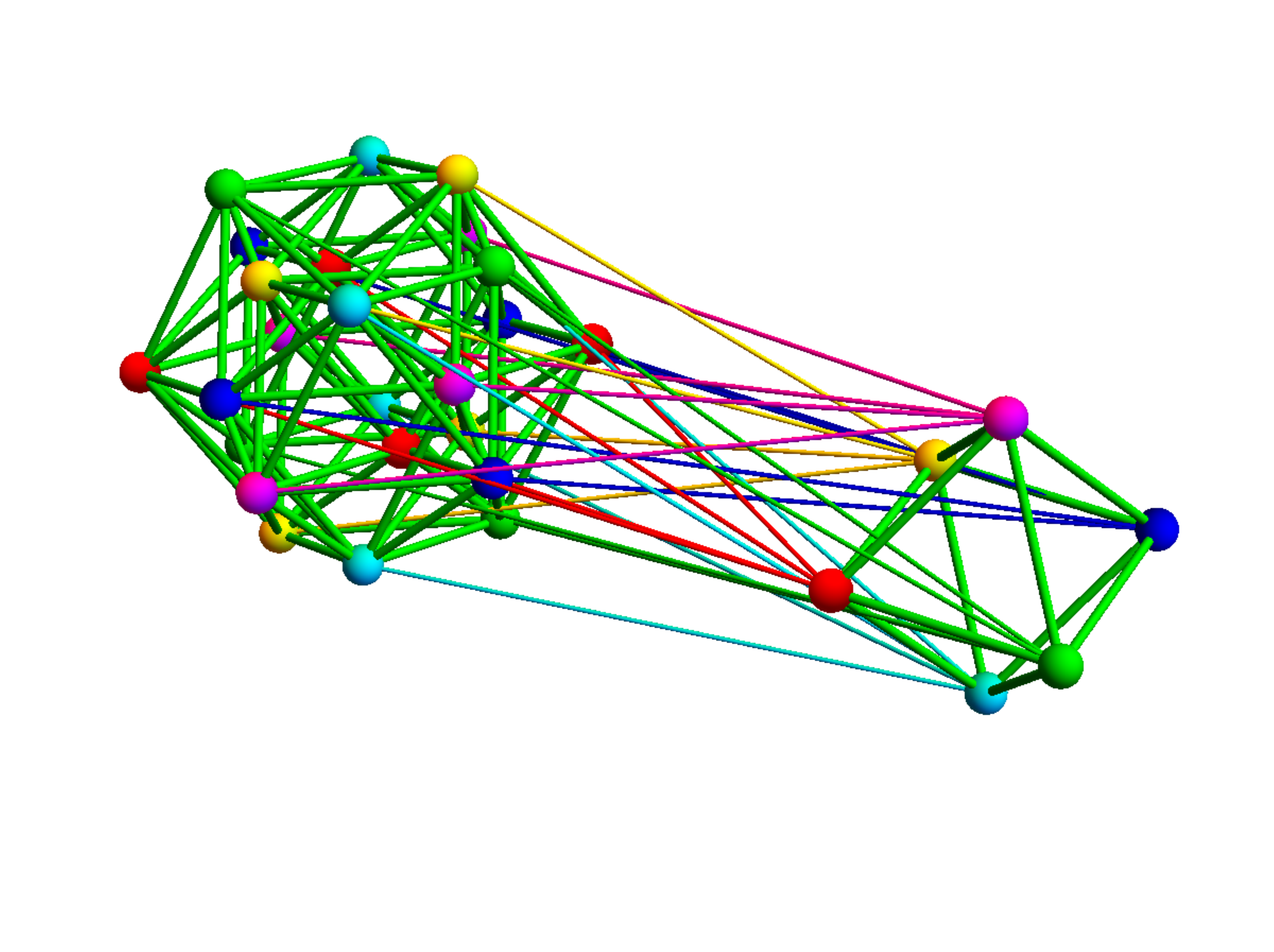}}
\caption{
A map $f$ from the $24$ cell to the octahedron, which
is a discrete 2 sphere. This is a discrete Hopf fibration.
Every inverse $f^{-1}(\{x\})$ consists of 4 points and plays the role of
a circle in the continuum. The inverse of a circle plays the role of the tori visible in
animations of the Hopf fibration. Note that the 24 cell $G$ is not 
geometric as ${\rm dim}(G)=2$. }
\end{figure}

\begin{figure}[h]
\parbox{13.8cm}{
\parbox{6.5cm}{\scalebox{0.19}{\includegraphics{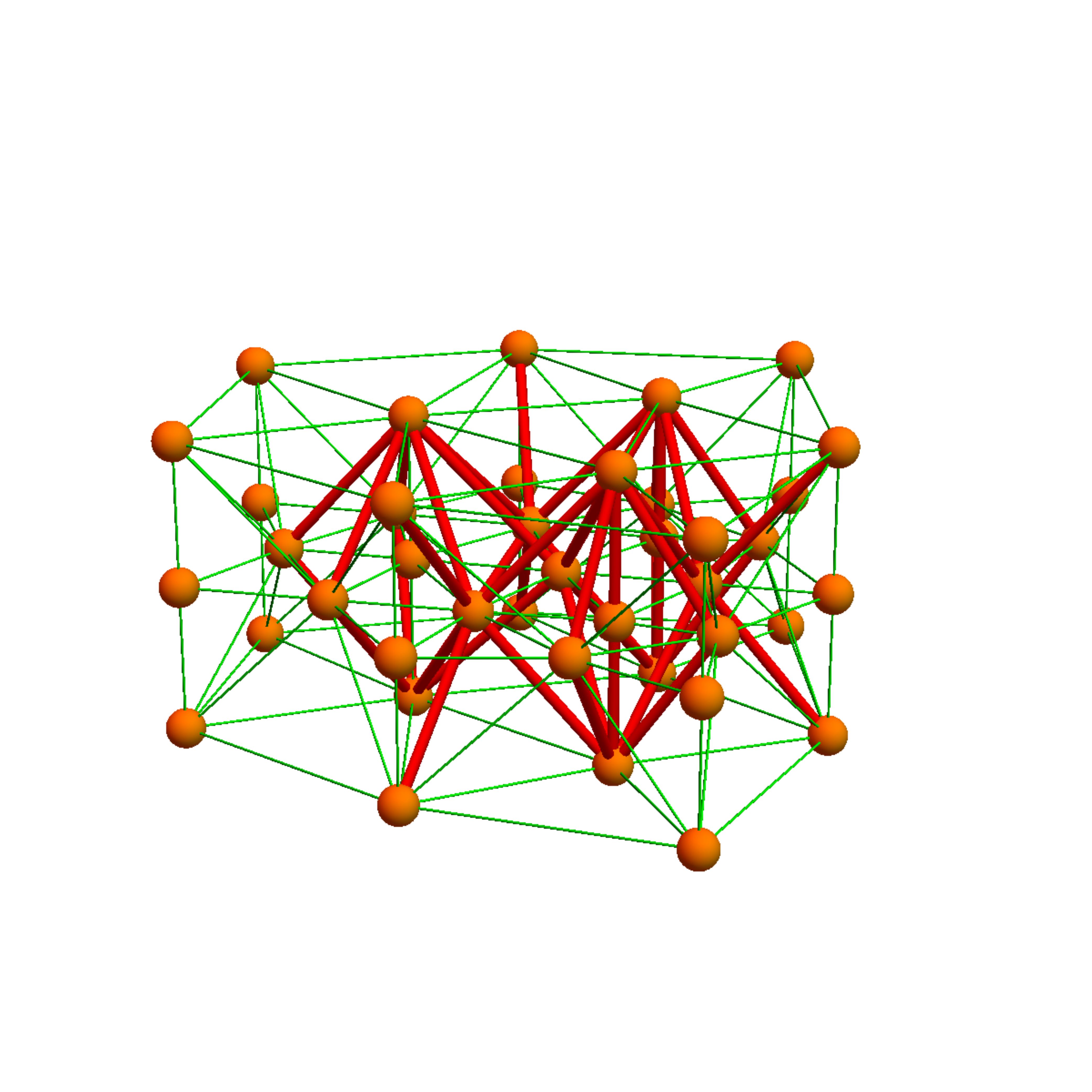}} }
\parbox{6.5cm}{\scalebox{0.19}{\includegraphics{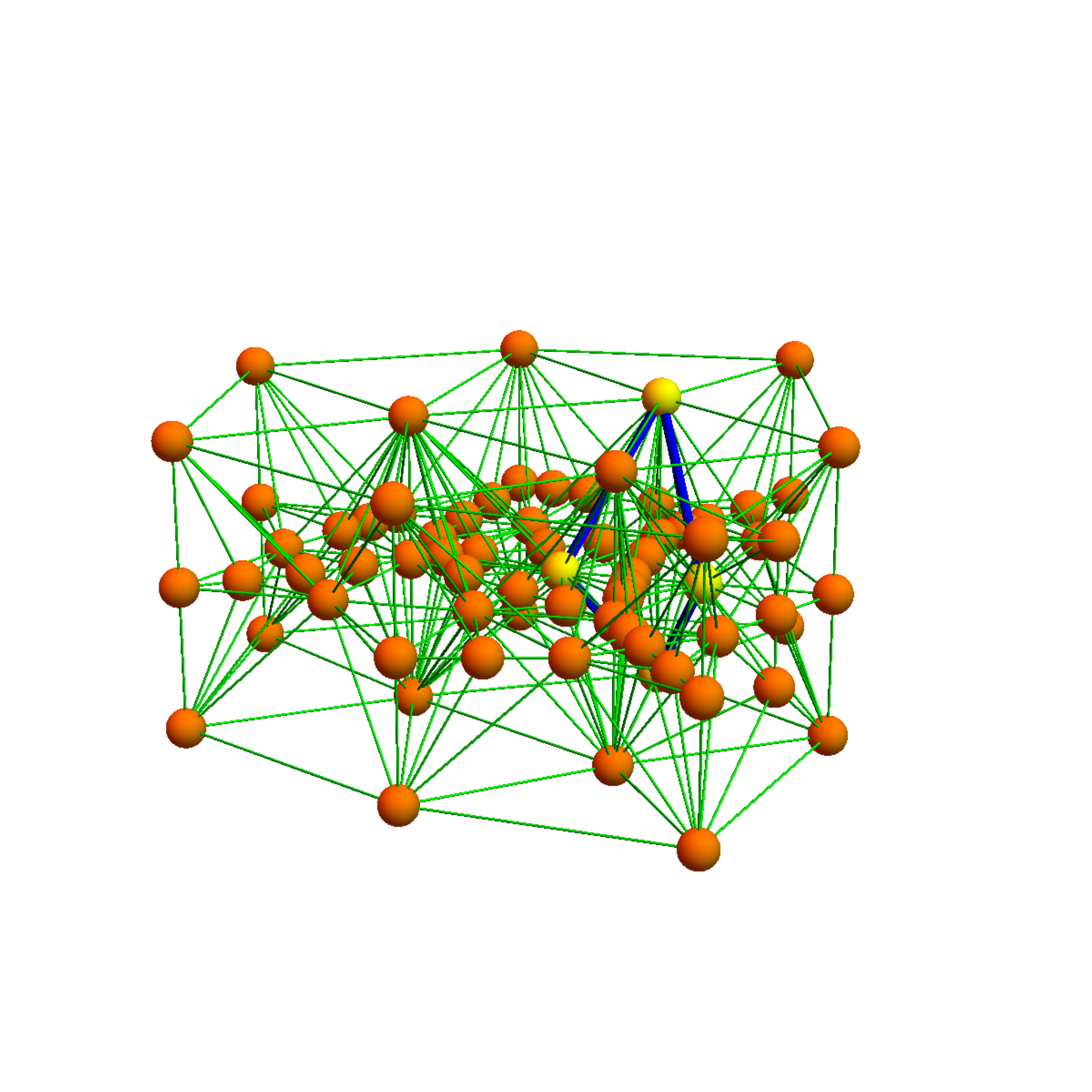}} }
}
\parbox{13.8cm}{
\parbox{6.5cm}{\scalebox{0.19}{\includegraphics{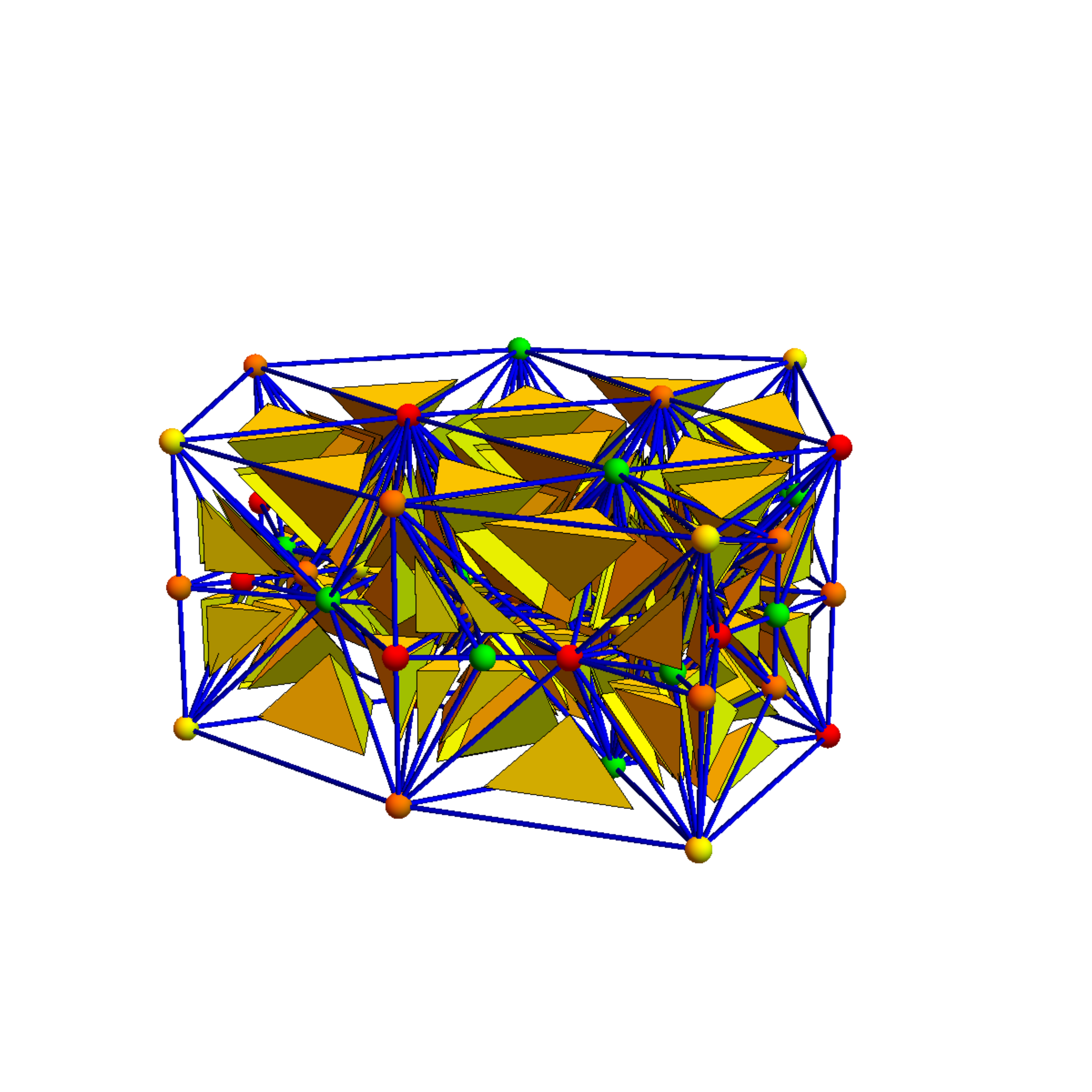}} }
\parbox{6.5cm}{\scalebox{0.19}{\includegraphics{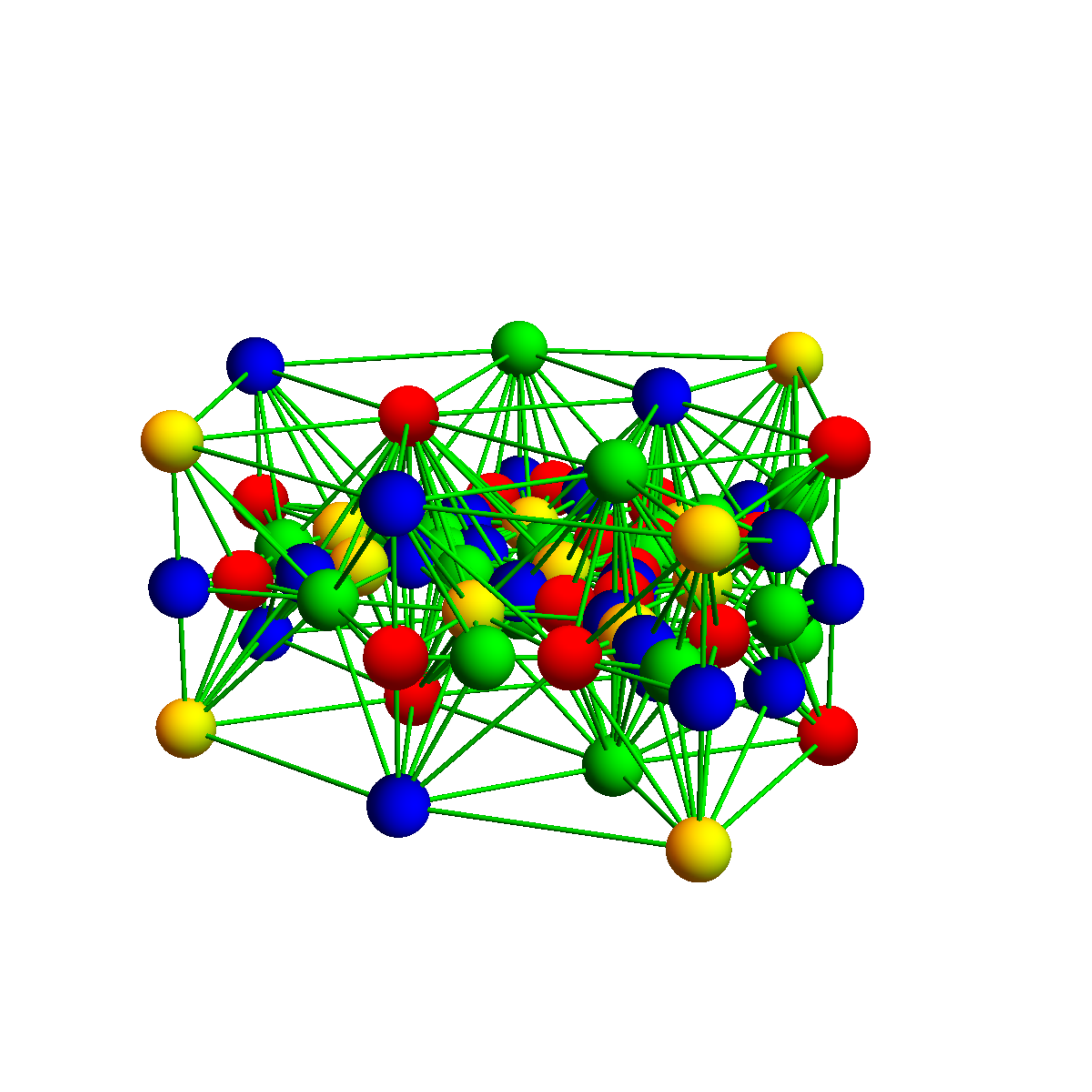}} }
}
\caption{
Coloring a $2$-dimensional graph using a self-cobordism. 
We could also just take half, but then there is no interior, but 
the later is only a cosmetic requirement.
}
\end{figure}

\begin{figure}[h]
\parbox{6.2cm}{\scalebox{0.19}{\includegraphics{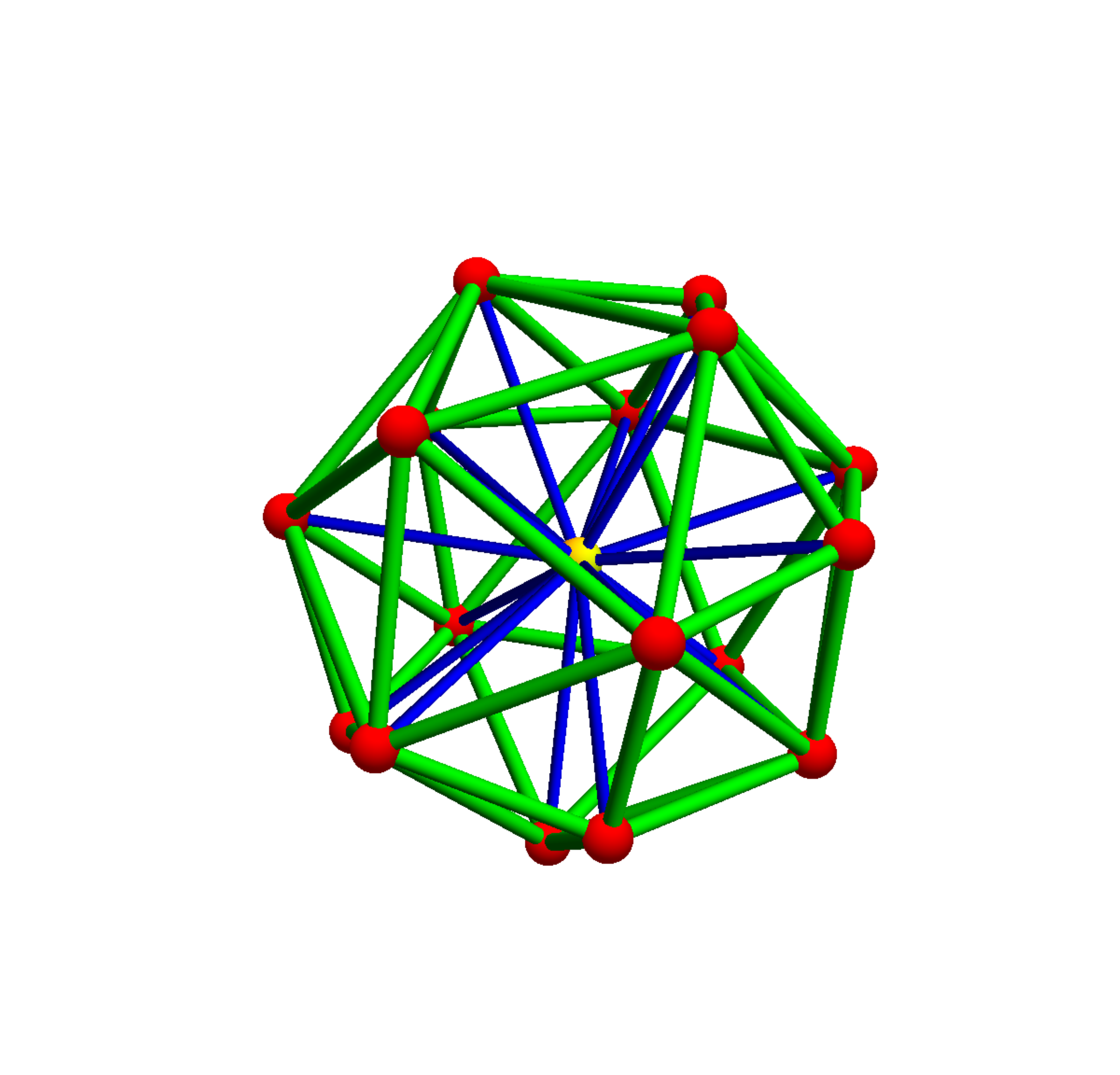}} }
\parbox{6.2cm}{\scalebox{0.19}{\includegraphics{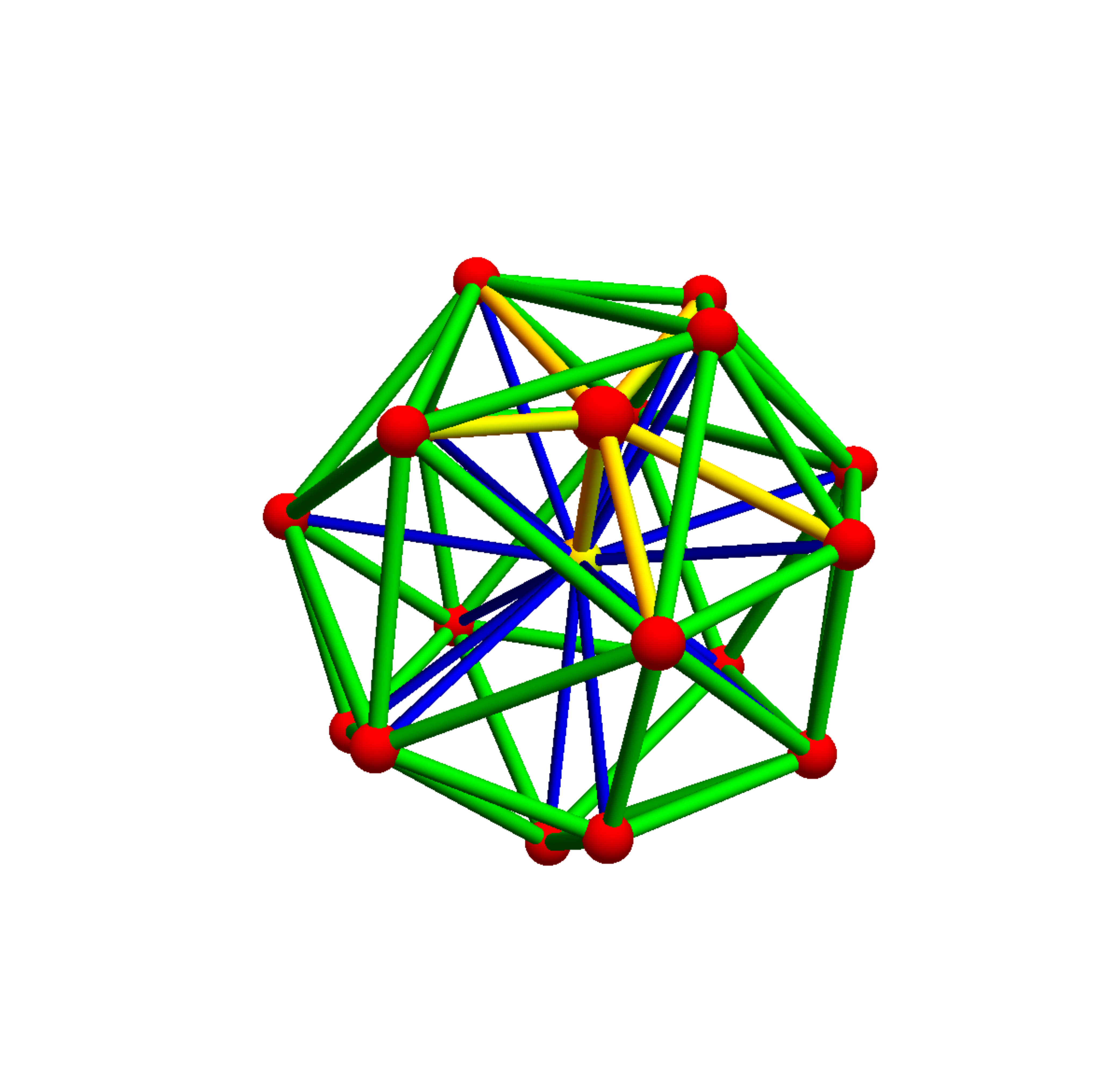}} }
\parbox{6.2cm}{\scalebox{0.19}{\includegraphics{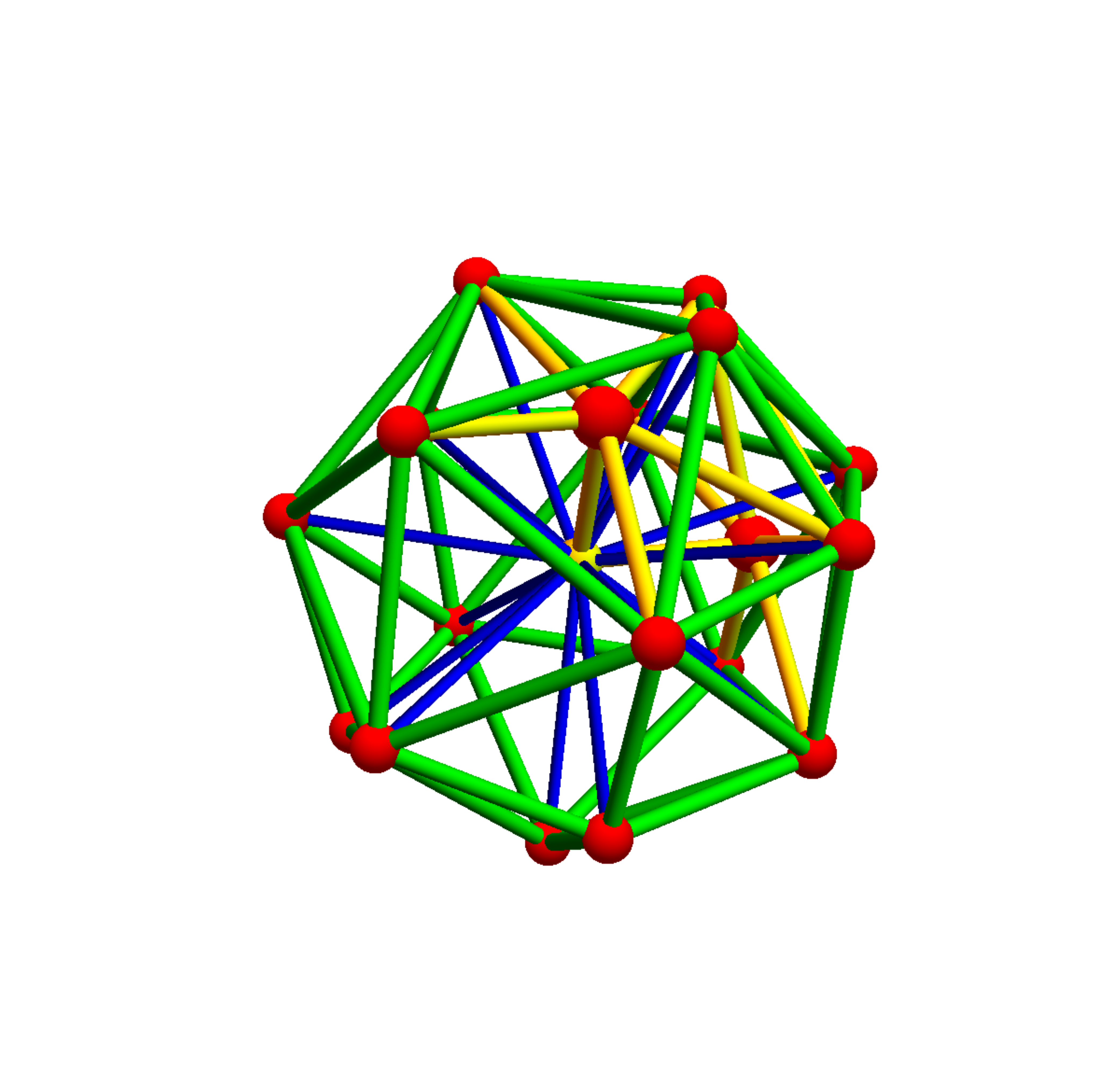}} }
\parbox{6.2cm}{\scalebox{0.19}{\includegraphics{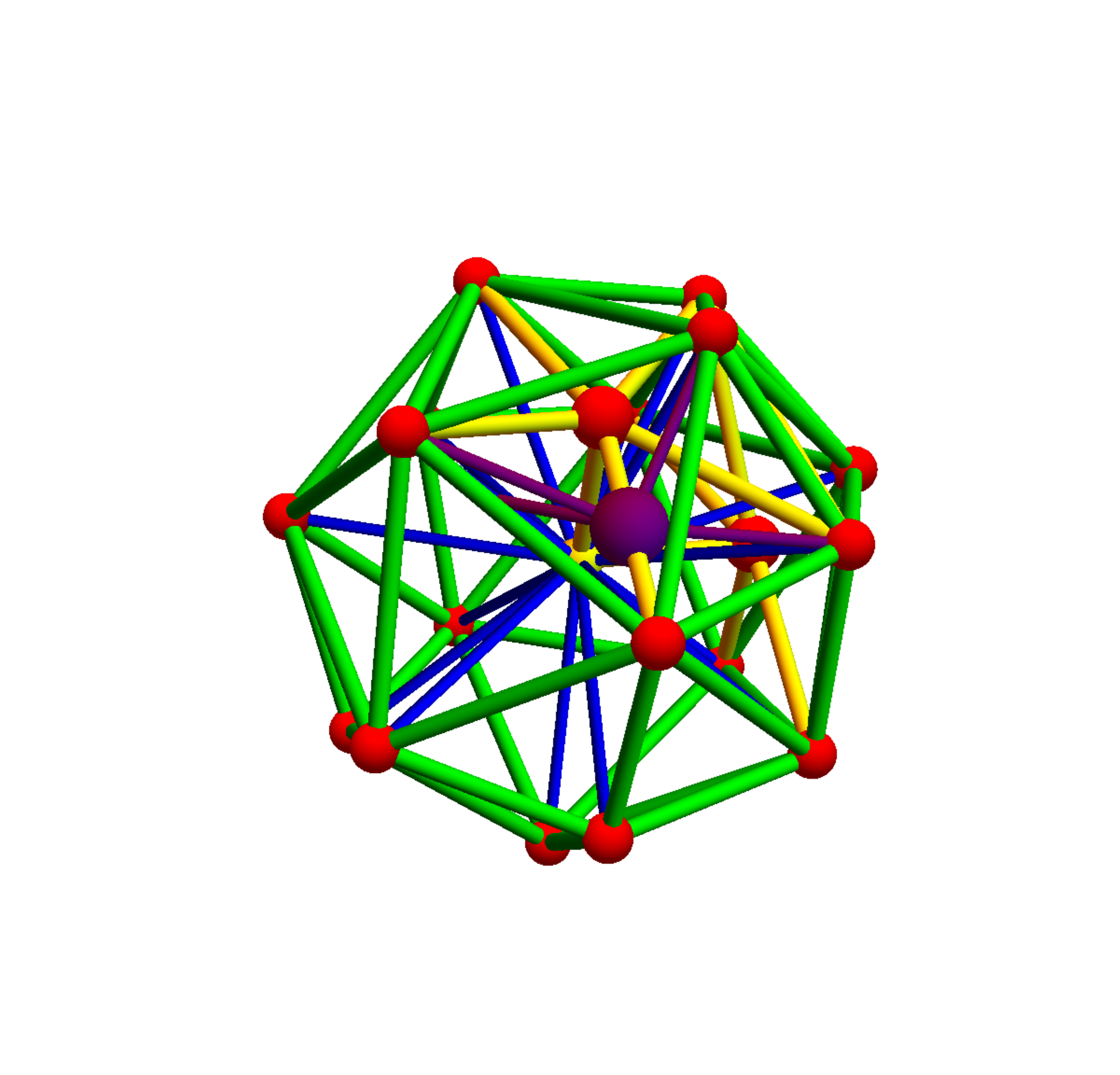}} }
\caption{
Refining a graph consists of building new $2$-dimensional surfaces in the 
inside: subdivide an edge $e=(a,b)$ with an additional vertex $v$ 
and will in a wheel graph centered at $v$ and boundary $S(a) \cap S(b)$. 
}
\end{figure}

\begin{figure}[h]
\parbox{6.2cm}{\scalebox{0.25}{\includegraphics{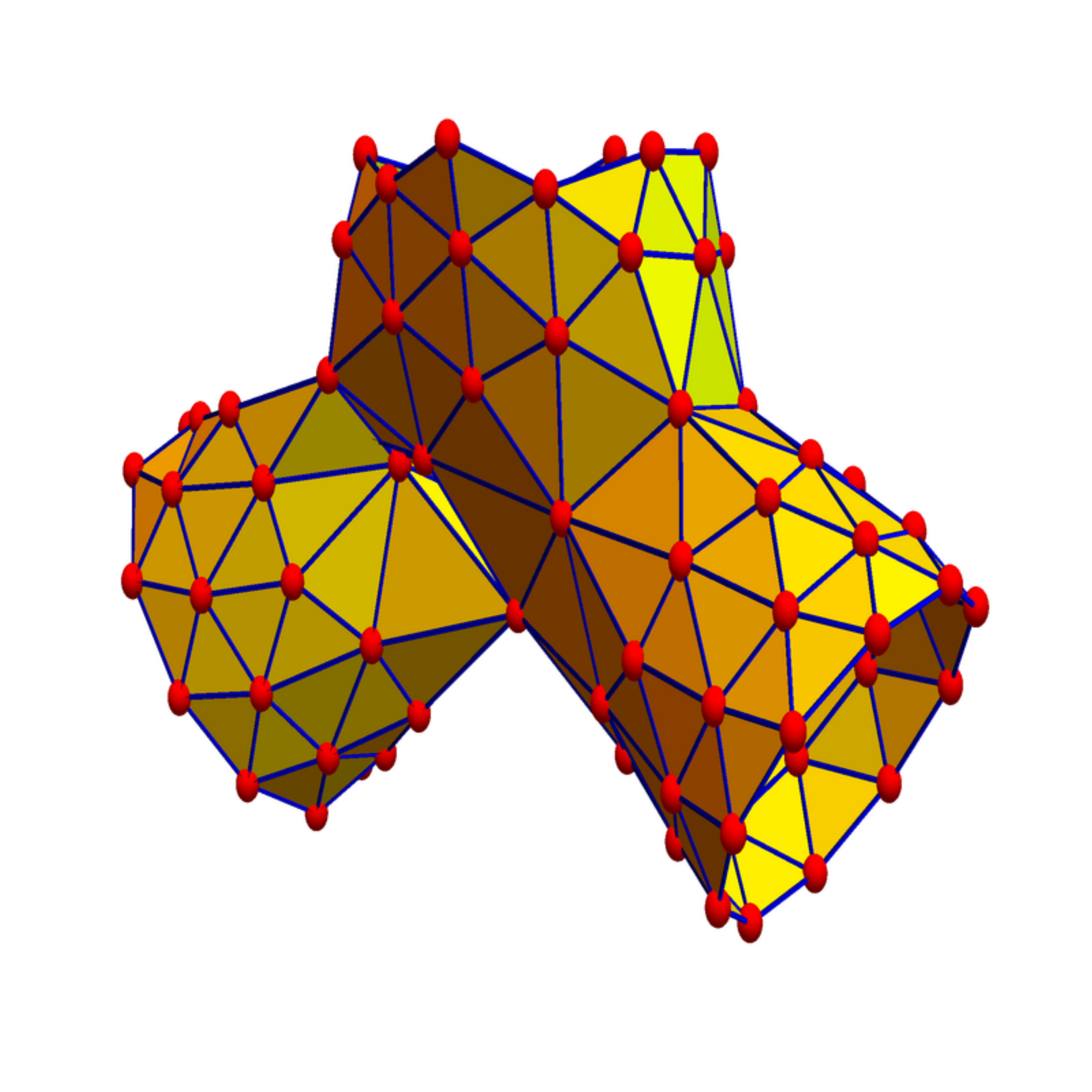}} }
\caption{
The notion of cobordism for graphs is the same as in the continuum. The earliest definition for
simplicial complexes is \cite{Fisk1980}. 
Two geometric graphs of dimension $d$ are cobordant if their disjoint union is the boundary of a
$(d+1)$-dimensional geometric graph. Strangely, the notion seems not have appeared in 
graph theory \cite{knillcalculus}.  }
\end{figure}

\begin{figure}
\parbox{6.2cm}{\scalebox{0.12}{\includegraphics{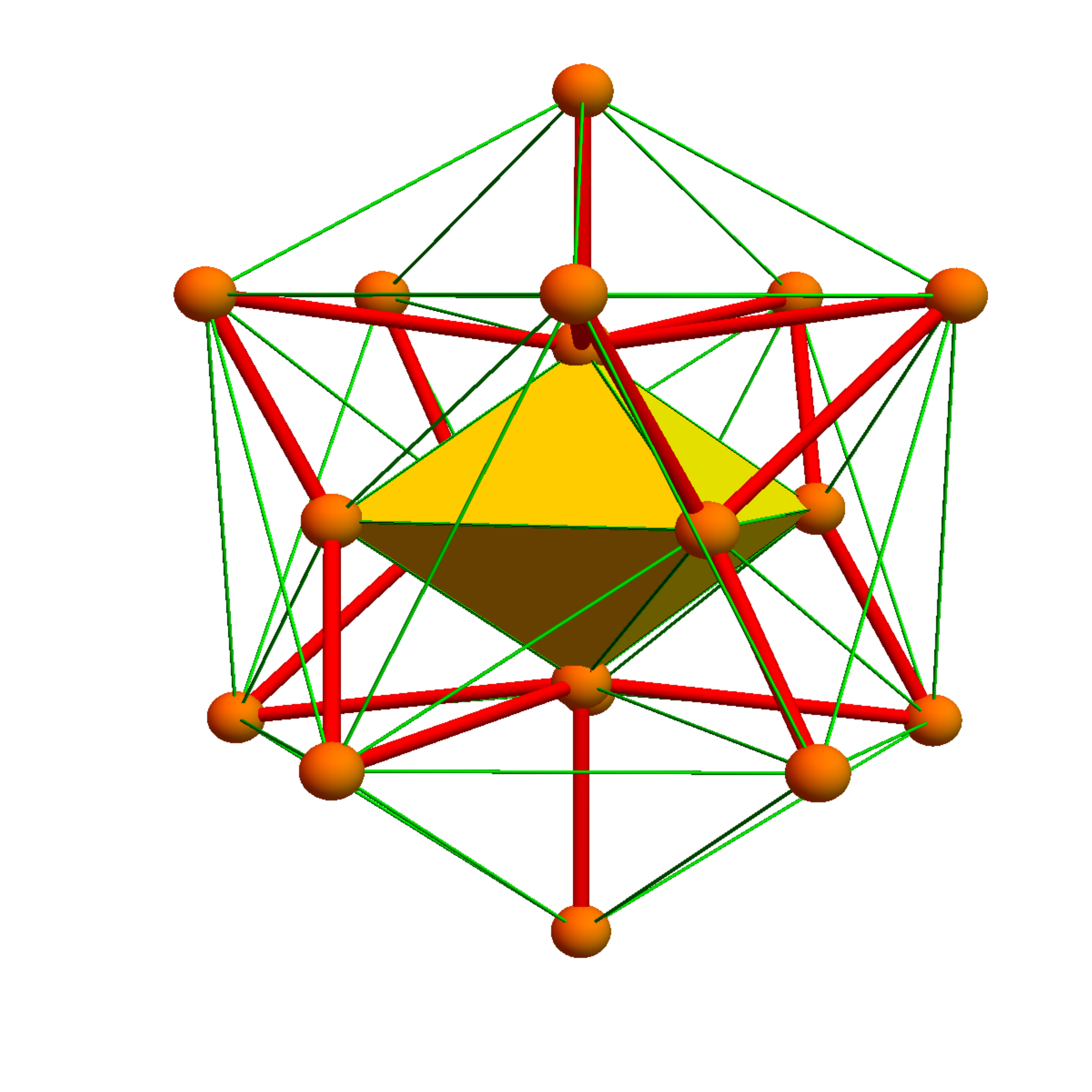}} }
\parbox{6.2cm}{\scalebox{0.12}{\includegraphics{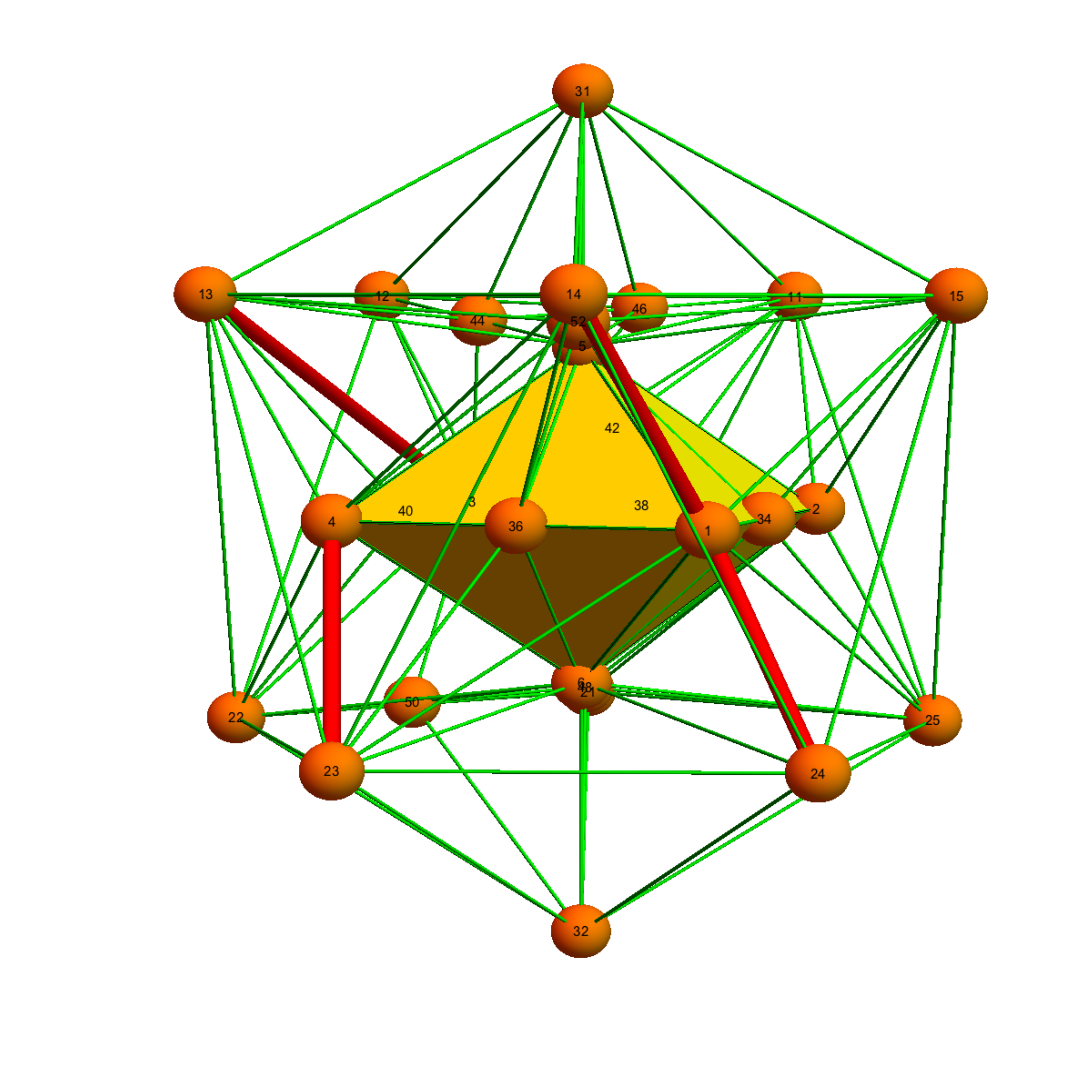}} }
\caption{
A cobordism between the octahedron $G_1$ and the icosahedron $G_2$ is
a $3$-dimensional annulus graph $H \in \Gcal_3$ with the topology of 
$S^2 \times [0,1]$.
Its Betti vector is  $(1,0,1)$ as it is simply connected but contains
a non-contractible sphere (the octahedron). This could be used to color the
icosahedron. Marked are the odd degree edges. The second picture shows
the situation after some subdivisions.
}
\end{figure}

\clearpage

\begin{figure}[h]
\parbox{14.5cm}{
\parbox{4.7cm}{\scalebox{0.12}{\includegraphics{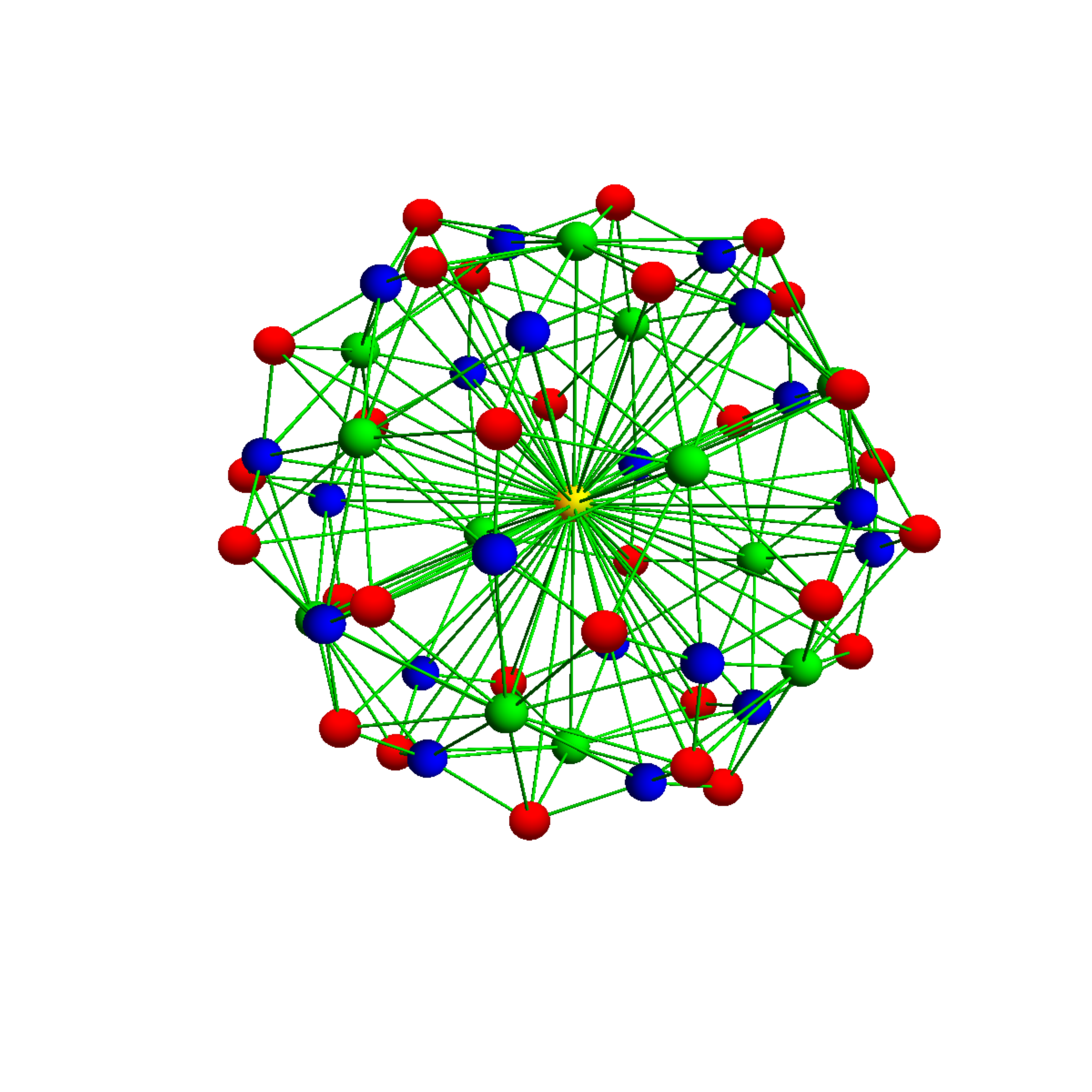}} }
\parbox{4.7cm}{\scalebox{0.12}{\includegraphics{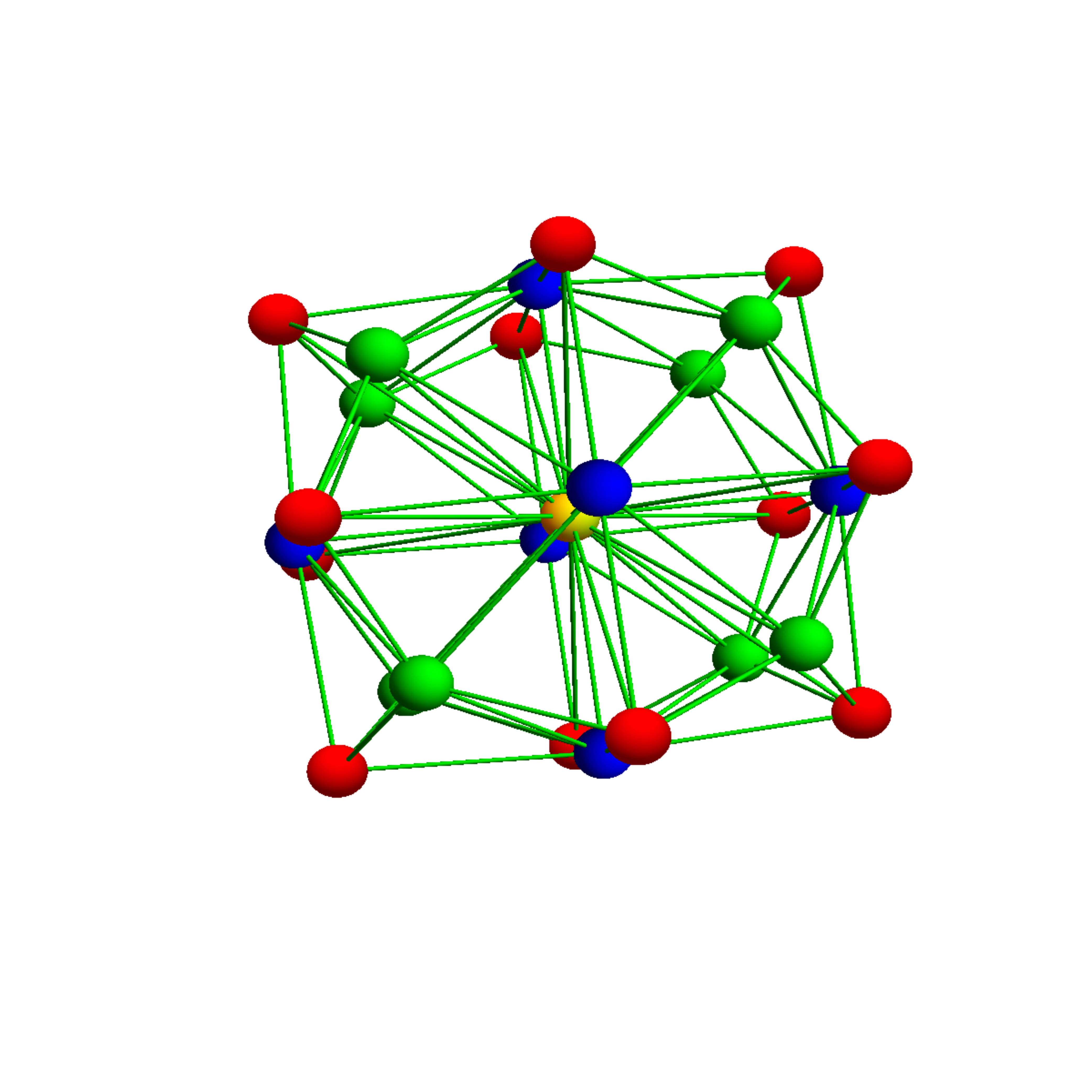}} }
\parbox{4.7cm}{\scalebox{0.12}{\includegraphics{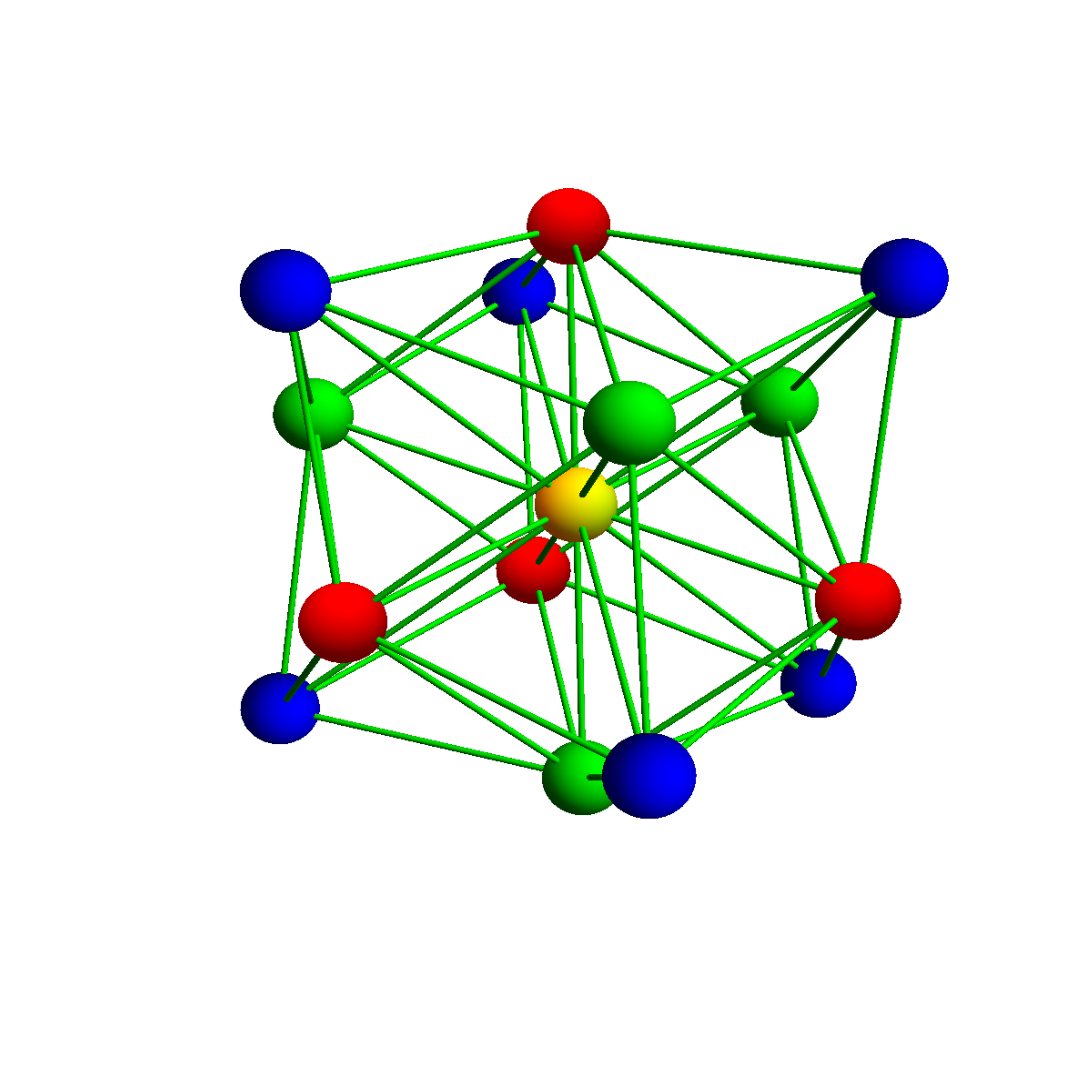}} }
}
\caption{
There are 4 Catalan polyhedra in $\Scal_2$. Three of them are in $\Ccal_3$.
We see them here as boundary of a ball $\Bcal_3 \cap \Ccal_4$. 
}
\end{figure}

\begin{figure}[h]
\parbox{13.8cm}{
\parbox{6.2cm}{\scalebox{0.14}{\includegraphics{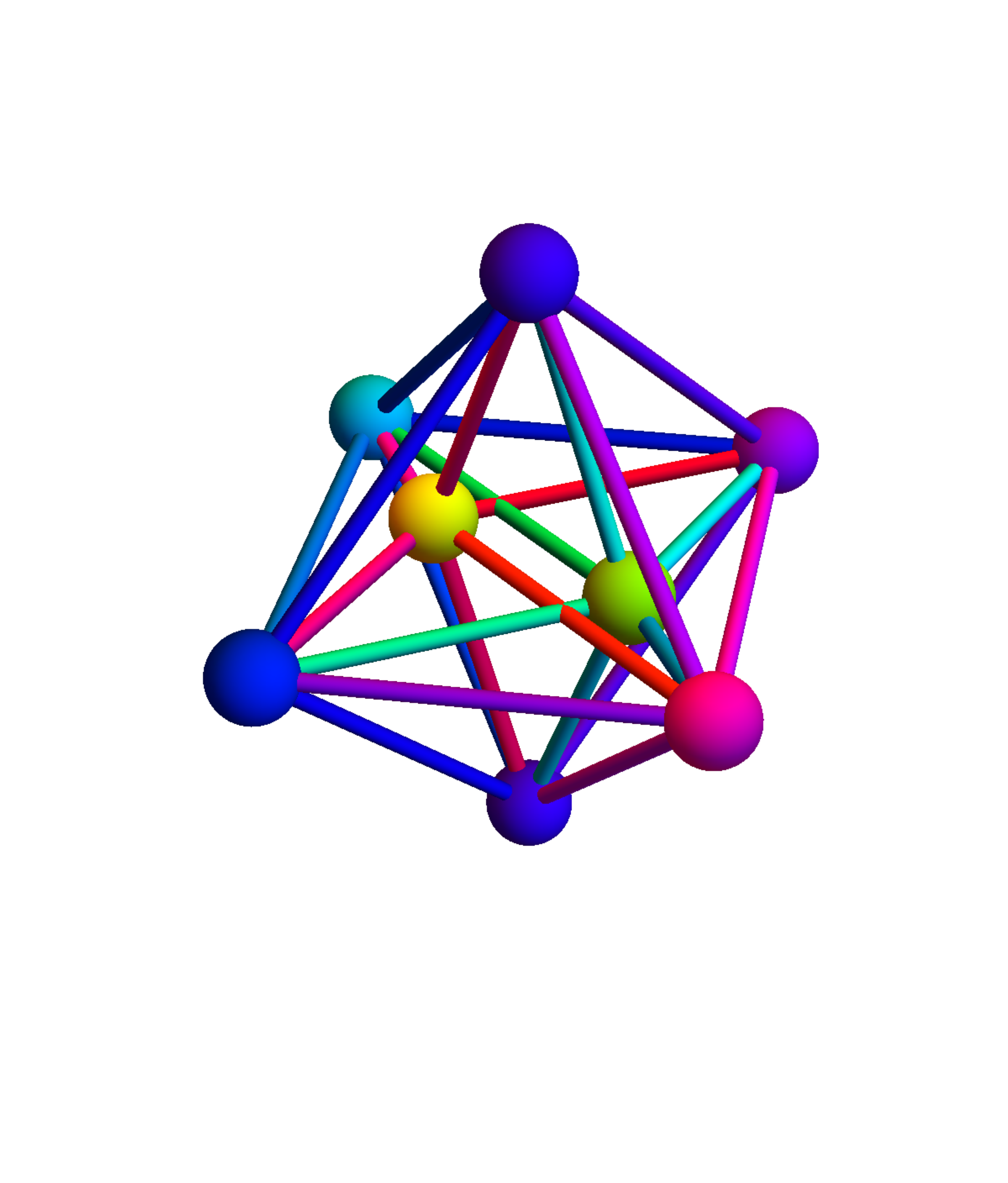}} }
\parbox{6.2cm}{\scalebox{0.14}{\includegraphics{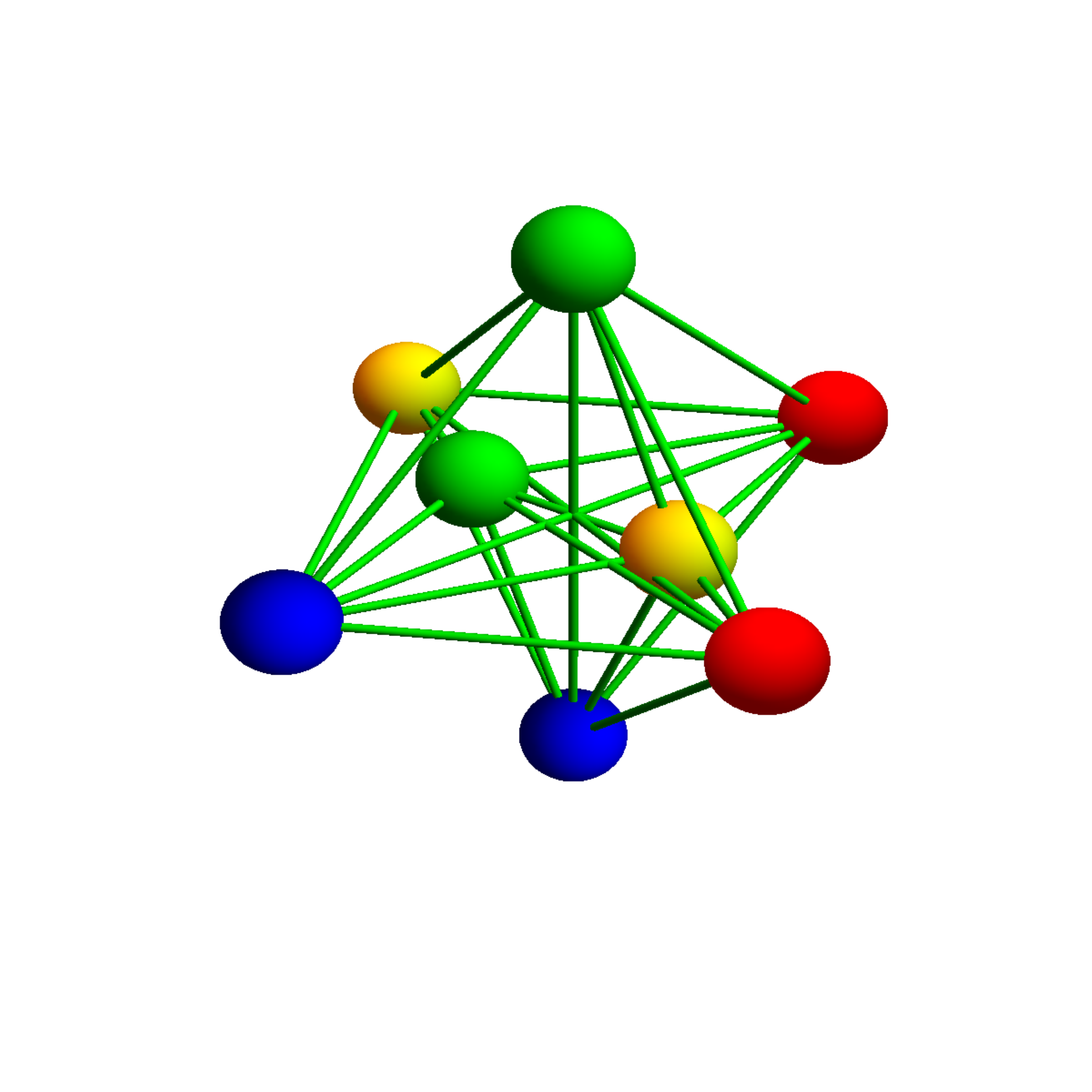}} }
}
\caption{
The 16 cell $G$ is the smallest element in $\Scal_3$. Its dual is the {\bf tesseract}, the 
$4$-dimensional cube $\hat{G}$. We see here the graph $G$ colored according to 
the 4th coordinate in the classical embedding of $R^4$, then colored with 4 colors,
as $G \in \Ccal_4$ follows from simply connectedness and the tesseract being Eulerian. 
}
\end{figure}

\begin{figure}[h]
\scalebox{0.22}{\includegraphics{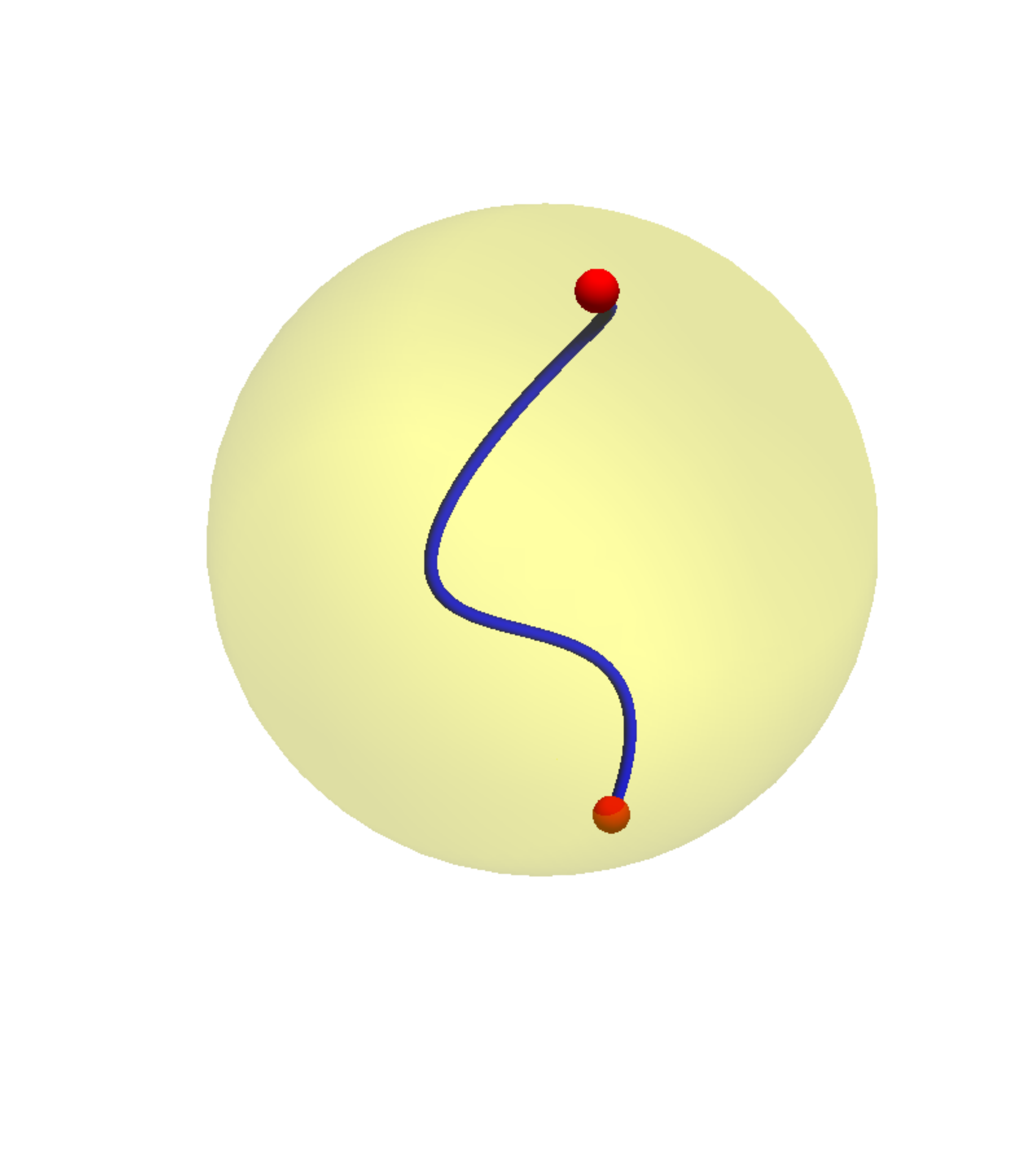}}
\scalebox{0.22}{\includegraphics{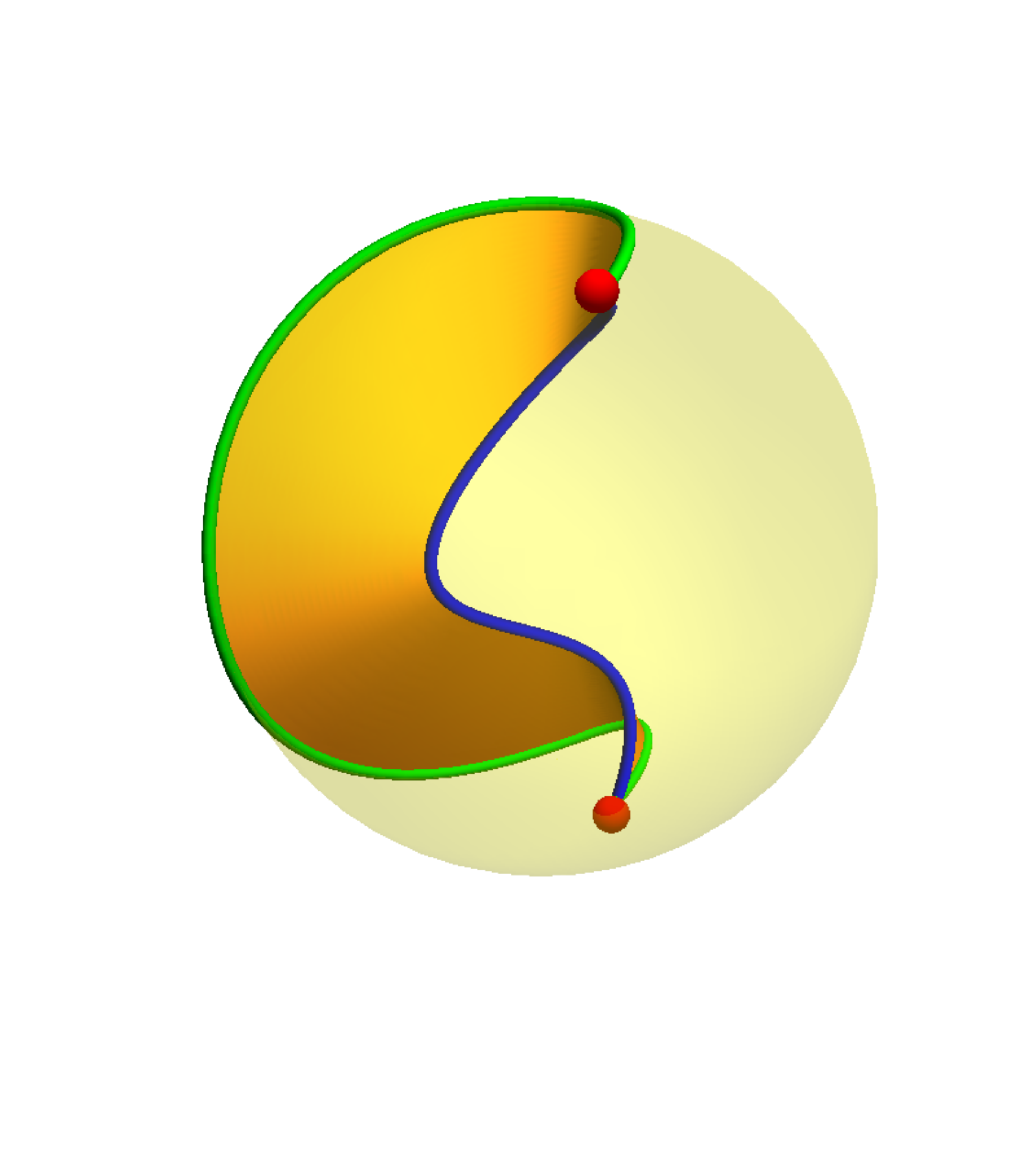}}
\scalebox{0.22}{\includegraphics{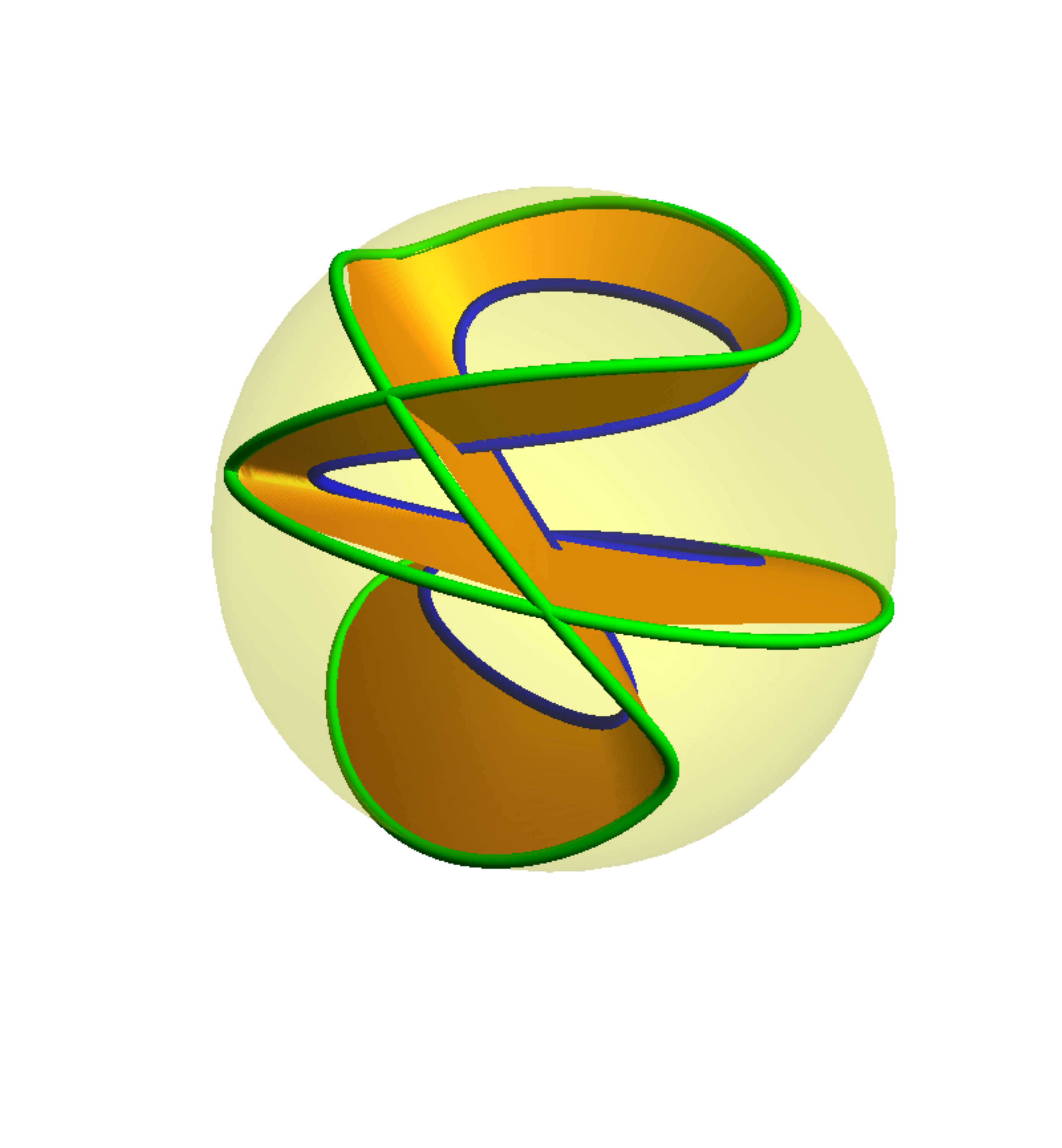}}
\scalebox{0.22}{\includegraphics{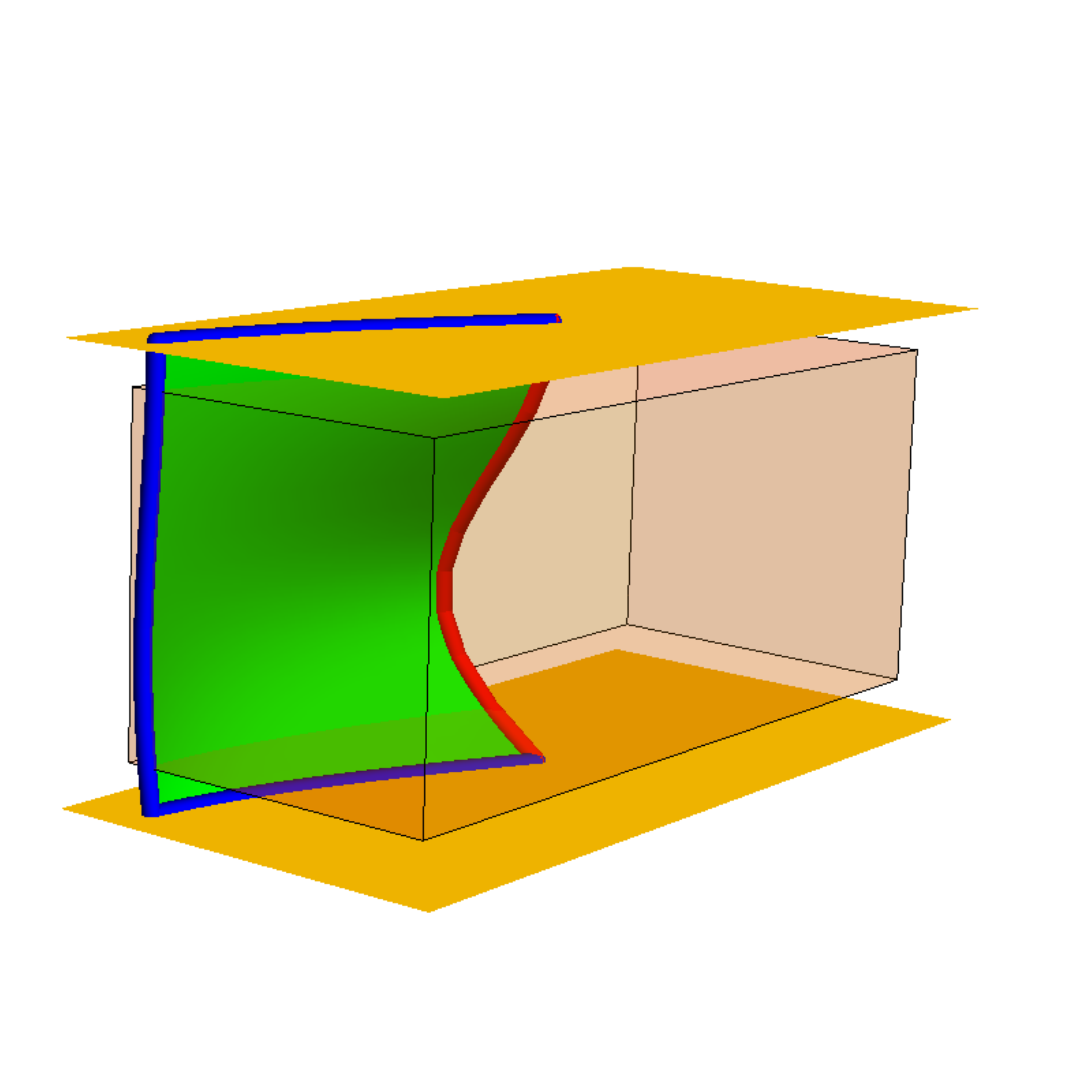}}
\caption{
Cutting along a surface containing part $\gamma$ of $O$ which is a closed curve
We have to refine the cobordism $H$ without affecting the boundary. 
These pictures illustrate the continuum limit.
}
\end{figure}

\begin{figure}[h]
\parbox{14.8cm}{
\parbox{7.2cm}{\scalebox{0.18}{\includegraphics{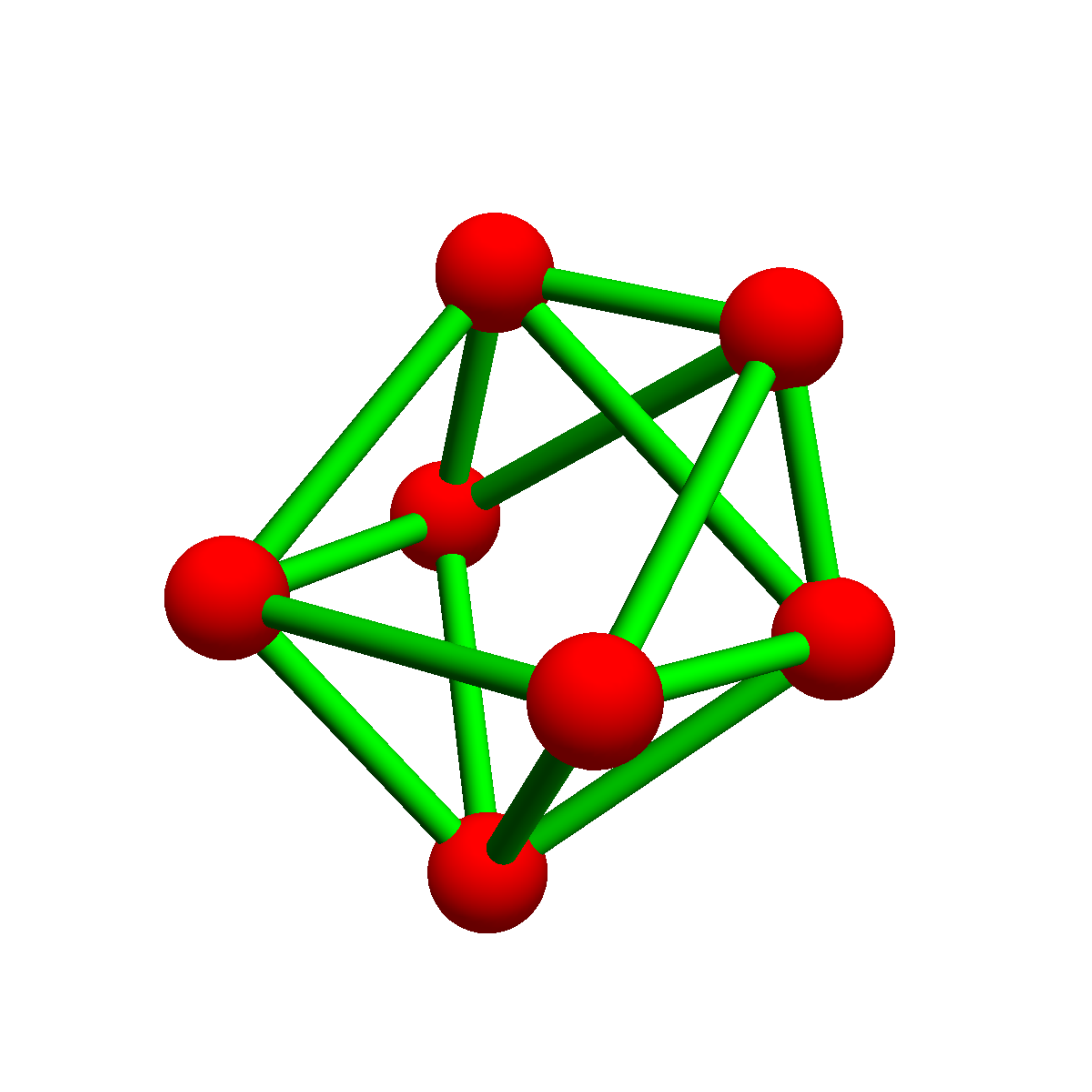}} }
\parbox{7.2cm}{\scalebox{0.18}{\includegraphics{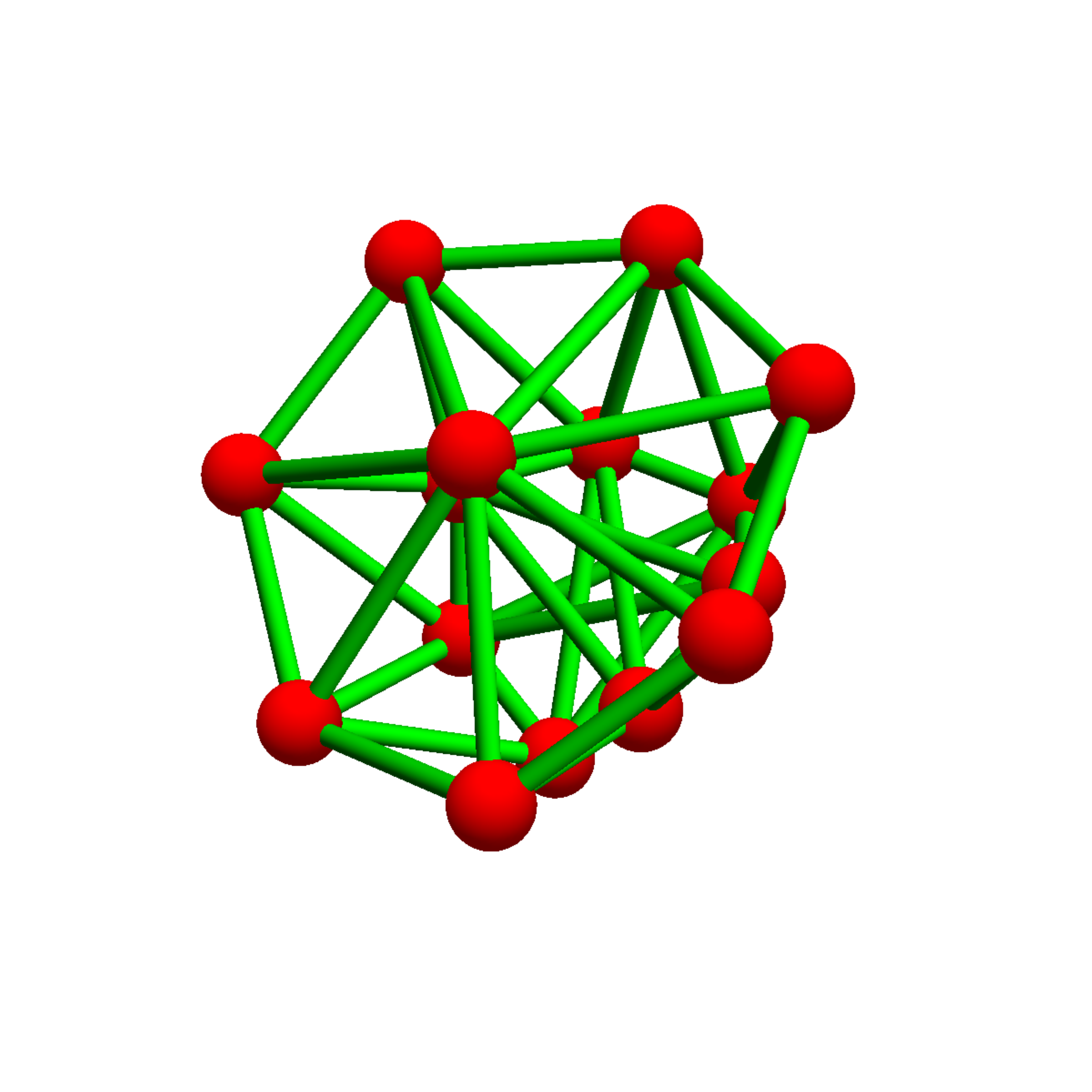}} }
}
\caption{
To the left we see a Moebius strip $M$ which is in $\Ccal_4 \setminus \Ccal_3$. 
The right picture shows a projective plane $P$ which is in $\Ccal_5 \setminus \Ccal_4$. It is not
the boundary of a $3$-dimensional graph.  }
\end{figure}

\begin{figure}[h]
\parbox{15.8cm}{
\parbox{5.1cm}{\scalebox{0.14}{\includegraphics{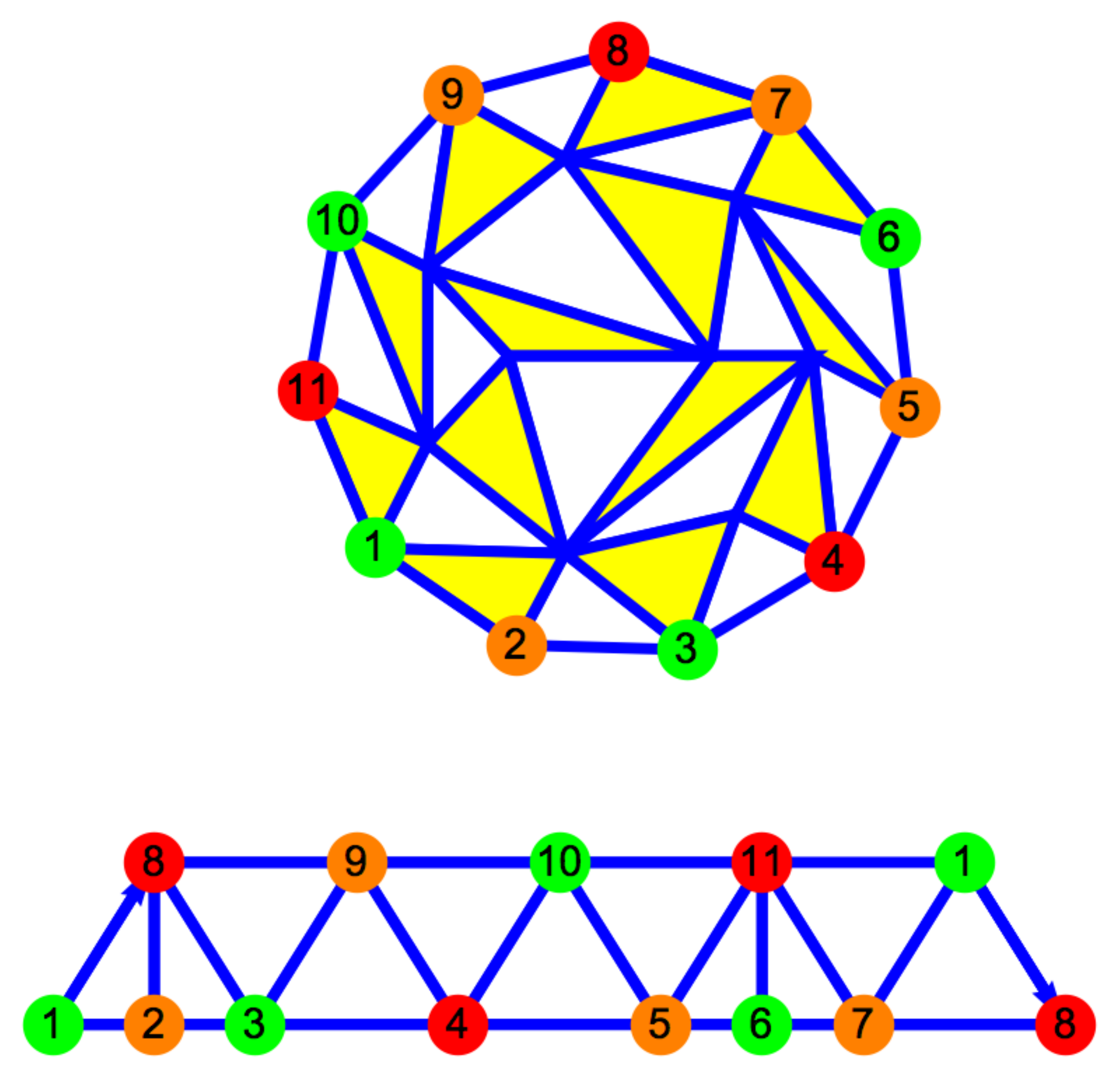}}}
\parbox{5.1cm}{\scalebox{0.14}{\includegraphics{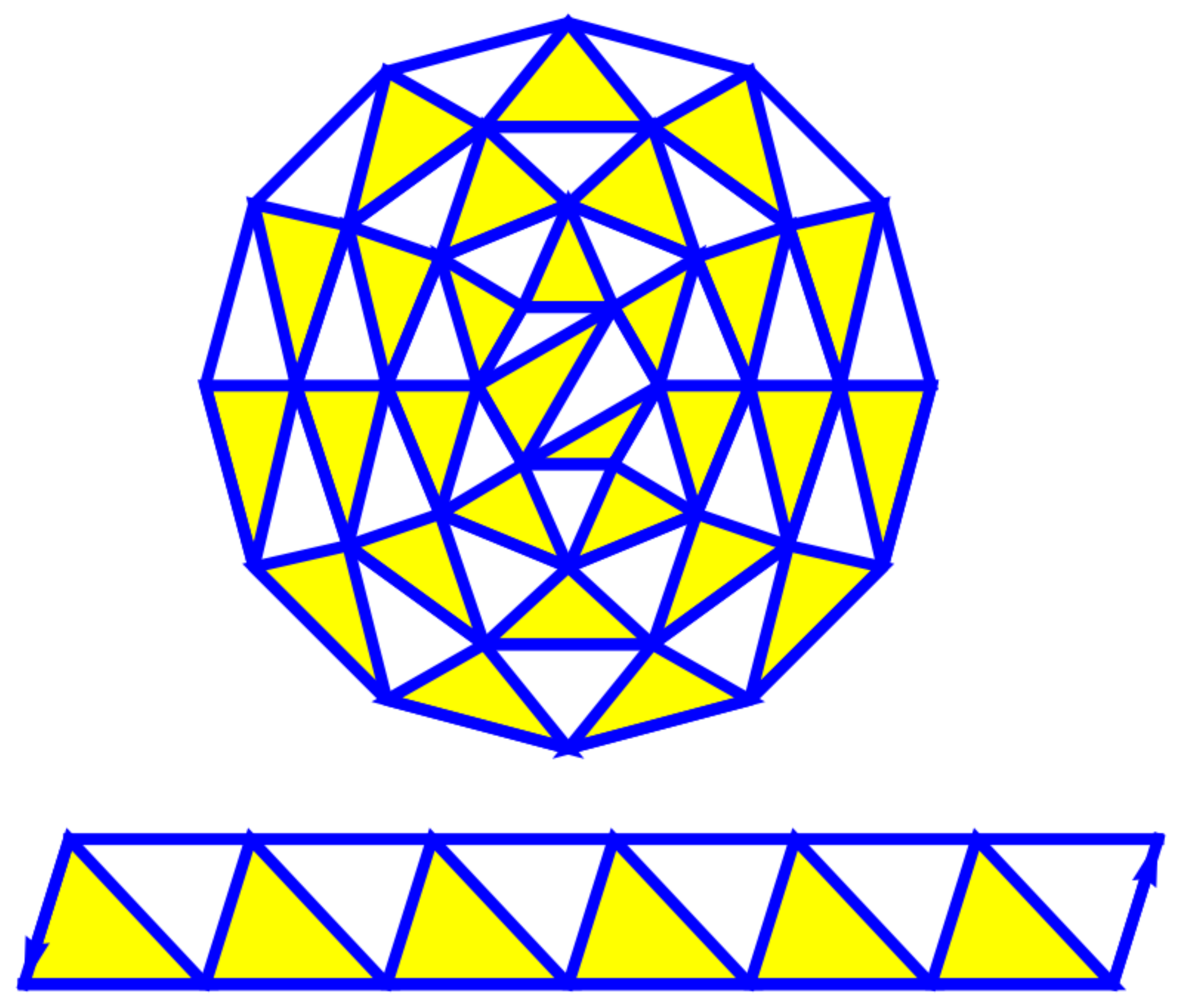}}}
\parbox{5.1cm}{\scalebox{0.14}{\includegraphics{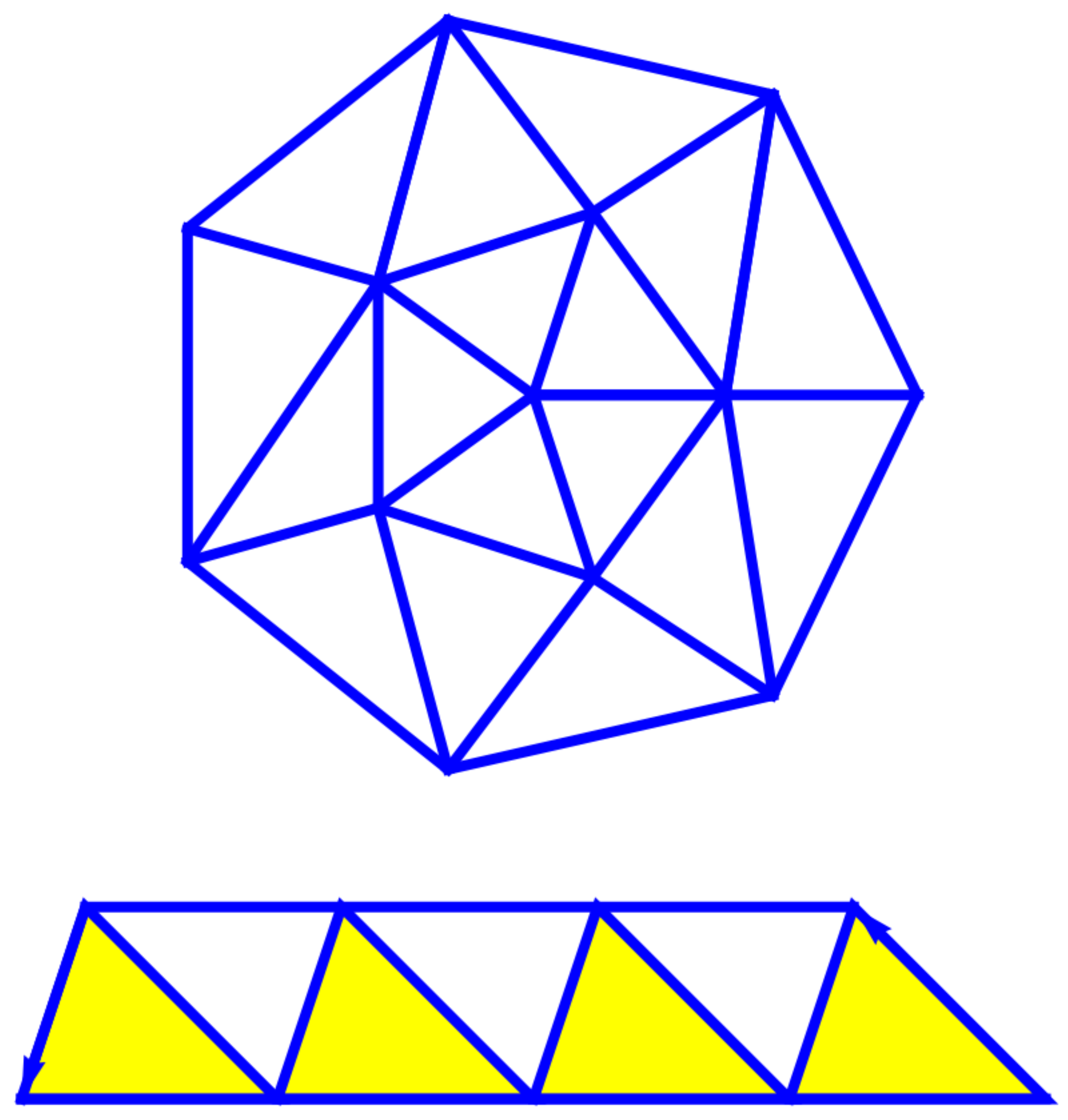}}}
}
\caption{
Discrete projective planes with chromatic number $3,4$ or $5$. 
In each case, we have a geometric graph obtained by taking a 
discrete ball and glue the lower M\"obius strip to the boundary.
(examples due to Jenny Nitishinskaya).
Literature search later revealed \cite{Fisk1977b} (page 231). 
}
\end{figure}

\begin{figure}[h]
\scalebox{0.3}{\includegraphics{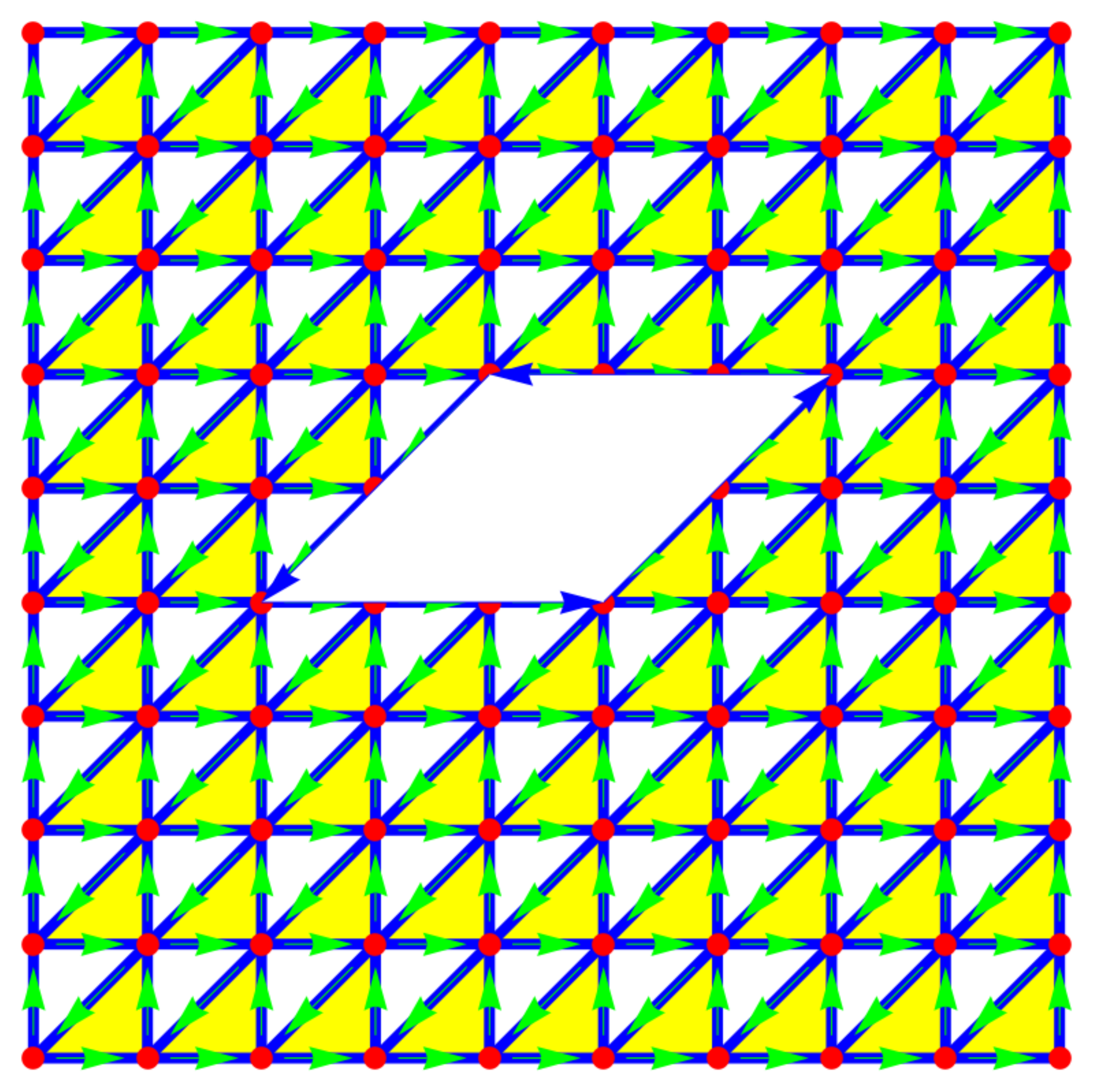}}
\scalebox{0.3}{\includegraphics{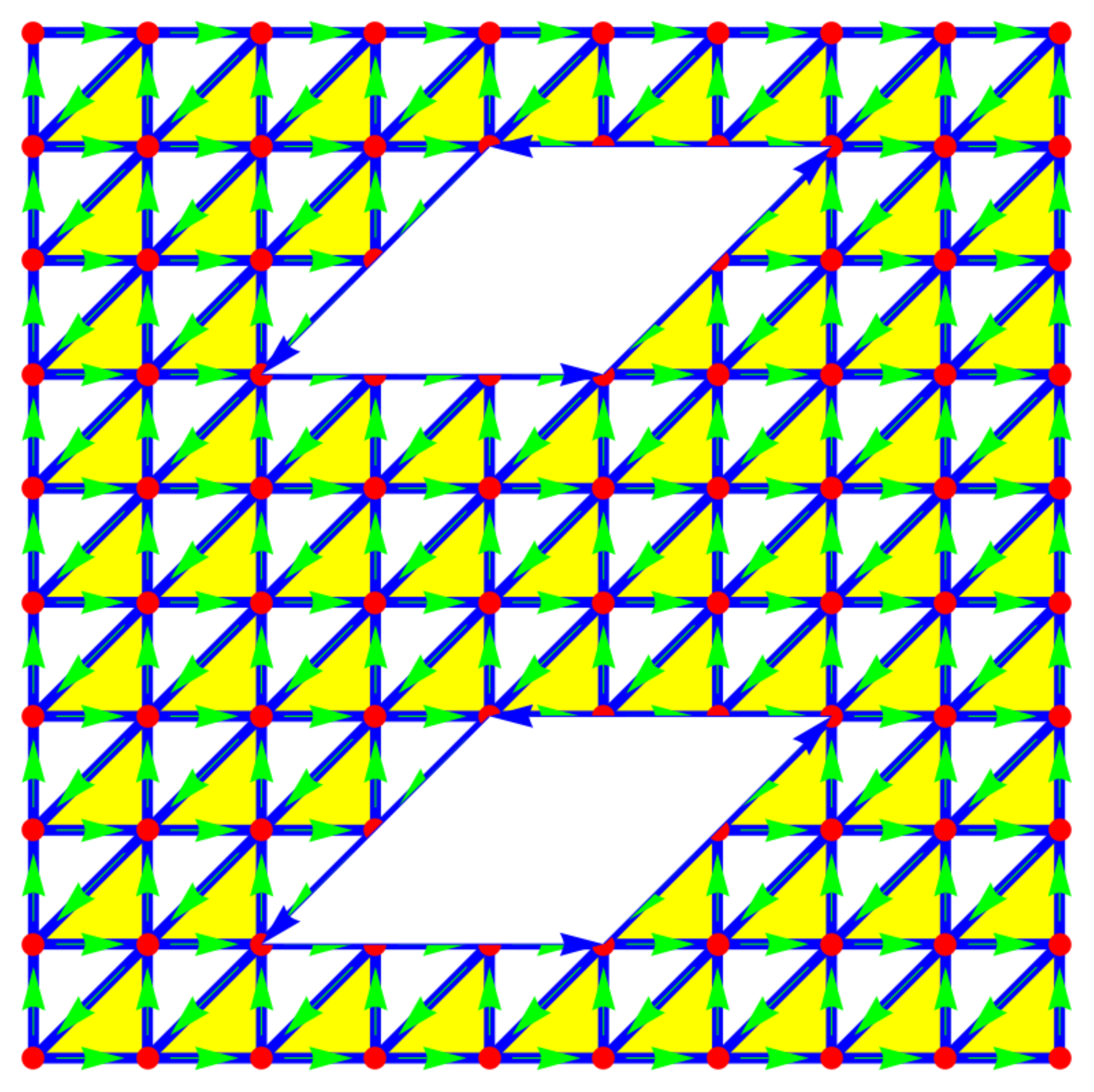}}
\caption{
In order to get higher genus $G \in \Gcal_2$, cut holes from a flat torus and identify the cyclic
boundaries. We can so get orientable higher genus graphs in $\Ccal_3$. 
to get a non-orientable surface cut a hole into a flat torus and glue in a  Moebius strip.
We see that $\Gcal_2 \cap \Ccal_3$ is not empty in every topological type. 
Doing a diagonal flip somewhere in the flat part produces examples of 
$\Gcal_2 \cap \Ccal_4$ in every topological class. We have only seen 
examples with chromatic number $5$ in the class of projective planes so far. 
\label{orientable}
}
\end{figure}

\begin{figure}[h]
\scalebox{0.30}{\includegraphics{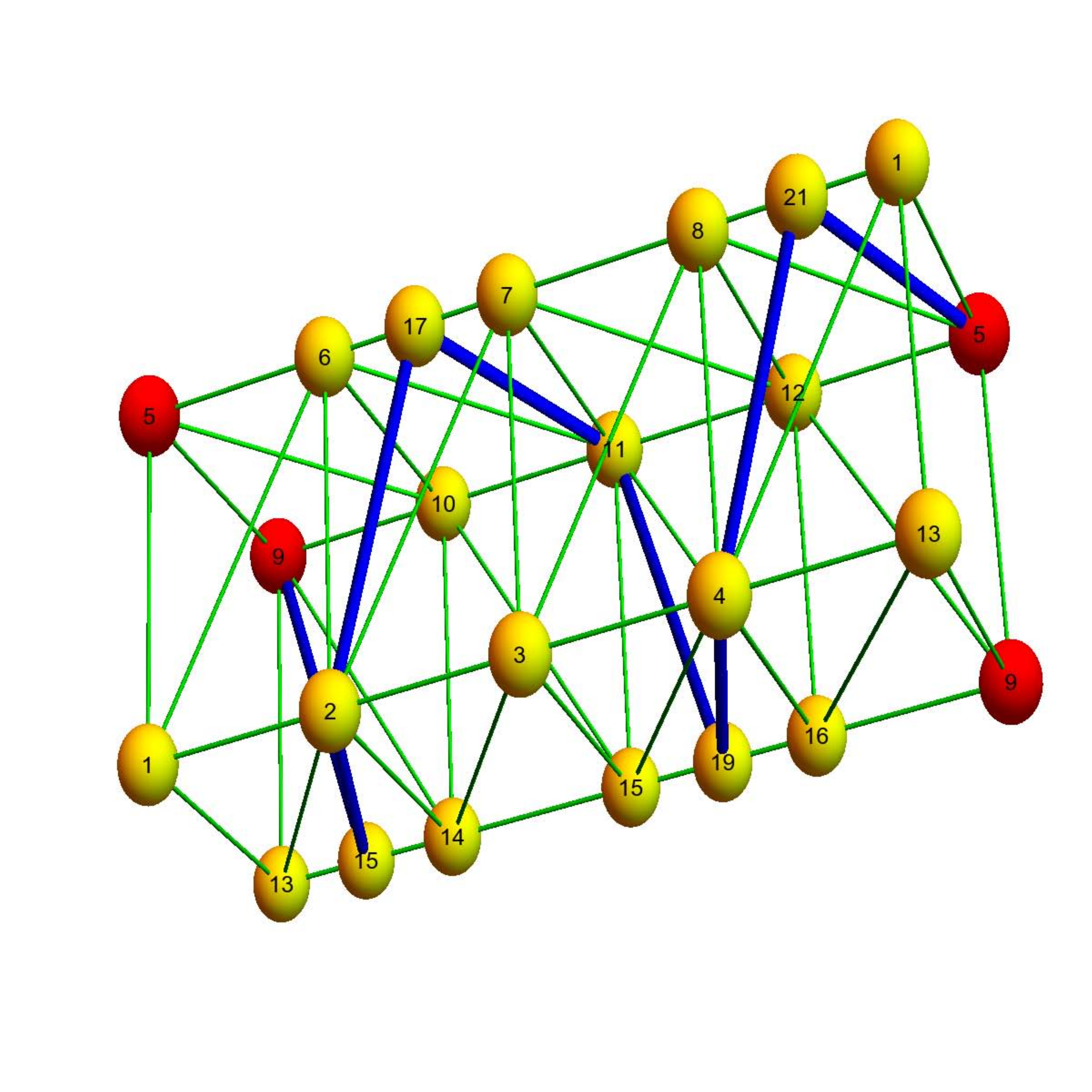}}
\caption{
A {\bf Fisk graph} is a toroidal graph in $\Gcal_2$. It is a triangularization of the torus
and chromatic number $5$. Receipt: take flat $4 \times 4$ torus, cut four edges so that
the graph has a closed geodesic $\gamma$ of winding type $(2,1)$ and length $8$,
Dehn twist the torus by $1$ along a shortest closed geodesic of type $(0,1)$
and length $4$ so that the now broken $\gamma$ has two end vertices $A="5",B="9"$ of odd degree
which are neighbors. The odd degree prevents already chromatic number
$3$ but since the coloring propagate uniquely along $\gamma$,
the 4 coloring of $S(A)$ and $S(B)$ are not compatible. We need a fifth color.
}
\end{figure}

\clearpage

\begin{figure}[h]
\scalebox{0.35}{\includegraphics{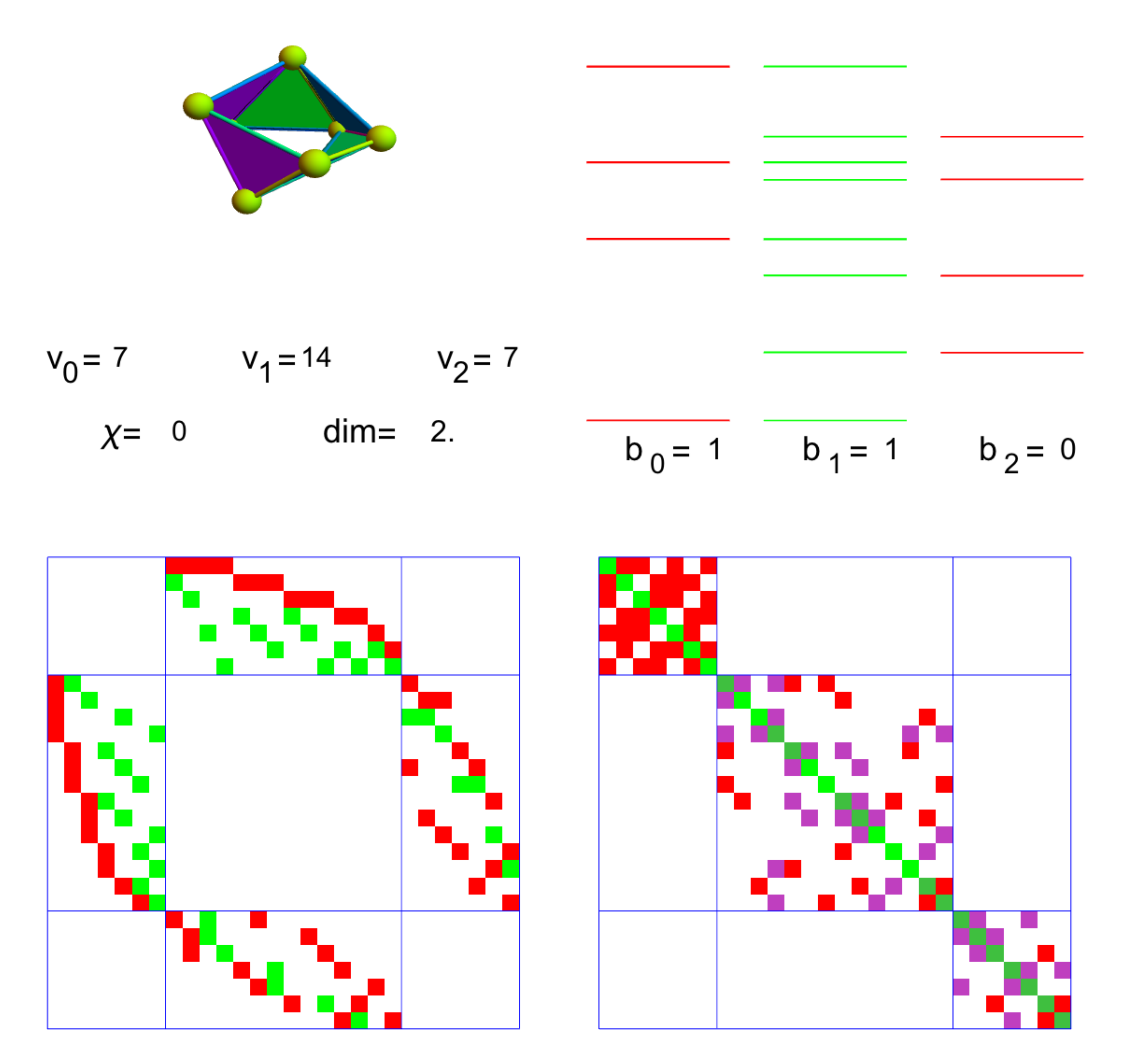}}
\caption{
The Dirac matrix $D$, the Laplacian $L$ and the spectrum $\sigma(L)$ of the 
Moebius graph $M$. These pictures appeared already (for other graphs however)
in  \cite{knillmckeansinger}. They illustrate the Euler-Poincar\'e formula
and the spectra of the form Laplacians. 
The union of the bosonic spectrum
$\sigma(L_0) \cup \sigma(L_2)$ is the same than the Fermionic spectrum
$\sigma(L_1)$. Isospectral deformations of the Dirac operator lead to an
expansion of geometry both in the Riemannian as well as in graph theory:
\cite{IsospectralDirac,IsospectralDirac2} and makes supersymmetry hard to
detect and also classical ``physics" is not affected as the Laplacian $L$
is unchanged. 
}
\end{figure}

\begin{figure}[h]
\scalebox{0.35}{\includegraphics{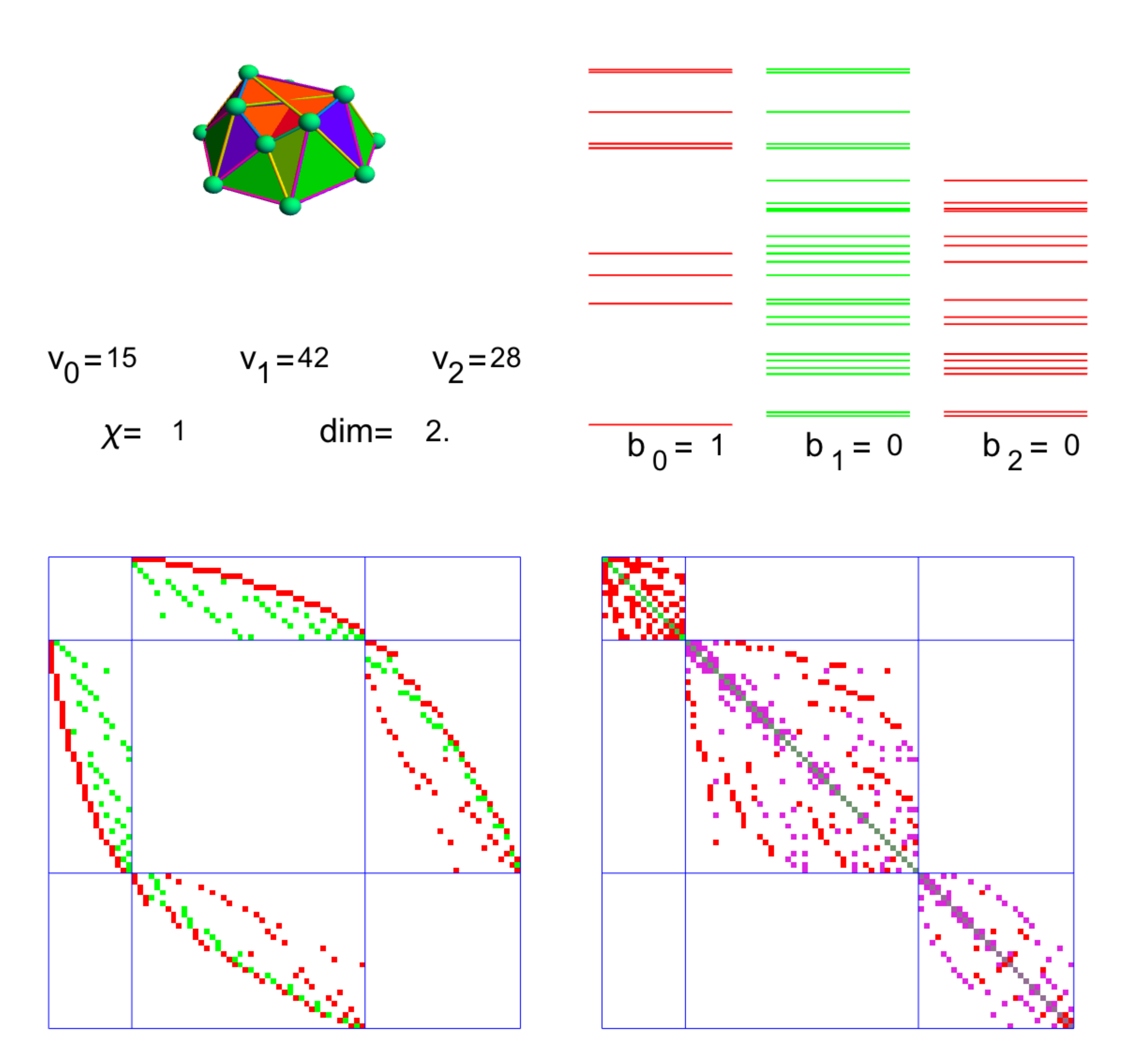}}
\caption{
Again the Dirac matrix $D$, the Laplacian $L$ and the spectrum $\sigma(L)$ of the
projective plane graph $P$. As for all graphs, supersymmetry holds in the sense
that the union of {\bf Bosonic eigenvalues} $\sigma(L_{2k})$ is the same than the
union of {\bf Fermionic eigenvalues} $\sigma(L_{2k+1})$. It is the reason for 
McKean-Singer. 
}
\end{figure}

\begin{figure}[h]
\scalebox{0.28}{\includegraphics{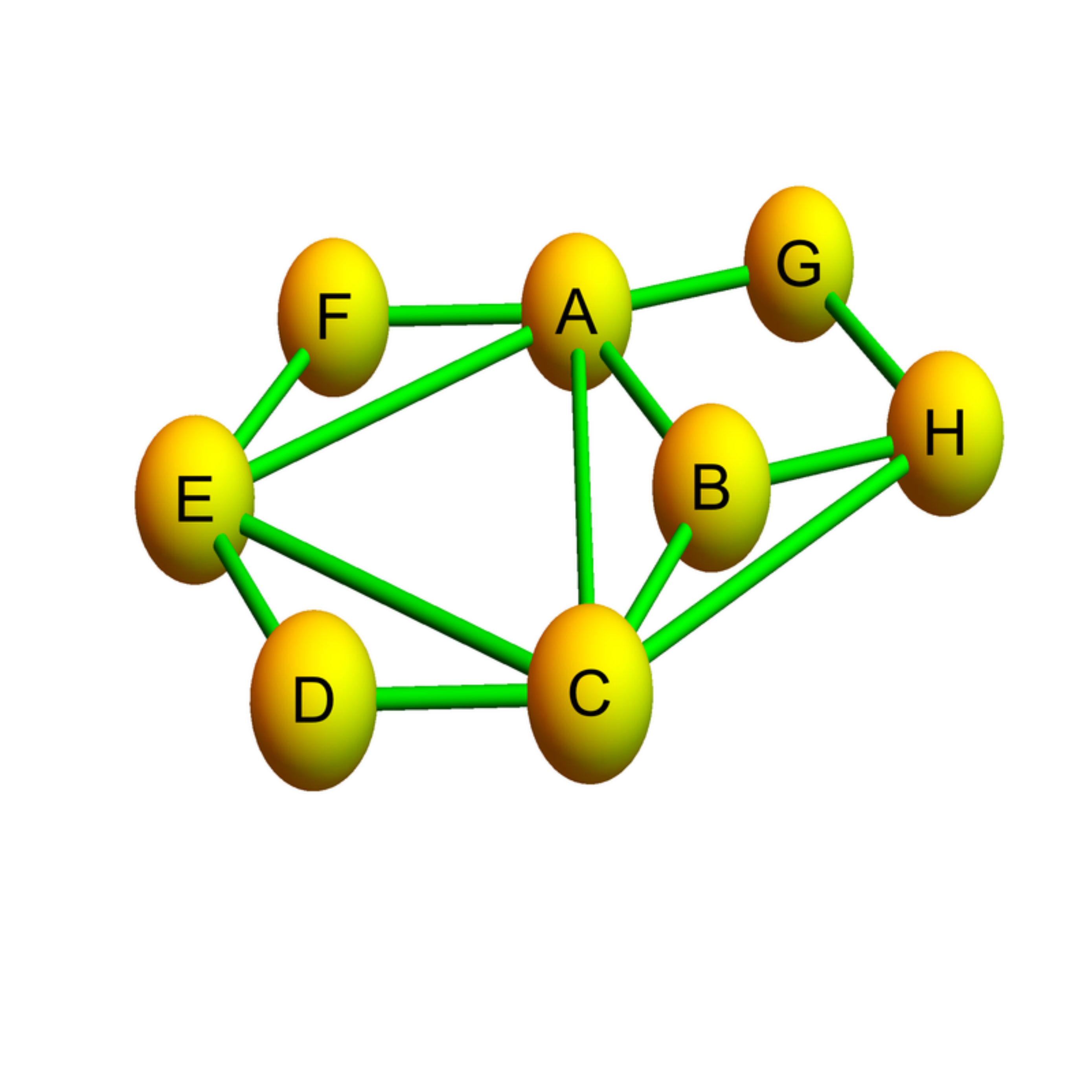}}
\scalebox{0.22}{\includegraphics{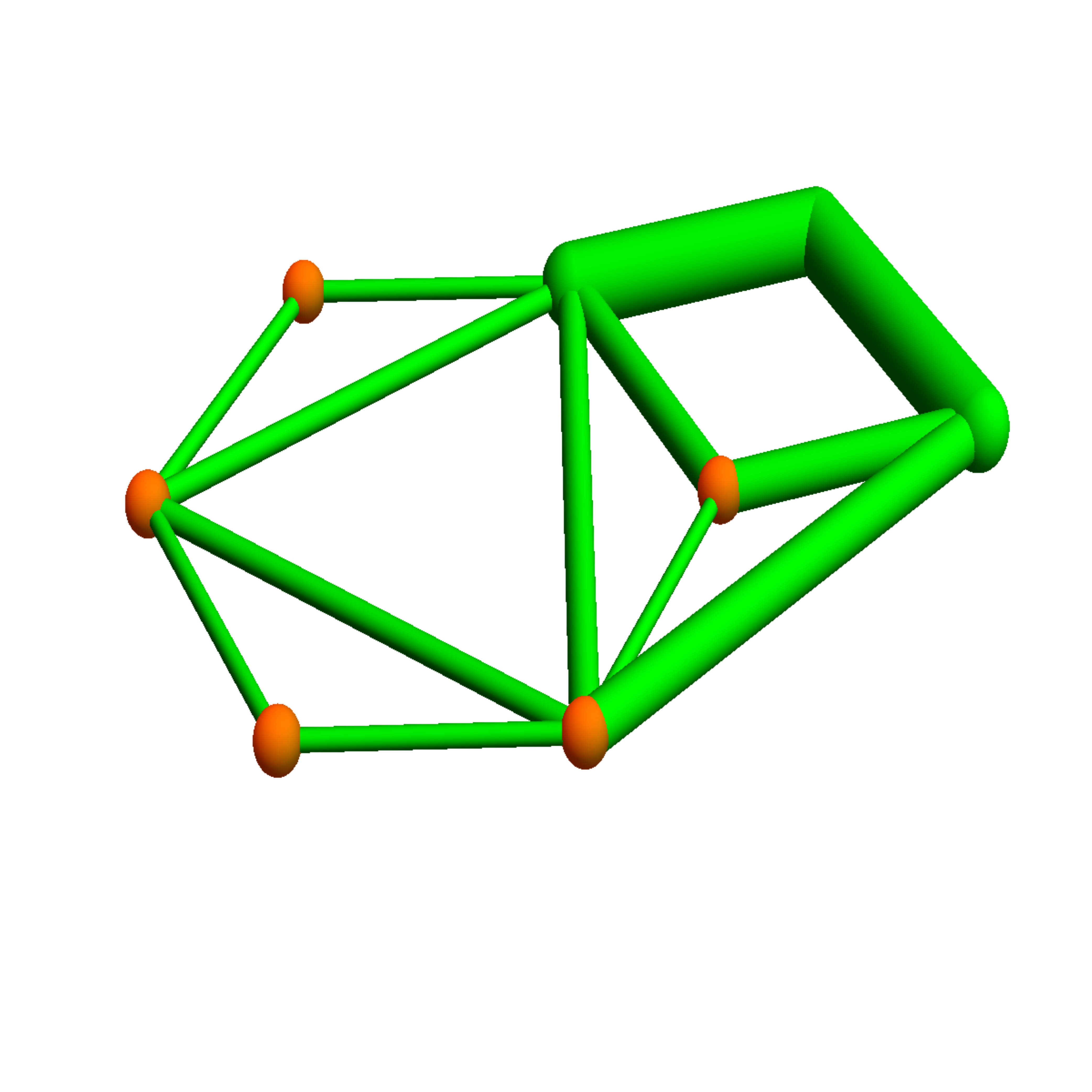}}
\caption{
The {\bf duke graph} $G$ from \cite{BergeOulipo5}
has $\pi_1(G)=Z, b_0=b_1=1$. The image of the Hurewicz map from $\pi_1(G)$ to $H^1(G)$
is obtained by taking the characteristic function of a path $\gamma$,
a one form $f$ and apply the heat flow $e^{-t L_1}$ to it.
It will converge for $t \to \infty$ to a harmonic element
$(4, 8, 16, 19, -47, 4, 8, 4, -4, 3, 25, 22, 47)$.
It is a representative of the cohomology class which 
spans the real vector space $H^1(G)$.
In the second picture the thickness of the cylinders depends on the value
which the harmonic $1$-form takes on the edges. We see that it is concentrated
around the noncontractible loop. 
}
\end{figure}

\begin{figure}[h]
\scalebox{0.19}{\includegraphics{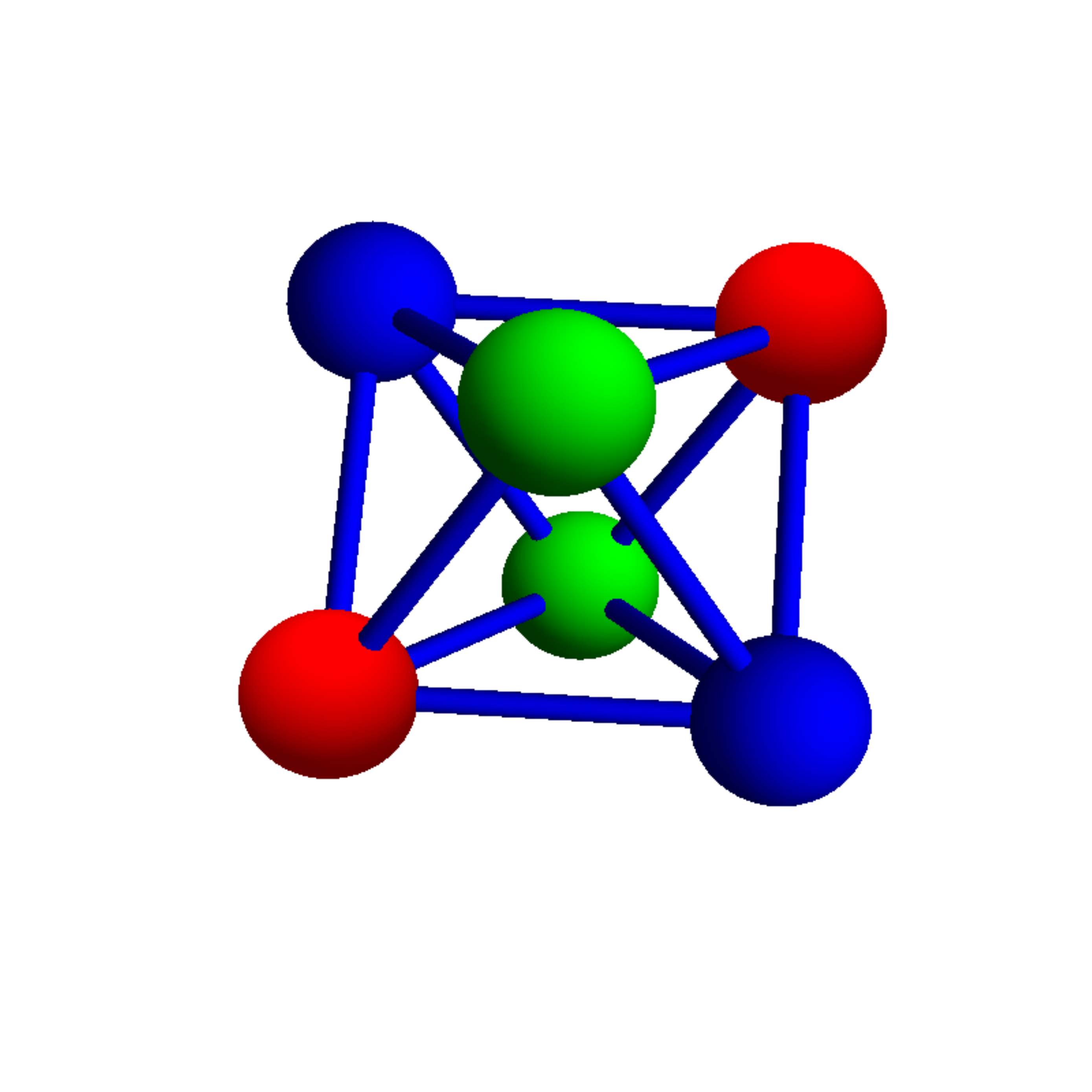}}
\scalebox{0.19}{\includegraphics{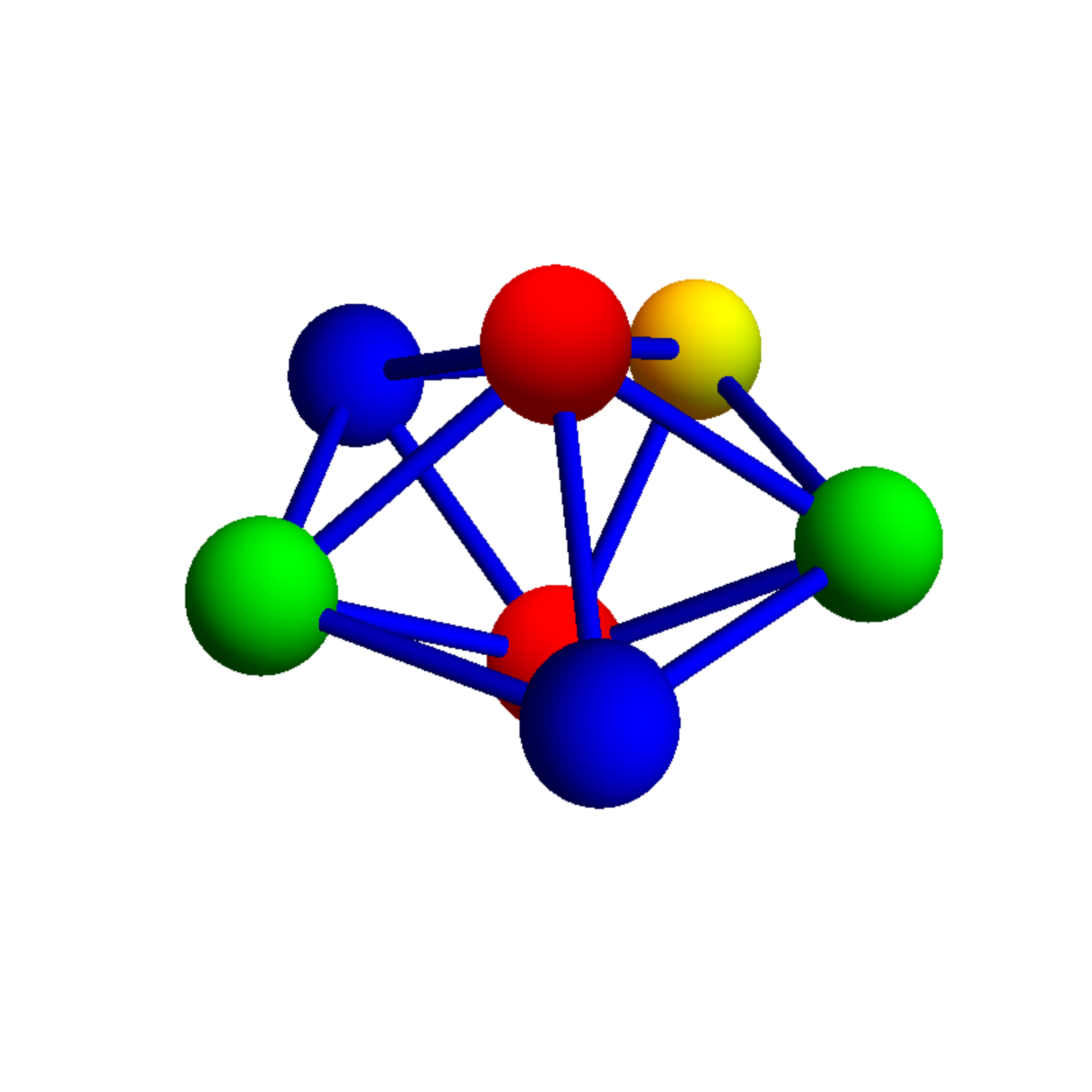}}
\scalebox{0.19}{\includegraphics{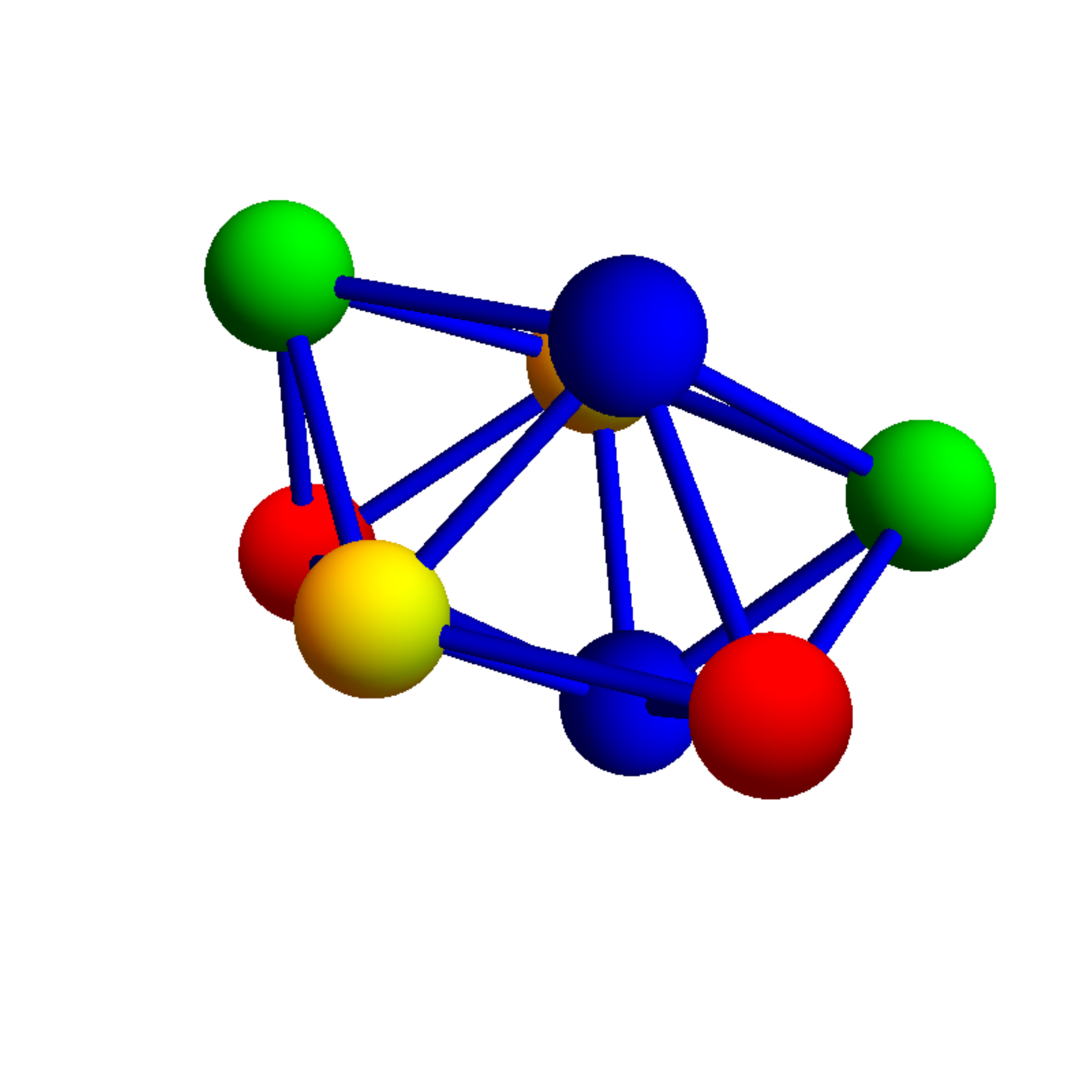}}
\scalebox{0.19}{\includegraphics{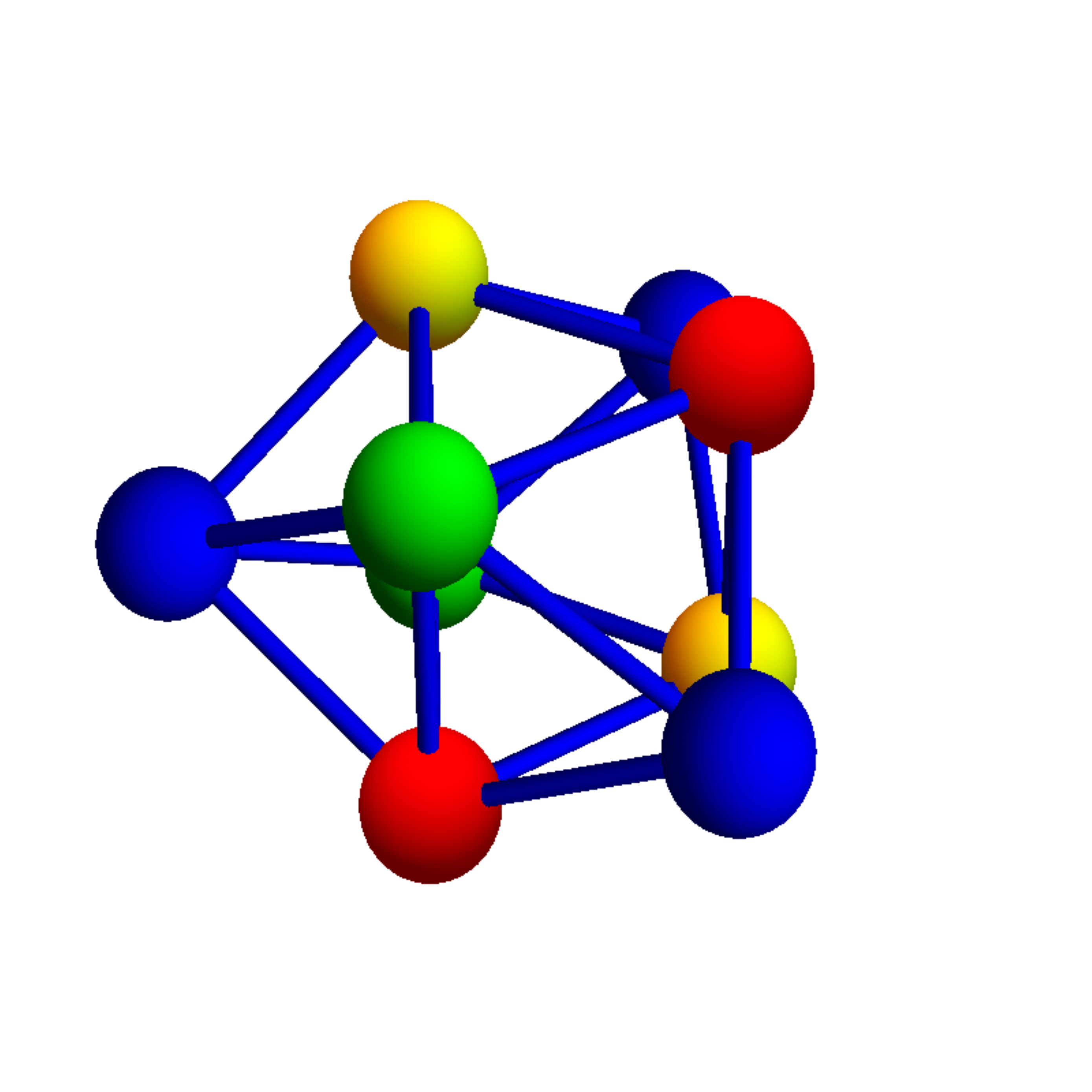}}
\scalebox{0.19}{\includegraphics{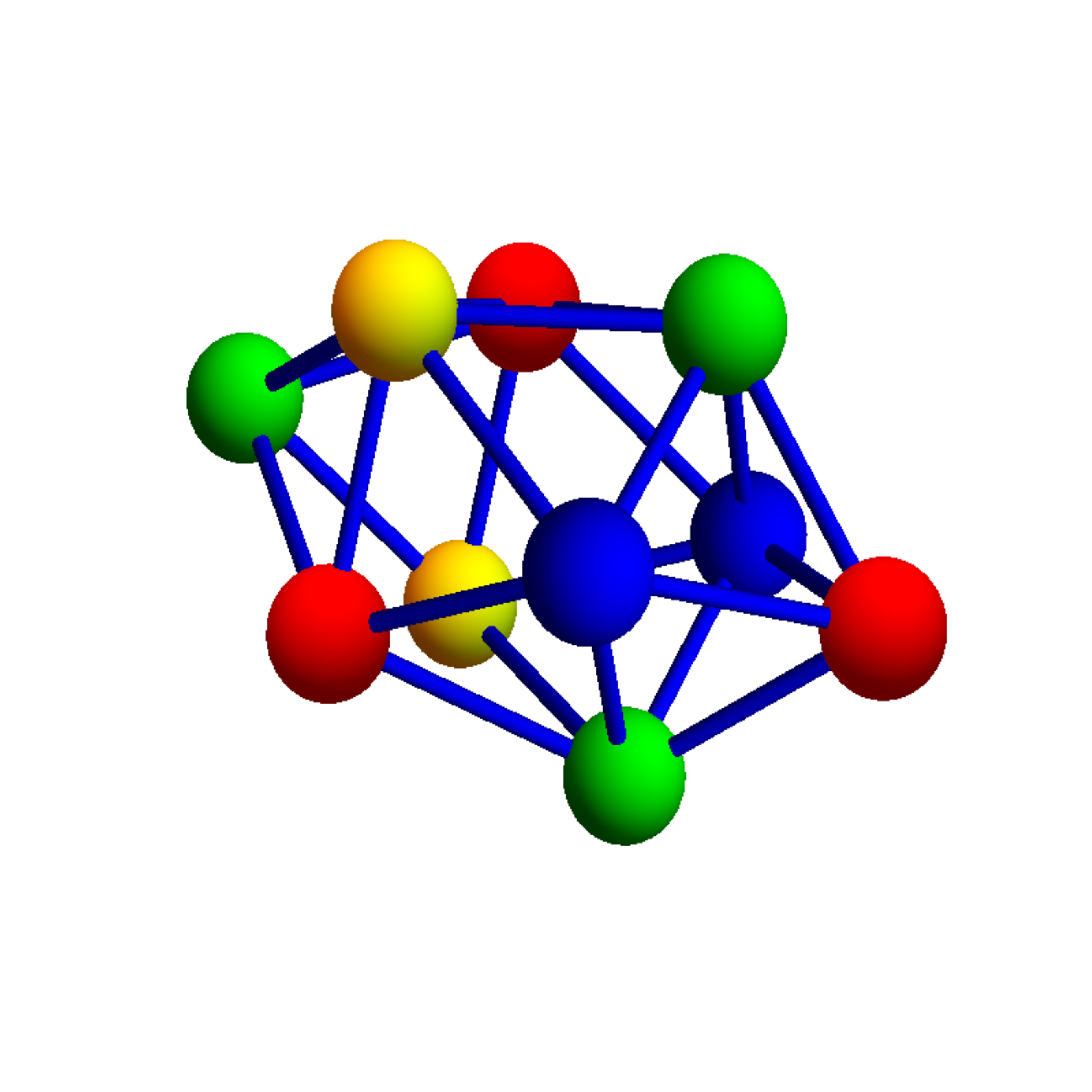}}
\scalebox{0.19}{\includegraphics{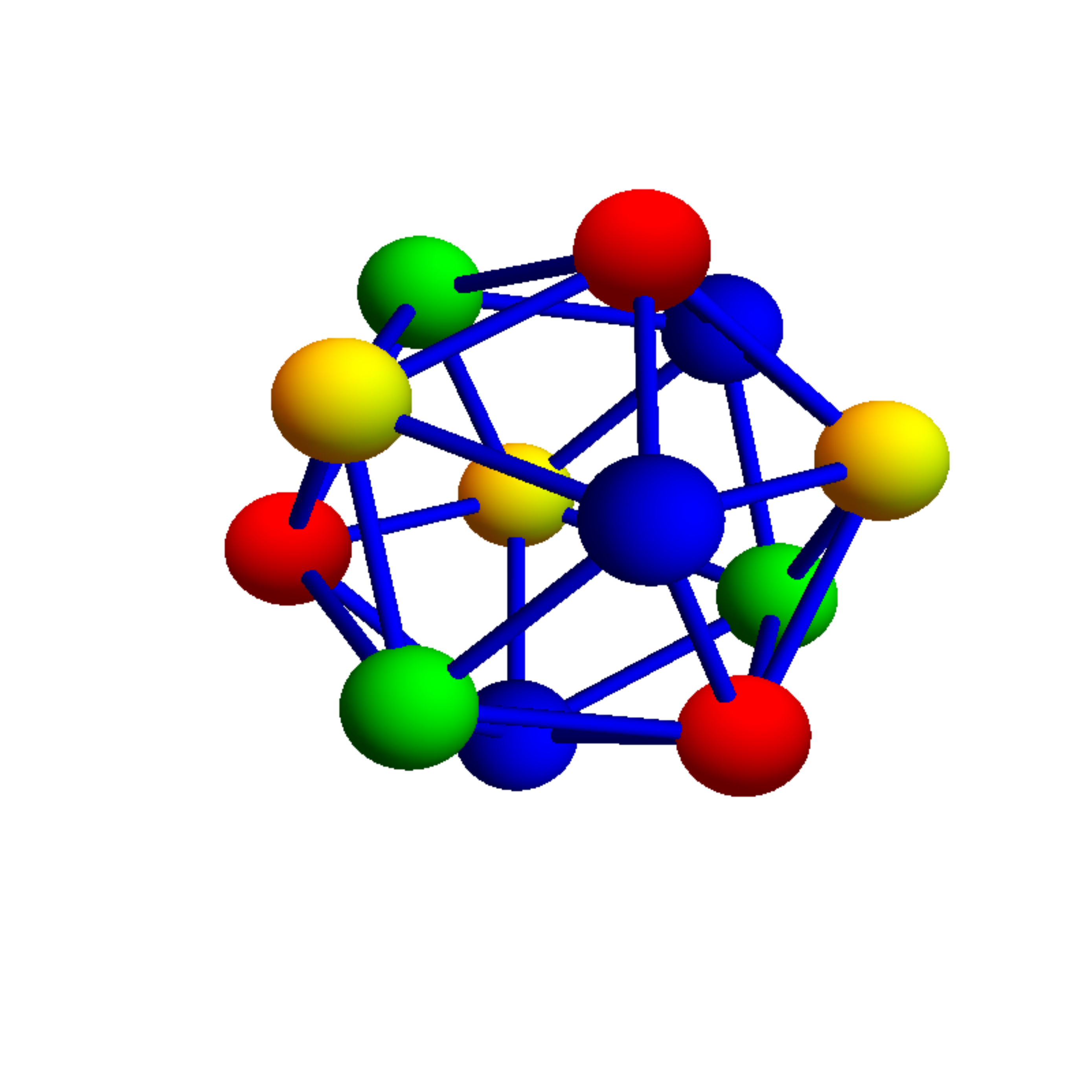}}
\caption{
Up to isomorphism, there are 6 positive curvature graphs in $\Gcal_2$. 
The smallest is the octahedron the largest the icosahedron. The graph with 
curvatures $1/3,1/3,1/3,1/3,1/6,1/6,1/6,/1/6$ is the {\bf Fritsch 
graph}, one of the smallest counter example to the Kempe proof. Only the octahedron
is in $\Ccal_3$. It has a unique coloring up to color permutation.
For the others, the color richness \cite{KnillFunctional}, (the 
chromatic polynomial evaluated at $4$ divided by the number of color
permutations $c(4)/4!$) is $5,3,2,8$ and $10$. 
}
\end{figure}

\clearpage

\begin{figure}[h]
\scalebox{0.1}{\includegraphics{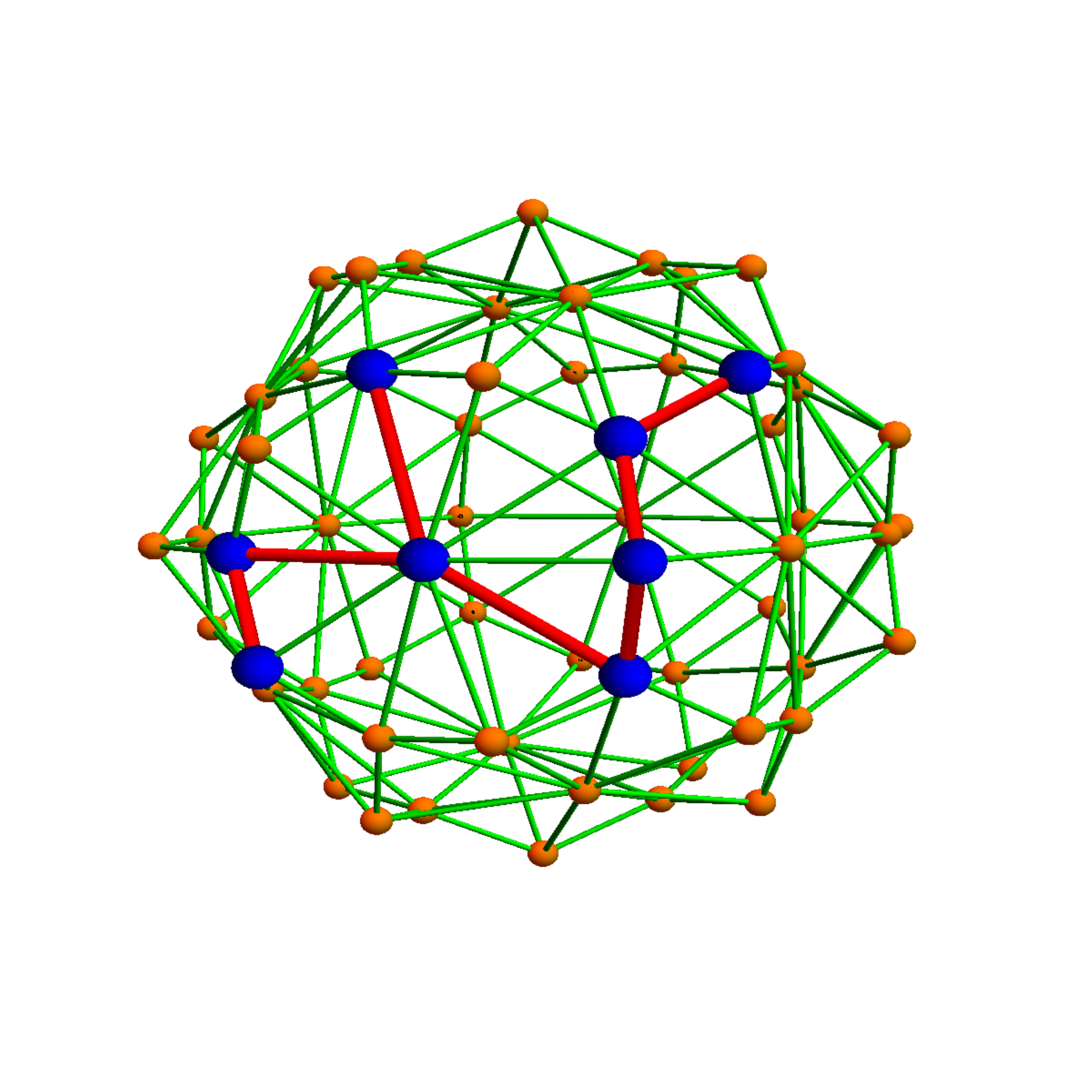}}
\scalebox{0.1}{\includegraphics{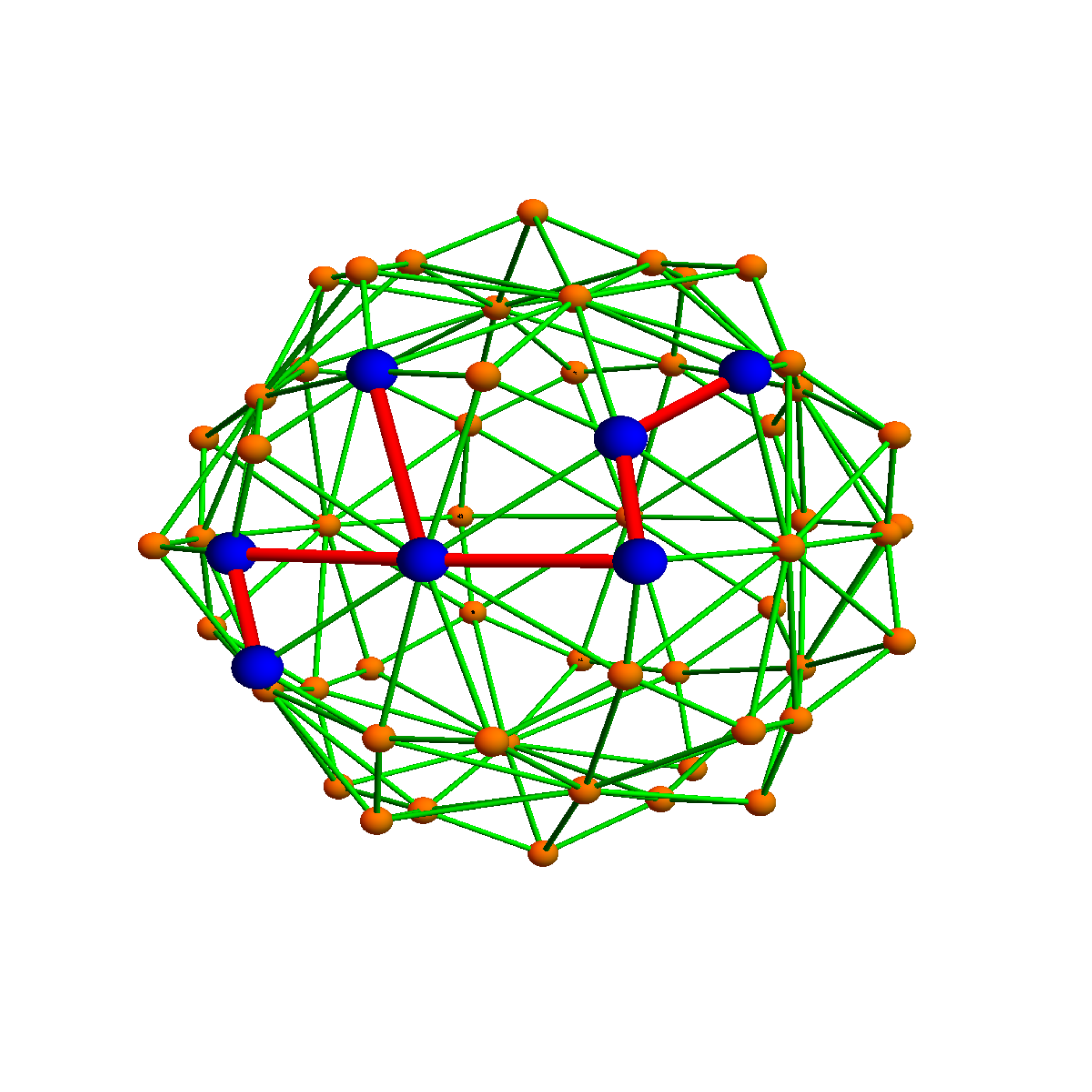}}
\scalebox{0.1}{\includegraphics{figures/fundamental2.pdf}}
\scalebox{0.1}{\includegraphics{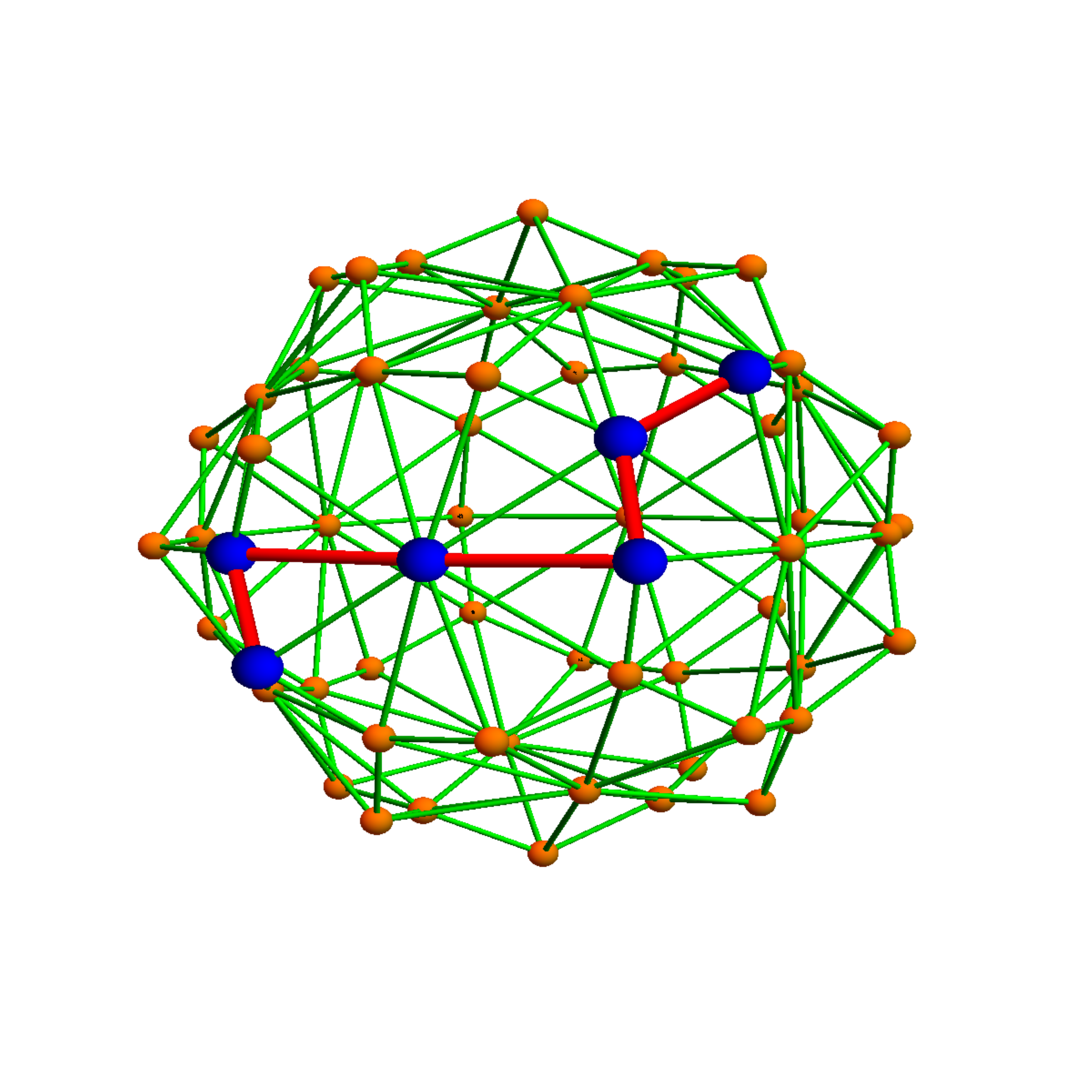}}
\caption{
Homotopy deformation steps $T$ and $S$ for curves in $G$.
The transformations $T,T^{-1},S,S^{-1}$ are applied to closed graphs
through a vertex $x_0$. }
\end{figure}

\begin{figure}[h]
\scalebox{0.13}{\includegraphics{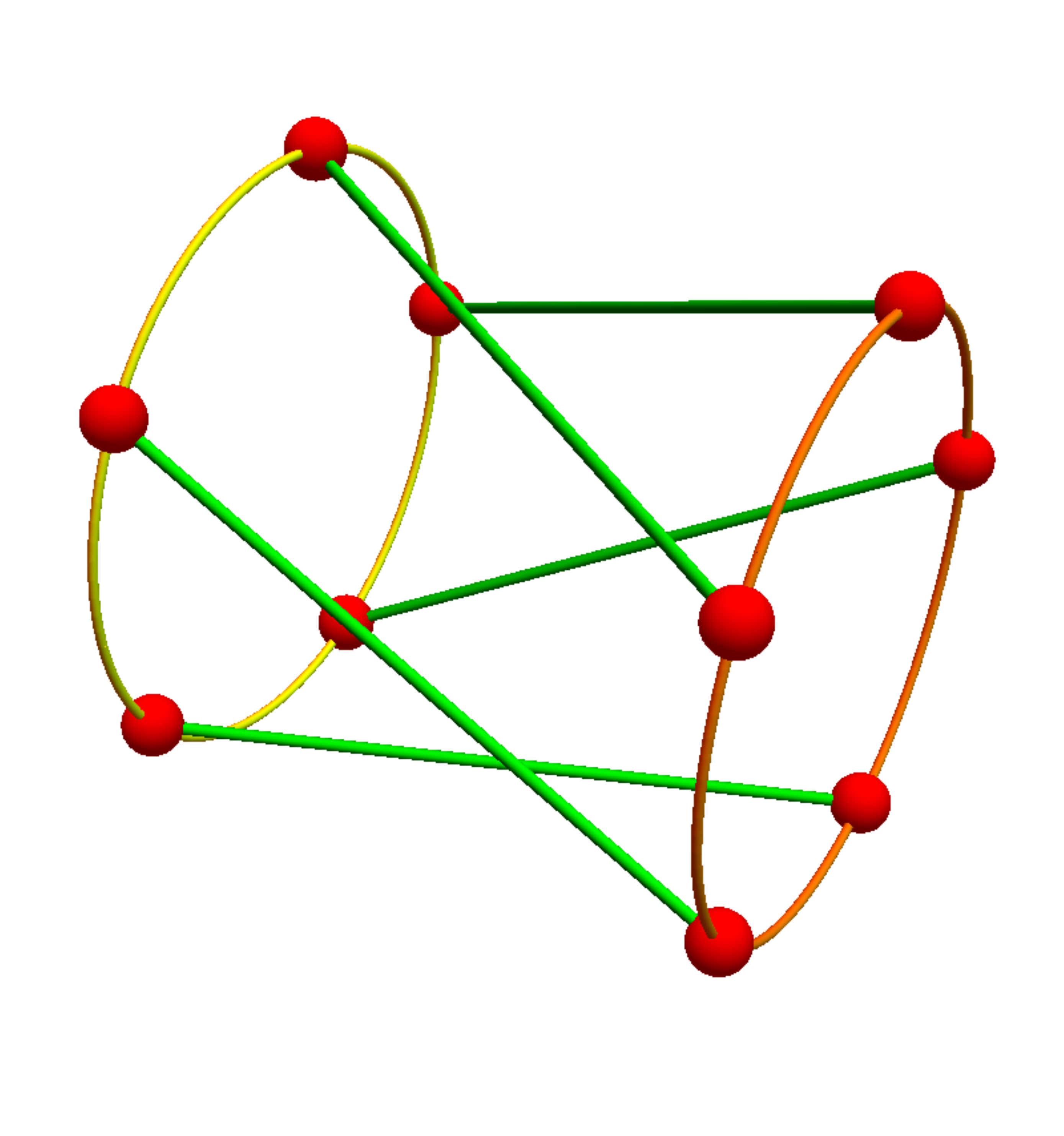}}
\scalebox{0.13}{\includegraphics{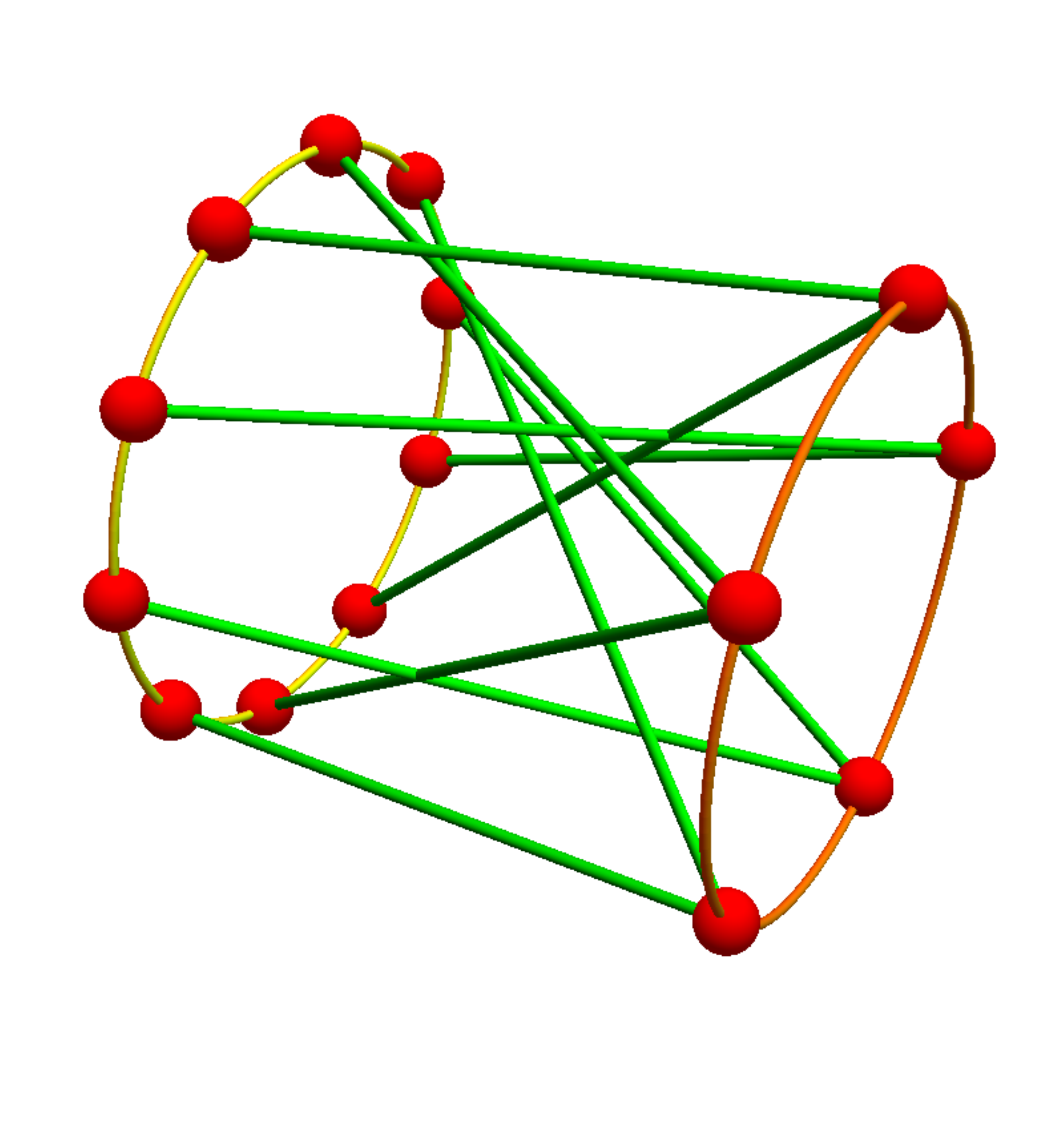}}
\scalebox{0.13}{\includegraphics{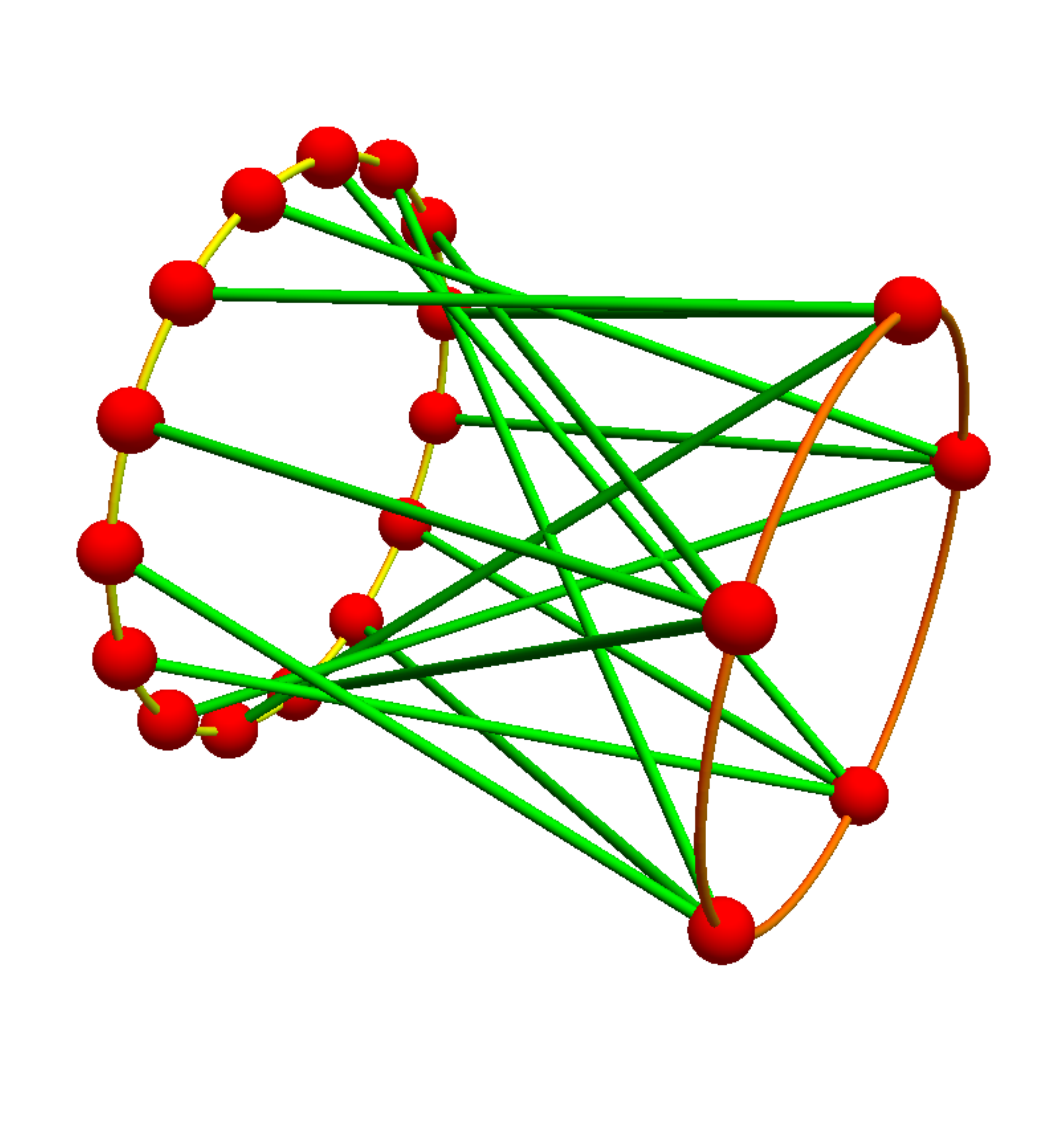}}
\caption{
A closed curve on a graph $G$ can be seen as a graph homomorphism $f$ from a circular
graph $H$ to $G$. We can encode this as a bipartite graph which,
is the disjoint union of $H$ and $G$ together with all the edges $(x,f(x))$.
We can identify the homotopy of the curve with the homotopy of the
{\bf immersion graphs}. }
\end{figure}

\begin{figure}[h]
\scalebox{0.3}{\includegraphics{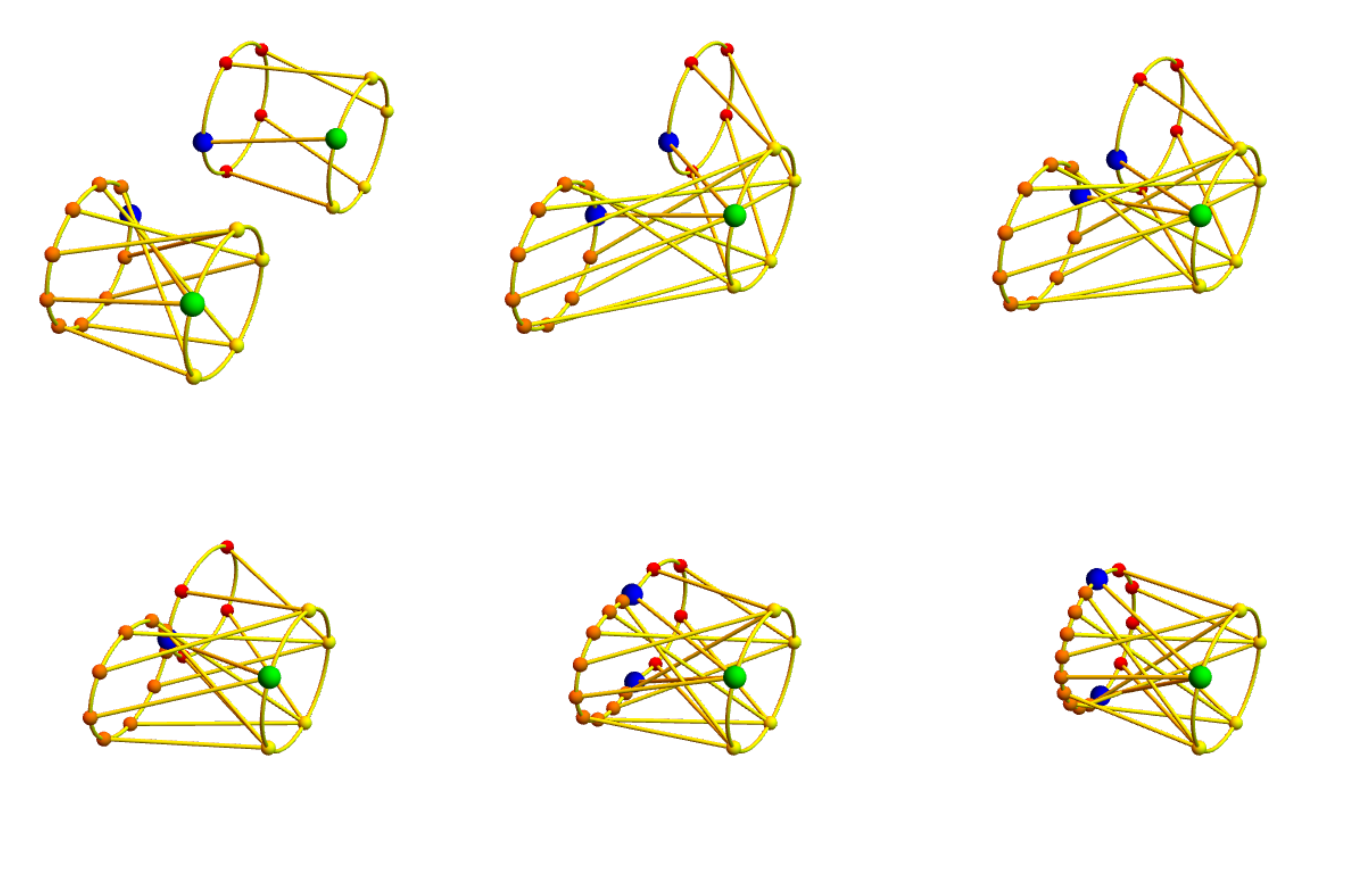}}
\caption{
The addition $3=1+2$ in the homotopy group $\pi_1(C_5)$.
We add two graph homomorphism $C_{5} \to C_5$
and $C_{10} \to C_{5}$ and get a graph homomorphism from
$C_{15} \to C_5$.
}
\end{figure}

\begin{figure}[h]
\scalebox{0.18}{\includegraphics{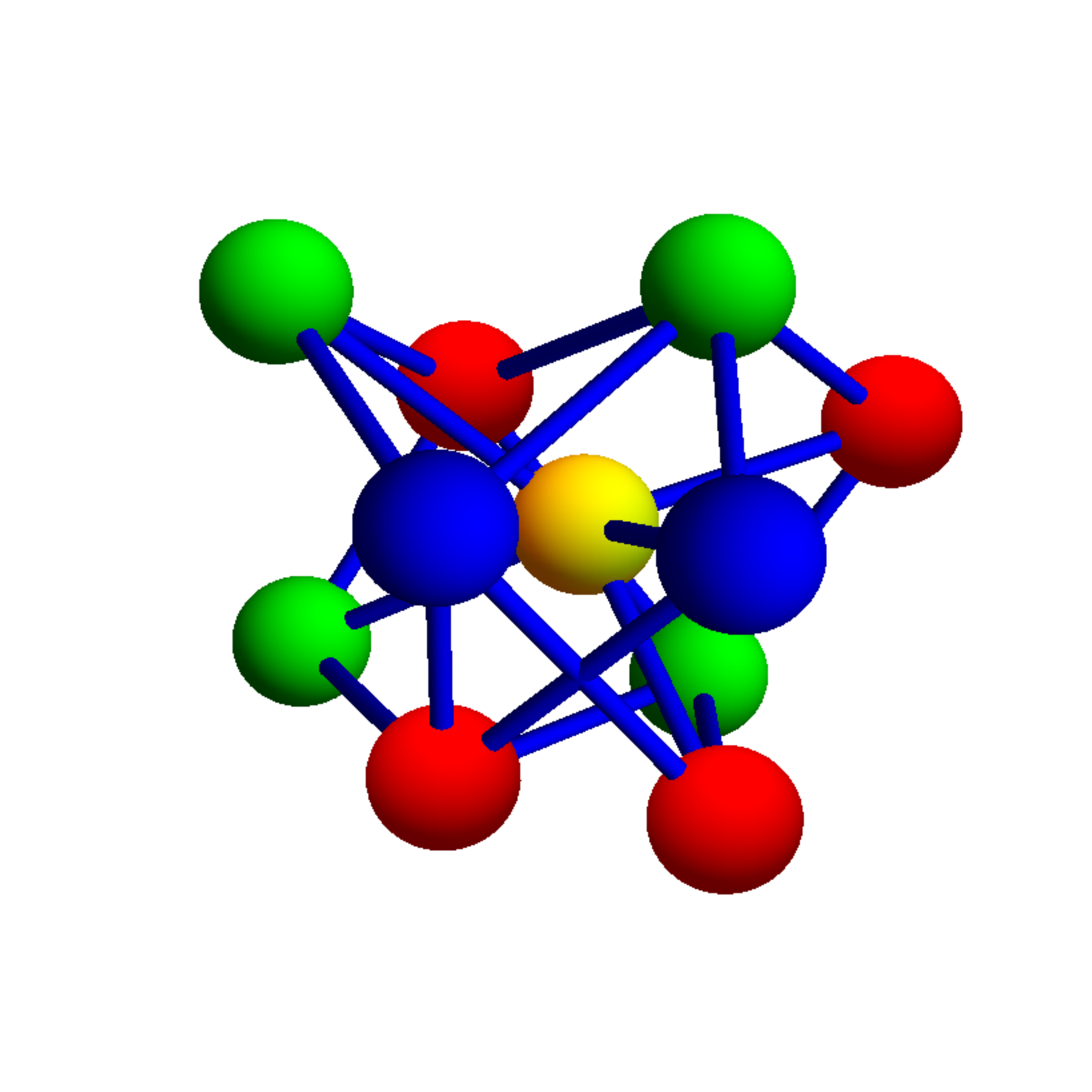}}
\scalebox{0.18}{\includegraphics{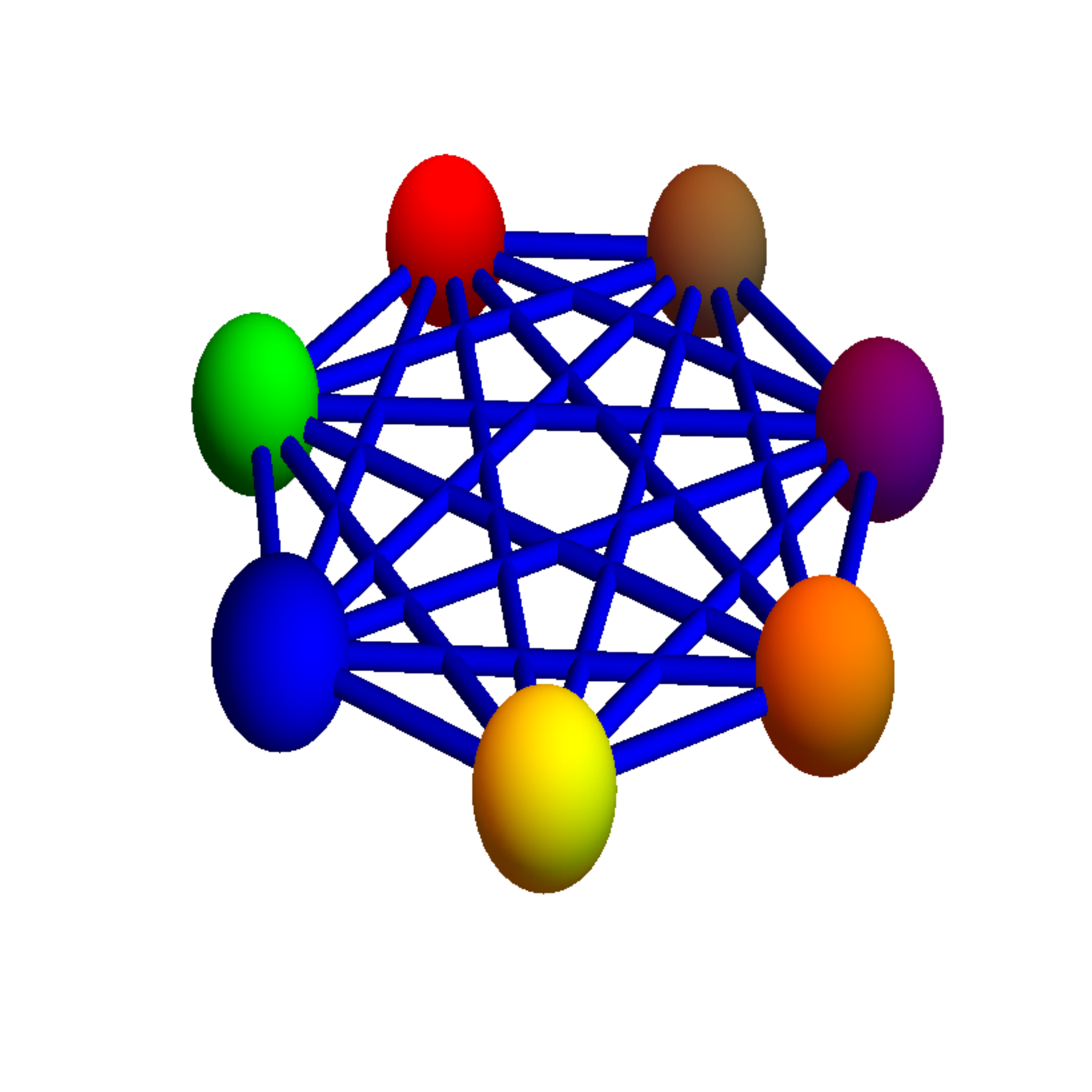}}
\caption{
The {\bf Groetsch graph} $G$ is a $1$-dimensional non-planar graph of Euler characteristic $-9$
and chromatic number $4$. It is not in $\Vcal_1$ as every vertex is a
singularity with the vertex degrees are $3,4,5$. 
$G$ is no counter example to the question whether $\Vcal_1 \subset \Ccal_3$.
To the right, $K_7$, a toroidal graph which is not in $\Ccal_6$. Even so it can be embedded on 
$T^2$ it is not in $\Gcal_2$.
\label{orientable}
}
\end{figure}

\begin{figure}[h]
\scalebox{0.20}{\includegraphics{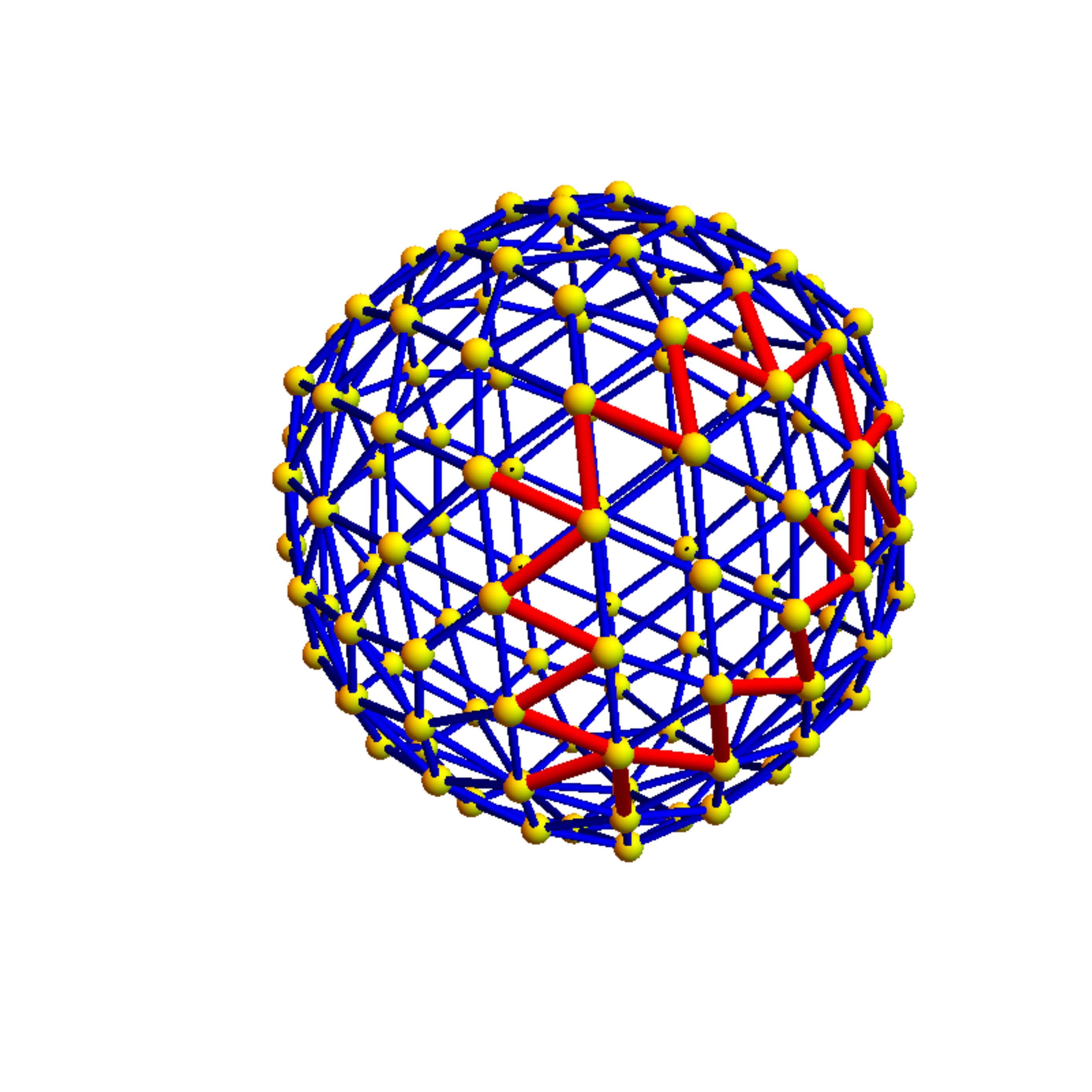}}
\scalebox{0.20}{\includegraphics{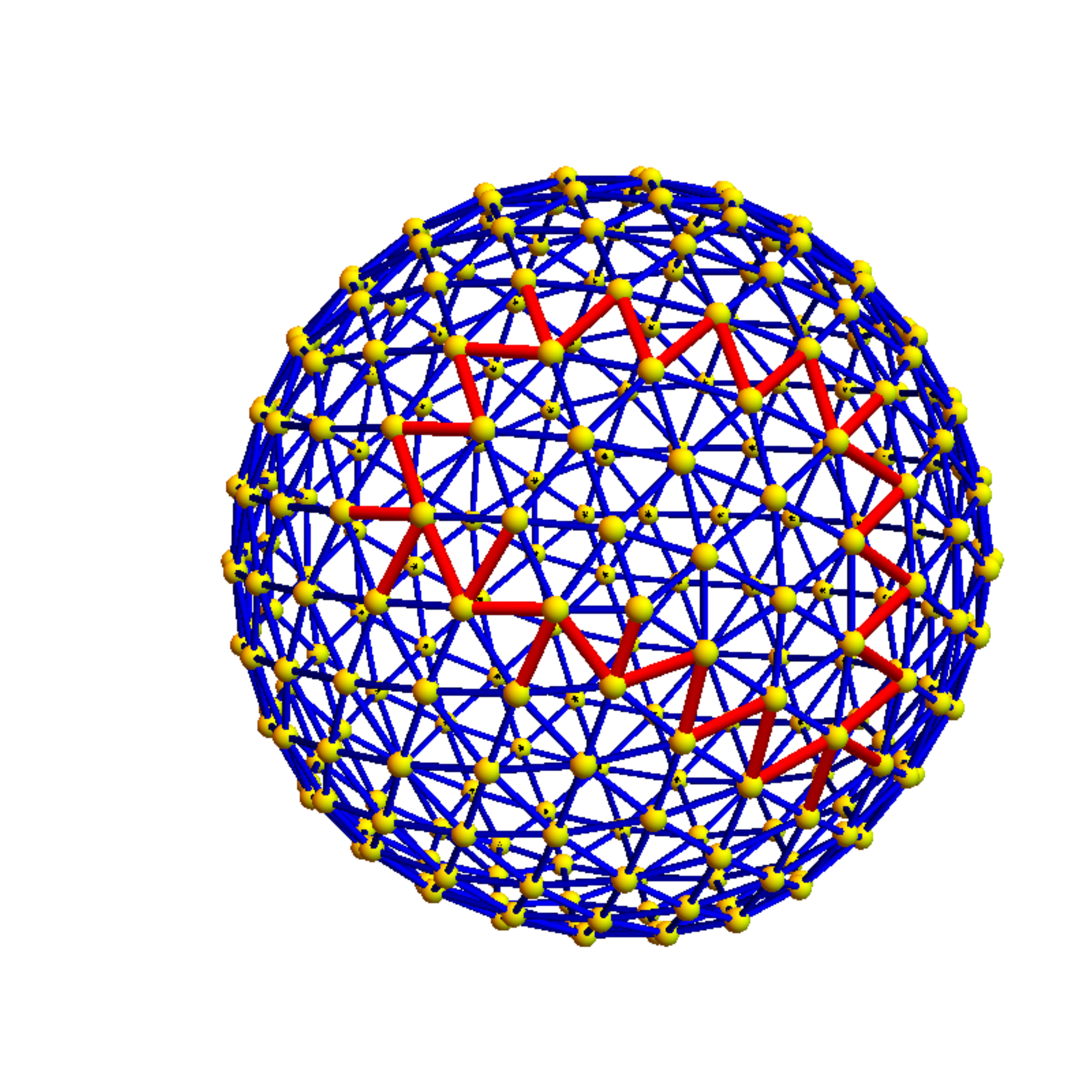}}
\caption{
Level curves $H = \{ f = c \}$ of a fullerene type graph $G$.
The graph $H$ is a subgraph of the upper line graph.
The edges contained in the level curve $f=c$ are high lighted.
}
\end{figure}

\begin{figure}[h]
\scalebox{0.23}{\includegraphics{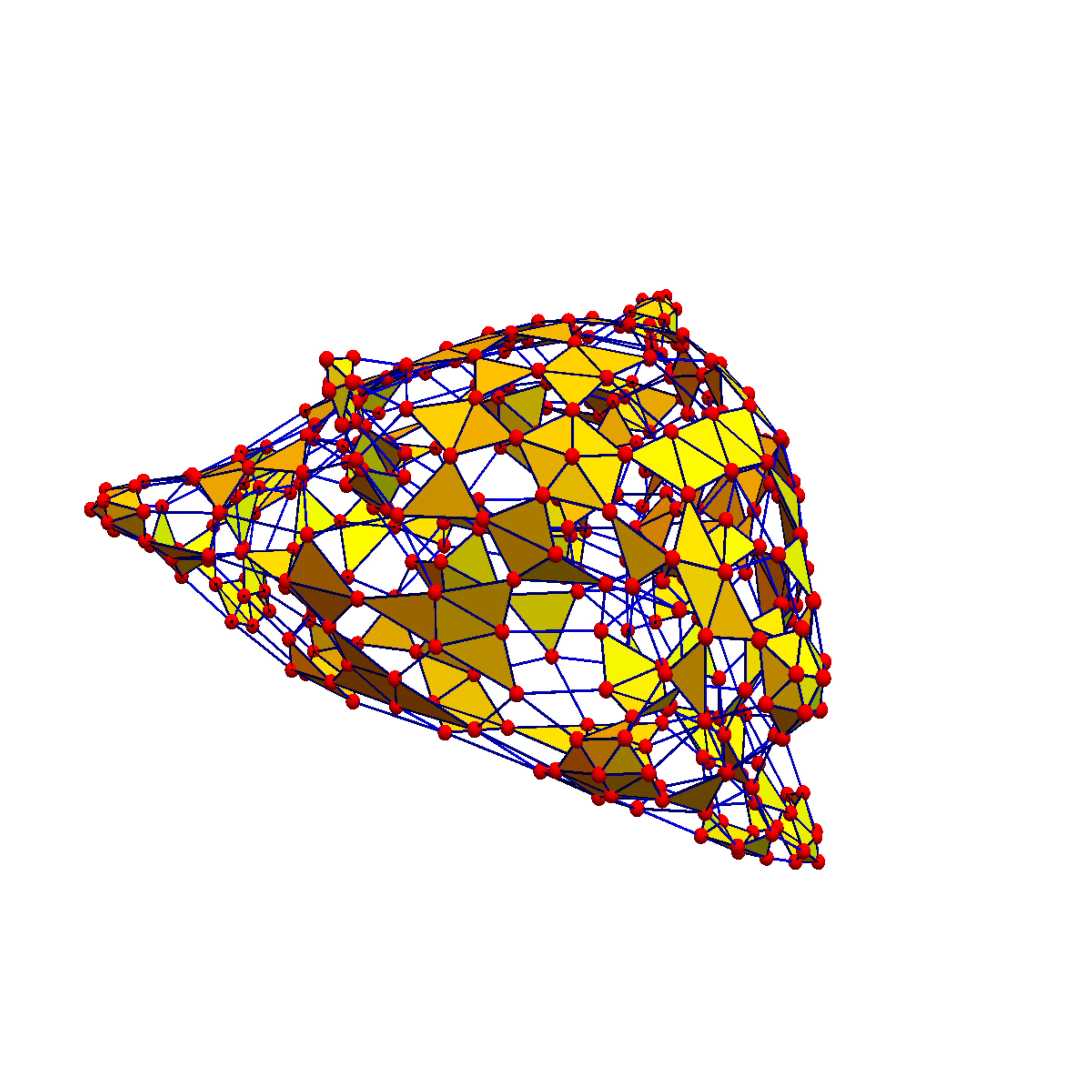}}
\caption{
A level surface $H=\{ f = c \}$ in the 600-cell $G \in \Scal_3$. 
Its faces are triangles or quadrangles. The later can be completed so
that $\overline{H} \in \Gcal_2$. 
}
\end{figure}

\begin{figure}[h]
\scalebox{0.43}{\includegraphics{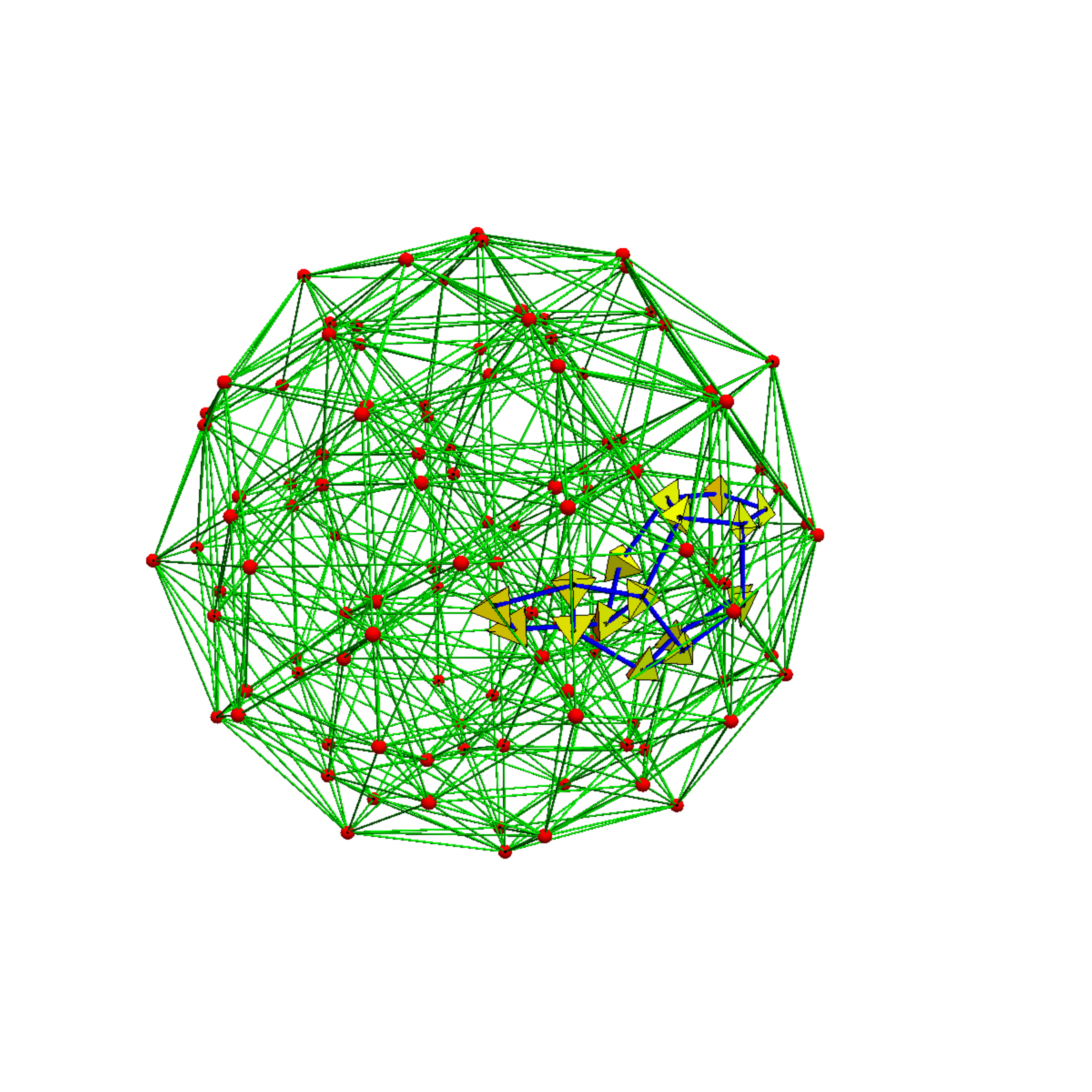}}
\caption{
The figure shows a $1$-dimensional
variety $K=\{ f = c, g= c \} \subset \Scal_1$ inside the same 600 cell host graph $G$.
It is the subgraph of the 120-cell $\hat{G}$ consisting of tetrahedra in $G$,
where both $f$ and $g$ change sign.
}
\end{figure}

\begin{figure}[h]
\scalebox{0.08}{\includegraphics{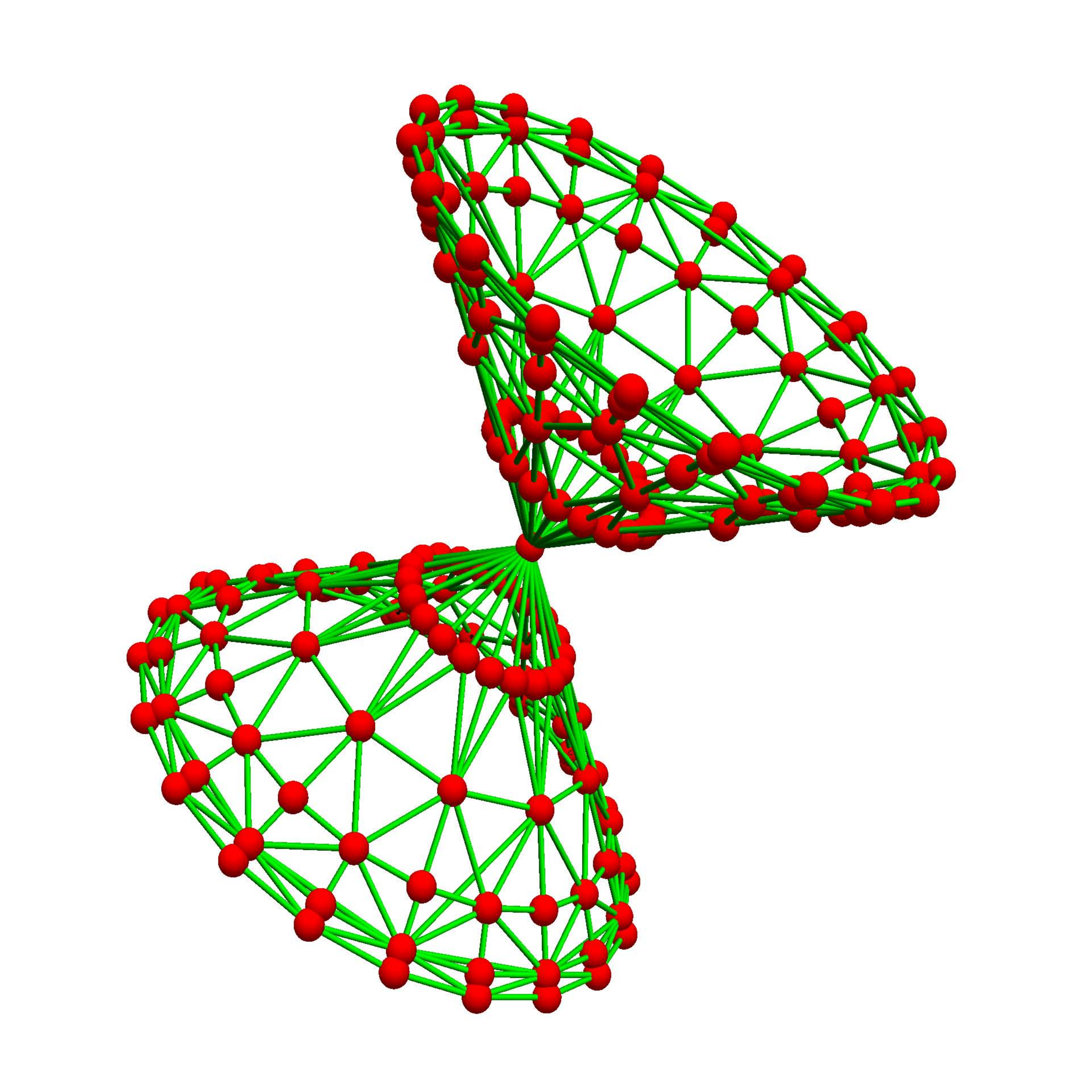}} 
\scalebox{0.28}{\includegraphics{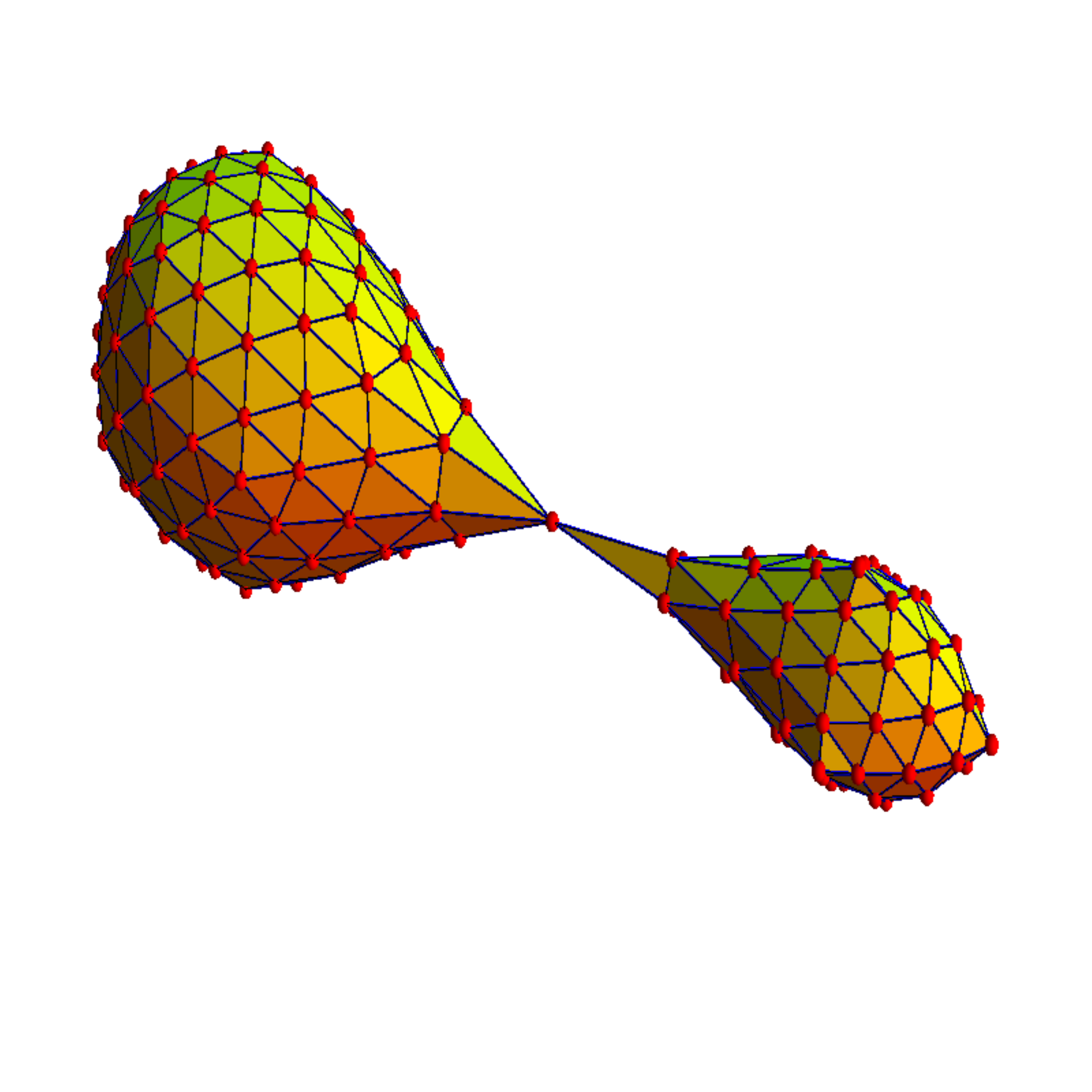}} 
\caption{
A graph $G \in \Vcal_2$ for which the singularity set $\sigma(G)$ consists of one vertex $x$
for which the unit sphere at that point is the union of two circles. All other points are regular:
their unit spheres are circular graphs. To the right a version without boundary. 
}
\end{figure}

\begin{figure}[h]
\scalebox{0.24}{\includegraphics{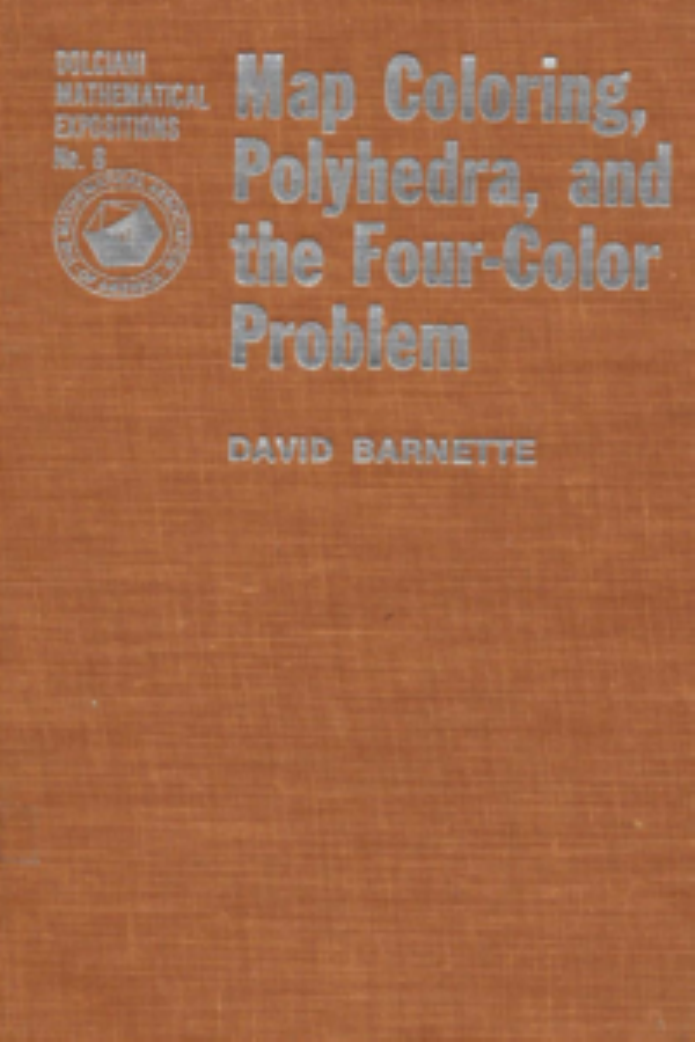}}
\scalebox{0.24}{\includegraphics{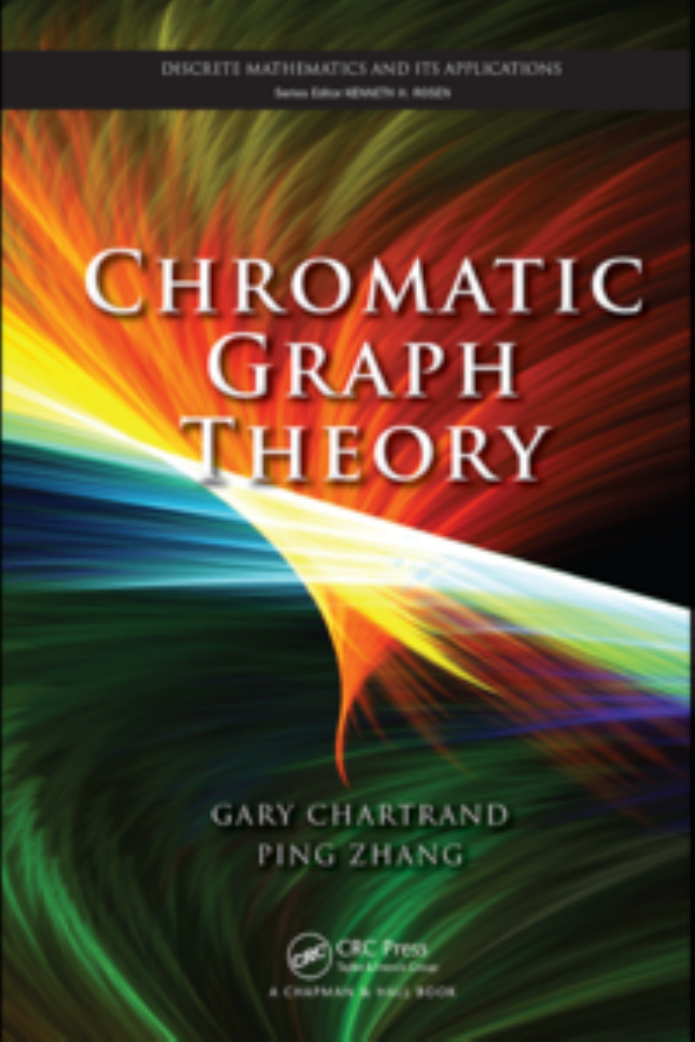}}
\scalebox{0.24}{\includegraphics{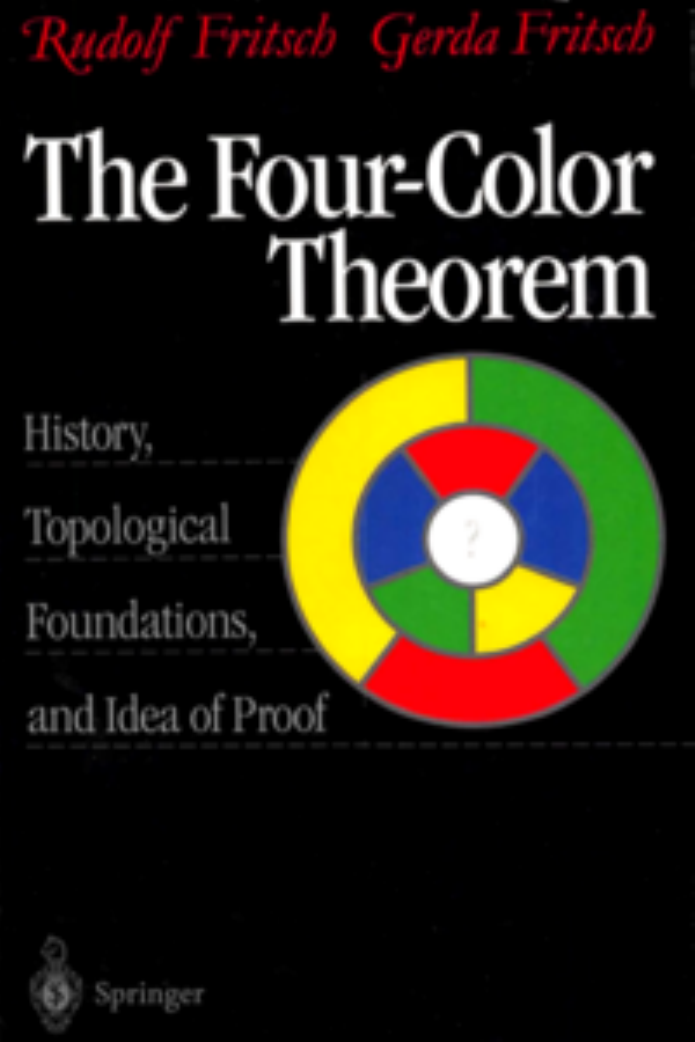}}
\scalebox{0.24}{\includegraphics{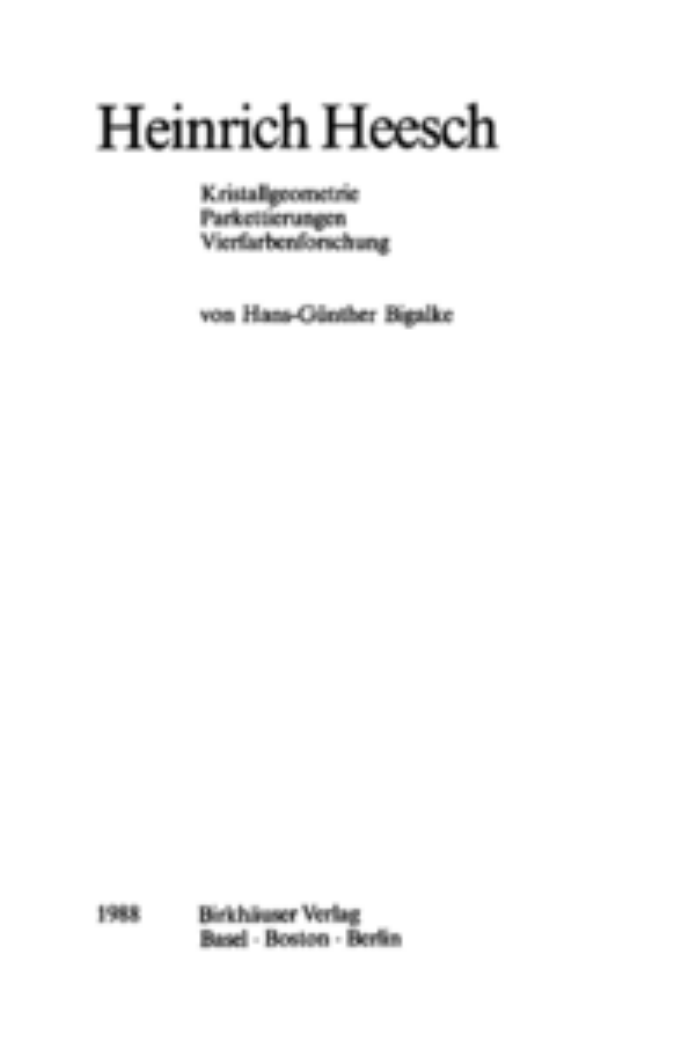}}
\scalebox{0.24}{\includegraphics{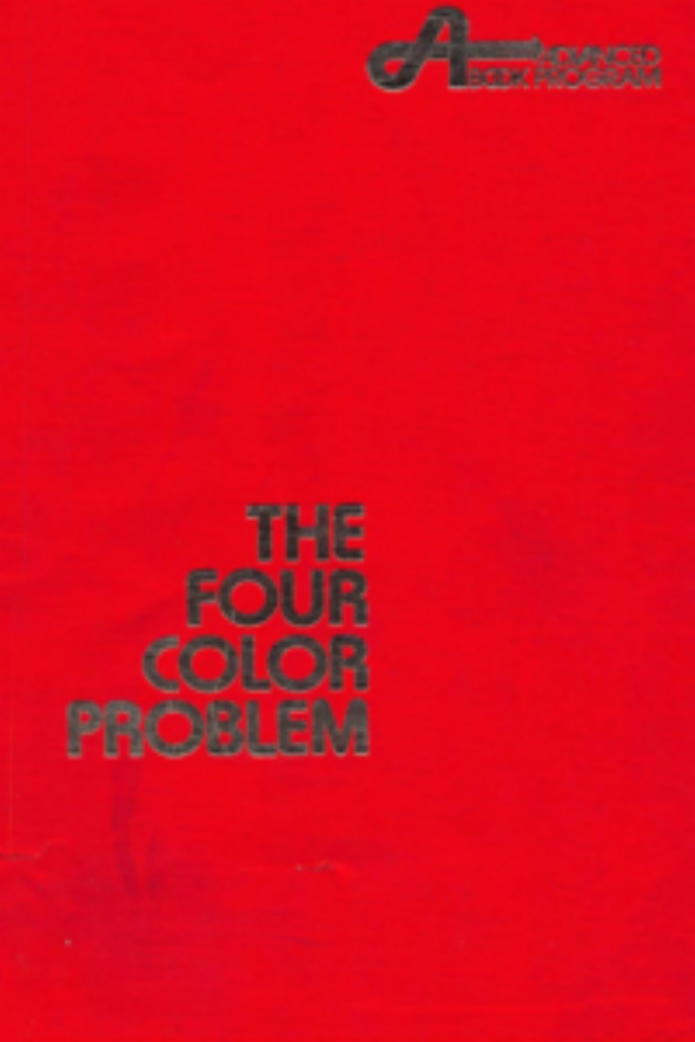}}
\scalebox{0.24}{\includegraphics{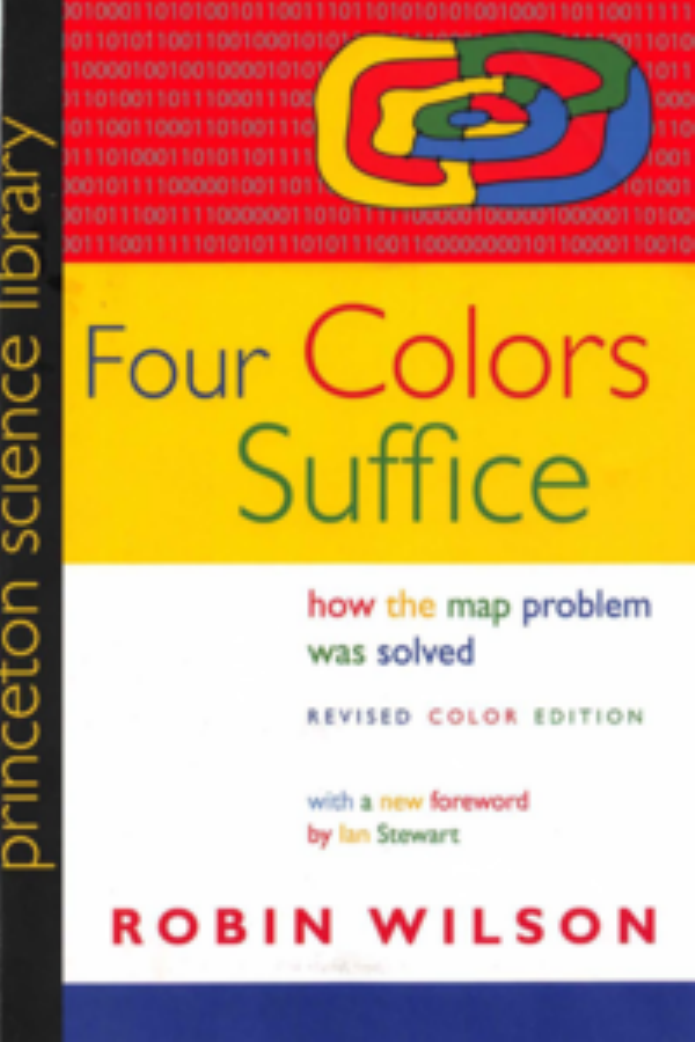}}
\scalebox{0.24}{\includegraphics{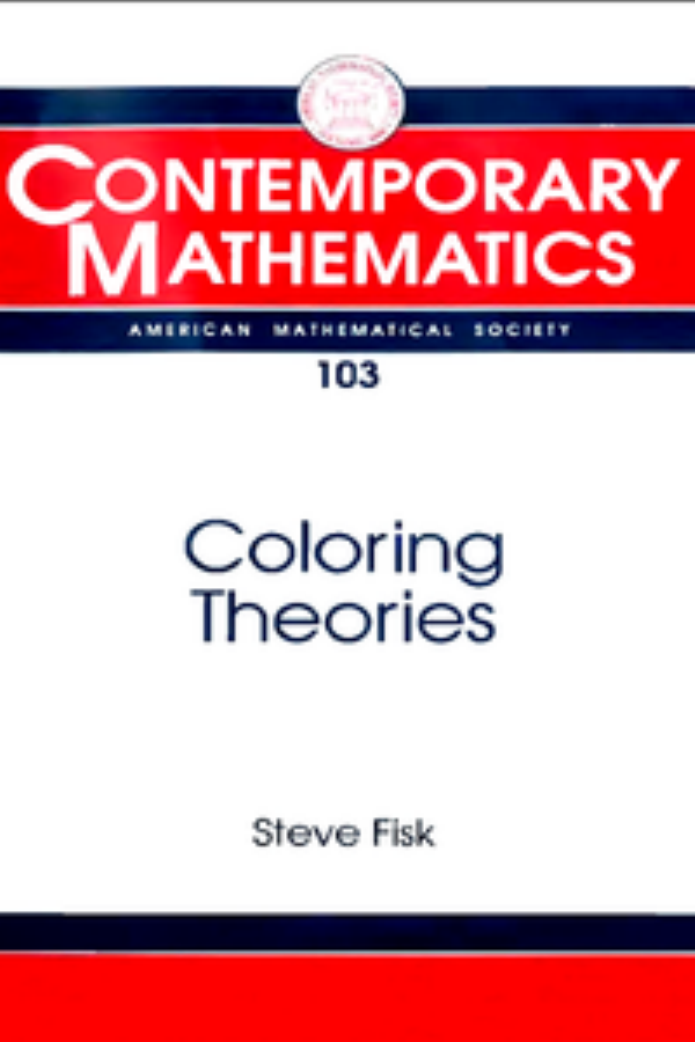}}
\caption{7 books devoted to graph coloring:
\cite{Barnette,ChartrandZhang2,FritschFritsch,Heesch,Ore,RobinWilson,Fisk1980}.
}
\end{figure}

\clearpage

\begin{figure}
\scalebox{0.13}{\includegraphics{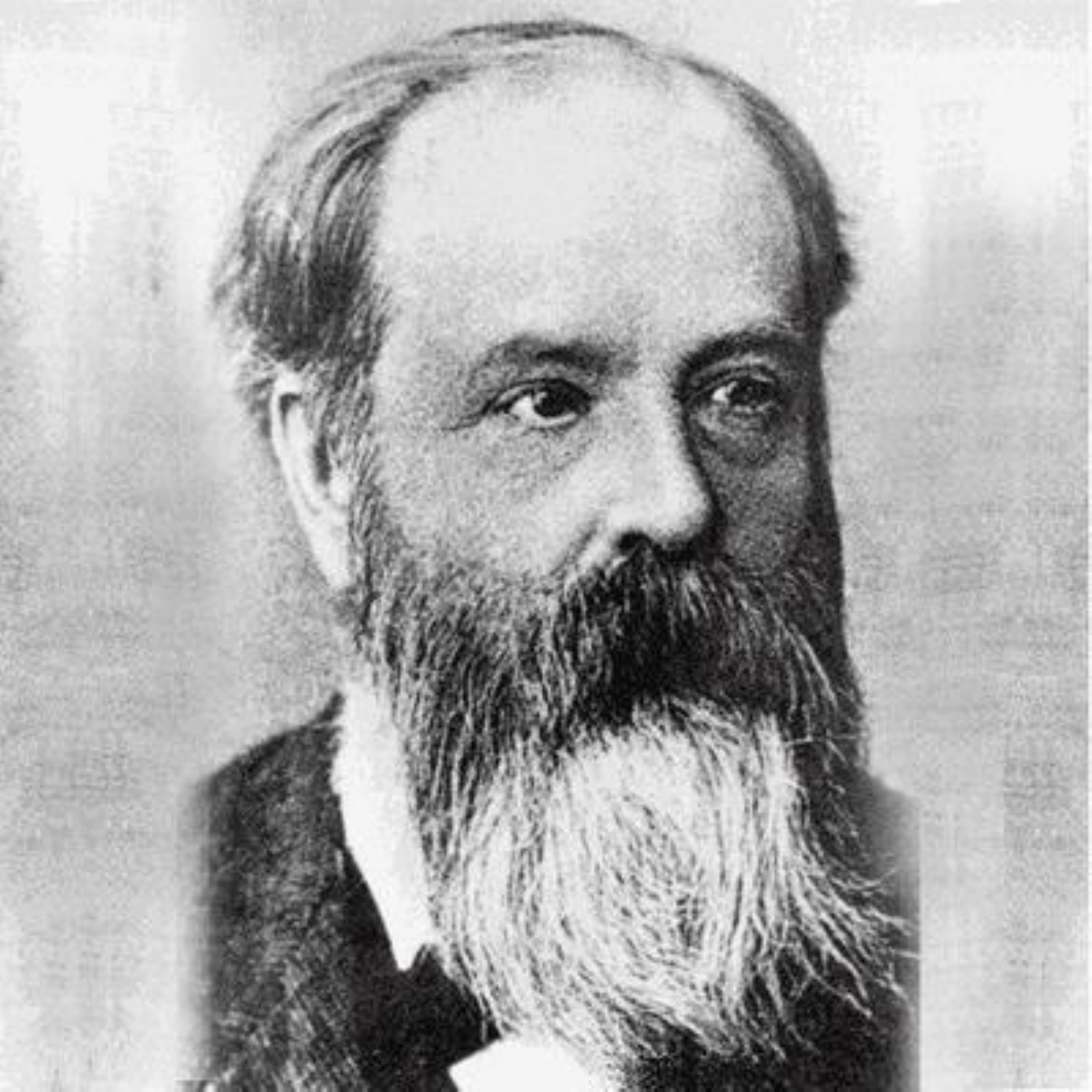}}
\scalebox{0.13}{\includegraphics{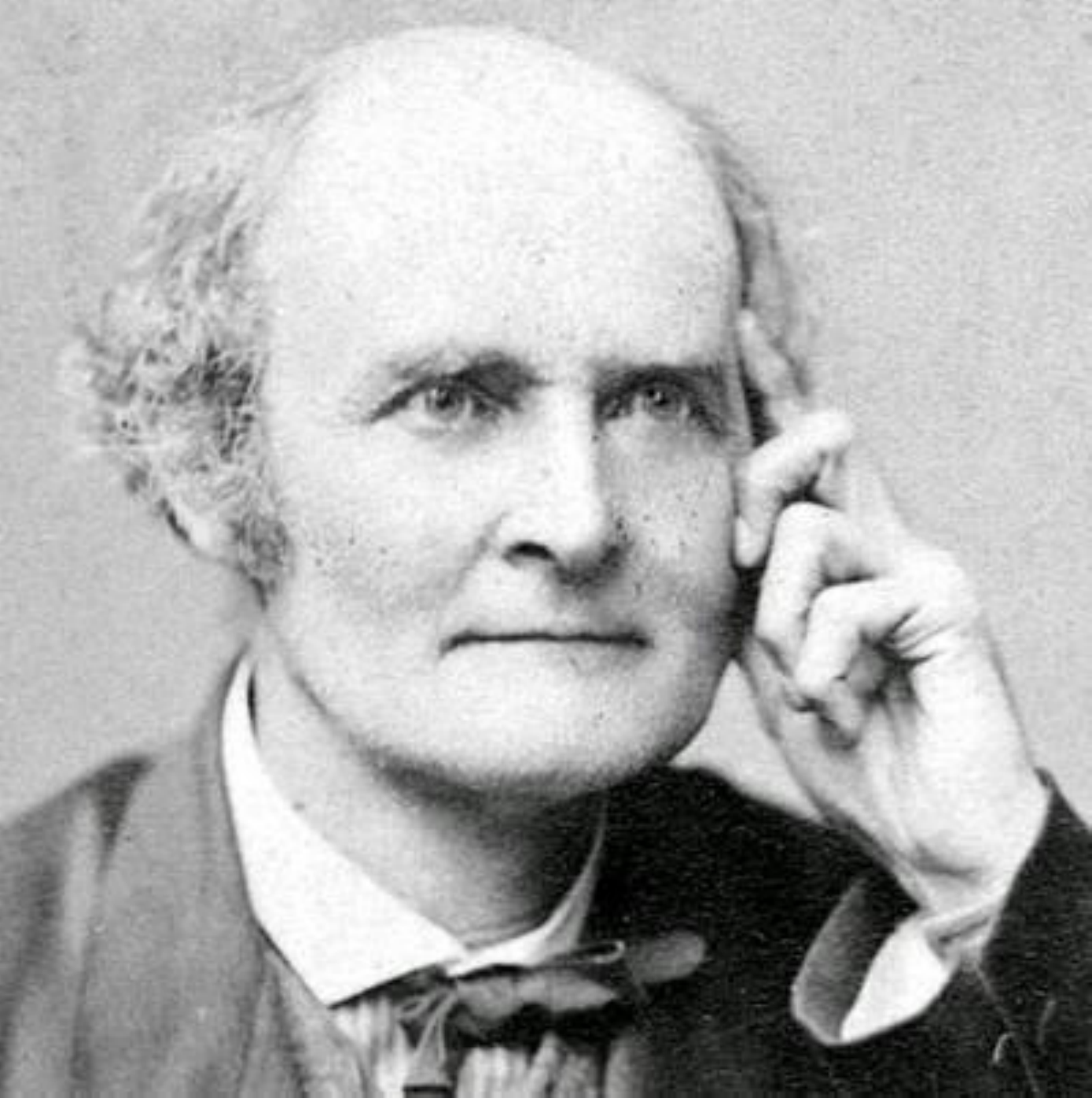}}
\scalebox{0.13}{\includegraphics{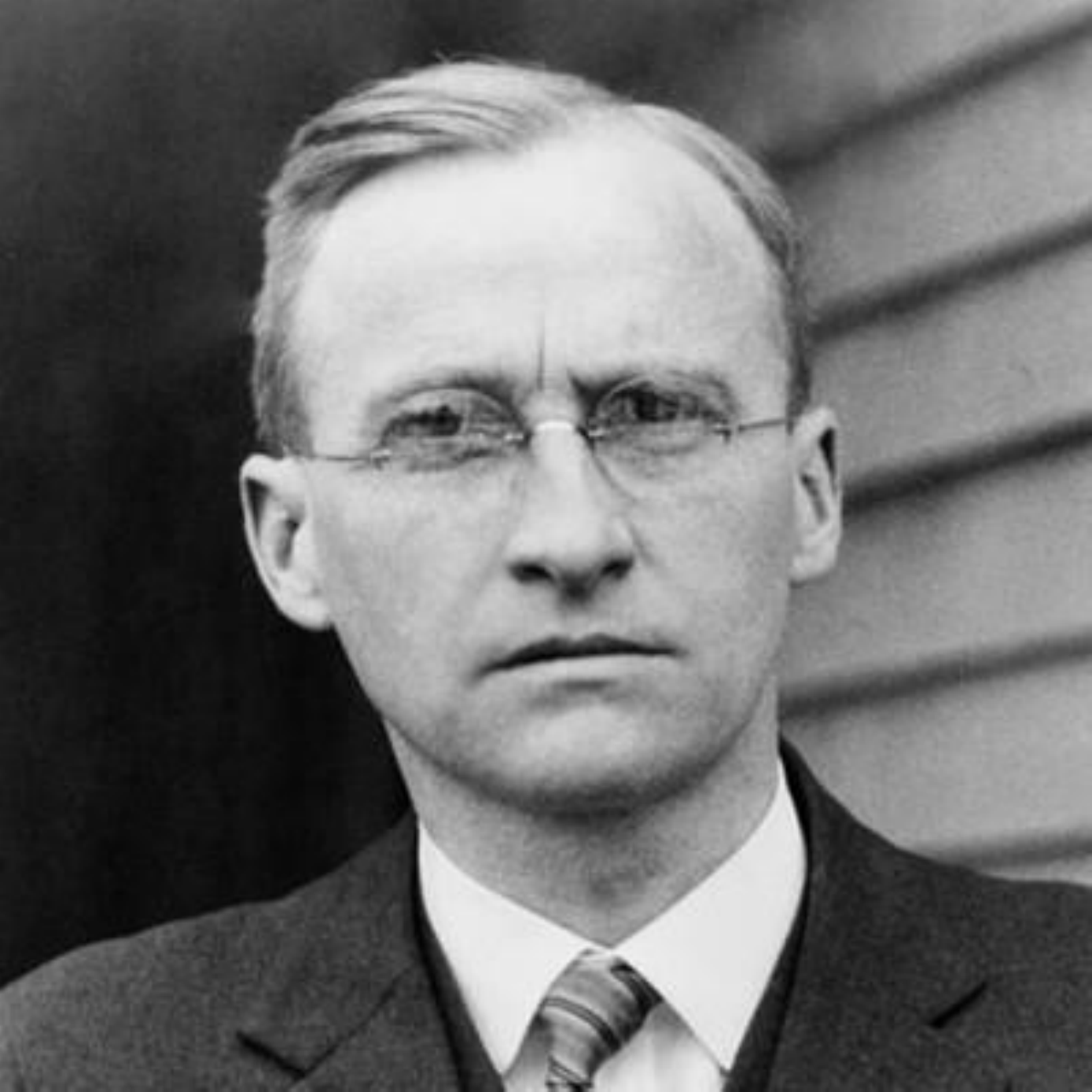}}
\scalebox{0.13}{\includegraphics{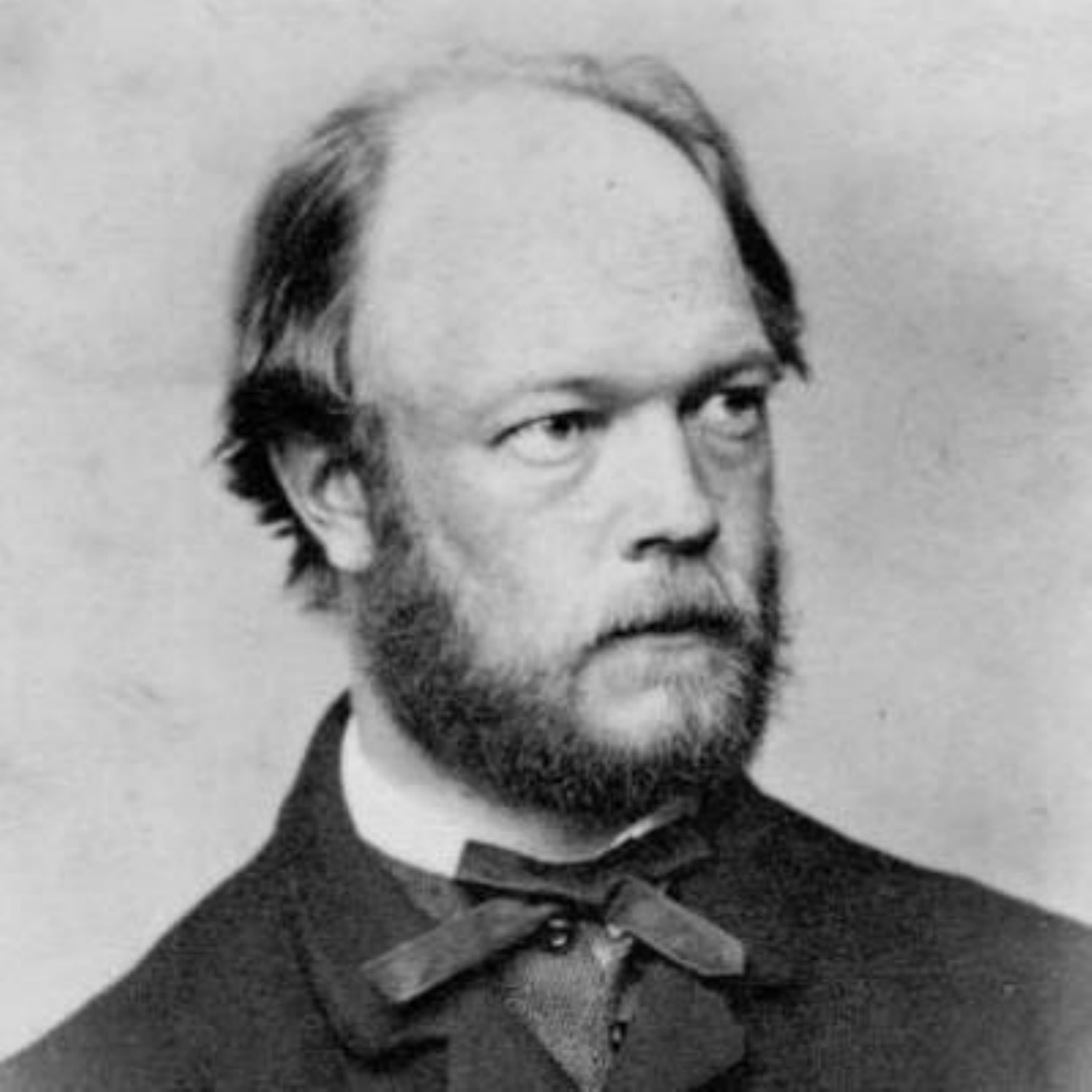}}
\scalebox{0.13}{\includegraphics{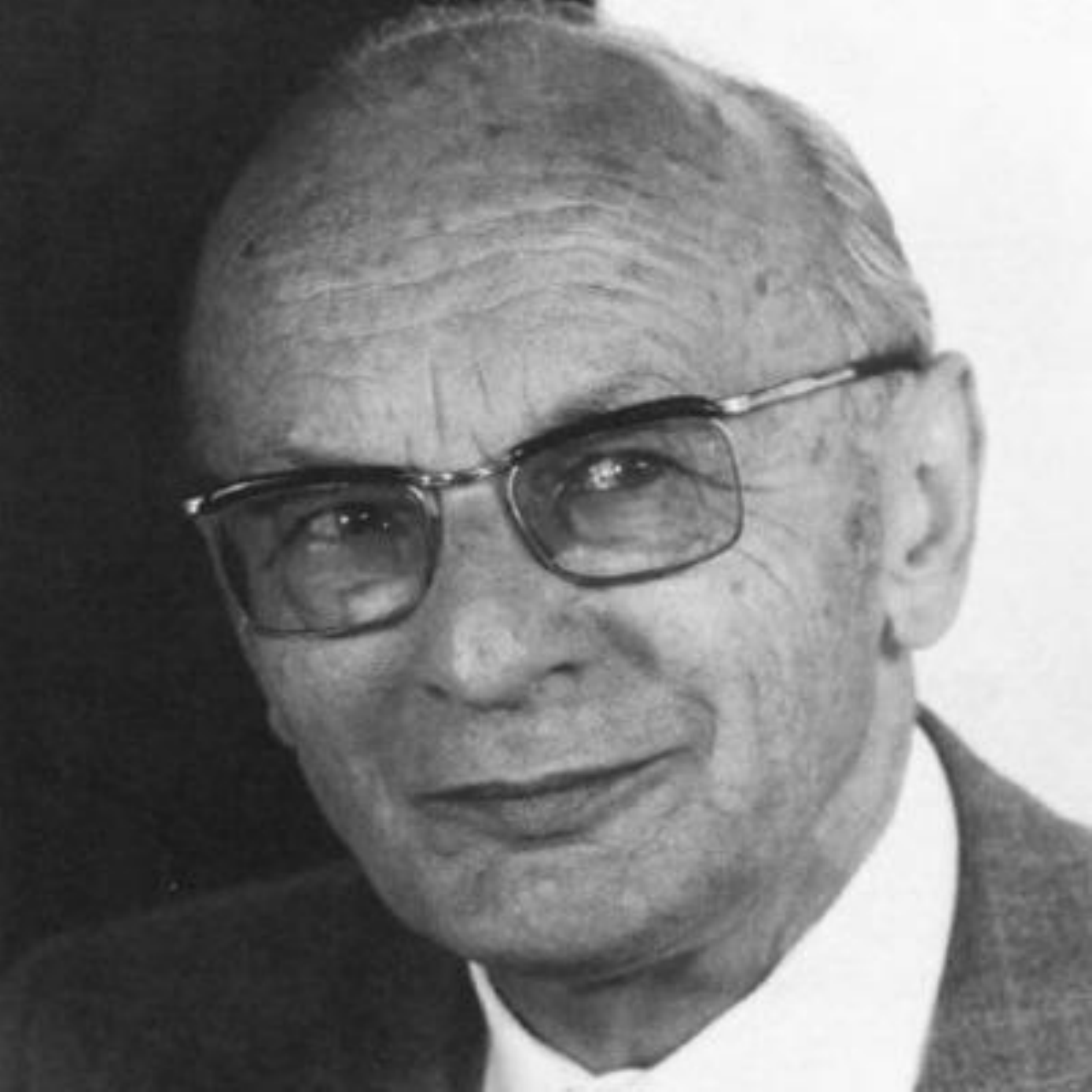}}
\scalebox{0.13}{\includegraphics{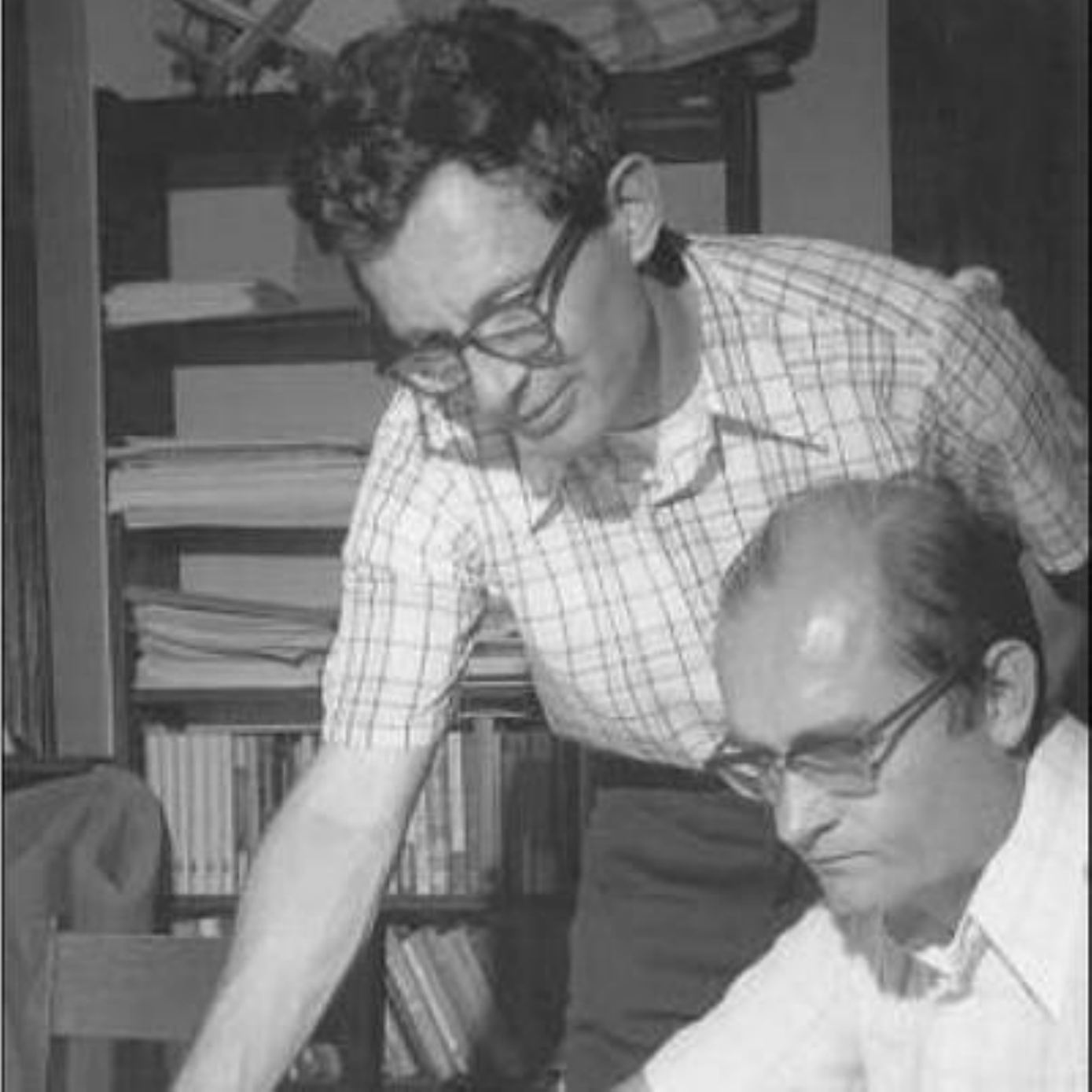}}
\caption{Guthrie, Cayley, Birkhoff, Tait, Heesch, Appel and Haken}
\end{figure}
\begin{figure}
\scalebox{0.13}{\includegraphics{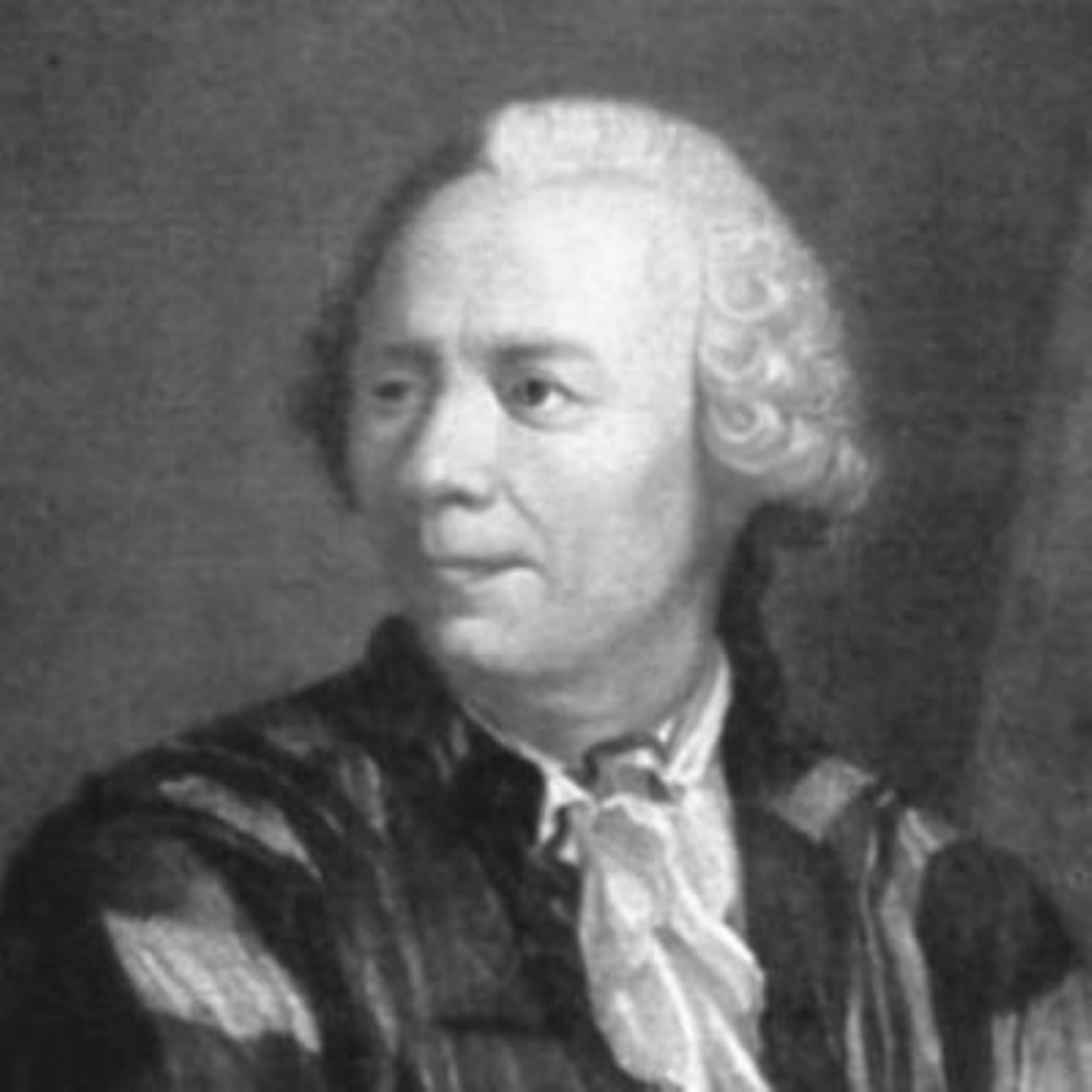}}
\scalebox{0.13}{\includegraphics{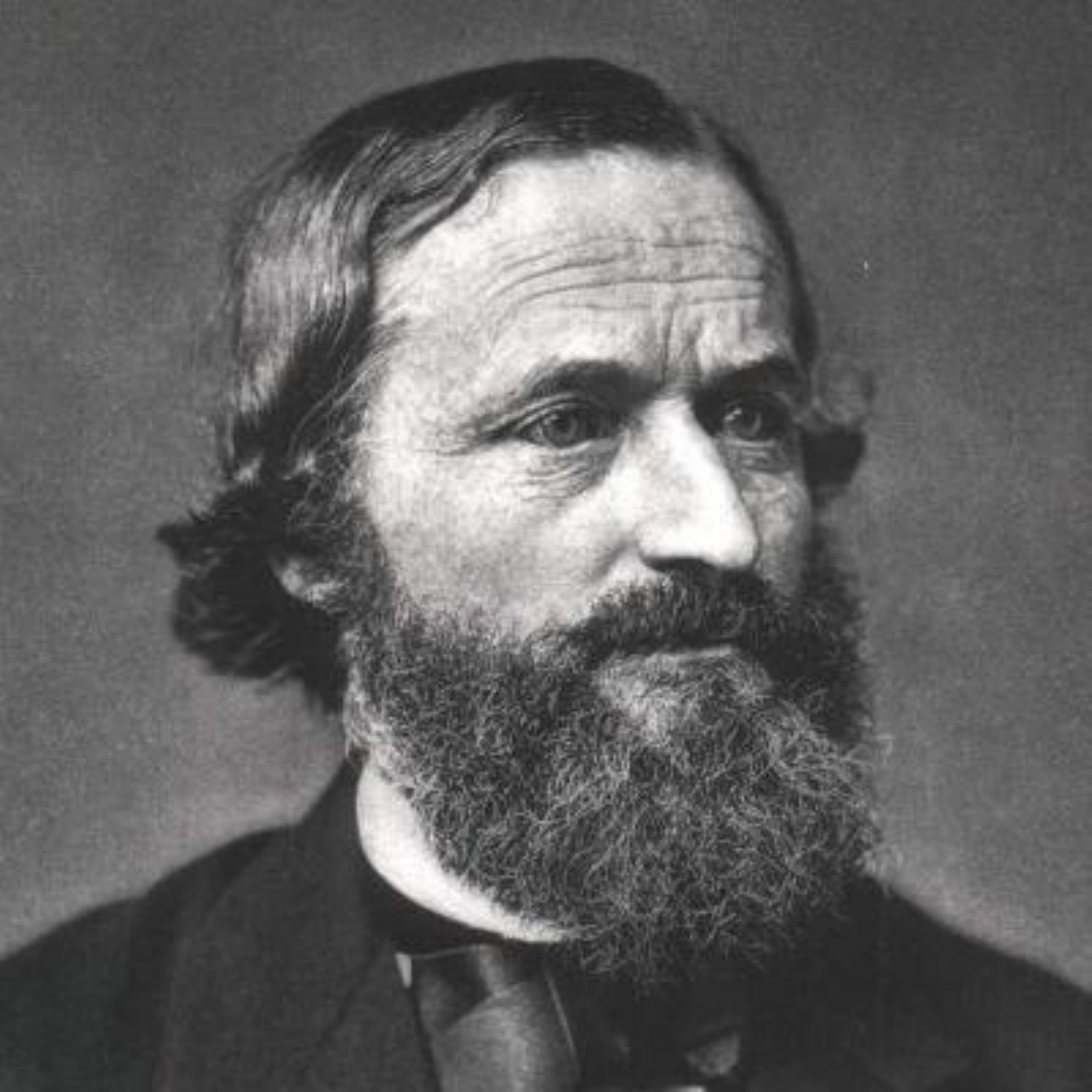}}
\scalebox{0.13}{\includegraphics{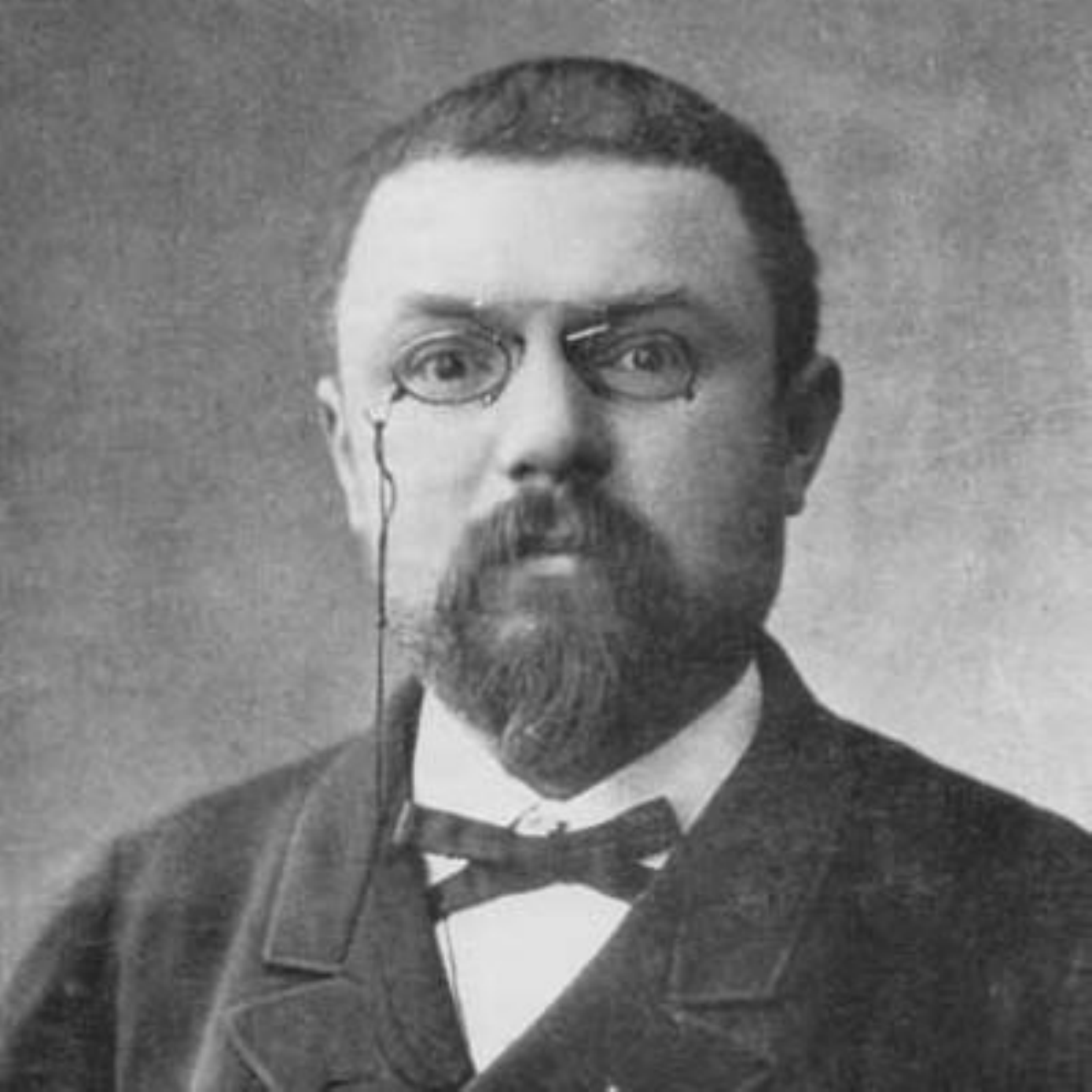}}
\scalebox{0.13}{\includegraphics{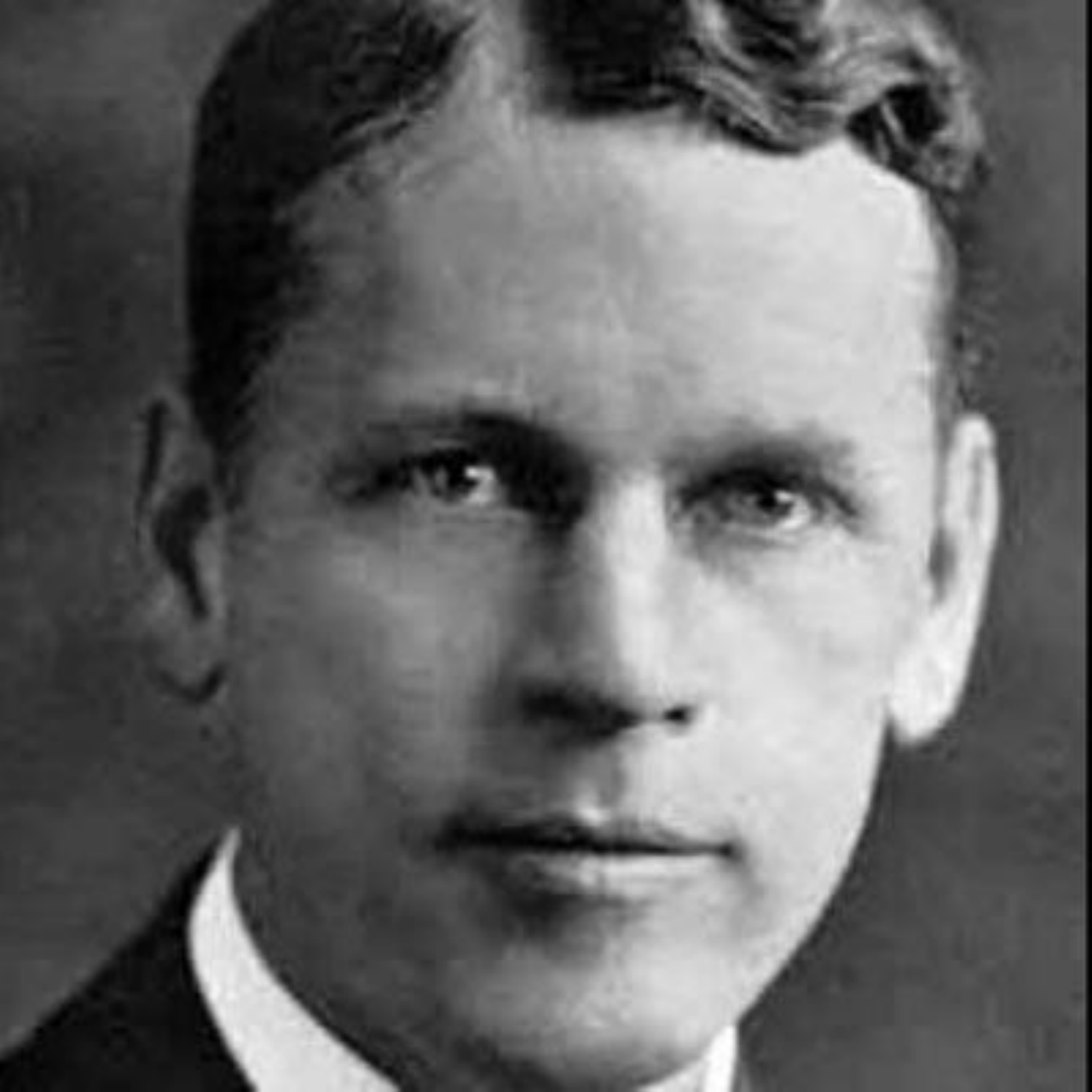}}
\scalebox{0.13}{\includegraphics{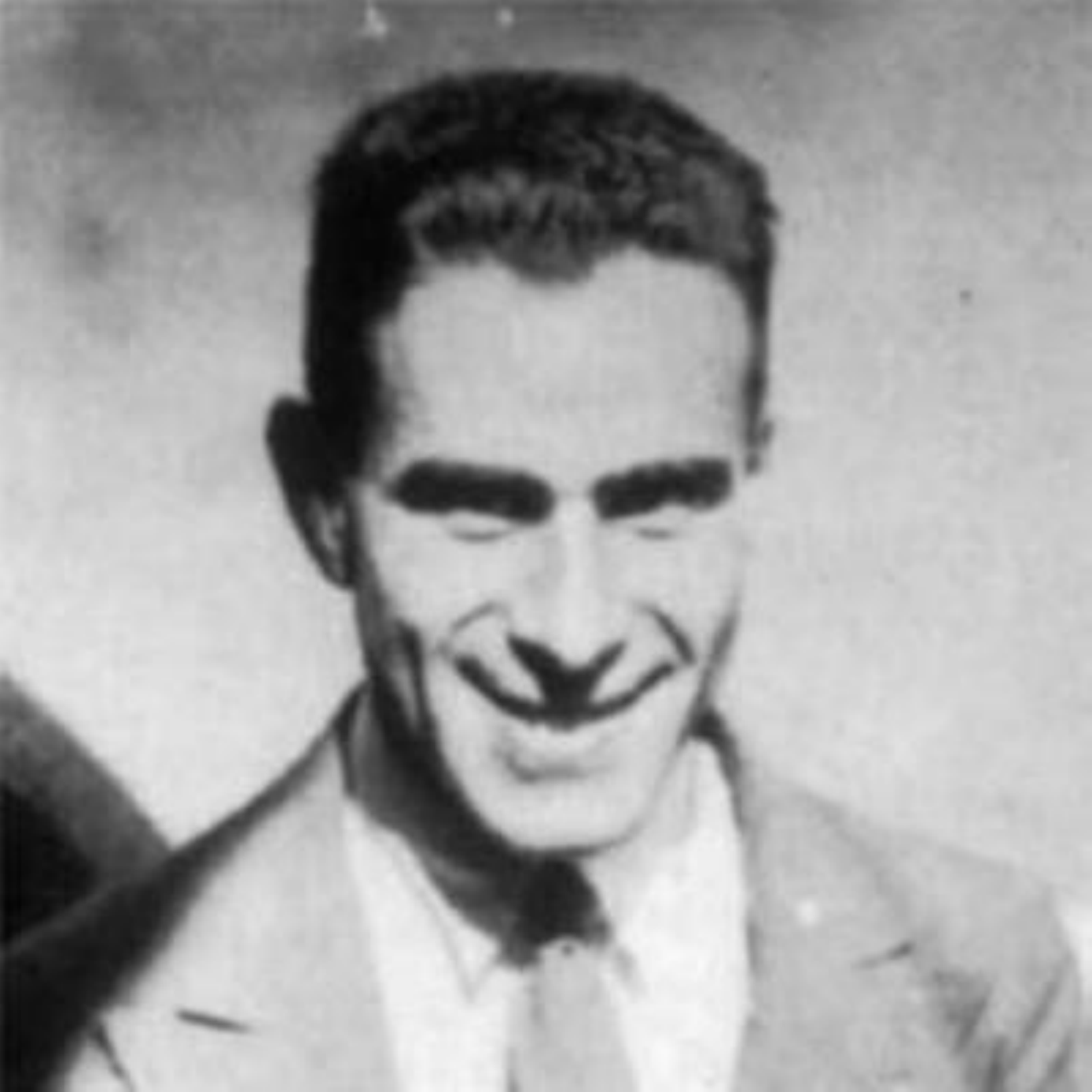}}
\scalebox{0.13}{\includegraphics{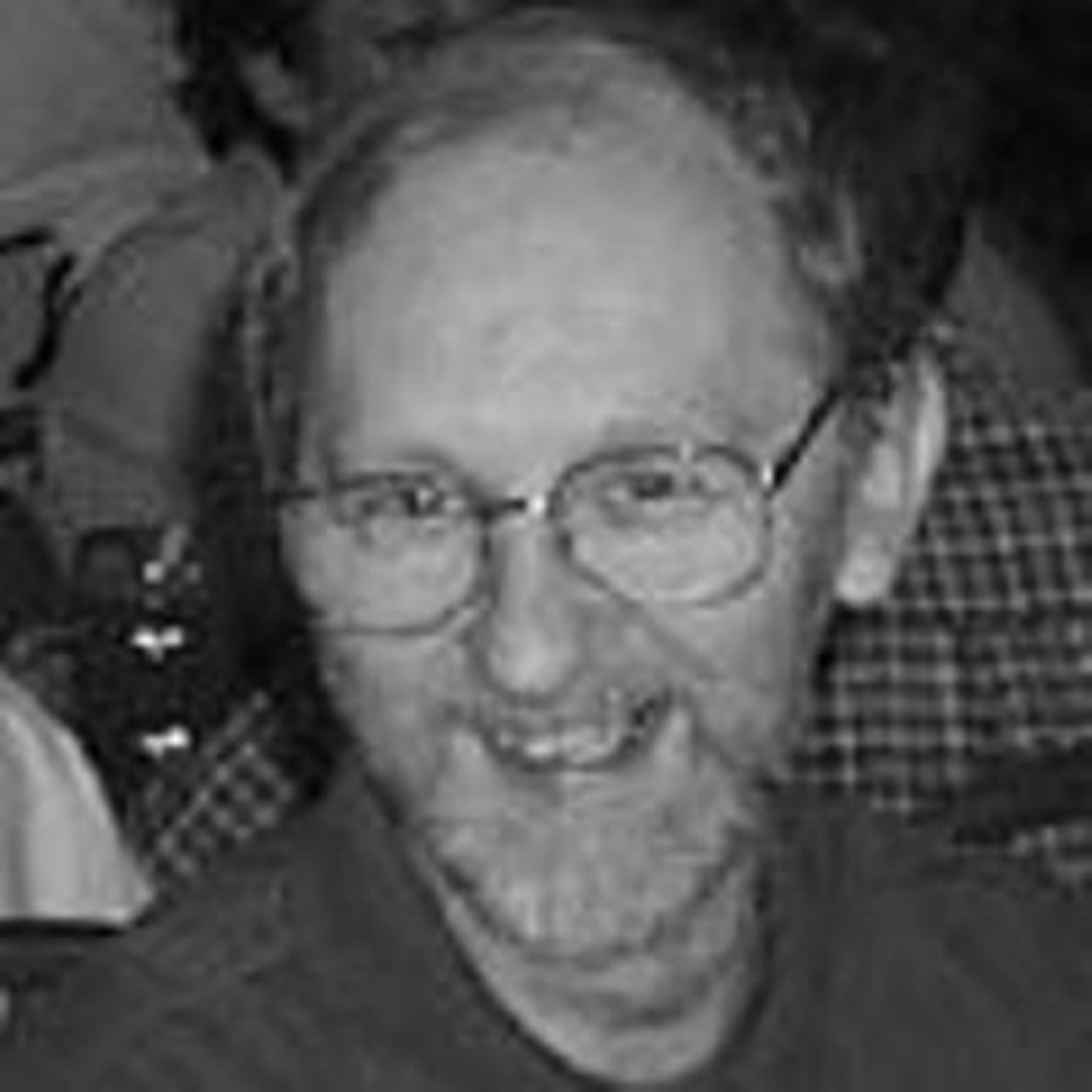}}
\caption{Euler, Kirchhoff, Poincar\'e, Veblen, Hurewicz, Barnette}
\end{figure}
\begin{figure}
\scalebox{0.13}{\includegraphics{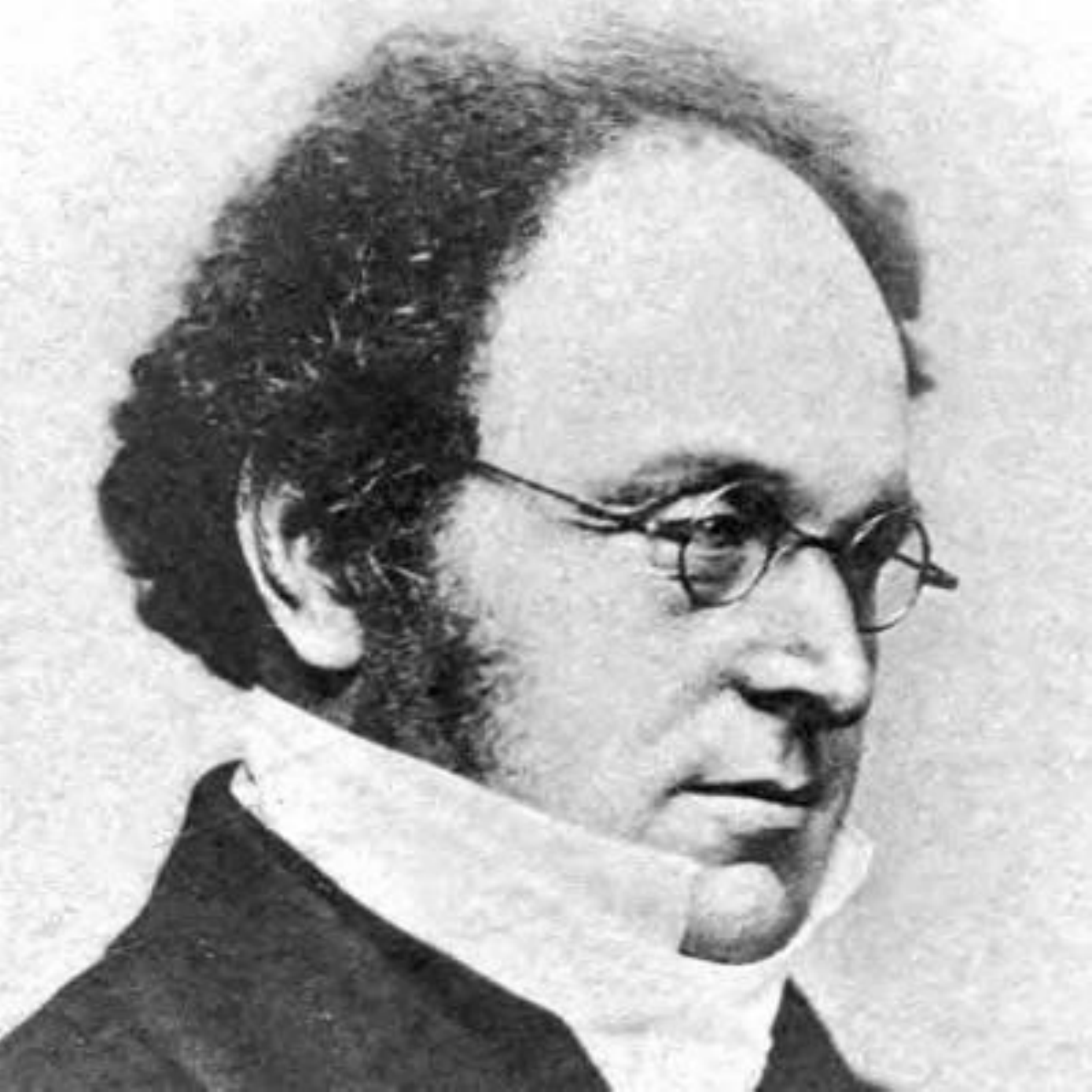}}
\scalebox{0.13}{\includegraphics{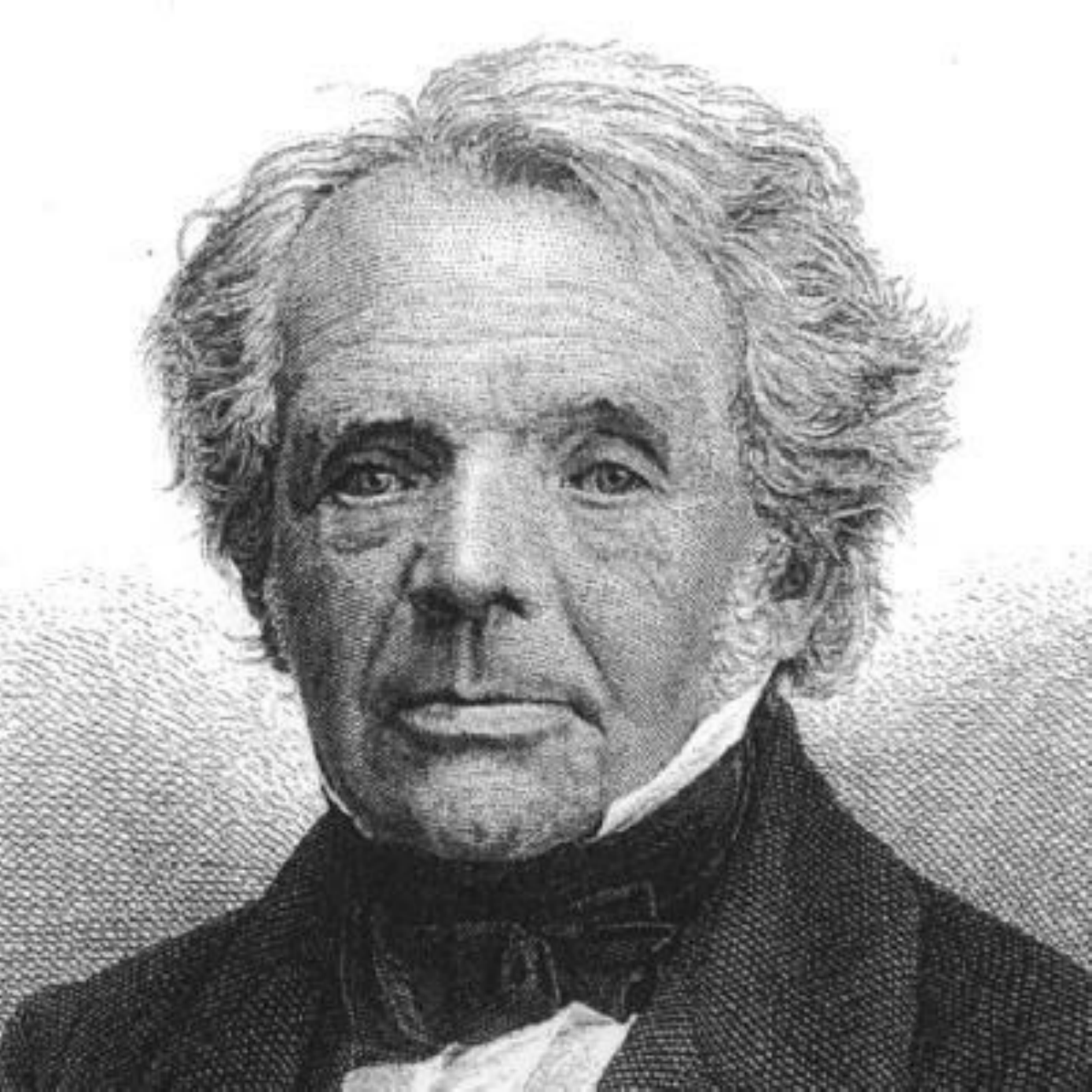}}
\scalebox{0.13}{\includegraphics{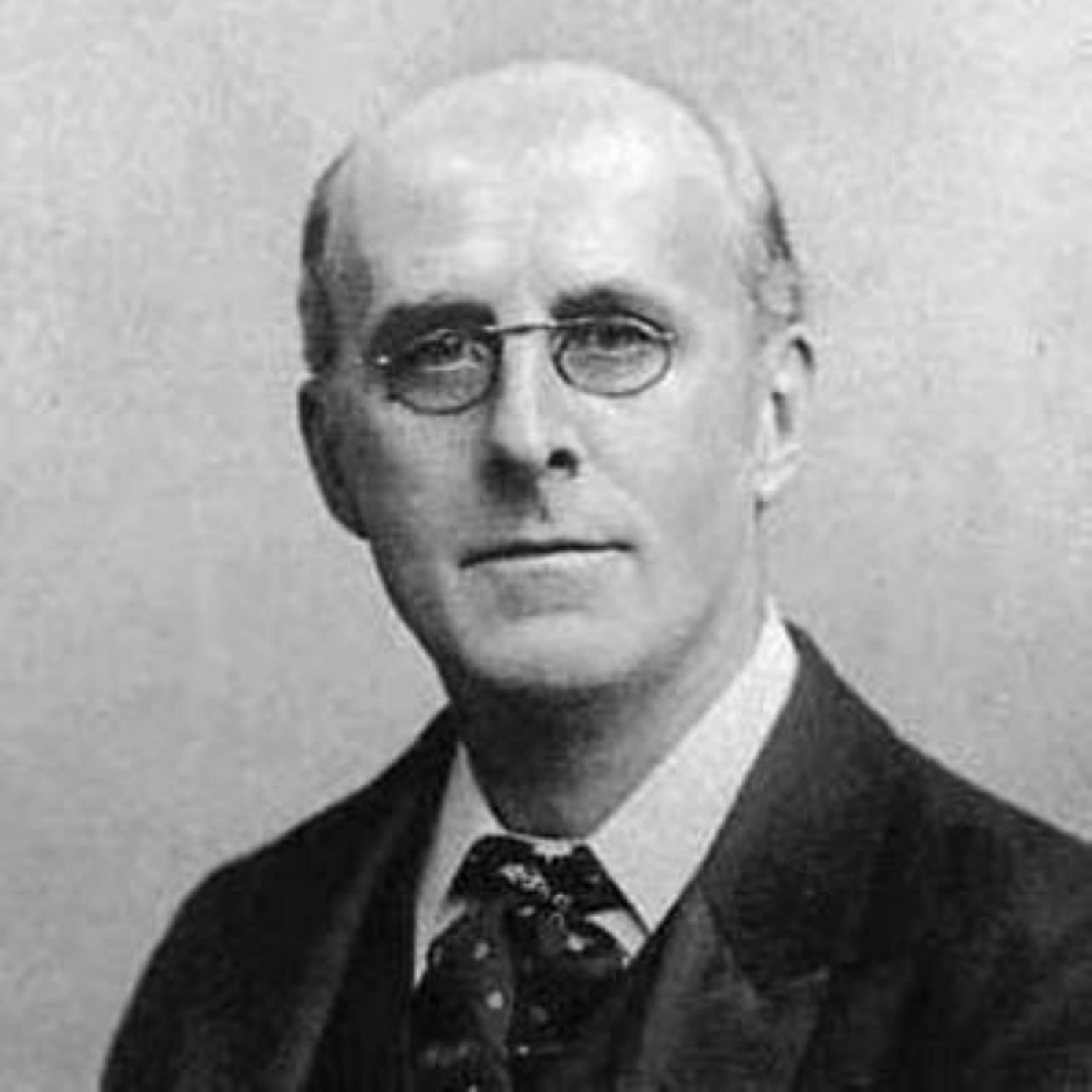}}
\scalebox{0.13}{\includegraphics{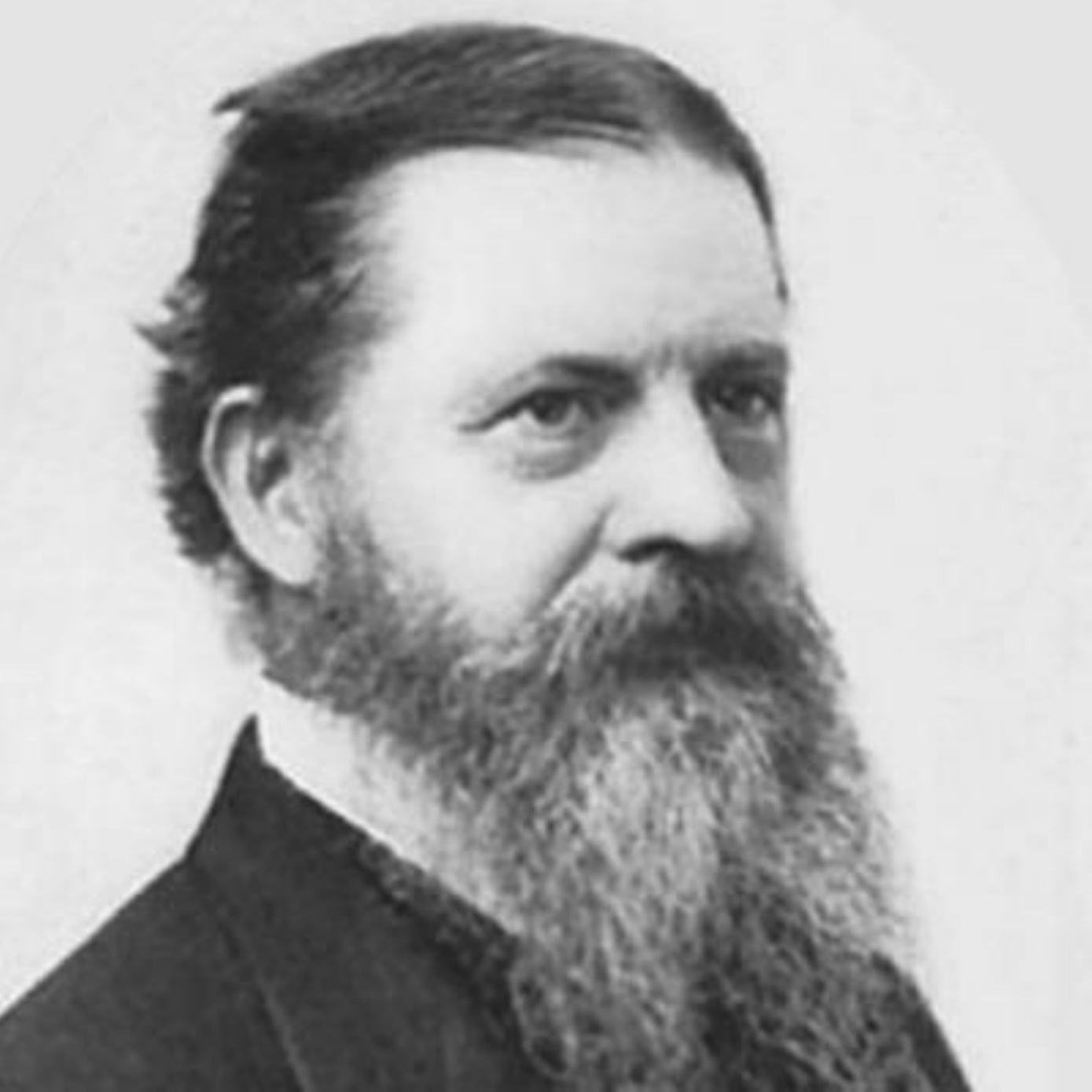}}
\scalebox{0.13}{\includegraphics{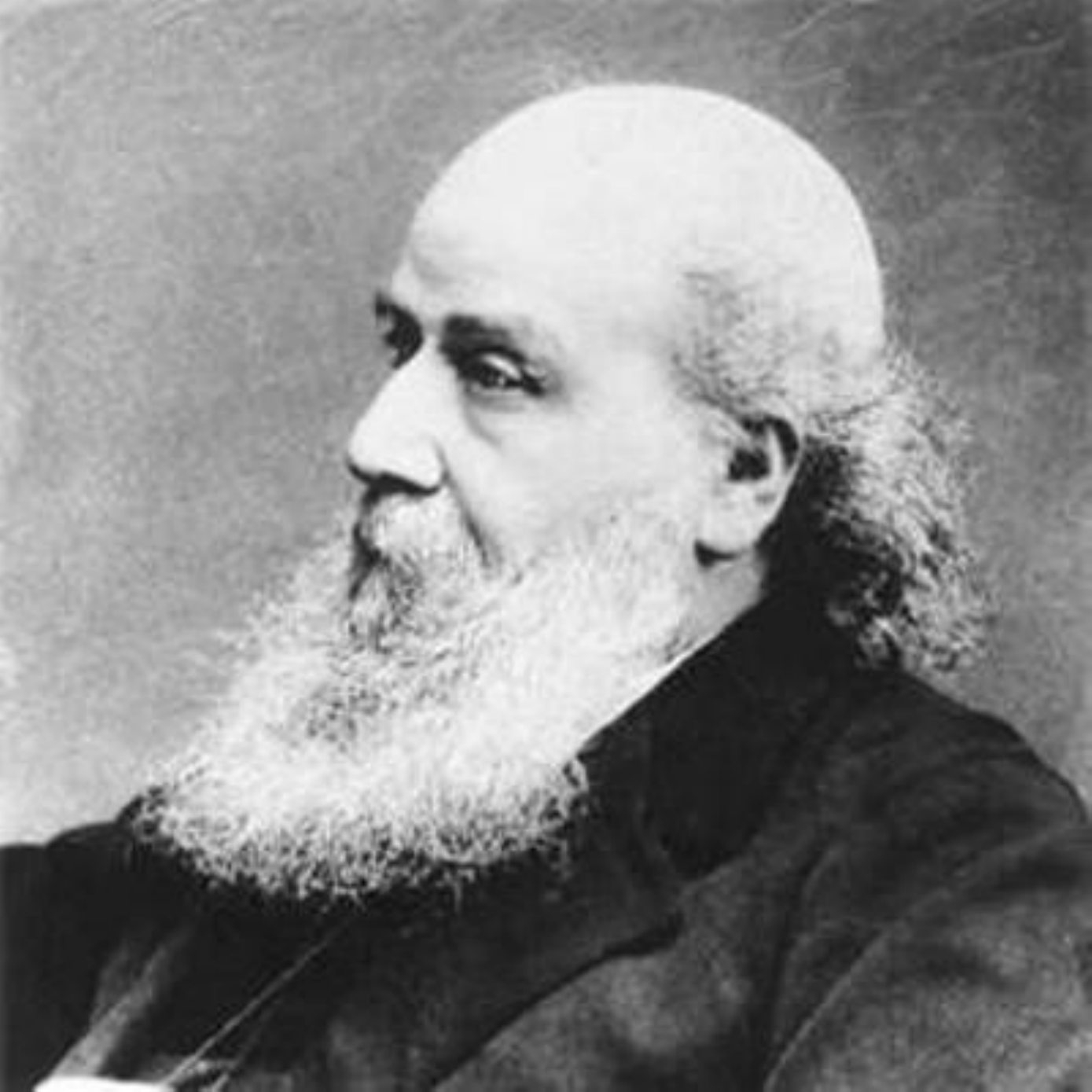}}
\scalebox{0.13}{\includegraphics{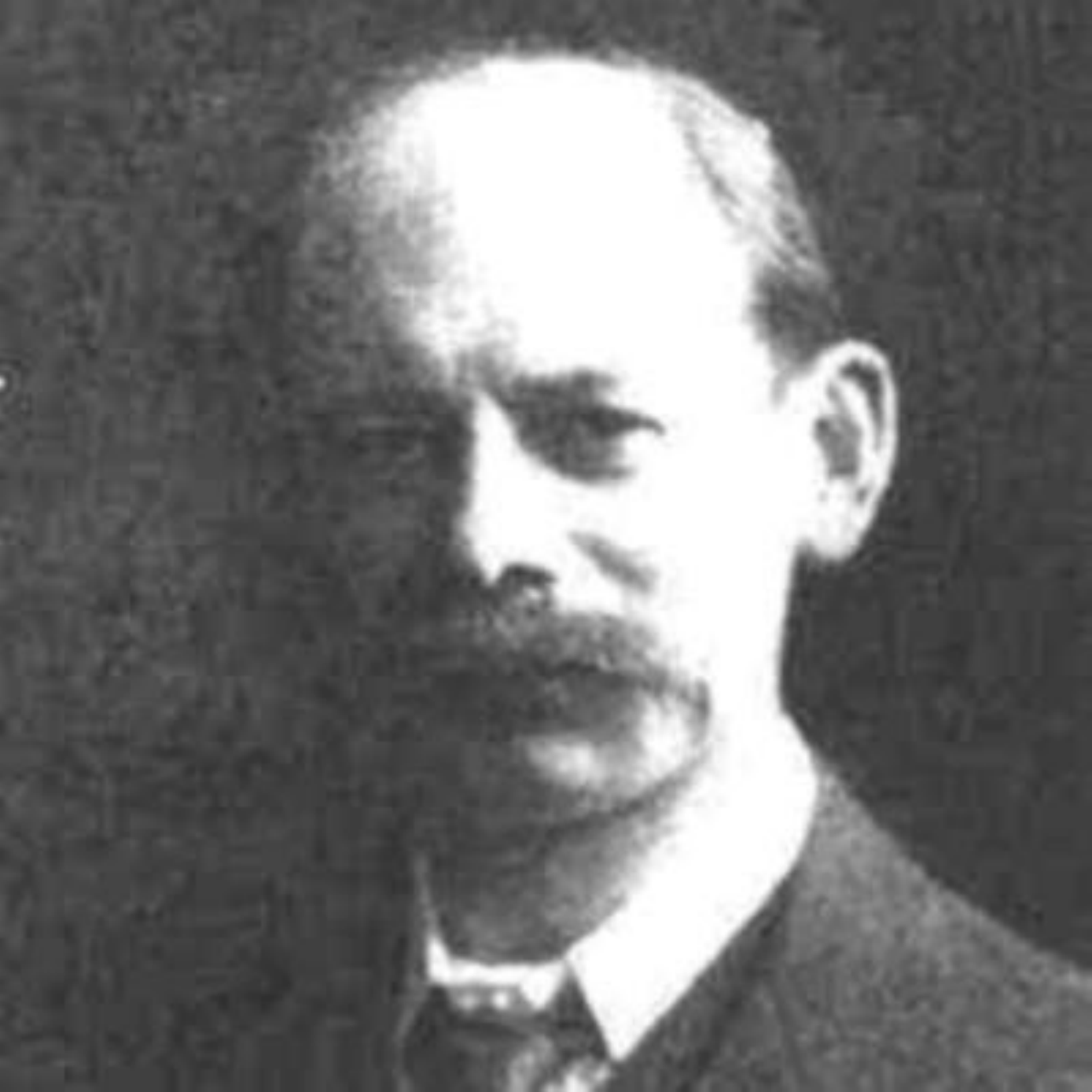}}
\caption{
Morgan, Moebius, Kempe, Pierce, Sylvester and Heawood.}
\end{figure}
\begin{figure}
\scalebox{0.13}{\includegraphics{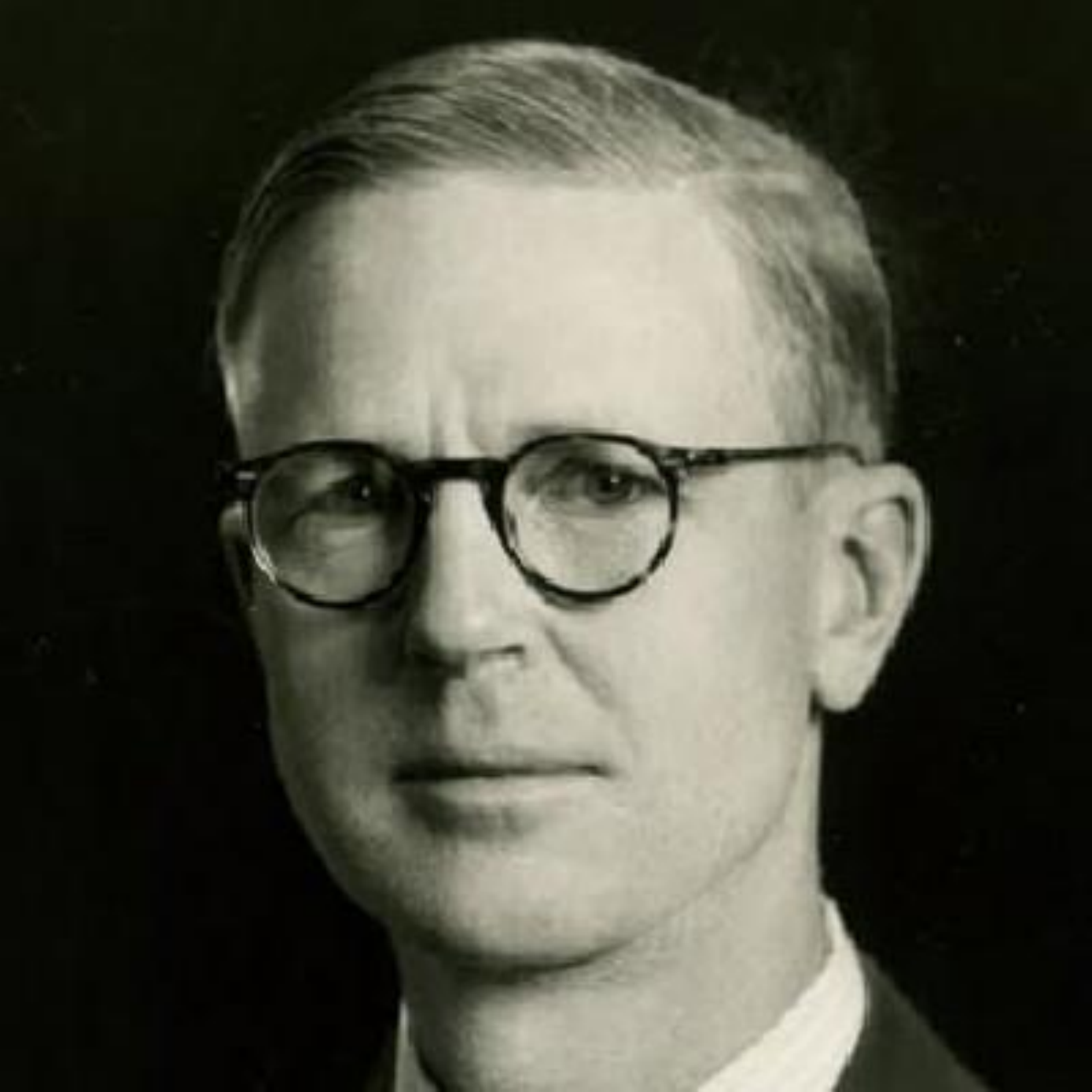}}
\scalebox{0.13}{\includegraphics{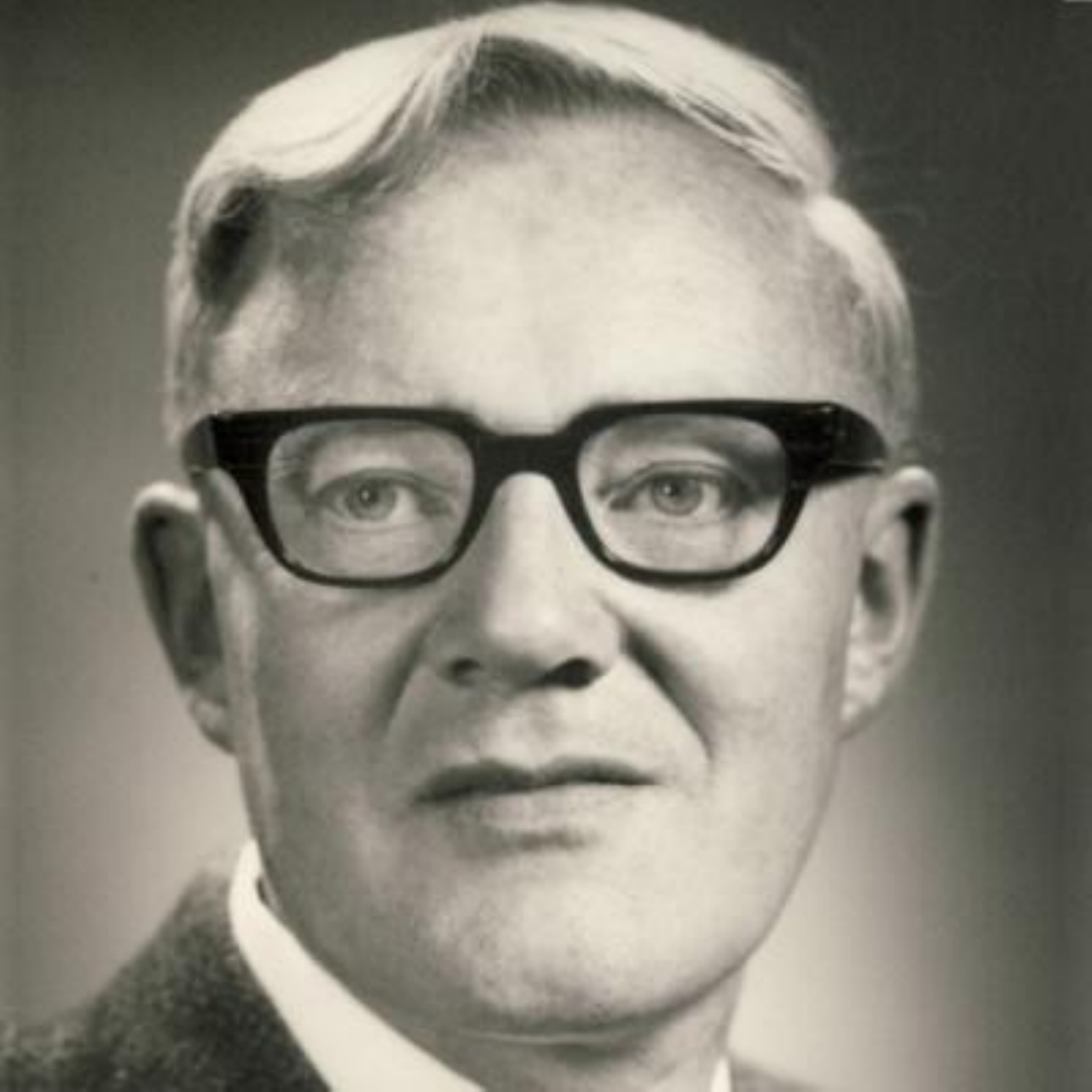}}
\scalebox{0.13}{\includegraphics{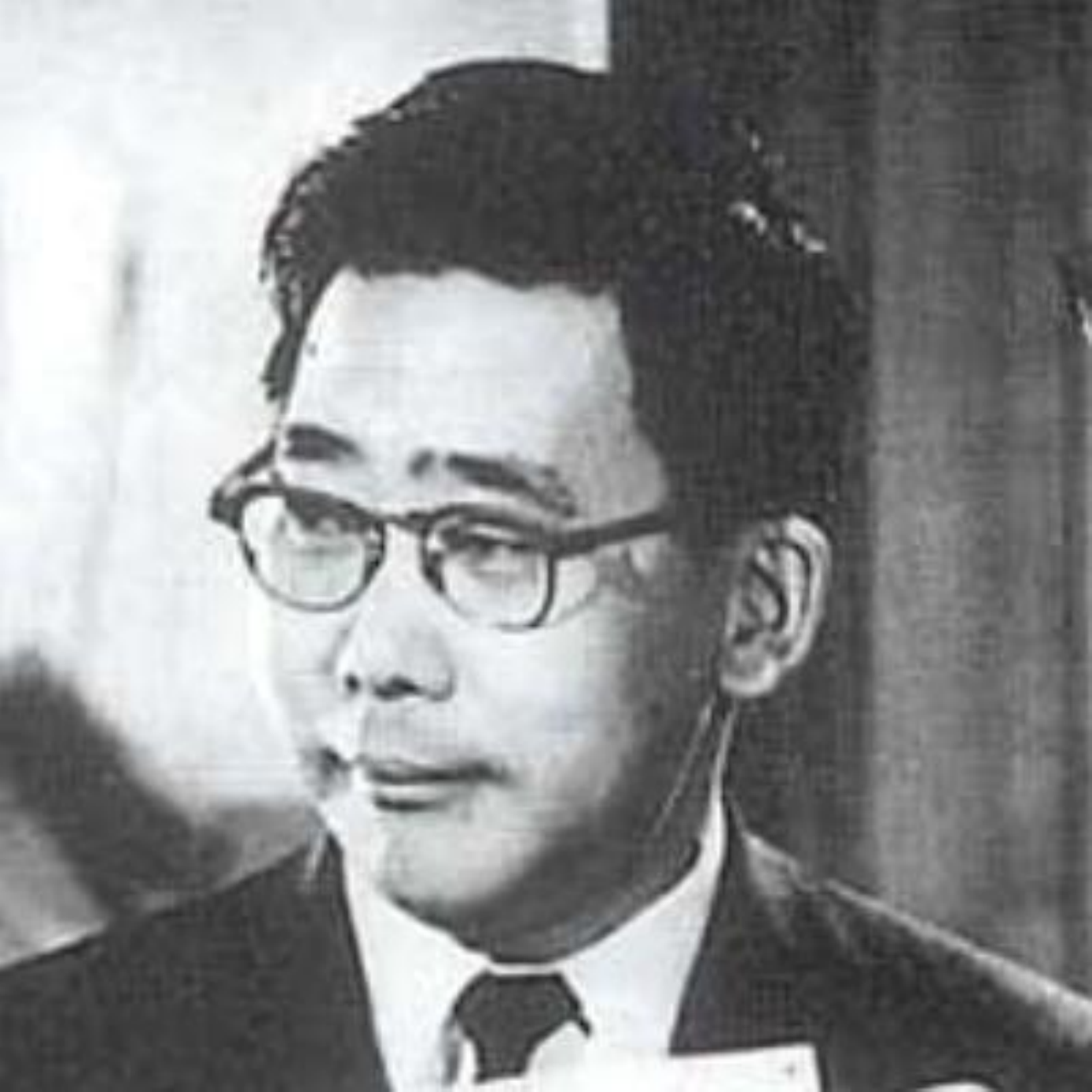}}
\scalebox{0.13}{\includegraphics{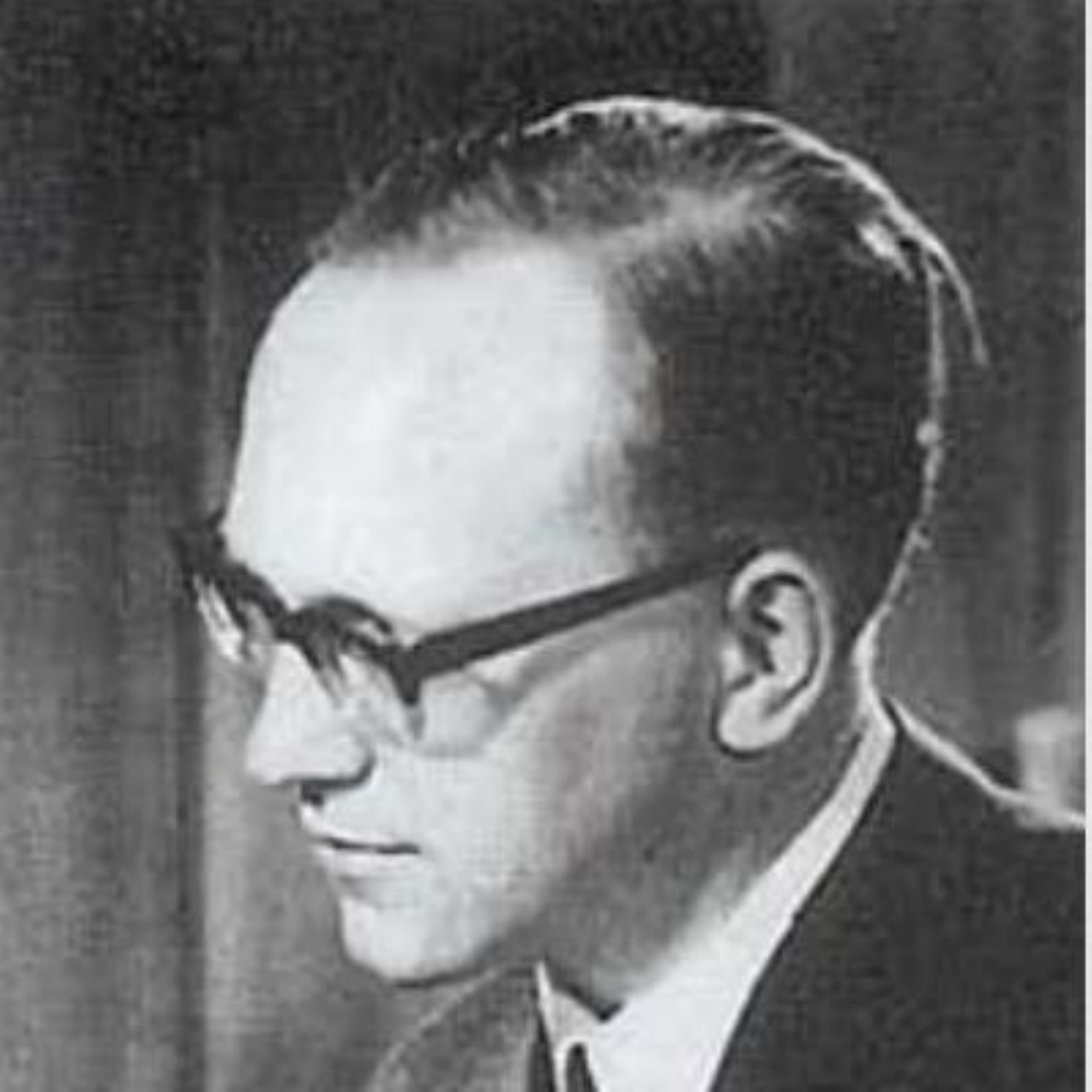}}
\scalebox{0.13}{\includegraphics{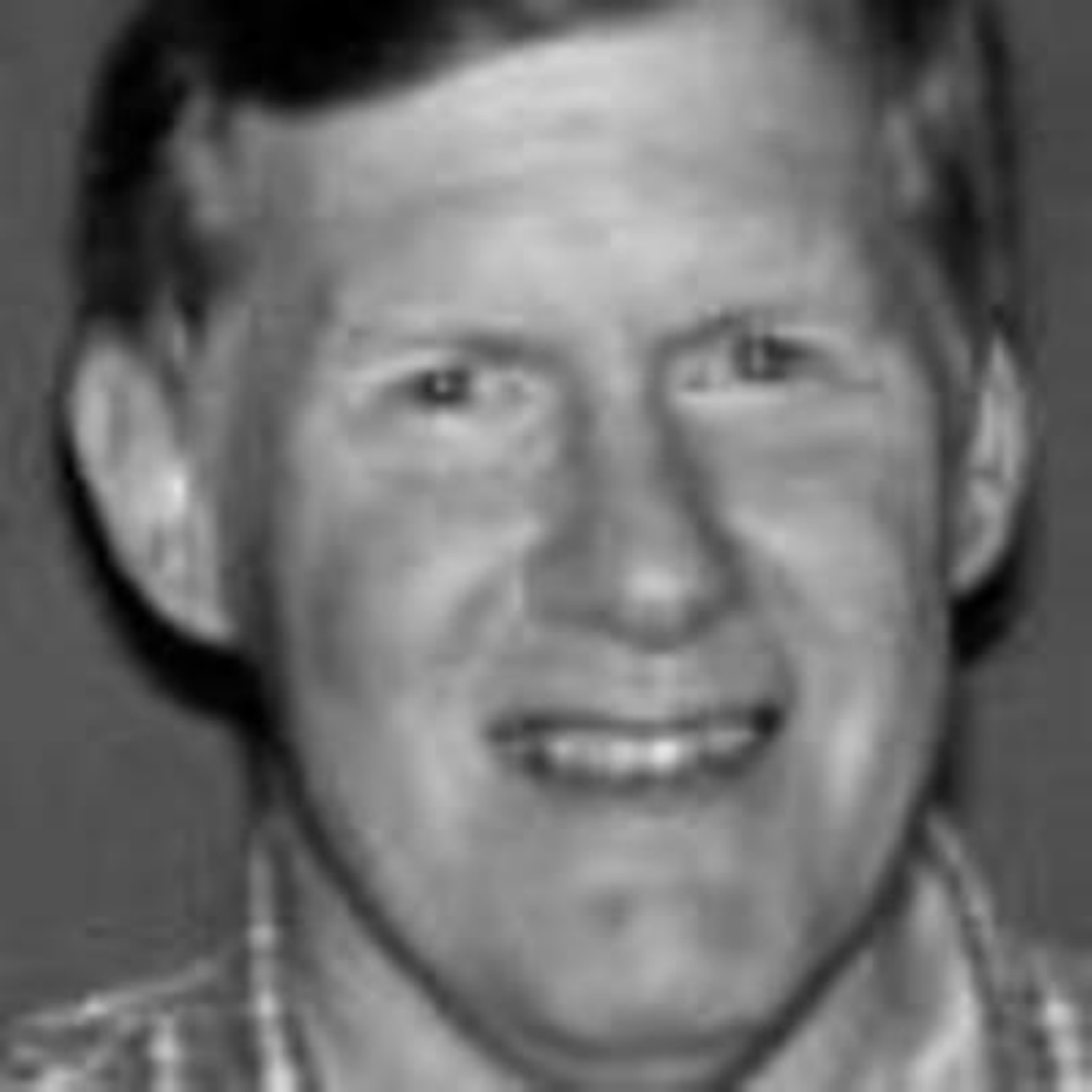}}
\scalebox{0.13}{\includegraphics{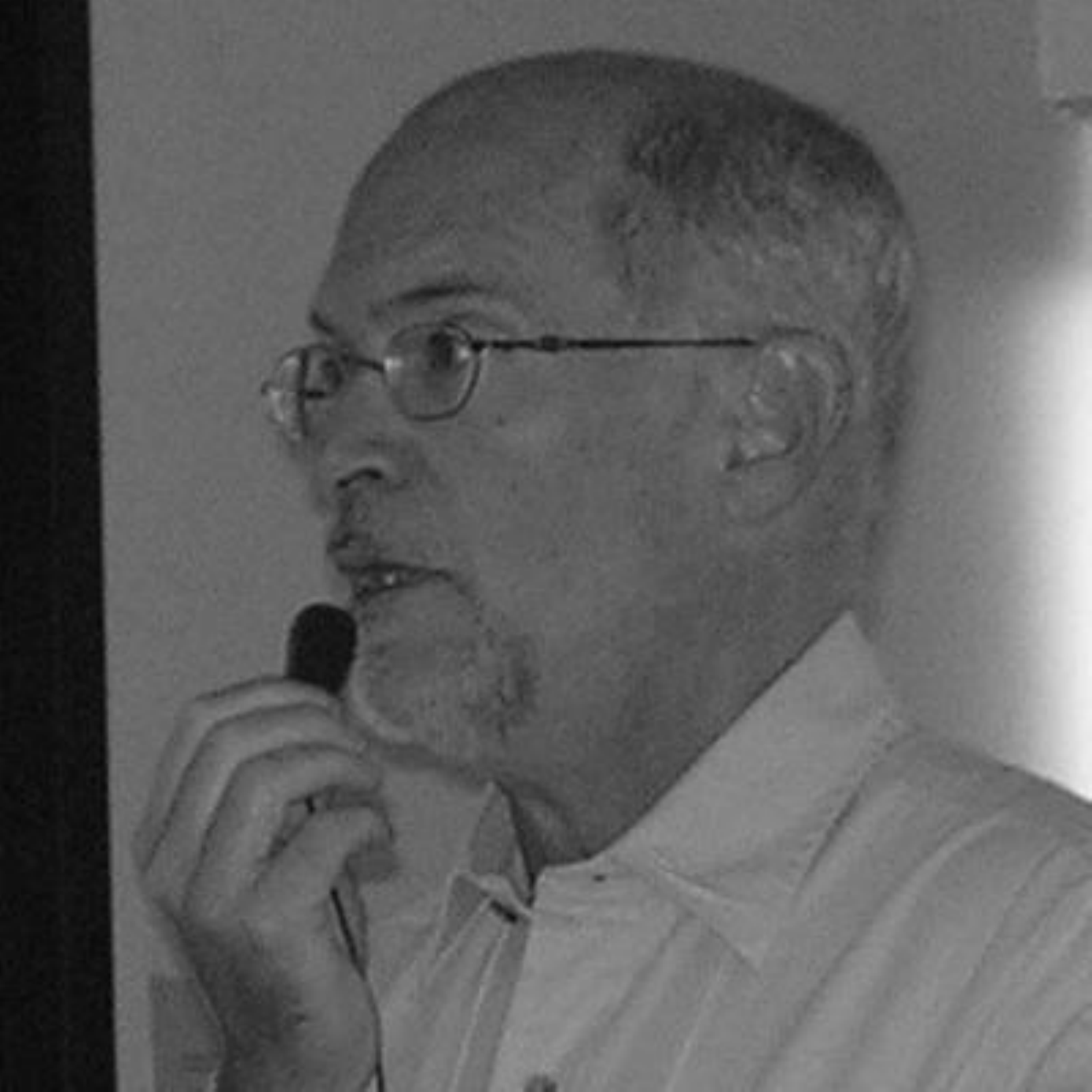}}
\caption{
Whitney, Tutte, Shimamoto, Duerre, Stromquist, Albertson}
\end{figure}
\begin{figure}
\scalebox{0.13}{\includegraphics{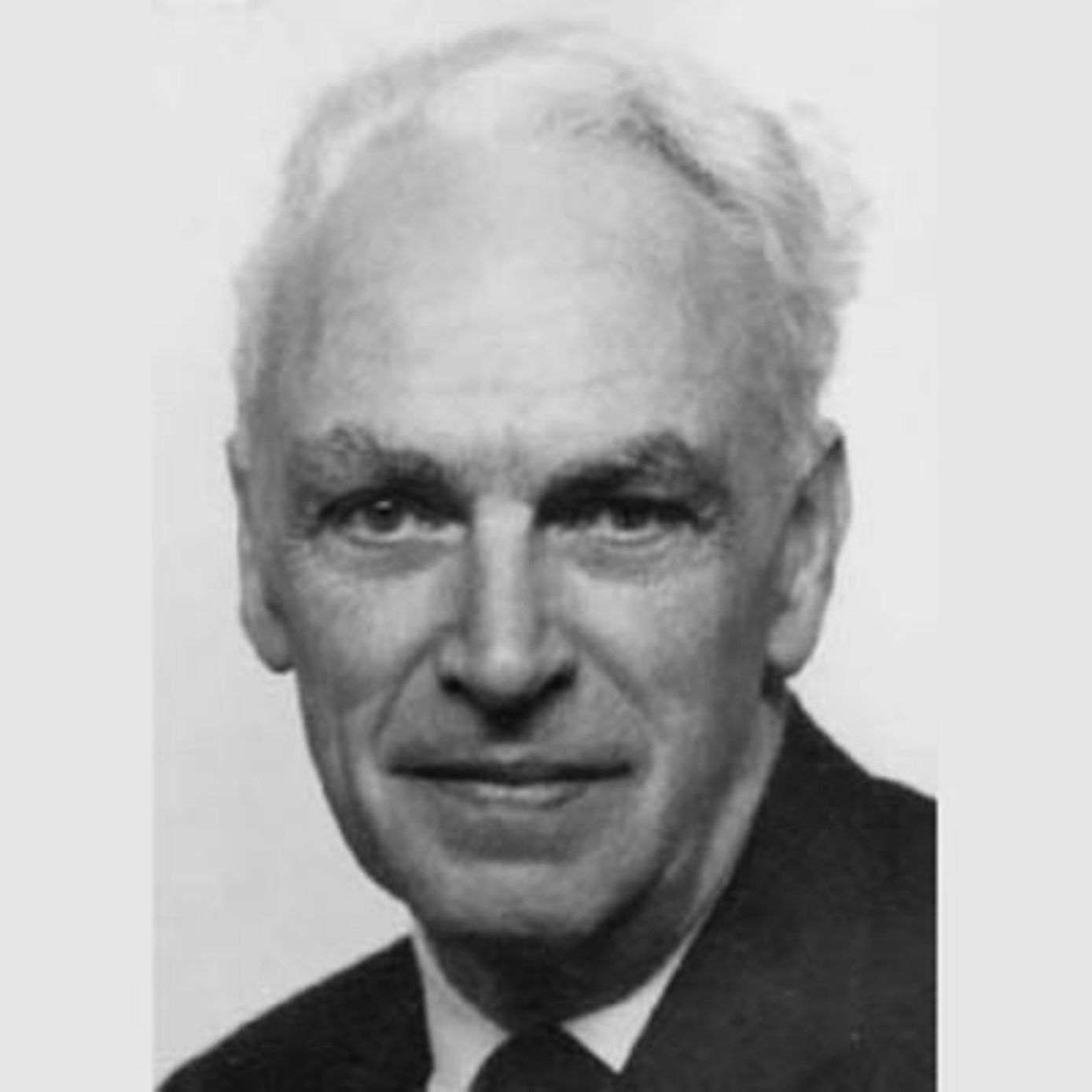}}
\scalebox{0.13}{\includegraphics{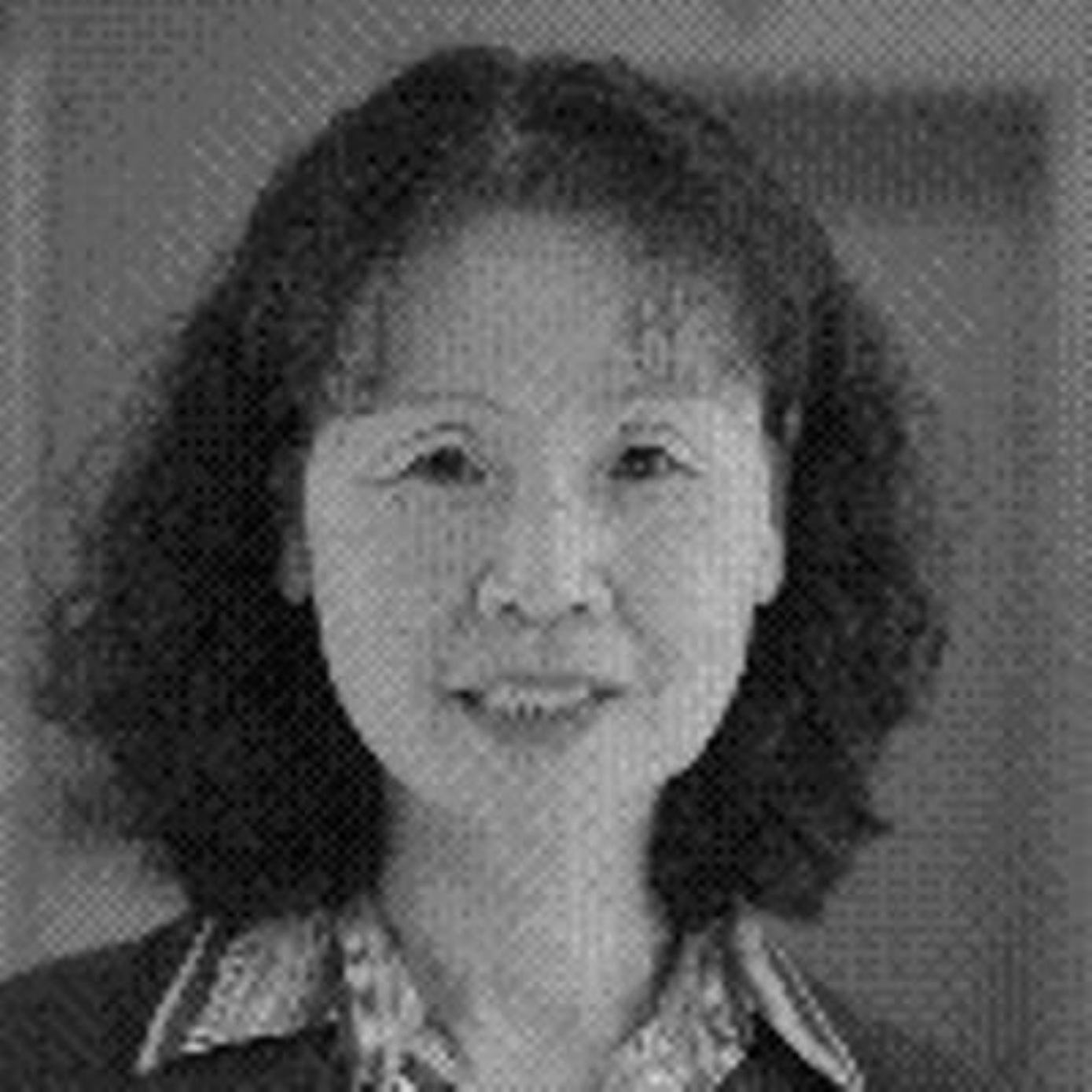}}
\scalebox{0.13}{\includegraphics{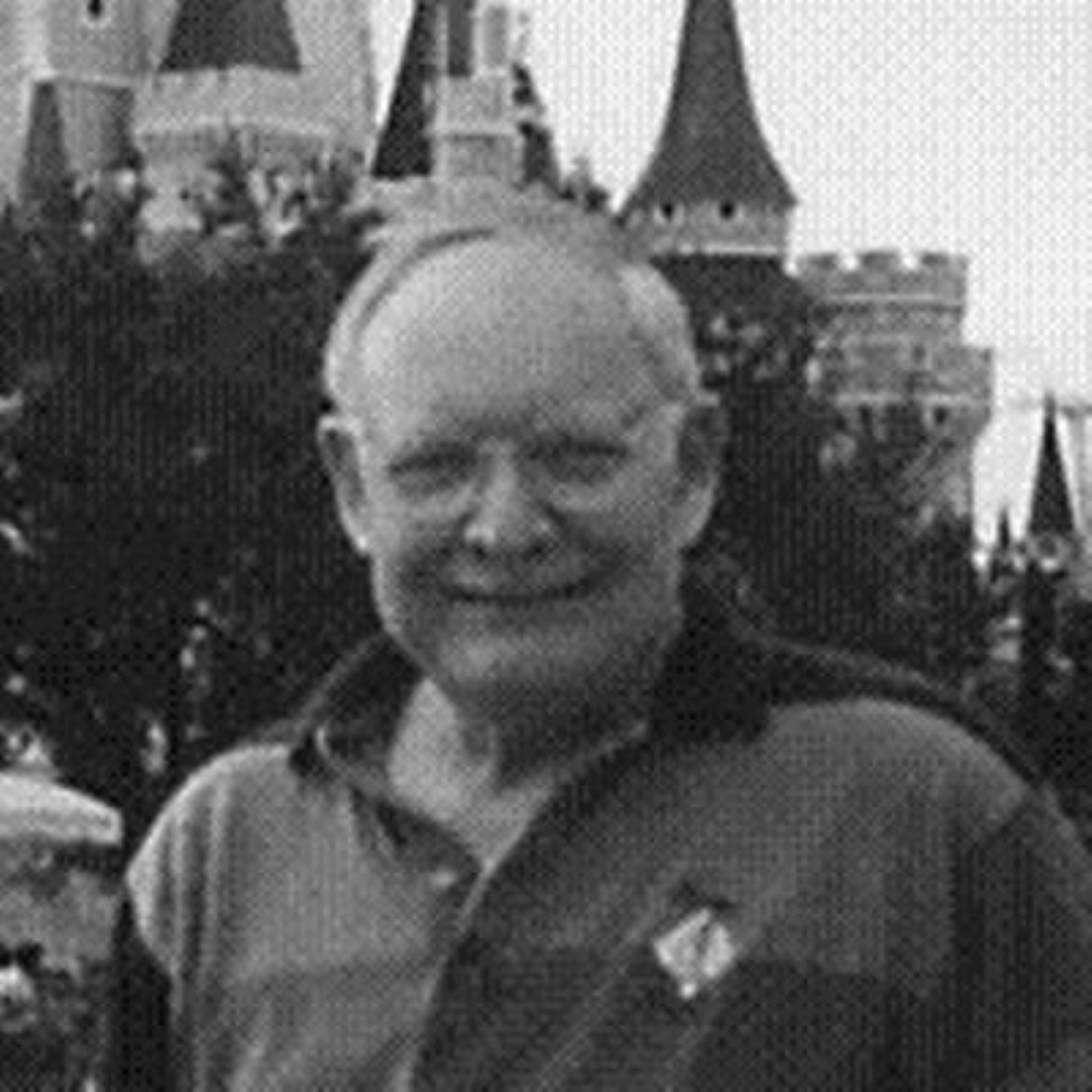}}
\scalebox{0.13}{\includegraphics{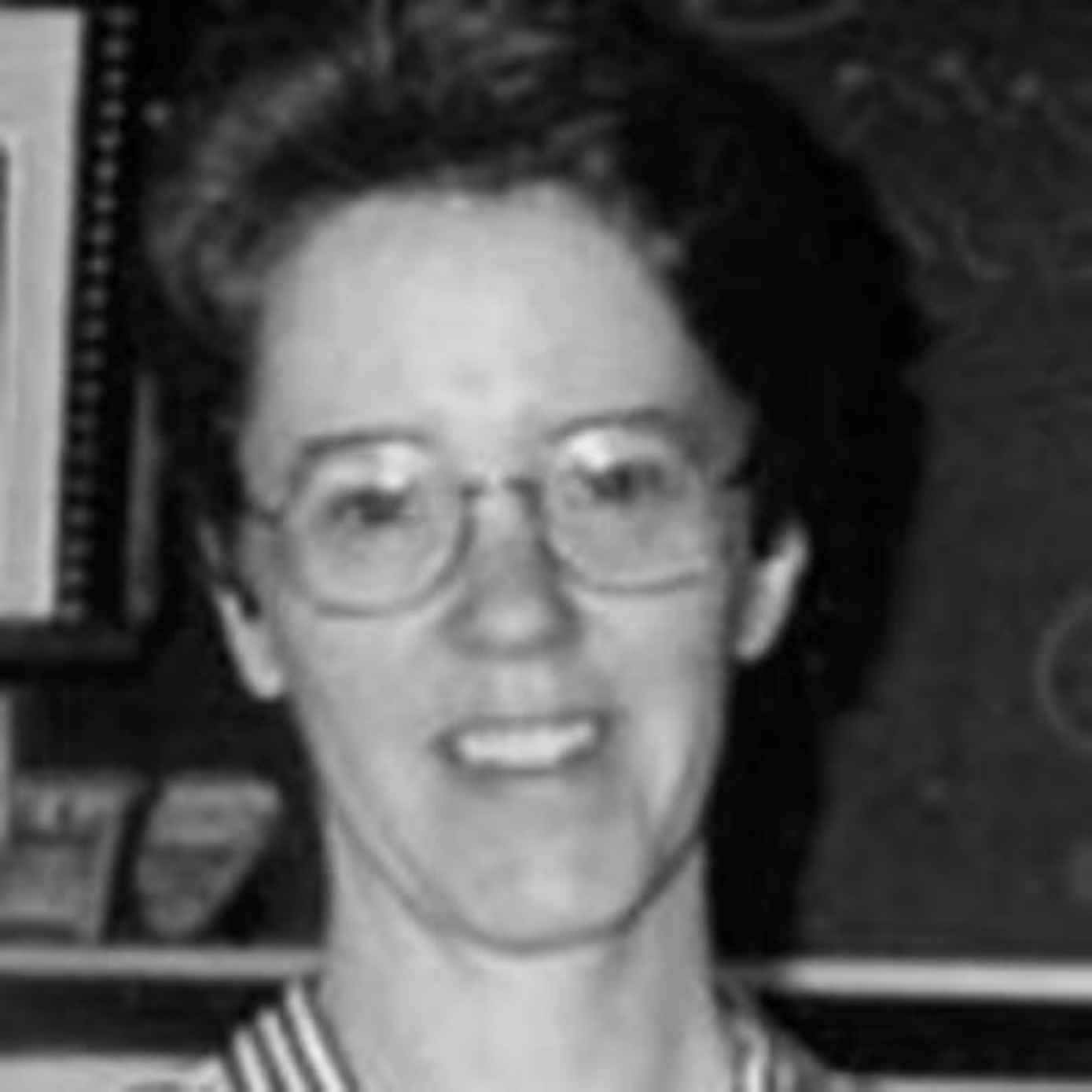}}
\scalebox{0.13}{\includegraphics{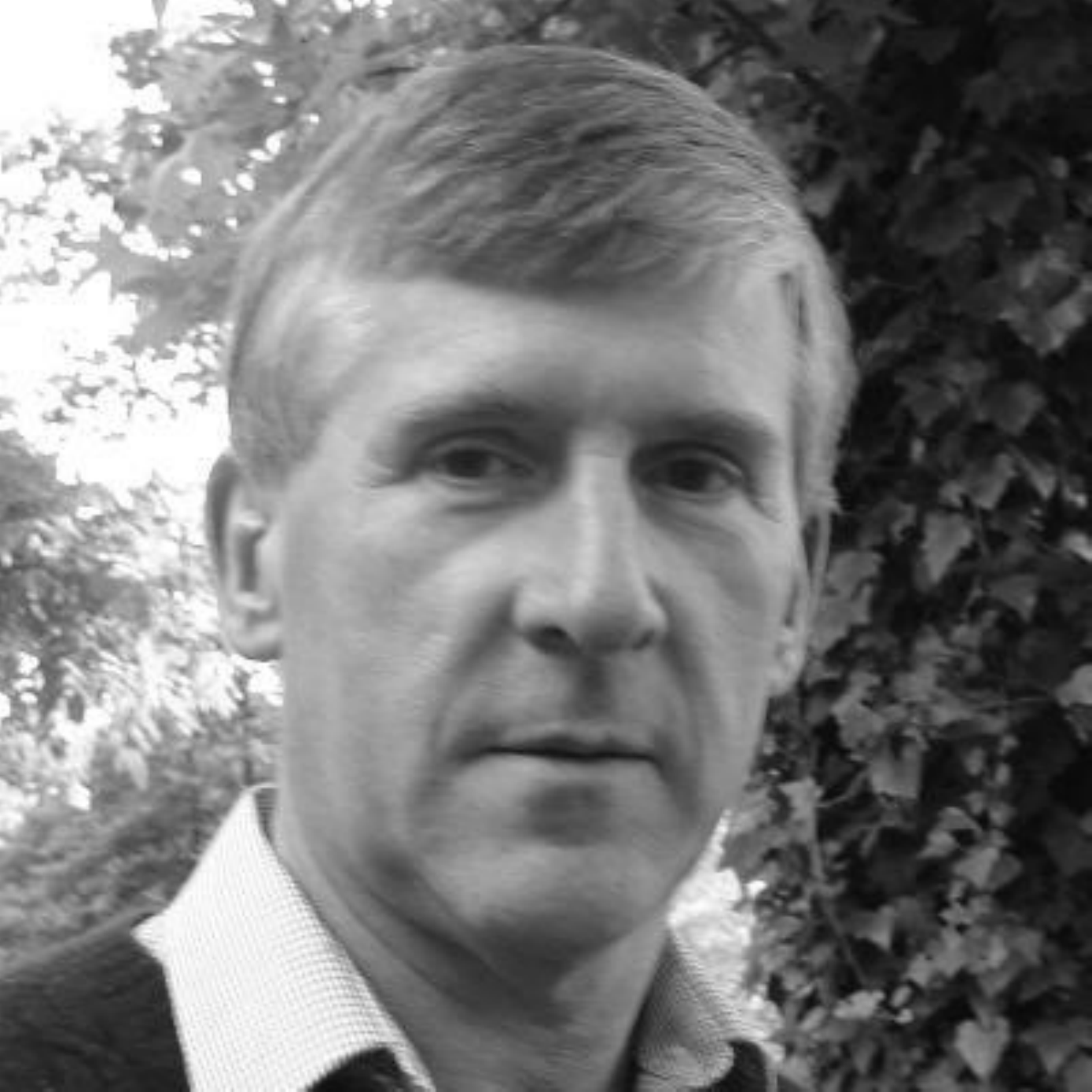}}
\scalebox{0.13}{\includegraphics{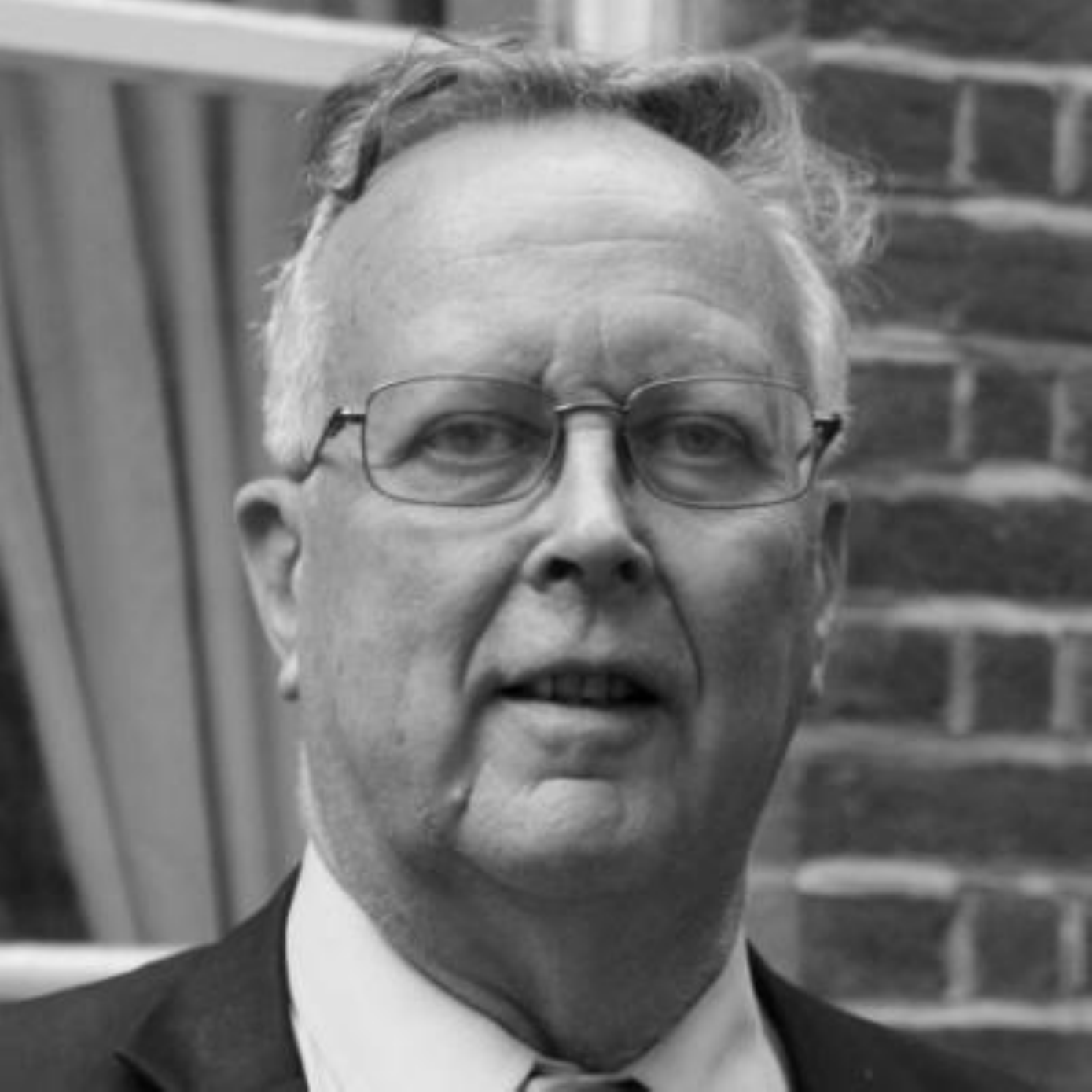}}
\caption{
Ore, Zhang, Chartrand, Hutchinson, Seymour, Wilson}
\end{figure}
\begin{figure}
\scalebox{0.13}{\includegraphics{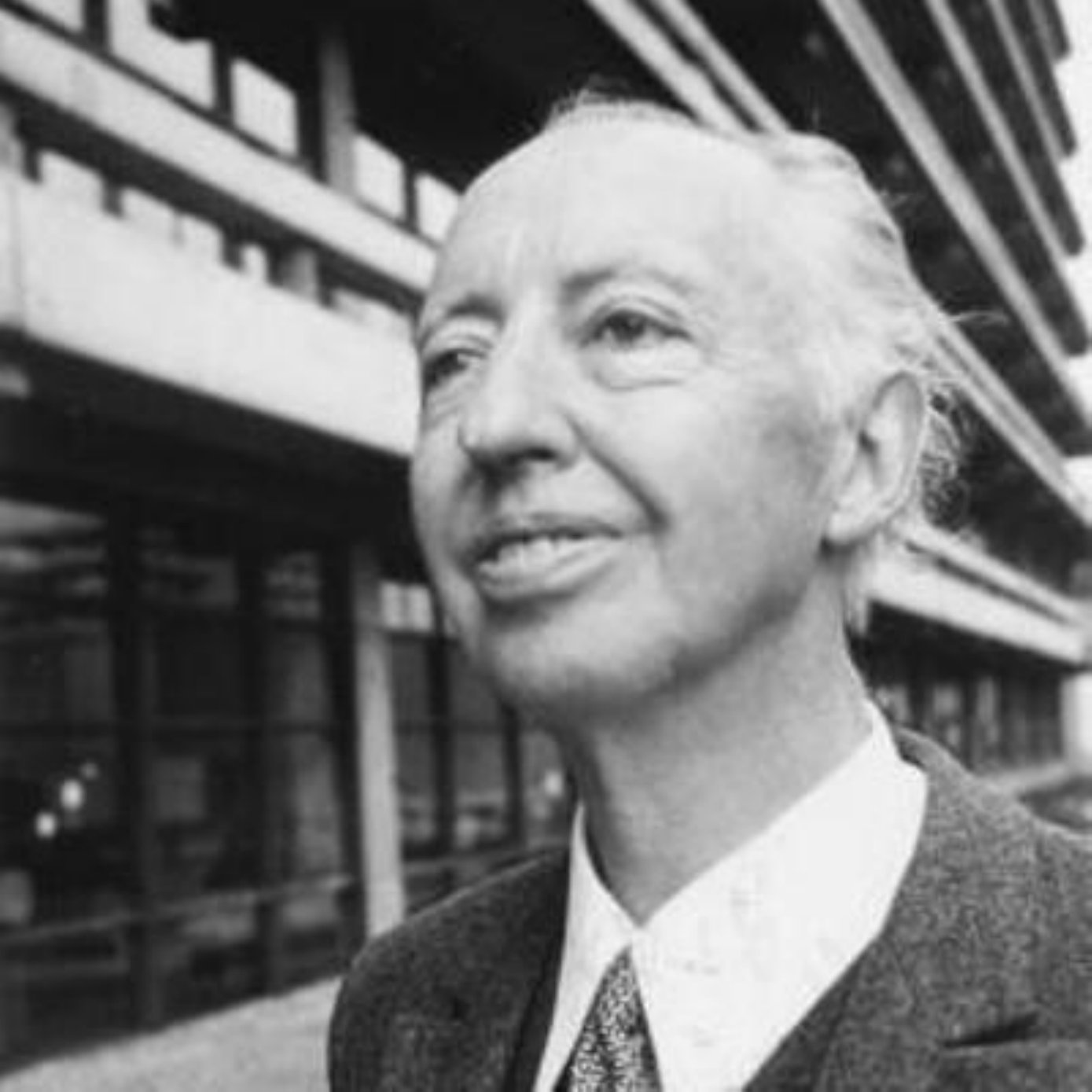}}
\scalebox{0.13}{\includegraphics{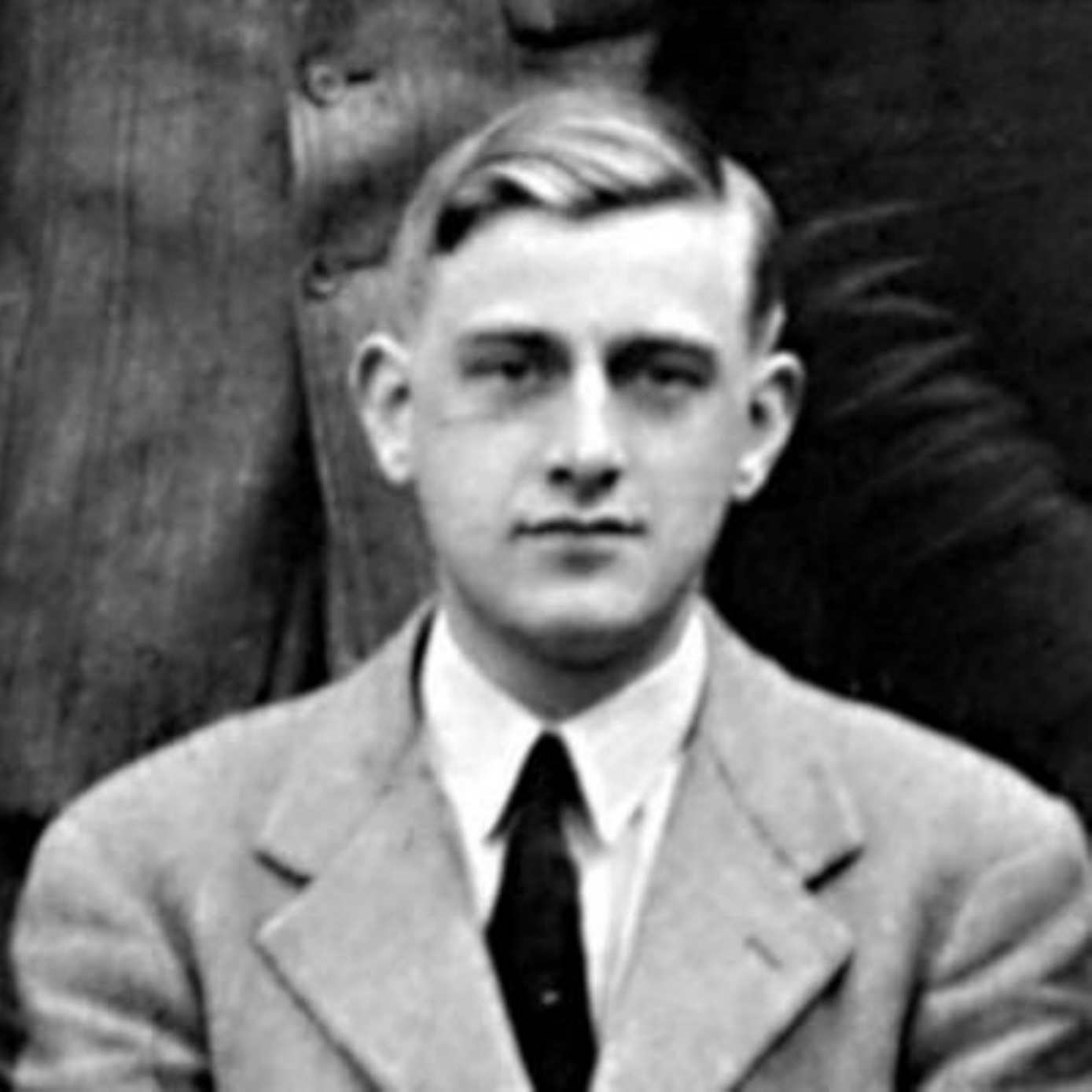}}
\scalebox{0.13}{\includegraphics{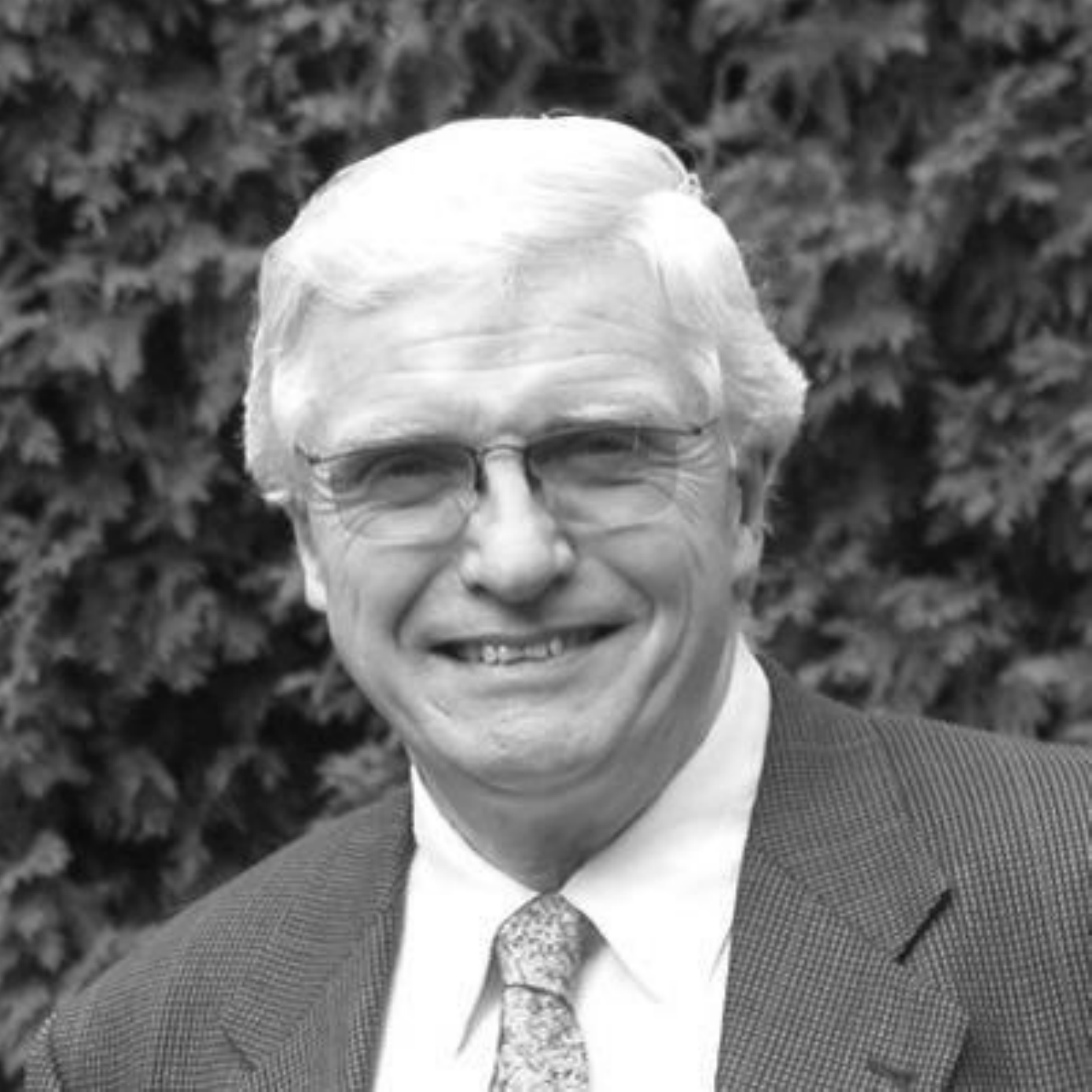}}
\scalebox{0.13}{\includegraphics{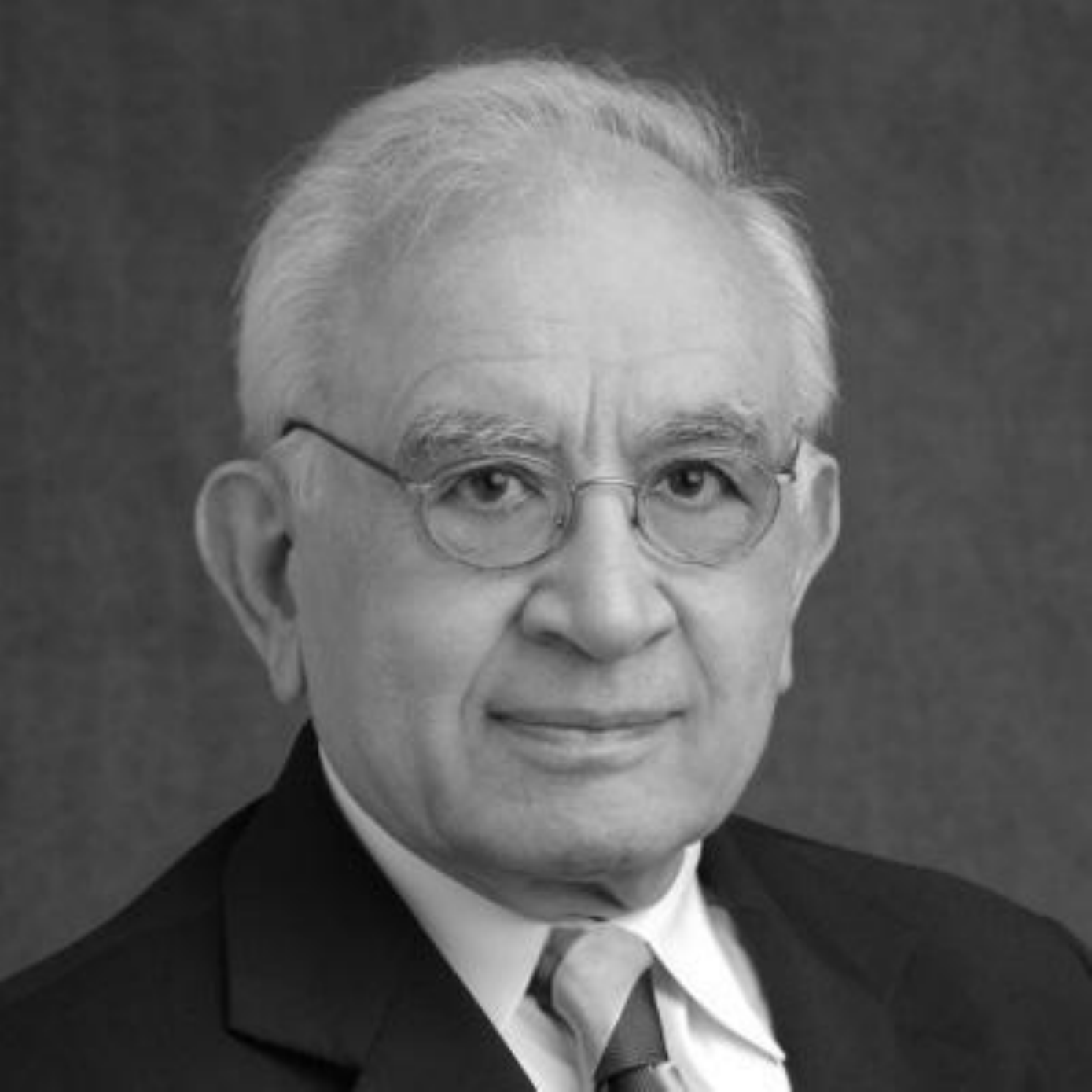}}
\scalebox{0.13}{\includegraphics{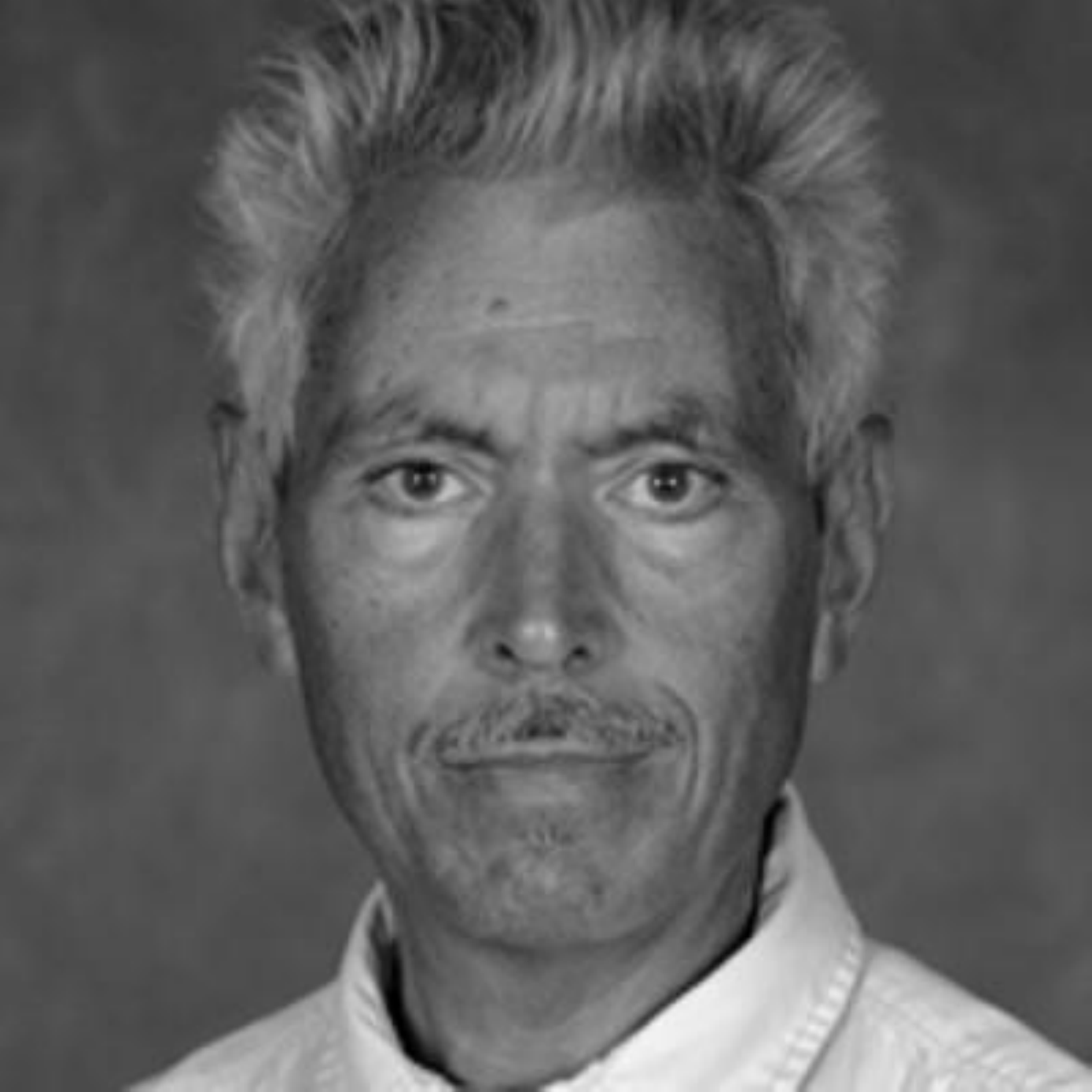}}
\scalebox{0.13}{\includegraphics{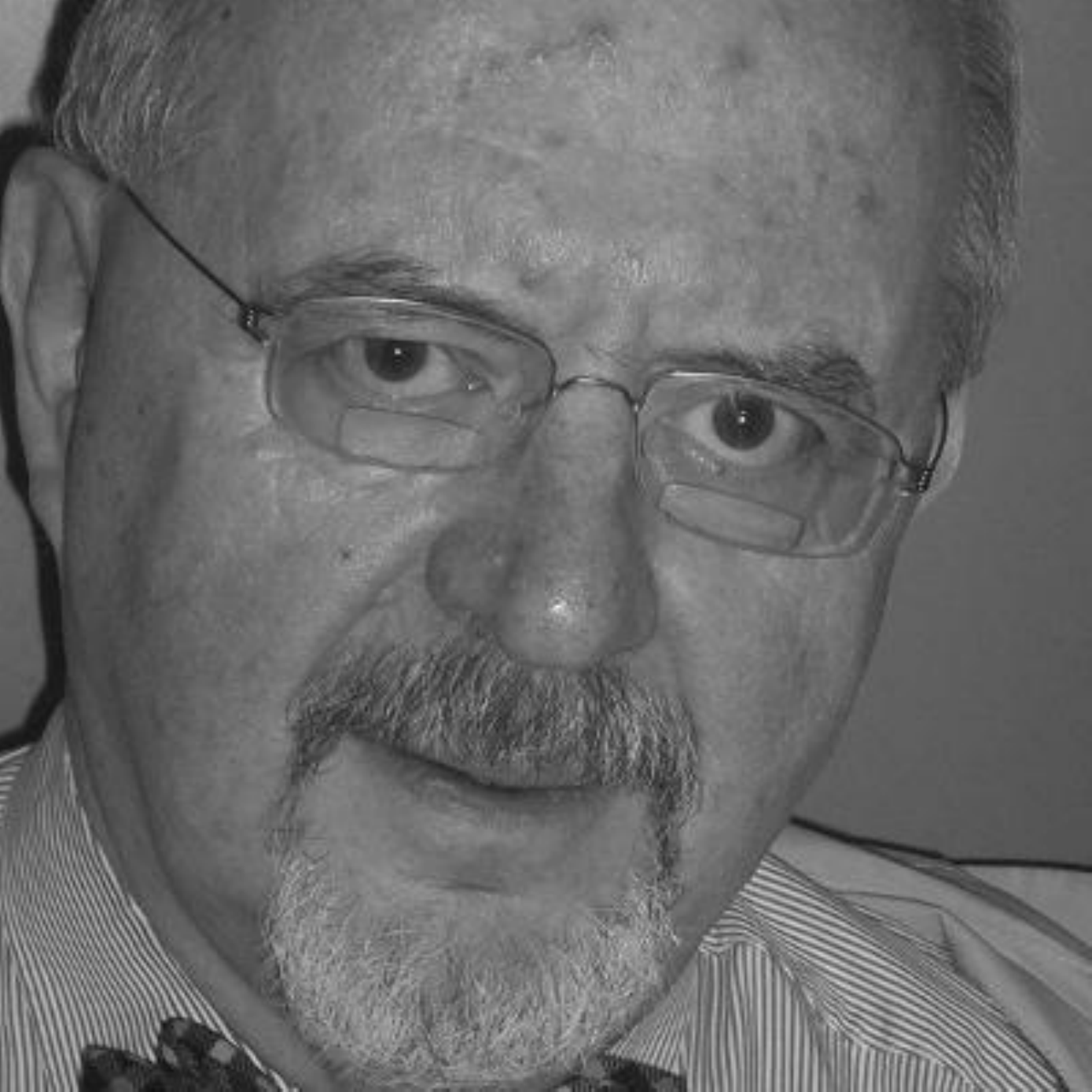}}
\caption{
Hadwiger, Brooks, Robertson, Saaty, Kainen,Fritsch,
}
\end{figure}

\section*{Appendix: some quotes}

The book \cite{Soifer} collects already many statements about the four color theorem.
What follows is a mostly complementary compilation intersecting with the quotes in \cite{Soifer} 
only in the Erdoes quote: \\

According to Haken \cite{Haken2013}, mathematicians
working on the 4 color problem were {\it ``either "optimists" or "modernists", where "optimists"
hoping for an enlightening proof"}. \\

In \cite{Barnette}: 
{\it "Many great mathematicians of the twentieth century have worked seriously on it
while almost every mathematician has given it at least a few idle thoughts.} and
on page 126: {\it "The four color conjecture has probably been worked on by more mathematicians
than any other problem.}. \\

In \cite{Soifer} is a quote of Jensen and Toft: 
{\it ``Does there exist a short proof of the four-color theorem ... in which all the details 
can be checked by hand by a competent mathematician in, say, two weeks?"} \\

Halmos \cite{Halmos1990}: {\it ``Except for that, however, feel that we humanity
learned mighty little from the proof; I am almost tempted to say that as mathematicians
learned nothing at all."} An other, longer quote of Halmos can be found in \cite{Soifer}. \\

Stewart \cite{StewartConcepts} warns:
{\it ``We have seen that in Hilbert's day the prevailing belief
among mathematicians was that every true theorem must have a proof, and how naive
Goedel showed this belief to be. It is probably just as naive to imagine that every
provable theorem has an intellectually satisfying proof".} \\

Erdoes writes in 1991 (this is in \cite{Soifer}) 
{\it ``I would be much happier with a computer-free proof
of the four color problem, but I am willing to accept
the Appel-Haken proof - beggars cannot be choosers."} \\

Davies in 2010: \cite{Davies2010}
{\it ``The proof of the four colour theorem by Appel and Haken in 1976
started a new era in which the proofs of some theorems depended on a
combination of human insight and laborious computer-based calculations.
Some mathematicians found this upsetting, but the next generation have
learned to accept it. As computers have grown more powerful, further
examples of a similar type have appeared, and there are now Fields in which
some of the key results depends on the use of computers.  Our own human
efforts are by no means error free, and it is fortunate that the kind of mathematical 
thinking that we are particularly good at complements the kind that computers are good at."} \\

On the same theme \cite{SaatyKainen} (p. 98):
{\it "To use the computer as an essential tool in their proofs, mathematicians
will be forced to give up hope to verifying proofs by hand, just as scientific
observations made with a microscope or telescope do not admit direct tactile
confirmation.}. \\

In \cite{Heesch}: 
{\it "An other aspect is tragic: for almost fifty years, 
Heinrich Heesch worked intensively on solving the four color theorem,
one of the great century problems in mathematics, sometimes at the limit
of a breakdown. He was the first to use computers to attack the problem, 
sometimes to an incredible amount."} \\

In \cite{KrantzPark}: 
{\it "The eminent mathematician Felix Klein in G\"ottingen heard of the problem 
and declared that the only reason the problem had never been solved is that no 
capable mathematician had ever worked on it. He, Felix Klein would offer a class,
the culmination of which would be a solution to the problem. He failed."} 
Klein seemed not have been the only one: In \cite{Wilson} (page 105) there is
similar story about Hermann Minkowski: {\it "This theorem has not yet been proved, but
that is because only mathematicians of the third rank have occupied themselves with it."}
with a similar attempt to prove the theorem in class, leading to an outburst "Heaven
is angered by my arrogance". \\

In the obituary to Appel, \cite{AppelObituary}, Edward Frenkel is cited to state about the 
proof of the four color theorem:
{\it "Like a landmark Supreme Court case, the proof's legacy is still felt and hotly debated,"} \\

And according to Noga Alon in 1993 \cite{ChartrandZhang2} p.147 
{\it Graph coloring is arguably the most popular subject in graph theory}. \\

\section*{Appendix: definitions}

Before giving an overview of definitions whose common theme is that they have 
an inductive setup, lets try to give some more references: \\

Geometric graph theory or discrete differential geometry has been developed in various forms, 
in the context of computational geometry \cite{GSD,Devadoss,ComputationalGeometry}, 
computational physics \cite{Regge}, integrable systems \cite{BobenkoSuris}, 
computer graphics or computational-numerical methods \cite{Bobenko,CompElectro2002}), 
discrete Morse theory \cite{forman95,forman98,Forman1999,Forman2003},
discrete algebraic geometry \cite{BakerNorine2007,PostnikovShapiro,Biggs} as well as in or own,
more naive work. Algebraic topology on graphs should be credited to
Poincar\'e, whereas Kirchhoff as the discoverer of Stokes theorem. See  
\cite{ArchimedesFunctions} for some history on calculus. 
While Poincar\'e defined the matrices $d_j$, he did not combine them yet to form the 
Dirac operator $D=d+d^*$ producing the form Laplacian $L=D^2$ decomposing to the individual form
Laplacians $L_k$ whose nullity is the $k$th Betti number. It is for us a simple tool to compute 
topological quantities quickly using the Hodge connection. The use of
the heat flow ${\rm exp}(-t L_k)$ to prove the Hurewicz homomorphism suggested here looks new, 
but might have appeared in the continuum using  Rham currents. We plan to elaborate on
Hurewicz more as the graph theoretical part elsewhere as the discrete connection could become 
quite approachable also in higher dimensions. \\

Our graph theoretical definition of geometric graphs repeated below avoids to specify what a 
``triangulation" is. Already in two dimensions, triangulations can contain tetrahedra. We will see in the
proof of the sphere characterization lemma that 4-disconnected triangularizations can lead to unit 
spheres which contain star graphs. 
In higher dimensions, there are even more difficulties. There can be exotic triangulations, 
where the unit sphere is not simply connected. 
In order to avoid any continuum theories or triangularization one has to have a notion of dimension
in the discrete. 
There is a plethora of definitions for dimension in graph theory. The inductive
dimension \cite{elemente11,randomgraph} is a discrete Menger-Uryson dimension.
This does everything right: a discretization of a $d$-dimensional manifold using geometric 
graphs in $\Gcal_d$ has discrete dimension $n$. The discrete spaces have the same homotopy, 
cohomology and cobordism groups as in the continuum. 

\definition{
Define ${\rm dim}(\emptyset)= -1$ and  inductively
${\rm dim}(G) = 1+\frac{1}{|V|} \sum_{v \in V} {\rm dim}(S(v))$.
}

After having developed homotopy notions for Poincar\'e-Hopf \cite{poincarehopf} and Ljusternik-Schnirelmann  \cite{josellisknill}, we found 
that it was introduced in \cite{I94} already and simplified in \cite{CYY} to become the version we have 
worked with in \cite{josellisknill}. The notion traces back to Whitehead. 

\definition{
The one point graph $K_1$ is {\bf contractible}. Inductively, $G=(V,E)$ is {\bf contractible}, 
if there is $x \in V$ such that both the unit sphere $S(x)$ and $V \setminus \{x \;\}$ 
are contractible. }

\definition{
A {\bf homotopy step} consists of removing a vertex $x$ from $G$ for which $S(x)$ is 
contractible or applying the inverse operation of adding a new vertex and connecting 
it to a contractible subgraph $G$. Two graphs are {\bf homotopic} if one can get from one to the other by 
applying a finite sequence of homotopy steps. }

\definition{
Let $\Gcal_0$ denote the set of graphs without edges,
let $\Scal_0 \subset \Gcal_0$ be the set of graphs with two vertices and no edges and let
$\Bcal_0 \subset \Gcal_0$ be the set of graphs with one vertex. Assume $\Gcal_{d-1},\Scal_{d-1},\Bcal_{d-1}$
are known, define {\bf geometric graphs} $\Gcal_d = \{ G \; | \;  S(x) \in \Scal_{d-1}$
or $S(x) \in \Bcal_{d-1} \; \}$. For $G \in \Gcal_d$, the {\bf boundary} is
$\delta G = \{ x \; | \; S(x) \in \Bcal_{d-1} \; \}$
and $G \setminus \delta G$ is the {\bf interior} which we ask to be nonempty. The set of {\bf balls}
$\Bcal_d \subset \Gcal_d$ is defined as the set of contractible graphs in $\Gcal_d$ for which the
boundary is in $\Scal_{d-1}$.
The set of {\bf spheres} $\Scal_d \subset \Gcal_d$ is the set of non-contractible graphs for which
removing an arbitrary vertex produces a graph in $\Bcal_d$. }

One can define subgraph homotopy deformations in the usual way. For closed curves in $G$ for example
one can specify two different homotopy deformation steps like adding a backtracking step or replacing
an edge with the two other edges of a triangle. The following
definition is more elegant, as it refers entirely to the already defined Ivashchenko homotopy: 

\definition{
A graph homomorphism $f:H \to G$ defines a graph $(f,H,G)$ containing the union of the graphs
$H$ and $G$ as well as the edges $(x,f(x))$. We call it the {\bf immersion graph}. 
Two immersions $H \to G$ and $H' \to G$ are called {\bf homotopic} if their
immersion graphs are homotopic as graphs. }

The notion of edge degree which is crucial in our story is certainly present in the Fisk theory.

\definition{
The {\bf degree} of a $(d-2)$-dimensional simplex $x$ in $G \in \Gcal_d$ with $d \geq 2$
is defined as the number of $d$-dimensional simplices in $G$ which contain $x$. }

For geometric three dimensional graphs one could define

\definition{
Call  $1-{\bf deg}(e)/6$ the {\bf geometric Ricci curvature} of $G$ at the edge $e$. 
}

We have then seen as a consequence of Gauss-Bonnet that the sum over all edge degrees at 
a vertex is even as it is $2 V_0$, where $V_0$ is the number of vertices of  $S(x)$. 
The sum over all Ricci curvatures at a vertex is then $2 V_0/3$ which is an even more naive
definition of scalar curvature when scaled and centered so that it is zero in the flat case. 
A definition of Ricci curvature is good if one can derive results like Schoenberg-Myers 
linking positive Ricci curvature with bounds on the diameter of the graph. \\

The notion of curvature in two dimensions can already be traced back to Wernicke \cite{Wernicke}.
It was extensively used in chromatic graph theory, especially by Heesch, who introduced a discharging method 
which very much looks algebro-geometric nowadays, as it can be seen as chip-firing of charge distributions 
which can be seen as divisors. Curvature and indices play an important role in Baker-Norine theory.
The later theory looks (like most graph theory literature) at graphs as one-dimensional objects. \\

Curvature $K(p)=6-V_1(p)$ for graphs was crucial for Heesch \cite{Presnov1990,Presnov1991} as a starting
point for discharging. It was realized early on that $\sum_p K(p)=12$ is a Gauss-Bonnet formula for planar graphs.
An other graph theoretical curvature defined in \cite{Gromov87} is
$K(p) = 1-\sum_{j \in S(p)} (1/2-1/d_j)$, where $d_j$ are the cardinalities of the neighboring face degrees
for vertices $j$ in the sphere $S(p)$. In the case of triangularizations $d_j=3$, it becomes $K(x)=1-V_1/6$. \\

Since curvature involves second derivatives, we have experimented with second neighborhood 
notions in \cite{elemente11} which shows that already in the simplest case of discrete flat planimetry,
some smoothness of the regions is required in order that Hopf-Gauss-Bonnet holds. \\

The following notion is a discrete Euler curvature form which satisfies in the continuum the Gauss-Bonnet-Chern
theorem, where it is only defined in even dimensions as it involves the Pfaffian of curvature form \cite{Cycon}.

\definition{
Assuming $V_{-1}=1$ and counting with $V_k$ the $K_{k+1}$ subgraphs in $S(x)$, 
the {\bf curvature} of a vertex in a finite simple graph is 
$K(x) = \sum_{k=0}^{\infty} (-1)^k \frac{V_{k-1}(x)}{k+1}$. }

Of course, there are discrete versions of sectional, Ricci and scalar curvature which are ``naive" 
in comparison with popular notions \cite{Ollivier,LinLuYau,JostLiu}. The following definition works
for any finite simple graph: 

\definition{
A {\bf two-dimensional section} in a graph $G \in \Gcal_d$ is a wheel subgraph $W_n \in \Gcal_2$. 
The sectional curvature is the curvature of the center point of this wheel graph. 
The {\bf Ricci curvature} for an edge $e$ is the average over all sectional curvatures for all
wheel graphs containing the edge. 
The {\bf scalar curvature} for a vertex $x$ is the average over all Ricci curvatures of all 
edges connected to $x$. }

To defend these simple notions, lets mention that integral geometry allows to deform the apparently
rigid notions pretty arbitrarily. Just change the probability space on functions (wave functions) 
and define the curvature as the expectation ${\rm E}[i_f(x)]$ of the index with respect to an other measure. 
This allows to deform curvature similarly and add flexibility similarly as changing the  metric $g$ does in 
Riemannian geometry. The advantage of the probabilistic link is that it works in the same way 
in the continuum and discrete case, that Gauss-Bonnet-Chern becomes an immediate consequence 
of Poincar\'e-Hopf.  \\

The integral geometric connection can be traced back to \cite{Banchoff67,Banchoff70} but again
in a piecewise linear manifold setup and not in graph theory.  \\

Finally, lets look at the generalization of geometric graphs which allows to describe 
graphs which are ``variety like". This is necessary for example when describing 
``algebraic varieties" in a host graph $G$ (where $G$ plays the role of the affine or projective space )
given as a completion of graphs formed by $k$-simplices for which $k$ functions $f_1, \dots , f_k$ 
change sign and two $k$ simplices are connected, if they share a common $(k-1)$ simplex.  \\

Lets first define the lass of {\bf uniformly $d$-dimensional graphs} $\Dcal_d$:

\definition{
Let $\Dcal_0 = \Gcal_0$ be the class of graphs without edges. 
Define inductively $G$ to be in $\Dcal_d$ if every unit sphere is $(d-1)$-dimensional and 
in $\Dcal_{d-1}$. }

A tree without seeds or the figure 8 graph are examples of a uniformly one-dimensional graph. 
Graphs in $\Gcal_d$ are examples of uniformly $d$ dimensional graphs but $\Dcal_d$ is 
much larger as it contains for example graphs for which some unit spheres are 
not in $\Scal_{d-1}$. 

\definition{
The set $\Vcal_d$ of $d$-dimensional varieties is defined inductively:
start with $0$-dimensional varieties as $\Vcal_0=\Gcal_0$ and define $\Vcal_d$ as a subclass of $\Dcal_d$ 
for which $S(x) \in \Vcal_{d-1}$ for all vertices $x$ and  $\sigma(G)$ and $\delta(G)$ are both in 
$\Vcal_{d-1} \cup \Vcal_{d-2} \cup \cdots \cup \Vcal_{0}$.
The set  $\sigma(G)$ of {\bf singularities} 
of $G \in \Vcal_d$ is defined as the graph generated by vertices $x$ for which $S(x)$ is in 
$\Vcal_{d-1} \setminus (\Scal_{d-1} \cup Bcal_{d-1})$. }

A tree without seeds for which no branch points are neighboring is an example of a one dimensional 
variety.  \\

In other words, we want a $d$-dimensional variety to be a graph which at most places looks like a 
geometric $d$-dimensional graph possibly with boundary. Some singularities or a boundary is allowed.
Both the boundary set $\delta(G)$ and the singularity set $\sigma(G)$ however should each should form
a smaller dimensional variety. 
This is motivated from the continuum, where singularities are defined by the zero locus of some equations
forming a lower dimensional variety. While this is not equivalent, we have to compensate in the
graph theoretical setup for the fact that we do not have traditional calculus available to define 
singularities. \\

The introduction and pictures and available theorems in graph theory which completely parallel the 
story in the continuum should indicate that we do not have to lament about the lack of calculus on
graphs. The same can be said about differential geometry and more recently about algebraic geometry. 
Things just looks a bit different, but what counts are the theorems. 

\section*{Appendix: characterization of the sphere}

Here is a proof of our elementary $\Scal_2$ characterization lemma given in the first section.
The lemma is important to us as it is used to show that the statement $\Scal_2 \subset \Ccal_4$ is equivalent 
to the four color theorem. We got dragged into working at colorings of $\Gcal_2$ because we were
initially under the impression that coloring graphs in $\Scal_2$ is much easier than coloring graphs in $\Pcal$.
The reason was that we believed that the geometric tool of level curves would make things possible
by foliating the line graph and coloring it. Without this naive assumption,
we probably would never have started. As level curves are $1$-dimensional and colored with three colors, 
we expected to be able to color the graph. This is a path which Tait followed. 
Lets denote by $\Wcal$ the class of $4$-connected 
maximal planar graphs. The letter $\Wcal$ is used because Whitney studied them in 1931 first, 
showing that all graphs in $\Wcal$ are Hamiltonian. The lemma now can be reformulated as:

\resultlemma{ $\Scal_2 = \Wcal$. }

\begin{proof}
{\bf (i)} Assume first that $G \in \Scal_2$. This means that every unit sphere $S(x)$ is a cyclic graph with $4$
or more elements and that removing one vertex of $G$ produces a contractible graph. \\

{\bf a) $G$ is tetrahedral-free:} otherwise, the unit sphere $S(x)$ would contain a triangle.  \\
{\bf b) $G$ is 4-connected:} assume it is not and removing three vertices makes the graph disconnected. 
Restricting to a unit disc $B(x)$ we see that the path $a-x-b$ must cut it into two.
As the vertices $a$ and $b$ are connected within $G$ differently by maximality
we conclude that $(a,x,b)$ is a triangle in $G$ which separates the graph. 
Having a separating triangle implies that all unit spheres $S(a),S(b)$ and $S(x)$ contain a star graph $S_3$ 
contradicting the requirement that $S(x)$ must be cyclic.\\
{\bf c) $G$ is planar:} removing one vertex in $\Scal_2$ renders the graph contractible. It is now a 
$2$-dimensional graph with boundary (such a ``disk" is sometimes also called a ``configuration"): 
the interior points have cyclic unit spheres, the boundary points have line graphs as unit spheres.
That such a graph can be drawn in the plane can 
be seen by induction, starting the smallest of this kind, the wheel graphs.
Every homotopy extension step can be drawn and at every stage the boundary of the graph is a cyclic graph. 
The extension step is to build a pyramid extension over a line subgraph of the boundary.\\
{\bf d) $G$ is maximal planar:} adding an other edge $(x,y)$ would have the effect that 
the unit sphere $S(x)$ would separate $x$ from any other vertex by the Jordan curve theorem. \\

{\bf (ii)} 
To show the reverse statement, assume now that $G$ is a 4-connected maximal planar graph $G$. (This proof could
be shortened by using the fact that maximal planar graphs are triangularizations. Lets try to make it self-contained.) 
A {\bf face} of $G$ is a closed path in the graph which encloses a connectivity component $A \subset \R^2$
of the complement of a planar embedding $\tilde{G}$ of $G$ as a $1$-dimensional simplicial complex. \\

{\bf a) Every face of $G$ is a triangle:} this follows from maximal 4-connected planarity: if we had $n$-gon face with $n>3$, 
   we could add additional diagonal connections without violating either 4-connectivity nor planarity. 
   This would contradict maximality. \\
{\bf b) The Euler characteristic $\chi(G)$ satisfies $\chi(G)=2$:} removing one vertex of degree $k$ also removes $k$
   edges and $k$ triangles. It reduces the Euler characteristic by $1-k+k=1$. As the result is contractible and contractible
   graphs have Euler characteristic $1$, the graph $G$ has Euler characteristic $2$. \\
{\bf c) The graph $G$ can not contain a tetrahedron $t$:} otherwise, removing a triangle in $t$
   would isolate the connectivity component of the fourth point in $t$ contradicting 4-connectivity. This again implicitly
   made use of the Jordan curve theorem because a triangle separates the inside from the outside. \\ 
{\bf d) A unit sphere $S(x)$ is triangle-free:} otherwise we had a tetrahedral unit ball $B(x)$ contradicting c). \\
{\bf e) The unit sphere $S(x)$ is connected:} if it were disconnected, then by maximality,
   an other connection in $S(x)$ could be added without violating maximal 4-connected planarity. \\
{\bf f) The degree of every $y \in S(x)$ within $S(x)$ is larger than $1$:} assume $y \in S(x)$ had only one neighbor $a$
   in $S(x)$. As removing $x,y$ keeps the graph connected, $a$ is connected via a path to an other point $b \in S(x)$
   and so directly connected to a neighboring point $y \in S(x)$
{\bf g) The degree of $y$ in $S(x)$ can not be $3$:} removing $x,y,b$ with not interconnected $(a,b,c)$ neighboring
   $y \in S(x)$ would render $G$ disconnected as a path $(a,\dots ,h, \dots,c)$ would lead to a homeomorphic
    copy of $K_{3,3}$ containing vertices $(h,y,x,a,b,c)$ inside $G$ contradicting the Kuratowski theorem. \\
{\bf h) Each unit sphere is cyclic:} from steps d)-g) follows that $S(x)$ is cyclic with $n \geq 4$ as every vertex $y$ in $S(x)$ has exactly $2$ neighbors. \\
{\bf i) The claim:} we know now from h) that $G \in \Gcal_2$. From step b) follows that $G \in \Scal_2$. 
\end{proof}

{\bf Remarks.} \\

a) {\bf Maximal planar graphs} play a role when studying
or {\bf minimal graphs} which are hypothetical smallest planar graphs needing 5 colors but for which the
removal of a vertex brings it to $\Ccal_4$. These graphs are sometimes dubbed ``irreducable configurations" or
``minimal criminals".  \\

b) Whitney characterized planar graphs using duality \cite{Ore}. The classical notion of duality requires planarity
since one takes the faces of the graph as the vertices of the dual graph. The notion of ``face" needs an Euclidean embedding.
We have defined the dual graph of a graph in $\Gcal_d$ using maximal simplices.
If the graph was $d$-dimensional, then the dual graph has as vertices the $d$ dimensional simplices in $G$. For $G \in \Gcal_d$
under some conditions like orientability, the dual graph $\hat{G}$ can be completed to be in $\Gcal_d$ again leading to
Poincar\'e duality in $\Gcal_d$. Whitney managed to define a notion of duality (called W-duality in \cite{Ore})
which does not need an Euclidean embedding: it can of algebro-topological nature. If a graph is treated as a one-dimensional
object (which is for us only the case when the graph is triangle free but is often assumed in the graph theory literature)
then $v_0-v_1 = b_0-b_1$ is the Euler-Poincar\'e formula. Now $b_1$ is called the {\bf circuit rank} and $v_1-b_1$ the
{\bf component rank}. The sum of these ranks is the number $v_1$ of edges. As in a duality situation, where faces and vertices
switch their role and edges stay the same, edges obviously play a central role and circuit and component rank switch.
Now Whitney gave an abstract notion of duality by asking (roughly speaking) that two graphs $G_1,G_2$ are W-dual, if they
have the same number of edges and for every subgraph $H_1$ of $G_1$ there is a corresponding subgraph $H_2$ of $G_2$ such that $H_1$ and $H_2$
have switched circuit and component ranks. Whitney proved that a graph is planar if and only if it has a W-dual. \\

In general, this notion is a bit unintuitive, like for a triangle already, where the dual graph is a triangle again,
the duality is much more geometric in maximal planar 4-connected graphs, aka graphs in $\Scal_2$, where according to Steinitz
graphs correspond to convex polyhedra for which the dual graph is a convex polyhedron again, what we usually call the dual
polyhedron. This generalizes to graphs in $\Scal_d$, where the completion of the dual graph is again a graph in $\Scal_d$.
Since the Betti vector $b=(b_0,\dots ,b_d) = (1,0, \dots, 1)$ for any graph in $\Scal_d$, the dual graph has the same
Betti vector $(b_d, \dots ,b_0)$. \\

c) \cite{SaatyKainen} show in Theorem 3-1 that
the class of {\bf triangularizations} (graphs for which every face is triangular) coincides with the class of maximal planar graphs. Then
in Theorem 3-6 and 3-7 that minimal graphs are maximal planar and 5-connected. The condition {\bf maximal planar and 5-connected}
comes close to the characterization of $\Scal_2$ as the class of maximally 4-connected planar graphs but there is no relation:
the tetrahedron is maximal planar and 5-connected but not in $\Scal_2$, the octahedron is in $\Scal_2$ but not 5-connected.  \\

d) A planar graph is called {\bf tetrahedral complete} if every triangle is part of a tetrahedron. Every planar graph
has a tetrahedral completion which is still planar. It is obtained by capping every triangle not part of a tetrahedron
with  an additional vertex. Given a graph $H \in \Scal$, let $\overline{H}$ be the
completion and let $\overline{\Scal}_2$ denote the set of all these completed spheres.
A graph in $\overline{\Scal}_2$ is no more 4-connected but it is a maximal planar graph. If $H \in \Ccal_4$
then $\overline{H} \in \Ccal_4$ because we can color each of the tetrahedron tips with the fourth color.
From the equivalence of $\Scal_2 \subset \Ccal_4$ with the four color theorem follows that the condition
$\overline{\Scal}_2 \subset \Ccal_4$ is equivalent to the four color theorem too.
Note that every $G \in \Pcal$ is a subgraph of a tetrahedral completion $\overline{H}$
of a $H \in \Scal_2$ but every 4-connected planar graph is a subgraph of a completed sphere graph.  \\

e) It is remarkable that {\bf dimension} and {\bf homotopy} together define what we call a locally
Euclidean structure or manifold. Both in the continuum as well as in the discrete, we can define a $d$-sphere as a
$d$-dimensional compact metric space for which the removal of a single point renders the space contractible.
This is how we have defined $\Scal_d$ in the category of graphs.
In \cite{KnillTopology} we noted that this characterizes manifolds among metric spaces:
define the dimension of a compact metric space $(X,d)$ inductively to be $k$ if for  sufficiently small positive radius (or
minimal positive radius in the discrete case), spheres $S_r(x)$ have
dimension $k-1$. A metric space is called a "sphere", if removing one point renders the metric space contractible, where
contractible is defined traditionally.
The inductive assumption is that the empty set has dimension $-1$ and that the one point set is contractible.
A $0$-dimensional sphere is a two point set. A $1$-dimensional sphere is a topological Euclidean circle,
a $2$-dimensional sphere homeomorphic to the standard two sphere etc.
The lemma proven here for graphs relates the topology of the unit sphere with some connectivity and planar embeddability.
An intuitionist would find the notion of "planar graph" unsatisfactory as it refers to Euclidean structure which
assumes more foundation and a stronger Axiom system.
The Kuratowski theorem is important for such a mathematician as it gets rid of Euclidean ballast
and describes planar graphs graph theoretically using finite sets. Dealing with graphs in $\Scal_2$ does the same: it
captures the planarity but is also more intuitive as it is based on dimension and deformation.
Unlike Kuratowski, where detecting homeomorphic images
of $K_{3,3}$ or $K_5$ can be tricky, one can see immediately whether a graph is in $\Scal_2$: check the local sphere condition
at points and a global condition for Euler characteristic which is constructive as we can punch a hole and deform
the rest to a point or simply count the degrees and add up the curvatures. \\

\pagebreak
 
\bibliographystyle{plain}

\end{document}